\documentclass{amsbook}
\usepackage{amssymb}

\usepackage{ifthen}
\usepackage{pstricks}

\input{diagrams.sty}

\newcommand{\trt}{2}
\newcommand{\lbl}[1]
{\ifthenelse{\trt=3}{\quad\label{#1}\ref{#1}}{\label{#1}}}

\newcommand{\bff}{\bf}

\begin{document}


\newcommand{\Arno}{ V. Arnold,
\emph{Dynamics of intersections},
        Proceedings of a Conference
 in Honour of
J.Moser, edited by
 P.Rabinowitz and R.Hill,
 Academic Press,
 1990
pp. 77--84.  }

\newcommand{\Arnoprob}{ V. Arnold,
\emph{Sur quelques probl\`emes de la th\'eorie 
des syst\`emes dynamiques},
Topological methods in nonlinear Analysis, 
1994, {\bff 4} (2), p. 209-225.}

\newcommand{\cerf}{ J.Cerf, 
\emph{La stratification naturelle des espaces de 
fonctions diff\'erentiables
r\'eelles et le th\'eor\`eme de la pseudo-isotopie},
        Publ. IHES, {\bff 39}, 1970
 }

\newcommand{\farrell}
{T.Farrell,
\emph{The obstruction to
fibering a manifold over a circle},
Indiana Univ.~J.~
{\bff 21}
(1971),
    315--346.
}

\newcommand{\fried}{  D.Fried,
\emph{Homological Identities for closed orbits},
Inv. Math.,
{\bff 71},
(1983)
419--442.
}

\newcommand{\latour}
{F.Latour,
\emph{Existence de 1-formes ferm\'ees non singuli\`eres dans une classe 
de cohomologie de de Rham}, Publ. IHES {\bff 80} (1994)}

\newcommand{\laudsiko}
{   F.Laudenbach, J.-C.Sikorav,
\emph{
Persistance d'intersection avec la section
nulle au cours d'une isotopic hamiltonienne
 dans un fibre cotangent
},
Invent.~Math.~
{\bff 82}
    (1985     ),
pp. 349--357.
 }

\newcommand{\milnWT}
{J.Milnor,
\emph{ Whitehead Torsion},
Bull. Amer. Math. Soc.
{\bff 72}
(1966),
358 - 426.
}

\newcommand{\noviquasi}
{S.P.Novikov,
\emph{Quasiperiodic Structures in topology},
in the  book:  Topological Methods in Modern Mathematics,
 Publish or Perish,
 1993,
pp. 223 -- 235.  }

\newcommand{\patou}
{ A.V.Pazhitnov,
\emph{
On the Novikov
complex for rational Morse forms},
 Annales de la Facult\'e de Sciences de Toulouse,
{\bff 4}      (1995), 297 -- 338.
 }

\newcommand{\pasur}
{ A.V.Pajitnov,
\emph{
Surgery on the Novikov Complex},
K-theory
{\bff 10} (1996),      323-412.
 }

\newcommand{\pamrl}
{  A.V.Pajitnov,
\emph{Rationality and exponential growth
properties of the boundary operators in the Novikov
Complex},
Mathematical Research Letters,
{\bff 3}
(1996),
  541-548.
 }

\newcommand{\peix}
{ M.M.Peixoto, \emph{On an approximation 
theorem of Kupka and Smale},
Journal of differential equations,
{\bff 3}, 214 -- 227 (1966)}

\newcommand{\rani}
{ A.Ranicki,
\emph{Finite domination and Novikov rings},
Topology, 
{\bff 34} $N^\circ 3$
(1995),
  619 -- 633.
 }

\newcommand{\sikothese}
{J.-Cl. Sikorav,
\emph{Th\`ese}}

\newcommand{\sikoens}
{ J.-Cl.~Sikorav,
\emph{Un probleme de
disjonction par isotopie symplectique
dans un fibr\'e cotangent},
 Ann.~Scient.
~Ecole~Norm.~Sup.,
{\bff 19}
 1986,
 543--552.
 }

\newcommand{\smale}
{  S.~Smale,
\emph{On the structure of manifolds},
 Am.~J.~Math.,
{\bff 84} (1962)
 387--399.}

\newcommand{\smapoi}
{  S.~Smale,
\emph{Generalized Poincare's conjecture in dimensions
greater than four},
 Ann.~Math.,
{\bff 74} (1961)
 391--406.}

\newcommand{\thom}
{ R.Thom,
\emph{Sur une partition en cellules associ\'ee
\`a une fonction sur une vari\'et\'e},
     Comptes Rendus de l'Acad\'emie de Sciences,
{\bff 228}
(1949),
 973 -- 975.
}

\newcommand{\witt}
{ E.Witten,
 \emph{Supersymmetry and Morse theory},
 Journal of Diff.~Geom.,
{\bff 17} (1985)
 no. 2.
}


\newcommand{\abrob}
{   R. Abraham, J.Robbin, 
\emph{Transversal mappings and flows
},
 Benjamin, New York,  1967.
}

\newcommand{\biro}
{
G.Birkhoff and G-C.Rota,
\emph{Ordinary differential equations},
Blaisdell Publishing Company, 1962.}

\newcommand{\bousp}
{N.Bourbaki
\emph{El\'ements de Math\'ematiques.
Th\'eories spectrales},
Hermann, 1967.}

\newcommand{\dold}
{   A.Dold 
\emph{Lectures on Algebraic Topology
},
Springer,  1972.
}

\newcommand{\hirsch}
{   M.Hirsch
\emph{Differential Topology
},
 Springer, 1976.
}

\newcommand{\kling}
{   W.~Klingenberg
\emph{Lectures on closed geodesics
},
 Springer, 1978.
}

\newcommand{\milnhcob}
{   J.~Milnor,
\emph{Lectures on the
h-cobordism theorem},
 Princeton University
Press,
 1965.
}

\newcommand{\milnstash}
{   J.Milnor and J.Stasheff,
\emph{ Characteristic Classes},
 Princeton University Press,
 1974.}

\newcommand{\morse}
{  M.Morse,
\emph{Calculus of Variations in the Large},
  American Mathematical Society Colloquium Publications,
Vol.18,
 1934.}

\newcommand{\cohn}
{ P.M.Cohn,
\emph{Free rings and their relations},
  Academic press
( 1971)}

\newcommand{\lam}
{    T.Y.Lam,
\emph{Serre's Conjecture,         }
Lecture Notes in Mathematics
{\bff 635},
(1978)
227 p.
}


\newcommand{\huli}{   M.Hutchings, Y.J.Lee
\emph{              Circle-valued Morse theory, Reidemeister torsion
and Seiberg-Witten invariants of 3-manifolds},
e-print dg-ga/9612004 3  Dec 1996.}

\newcommand{\mengt}
{ G.Meng, C.H.Taubes
\emph{SW=Milnor Torsion }
Preprint
(1996)
}

\newcommand{\paodense}
{ A.V.Pazhitnov,
\emph{
On the Novikov
complex for rational Morse forms
},
Preprint of Odense University,
 Odense,
  1991.
 }

\newcommand{\paura}
{  A.V.Pajitnov,
\emph{Surgery on the Novikov Complex},
Rapport de Recherche CNRS URA 758,
 Nantes,
  1993.
}

\def\paepr{  A.V.Pajitnov,
\emph{The incidence coefficients in the Novikov Complex
are generically rational functions},
e-print dg-ga/9603006 14 Mar 96.
 }

\def\paepri{ A.V.Pajitnov,
\emph{Incidence coefficients in the Novikov Complex
for Morse forms: rationality and exponential growth properties},
e-print dg-ga/9604004 20 Apr 96.}

\def\pator{ A.V.Pajitnov,
\emph{Simple homotopy type of Novikov Complex for closed 1-forms and
Lefschetz $\zeta$-functions of the gradient flow},
e-print dg-ga/9706014 26 Jun 97.}

\newcommand{\pozniak}
{M.Pozniak, \emph{ Triangulation of Compact 
Smooth Manifolds and Morse theory},
Warwick preprint: 11/1990, November 1990}

\newcommand{\ranprepr}
{  A.Ranicki,
      \emph{
Finite domination and Novikov rings},
preprint,
 1993 }


\def\novidok{S.P.Novikov,
\emph{Mnogoznachnye funktsii i funktsionaly.
Analog teorii Morsa}, Doklady AN SSSR
{\bf 260}
(1981),
 31-35.
   English translation:
 S.P.Novikov. \emph{ Many-valued functions
 and functionals. An analogue of Morse theory},
 Sov.Math.Dokl.
{\bff 24}
(1981),
       222-226. }

\newcommand{\farber}
{M.Farber, \emph{Tochnost' neravenstv Novikova},
Funktsional'nyi analiz i ego prilozheniya
 {\bff 19}
 (1985),
     49--59
(in   Russian).
English translation:M.Farber, \emph{
 Exactness of Novikov inequalities },
Functional Analysis and Applications  {\bff
 19},
 1985.
 }

\newcommand{\novshme}
{  S.P.Novikov, I.Schmeltzer,
\emph{ Periodicheskie resheniya uravneniya tipa Kirhgofa dlya dvizheniya
tverdogo tela v zhidkosti i rasshirennaya teoriya
Lusternika-Schirel'mana-Morsa LSchM 1},
Funktsional'nyi analiz i ego prilozheniya
{\bff 15}
no. 3
 (1981).
English translation: S.P.Novikov, I.Shmel'tser, \emph{
 Periodic Solutions of
Kirchhof Equation for the free motion of a rigid body
in a fluid and the extended Theory of Lyusternik - Shnirelman -
Morse (LSM) 1                                        },
 Functional Analysis and Applications
{\bff 15} (1981),
 197 -- 207.
}

\newcommand{\novikirh}
{  S.P.Novikov, \emph{
Variatsionnye metody i periodicheskie resheniya uravneniya tipa
Kirhgofa},
Funktsional'nyi analiz i ego prilozheniya,
{\bff 15},
no. 4
 (1981).
English translation: S.P.Novikov,
\emph{ Variational Methods and Periodic Solutions
of Kirchhof type equations 2},
   Functional Analysis and applications,
{\bff 15}
no. 4
   (1981)
}

\newcommand{\noviuspe}
{ S.P.Novikov,
\emph{ Gamil'tonov formalizm i mnogoznchnyi
analog teorii Morsa}
{\bff 37}
   (1982),
 3-49.
English translation: S.P.Novikov,\emph{
The hamiltonian formalism and a
multivalued analogue of
Morse theory                          },
 Russ. Math. Surveys,
{\bff 37} (1982),
 1 -- 56.}

\def\padok{  A.Pajitnov
\emph{O tochnosti neravenstv tipa Novikova
 dlya sluchaya $\pi_1(M)=\ZZZ^m$
i klassov kogomologii v obschem polozhenii}
Doklady ANSSSR
{\bff 306},
 1989
no. 3.
English translation: A.V.
Pazhitnov, \emph{On the sharpness
of the inequalities of Novikov
type in the case
 $\pi_1(M)={\bold Z}^m$
for Morse forms whose cohomology classes are in general
position           },
 Soviet Math. Dokl.
{\bff 39}
 (1989),
no. 3.}

\newcommand{\pasbor}
{  A.Pajitnov
\emph{ O tochnosti neravenstv tipa Novikova dlya sluchaya
$\pi_1(M)=\ZZZ^m$
dlya mnogoobrazii so svobodnoi abelevoi fundamental'noi gruppoi}
(1989)
no. 11
(in   Russian).
English translation: A.V.Pazhitnov,\emph{
 On the sharpness of
 Novikov-type inequalities for
manifolds with free abelian fundamental group.},
 Math. USSR Sbornik,
{\bff 68}
 (1991),
 351 - 389.
}

\newcommand{\pastpet}
{ A.Pajitnov
\emph{
      Ratsional'nost' granichnyh operatorov v komplekse Novikova
v sluchae obschego polozheniya}
{\bff 9}, no.5 (1997), бва. 92--139
 }

         \newcommand{\tur}
{ V.Turaev
\emph{Kruchenie Raidemaistera v teorii uzlov}
{\bff 41.1}
(1986),
  119 - 182.
(in Russian).
English Translation:
V.G.Turaev,
\emph{
Reidemeister Torsion in knot theory,}
Russian Math. Surveys, 41:1 (1986), 119 - 182.}

\newcommand{\turtur}
{V.Turaev,
\emph{ Eilerovy structury, neosobye vektornye polya
i krucheniya tipa Raidemaistera}
Izv ANSSSR , 53:3
(1989),
(in Russian)
English translation:
V.G.Turaev,
\emph{ Euler structures, nonsingular vector fields
and torsions of Reidemeister type},
Math. USSR Izvestia 34:3 (1990), 627 - 662.}


\newcounter{e}[chapter]

\renewcommand{\thee}{\thechapter.\arabic{e} }
\renewcommand{\thesection}{\thechapter.\arabic{section}}

\newenvironment{prop}{\vskip0.1in     \noindent
\par\refstepcounter{e}%
   {\bf Proposition
                   \arabic{chapter}.\arabic{e}. }
\quad\it\vskip0.01in
  }{\vskip0.1in}

\newenvironment{propcit}[1]{\vskip0.1in     \noindent
\par\refstepcounter{e}%
   {\bf Proposition
                   \arabic{chapter}.\arabic{e}. {#1}}
\quad\it\vskip0.01in
  }{\vskip0.1in}

\newenvironment{theocitbeznom}[1]{\vskip0.1in     \noindent
   {\bf Theorem
                   {#1}}
\quad\it\vskip0.01in
  }{\vskip0.1in}

\newenvironment{thbeznom}{\vskip0.1in     \noindent
   {\bf Theorem
                   }
\quad\it\vskip0.01in
 }{\vskip0.1in}

\newenvironment{theo}{\vskip0.1in     \noindent
\par\refstepcounter{e}%
   {\bf Theorem
                   \arabic{chapter}.\arabic{e}. }
\quad\it\vskip0.01in
 }{\vskip0.1in}

\newenvironment{coro}{\vskip0.1in     \noindent
\par\refstepcounter{e}%
   {\bf Corollary
                   \arabic{chapter}.\arabic{e}. }
\quad\it\vskip0.01in
  }{\vskip0.1in}

\newenvironment{lemm}{\vskip0.1in     \noindent
\par\refstepcounter{e}%
   {\bf Lemma
                   \arabic{chapter}.\arabic{e}. }
\quad\it\vskip0.01in
  }{\vskip0.1in}

\newenvironment{defi}{\vskip0.1in     \noindent
\par\refstepcounter{e}%
   {\bf Definition
                   \arabic{chapter}.\arabic{e}. }
\quad\vskip0.01in
  }{\hfill$\qt$\vskip0.1in}

\newenvironment{defis}{\vskip0.1in     \noindent
\par\refstepcounter{e}%
   {\bf Definitions
                   \arabic{chapter}.\arabic{e}. }
\quad\vskip0.01in
  }{\vskip0.1in}

\newenvironment{rema}{\vskip0.1in     \noindent
\par\refstepcounter{e}%
   {\bf Remark
                   \arabic{chapter}.\arabic{e}. }
\quad\vskip0.01in
  }{\hfill$\qt$\vskip0.1in}

\renewcommand{\a}{\alpha}
\renewcommand{\b}{\beta}
\newcommand{\g}{\gamma}
\renewcommand{\d}{\delta}
\newcommand{\e}{\epsilon}
\newcommand{\ve}{\varepsilon}
\newcommand{\z}{\zeta}
\renewcommand{\t}{\theta}
\renewcommand{\l}{\lambda}
\renewcommand{\k}{\varkappa}
\newcommand{\m}{\mu}
\newcommand{\n}{\nu}
\renewcommand{\r}{\rho}
\newcommand{\vr}{\varrho}
\newcommand{\s}{\sigma}
\newcommand{\vp}{\varphi}
\renewcommand{\o}{\omega}


\newcommand{\G}{\Gamma}
\newcommand{\D}{\Delta}
\newcommand{\T}{\Theta}
\renewcommand{\L}{\Lambda}
\renewcommand{\P}{\Pi}
\newcommand{\Si}{\Sigma}
\renewcommand{\O}{\Omega}

\renewcommand{\AA}{{\mathcal A}}
\newcommand{\BB}{{\mathcal B}}
\newcommand{\CC}{{\mathcal C}}
\newcommand{\DD}{{\mathcal D}}
\newcommand{\EE}{{\mathcal E}}
\newcommand{\FF}{{\mathcal F}}
\newcommand{\GG}{{\mathcal G}}
\newcommand{\HH}{{\mathcal H}}
\newcommand{\II}{{\mathcal I}}
\newcommand{\JJ}{{\mathcal J}}
\newcommand{\KK}{{\mathcal K}}
\newcommand{\LL}{{\mathcal L}}
\newcommand{\MM}{{\mathcal M}}
\newcommand{\NN}{{\mathcal N}}
\newcommand{\OO}{{\mathcal O}}
\newcommand{\PP}{{\mathcal P}}
\newcommand{\QQ}{{\mathcal Q}}
\newcommand{\RR}{{\mathcal R}}
\renewcommand{\SS}{{\mathcal S}}
\newcommand{\TT}{{\mathcal T}}
\newcommand{\UU}{{\mathcal U}}
\newcommand{\VV}{{\mathcal V}}
\newcommand{\WW}{{\mathcal W}}
\newcommand{\XX}{{\mathcal X}}
\newcommand{\YY}{{\mathcal Y}}
\newcommand{\ZZ}{{\mathcal Z}}

\renewcommand{\aa}{{\mathbb{A}}}
\newcommand{\bb}{{\mathbb{B}}}
\newcommand{\cc}{{\mathbb{C}}}
\newcommand{\dd}{{\mathbb{D}}}
\newcommand{\ee}{{\mathbb{E}}}
\newcommand{\ff}{{\mathbb{F}}}
\renewcommand{\gg}{{\mathbb{G}}}
\newcommand{\hh}{{\mathbb{H}}}
\newcommand{\ii}{{\mathbb{I}}}
\newcommand{\jj}{{\mathbb{J}}}
\newcommand{\kk}{{\mathbb{K}}}
\renewcommand{\ll}{{\mathbb{L}}}
\newcommand{\mm}{{\mathbb{M}}}
\newcommand{\nn}{{\mathbb{N}}}
\newcommand{\oo}{{\mathbb{O}}}
\newcommand{\pp}{{\mathbb{P}}}
\newcommand{\qq}{{\mathbb{Q}}}
\newcommand{\rr}{{\mathbb{R}}}
\renewcommand{\ss}{{\mathbb{S}}}
\newcommand{\ttt}{{\mathbb{T}}}
\newcommand{\uu}{{\mathbb{U}}}
\newcommand{\vv}{{\mathbb{V}}}
\newcommand{\ww}{{\mathbb{W}}}
\newcommand{\xx}{{\mathbb{X}}}
\newcommand{\yy}{{\mathbb{Y}}}
\newcommand{\zz}{{\mathbb{Z}}}

\newcommand{\AAA}{{\mathbf{A}}}
\newcommand{\BBB}{{\mathbf{B}} }
\newcommand{\CCC}{{\mathbf{C}} }
\newcommand{\DDD}{{\mathbf{D}} }
\newcommand{\EEE}{{\mathbf{E}} }
\newcommand{\FFF}{{\mathbf{F}} }
\newcommand{\GGG}{{\mathbf{G}}}
\newcommand{\HHH}{{\mathbf{H}}}
\newcommand{\III}{{\mathbf{I}}}
\newcommand{\JJJ}{{\mathbf{J}}}
\newcommand{\KKK}{{\mathbf{K}}}
\newcommand{\LLL}{{\mathbf{L}}}
\newcommand{\MMM}{{\mathbf{M}}}
\newcommand{\NNN}{{\mathbf{N}}}
\newcommand{\OOO}{{\mathbf{O}}}
\newcommand{\PPP}{{\mathbf{P}}}
\newcommand{\QQQ}{{\mathbf{Q}}}
\newcommand{\RRR}{{\mathbf{R}}}
\newcommand{\SSS}{{\mathbf{S}}}
\newcommand{\TTT}{{\mathbf{T}}}
\newcommand{\UUU}{{\mathbf{U}}}
\newcommand{\VVV}{{\mathbf{V}}}
\newcommand{\WWW}{{\mathbf{W}}}
\newcommand{\XXX}{{\mathbf{X}}}
\newcommand{\YYY}{{\mathbf{Y}}}
\newcommand{\ZZZ}{{\mathbf{Z}}}

\newcommand{\gA}{{\mathfrak{A}}}
\newcommand{\gB}{{\mathfrak{B}}}
\newcommand{\gC}{{\mathfrak{C}}}
\newcommand{\gD}{{\mathfrak{D}}}
\newcommand{\gE}{{\mathfrak{E}}}
\newcommand{\gF}{{\mathfrak{F}}}
\newcommand{\gG}{{\mathfrak{G}}}
\newcommand{\gH}{{\mathfrak{H}}}
\newcommand{\gI}{{\mathfrak{I}}}
\newcommand{\gJ}{{\mathfrak{J}}}
\newcommand{\gK}{{\mathfrak{K}}}
\newcommand{\gL}{{\mathfrak{L}}}
\newcommand{\gM}{{\mathfrak{M}}}
\newcommand{\gN}{{\mathfrak{N}}}
\newcommand{\gO}{{\mathfrak{O}}}
\newcommand{\gP}{{\mathfrak{P}}}
\newcommand{\gQ}{{\mathfrak{Q}}}
\newcommand{\gR}{{\mathfrak{R}}}
\newcommand{\gS}{{\mathfrak{S}}}
\newcommand{\gT}{{\mathfrak{T}}}
\newcommand{\gU}{{\mathfrak{U}}}
\newcommand{\gV}{{\mathfrak{V}}}
\newcommand{\gW}{{\mathfrak{W}}}
\newcommand{\gX}{{\mathfrak{X}}}
\newcommand{\gY}{{\mathfrak{Y}}}
\newcommand{\gZ}{{\mathfrak{Z}}}

\newcommand{\gota}{{\mathfrak{a}}}
\newcommand{\gotb}{{\mathfrak{b}}}
\newcommand{\gotc}{{\mathfrak{c}}}
\newcommand{\gotd}{{\mathfrak{d}}}
\newcommand{\gote}{{\mathfrak{e}}}
\newcommand{\gotf}{{\mathfrak{f}}}
\newcommand{\gotg}{{\mathfrak{g}}}
\newcommand{\goth}{{\mathfrak{h}}}
\newcommand{\goti}{{\mathfrak{i}}}
\newcommand{\gotj}{{\mathfrak{j}}}
\newcommand{\gotk}{{\mathfrak{k}}}
\newcommand{\gotl}{{\mathfrak{l}}}
\newcommand{\gotm}{{\mathfrak{m}}}
\newcommand{\gotn}{{\mathfrak{n}}}
\newcommand{\goto}{{\mathfrak{o}}}
\newcommand{\gotp}{{\mathfrak{p}}}
\newcommand{\gotq}{{\mathfrak{q}}}
\newcommand{\gotr}{{\mathfrak{r}}}
\newcommand{\gots}{{\mathfrak{s}}}
\newcommand{\gott}{{\mathfrak{t}}}
\newcommand{\gotu}{{\mathfrak{u}}}
\newcommand{\gotv}{{\mathfrak{v}}}
\newcommand{\gotw}{{\mathfrak{w}}}
\newcommand{\gotx}{{\mathfrak{x}}}
\newcommand{\goty}{{\mathfrak{y}}}
\newcommand{\gotz}{{\mathfrak{z}}}

\newcommand{\krest}{\begin{picture}(14,14)
\put(00,04){\line(1,0){14}}
\put(00,02){\line(1,0){14}}
\put(06,-4){\line(0,1){14}}
\put(08,-4){\line(0,1){14}}
\end{picture}     }

\newcommand{\tret}{{\frac 13}}
\newcommand{\dvet}{{\frac 23}}
\newcommand{\polt}{{\frac 32}}
\newcommand{\polo}{{\frac 12}}

\renewcommand{\leq}{\leqslant}
\renewcommand{\geq}{\geqslant}

\newcommand{\bv}{B(v,\d)}

\newcommand{\ti}{\times}

\newcommand{\FR}{{\mathcal{F}}r}
\newcommand{\gt}{{\mathcal{G}}t}

\newcommand{\aand}{\quad\text{and}\quad}
\newcommand{\wwhere}{\quad\text{where}\quad}
\newcommand{\ffor}{\quad\text{for}\quad}
\newcommand{\for}{~\text{for}~}
\newcommand{\iif}{\quad\text{if}\quad}
\newcommand{\iiif}{~\text{if}~}

\newcommand{\en}{enumerate}

\newcommand{\Prf}{{\it Proof.\quad}}

\newcommand{\Wkr}{W^{\circ}}

\newcommand{\Ker}{\text{\rm Ker }}
\newcommand{\ind}{\text{\rm ind}}
\newcommand{\rk}{\text{\rm rk }}
\renewcommand{\Im}{\text{\rm Im }}
\newcommand{\supp}{\text{\rm supp }}
\newcommand{\Int}{\text{\rm Int }}
\newcommand{\grad}{\text{\rm grad }}

\newcommand{\nr}{\Vert}
\newcommand{\smo}{C^{\infty}}

\newcommand{\fpr}[2]{{#1}^{-1}({#2})}
\newcommand{\sdvg}[3]{\widehat{#1}_{[{#2},{#3}]}}
\newcommand{\disc}[3]{B^{({#1})}_{#2}({#3})}
\newcommand{\Disc}[3]{D^{({#1})}_{#2}({#3})}
\newcommand{\desc}[3]{B_{#1}(\leq{#2},{#3})}
\newcommand{\Desc}[3]{D_{#1}(\leq{#2},{#3})}
\newcommand{\komp}[3]{{\bold K}({#1})^{({#2})}({#3})}
\newcommand{\Komp}[3]{\big({\bold K}({#1})\big)^{({#2})}({#3})}

\newcommand{\ran}{\{(A_\lambda , B_\lambda)\}_{\lambda\in\Lambda}}
\newcommand{\rran}{\{(A_\lambda^{(s)},
 B_\lambda^{(s)}  )\}_{\lambda\in\Lambda, 0\leq s\leq n }}
\newcommand{\brs}{\rran}
\newcommand{\rans}{\{(A_\sigma , B_\sigma)\}_{\sigma\in\Sigma}}

\newcommand{\fmin}{F^{-1}}
\newcommand{\vh}{\widehat{(-v)}}

\newcommand{\chart}{\Phi_p:U_p\to B^n(0,r_p)}
\newcommand{\atlas}{\{\Phi_p:U_p\to B^n(0,r_p)\}_{p\in S(f)}}
\newcommand{\flow}{{\VV}=(f,v, \UU)}

\newcommand{\Rn}{\bold R^n}
\newcommand{\Rk}{\bold R^k}

\newcommand{\fcob}{f:W\to[a,b]}

\newcommand{\phicob}{\phi:W\to[a,b]}

\newcommand{\vphi}{\varphi}

\newcommand{\crr}{p\in S(f)}
\newcommand{\nrv}{\Vert v \Vert}
\newcommand{\nrw}{\Vert w \Vert}
\newcommand{\nru}{\Vert u \Vert}

\newcommand{\obb}{\cup_{p\in S(f)} U_p}
\newcommand{\proob}{\Phi_p^{-1}(B^n(0,}

\newcommand{\stind}[3]{{#1}^{\displaystyle\rightsquigarrow}_
{[{#2},{#3}]}}

\newcommand{\indl}[1]{{\scriptstyle{\text{\rm ind}\leqslant {#1}~}}}
\newcommand{\inde}[1]{{\scriptstyle{\text{\rm ind}      =   {#1}~}}}
\newcommand{\indg}[1]{{\scriptstyle{\text{\rm ind}\geqslant {#1}~}}}

\newcommand{\obbi}{\cup_{p\in S_i(f)}}
\newcommand{\vem}{\text{Vectt}(M)}
\newcommand{\pr}{\partial}

\newcommand{\id}{\text{id}}

\newcommand{\lau}[1]{{\xleftarrow{#1}}}
\newcommand{\rau}[1]{{\xrightarrow{#1}}}
\newcommand{\rad}[1]{ {\xrightarrow[#1]{}} }

\newcommand{\xit}{\tilde\xi_t}

\newcommand{\st}[1]{\overset{\rightsquigarrow}{#1}}
\newcommand{\bst}[1]{\overset{\displaystyle\rightsquigarrow}
\to{\boldkey{#1}}}

\newcommand{\stexp}[1]{{#1}^{\rightsquigarrow}}
\newcommand{\bstexp}[1]{{#1}^{\displaystyle\rightsquigarrow}}

\newcommand{\bstind}[3]{{\boldkey{#1}}^{\displaystyle\rightsquigarrow}_
{[{#2},{#3}]}}
\newcommand{\bminstind}[3]{\stind{({\boldkey{-}\boldkey{#1}})}{#2}{#3}}

\newcommand{\Tb}{\text{ \rm Tb}}

\newcommand{\VODIN}{V_{1/3}}
\newcommand{\VDVA}{V_{2/3}}
\newcommand{\VM}{V_{1/2}}
\newcommand{\ddd}{\cup_{p\in S_i(F_1)} D_p(u)}
\newcommand{\dddmin}{\cup_{p\in S_i(F_1)} D_p(-u)}
\newcommand{\where}{\quad\text{\rm where}\quad}

\newcommand{\kr}[1]{{#1}^{\circ}}

\newcommand{\vew}{\text{\rm Vect}^1 (W,\bot)}

\newcommand{\veww}{\text{\rm Vect}^1 (W)}

\newcommand{\Imm}{\text{\rm Im}}
\newcommand{\hrrr}{\text{\rm Vect}^1(M)}
\newcommand{\vemm}{\text{\rm Vect}^1_0(M)}
\newcommand{\ver}{\text{\rm Vect}^1(\RRR^ n)}
\newcommand{\verr}{\text{\rm Vect}^1_0(\RRR^ n)}

\newcommand{\mods}{\vert s(t)\vert}
\newcommand{\exd}{e^{2(D+\alpha)t}}
\newcommand{\exmin}{e^{-2(D+\alpha)t}}

\newcommand{\intt}{[-\theta,\theta]}

\newcommand{\ffmin}{f^{-1}}
\newcommand{\vvol}{\text{\rm vol}}
\newcommand{\Mat}{\text\rm Mat}
\newcommand{\Tub}{\text{\rm Tub}}
\newcommand{\vxi}{v\langle\vec\xi\rangle}
\newcommand{
\ST}
{\stexp}

\newcommand{\qt}{\hfill\triangle}
\newcommand{\qs}{\hfill\square}

\newcommand{\Vect}{\text{\rm Vect}}

\newcommand{\pa}{\vskip0.1in}

\newcommand{\wi}{\widetilde}

\newcommand{\ove}{\overline}
\newcommand{\unde}{\underline}
\newcommand{\ptf}{\pitchfork}

\renewcommand{\(}{\big(}
\renewcommand{\)}{\big)}

\newcommand{\grd}{{\text{\rm grd}}}

\newcommand{\RA}{\Rightarrow}
\newcommand{\LA}{\Leftarrow}

\newcommand{\emp}{\emptyset}
\newcommand{\wh}{\widehat}

\newcommand{\GC}{\GG\CC}
\newcommand{\GCT}{\GG\CC\TT}
\newcommand{\GT}{\GG\TT}

\newcommand{\GA}{\GG\AA}
\newcommand{\GRP}{\GG\RR\PP}

\newcommand{\GgC}{\GG\gC}
\newcommand{\GgCC}{\GG\gC\CC}

\newcommand{\GgCT}{\GG\gC\TT}

\newcommand{\GgCY}{\GG\gC\YY}
\newcommand{\GgCYT}{\GG\gC\YY\TT}

\newcommand{\GCCT}{\GG\gC\CC\TT}
\newcommand{\GCC}{\GG\gC\CC}

\newcommand{\stv}{\stexp {(-v)}}
\newcommand{\stu}{\stexp {(-u)}}

\newcommand{\strv}[2]{\stind {(-v)}{#1}{#2}}
\newcommand{\strw}[2]{\stind {(-w)}{#1}{#2}}
\newcommand{\stru}[2]{\stind {(-u)}{#1}{#2}}

\newcommand{\stvv}[2]{\stind {v}{#1}{#2}}
\newcommand{\stuu}[2]{\stind {u}{#1}{#2}}
\newcommand{\stww}[2]{\stind {w}{#1}{#2}}

\newcommand{\mx}{\mbox}

\newcommand{\vbsm}{V_b^{\{\leq s-1\}}    }
\newcommand{\vasm}{V_a^{\{\leq s-1\}}    }
\newcommand{\vbs}{V_b^{\{\leq s\}}    }
\newcommand{\vas}{V_a^{\{\leq s\}}    }

\newcommand{\Vbsm}{V_b^{[\leq s-1]}(\d)    }
\newcommand{\Vasm}{V_a^{[\leq s-1]}(\d)    }
\newcommand{\Vbs}{V_b^{[\leq s]}(\d)    }
\newcommand{\Vas}{V_a^{[\leq s]}(\d)    }

\newcommand{\vass}{V_{a_{s+1}}}

\newcommand{\vovo}{\stexp {(-v0)}}
\newcommand{\vov}{\stexp {(-v1)}}

\newcommand{\vbkm}{V_b^{\{\leq k-1\}}    }
\newcommand{\vakm}{V_a^{\{\leq k-1\}}    }
\newcommand{\vbk}{V_b^{\{\leq k\}}    }
\newcommand{\vak}{V_a^{\{\leq k\}}    }

\newcommand{\Vbkm}{V_b^{[\leq k-1]}(\d)    }
\newcommand{\Vakm}{V_a^{[\leq k-1]}(\d)    }
\newcommand{\Vbk}{V_b^{[\leq k]}(\d)    }
\newcommand{\Vak}{V_a^{[\leq s]}(\d)    }

\newcommand{\lc}{\lceil}
\newcommand{\rc}{\rceil}

\newcommand{\sps}{\supset}

\newcommand{\bere}{\begin{rema}}
\newcommand{\bede}{\begin{defi}}

\renewcommand{\beth}{\begin{theo}}
\newcommand{\bele}{\begin{lemm}}
\newcommand{\bepr}{\begin{prop}}
\newcommand{\beeq}{\begin{equation}}
\newcommand{\bega}{\begin{gather}}
\newcommand{\been}{\begin{enumerate}}

\newcommand{\beal}{\begin{aligned}}

\newcommand{\enre}{\end{rema}}

\newcommand{\enpr}{\end{prop}}
\newcommand{\enth}{\end{theo}}
\newcommand{\enle}{\end{lemm}}
\newcommand{\enen}{\end{enumerate}}
\newcommand{\enga}{\end{gather}}
\newcommand{\eneq}{\end{equation}}
\newcommand{\enal}{\end{aligned}}

\newcommand{\subs}{\subsection}
\newcommand{\ity}{\infty}

\newcommand{\lb}{\label}

\newcommand{\Vv}{{\boldsymbol{v}}}

\newcommand{\Lxi}{{\wh \L}_\xi}
\newcommand{\lL}{\wh{\wh L}}

\newcommand{\sub}{ Subsection}

\newcommand{\tens}[1]{\underset{#1}{\otimes}}

\newcommand{\LLxi}{\wi \L_\xi}

\newcommand{\hV}{\wh V}
\newcommand{\hHH}{\wh \HH}

\newcommand{\gama}[2]{\g({#1}, \tau_a({#2},{#1}); w )}

\newcommand{\gam}[2]{\g({#1}, \tau_0({#2},{#1}); w )}
\newcommand{\ga}[2]{\g({#1}, \tau({#2},{#1}); w )}

\newcommand{\ifff}{if and only if}

\renewcommand{\th}{therefore}
\newcommand{\ata}{almost~ transversality~ assumption}
\newcommand{\gr}{gradient}
\newcommand{\Mf}{Morse function}
\newcommand{\iis}{it is sufficient}
\newcommand{\sut}{~such~that~}
\newcommand{\sufsm}{~sufficiently~ small~}
\newcommand{\wrt}{~with respect to}
\newcommand{\sm}{\setminus}
\newcommand{\ems}{\emptyset}
\newcommand{\sbs}{\subset}
\newcommand{\ho}{homomorphism}
\newcommand{\ma}{manifold}
\newcommand{\nei}{neighborhood}
\newcommand{\dfm}{diffeomorphism}

\newcommand{\vf}{vector field}

\newcommand{\vfs}{vector fields}

\newcommand{\fe}{for every}

\renewcommand{\top}{topology}

\newcommand{\tr}{~trajectory }

\newcommand{\grs}{~gradients}
\newcommand{\trs}{~trajectories}

\newcommand{ \co}{~cobordism}
\newcommand{
\sma}{submanifold}
\newcommand{
\hos}{~homomorphisms}
\newcommand{
\Th}{~Therefore}

\newcommand{
\tthen}{\text \rm ~then}

\newcommand{
\wwe}{\text \rm ~we  }
\newcommand{
\hhave}{\text \rm ~have}
\newcommand{
\eevery}{\text \rm ~every}

\newcommand{\noconf}{~there~is~no~possibility~of~confusion}

\newcommand{\ATA}{Almost~ Transversality~ Condition}
\newcommand{\cob}{~cobordism}

\newcommand{\hot}{homotopy}

\newcommand{\TA}{Transversality Condition}

\newcommand{\hog}{homology}

\newcommand{\cog}{cohomology}

\newcommand{\wat}{ We shall assume that}

\newcommand{\sclv}{sufficiently close to $v$ in $C^0$-topology}

\newcommand{\cf}{continuous function }

\newcommand{\eg}{exponential growth}

\newcommand{\nics}{Novikov incidence coefficients}
\newcommand{\nic}{Novikov incidence coefficient}

\newcommand{\negc}{Novikov exponential growth conjecture}

\newcommand{\mc}{Morse Complex   }

\newcommand{\mas}{manifolds   }

\newcommand{\nc}{Novikov Complex   }

\newcommand{\glvf}{gradient-like vector field}

\newcommand{\glvfs}{gradient-like vector fields}

\newcommand{\fg}{finitely generated   }

\newcommand{\ta}{transversality condition}

\newcommand{\mnp}{Morse-Novikov pair}

\newcommand{\rp}{rationality property}

\newcommand{\egp}{exponential growth property}

\newcommand{\mi}[3]{{#1}^{-1}\([{#2},{#3}]\)}

\newcommand{\fii}[2]{\mi {\phi}{a_{#1}}{a_{#2}} }

\newcommand{\fifi}[2]{\mi {\phi}{#1}{#2} }

\newcommand{\pf}[2]{\mi {\phi_1}{\a_{#1}}{\a_{#2}} }

\newcommand{\mf}[2]{\mi {\phi_0}{\b_{#1}}{\b_{#2}}}

\newcommand{\wa}[2]{ W_{[a_{#1}, a_{#2}]}}

\newcommand{\waa}[1]{ W_{[a, a_{#1}]}}

\newcommand{\Wa}[2]{ W_{[{#1}, {#2}]}}

\newcommand{\WS}[1]{ W^{\{\leq {#1}\}}}

\newcommand{\ws}{\WS {s}}

\newcommand{\wsm}{\WS {s-1}}

\newcommand{\wsmm}{\WS {s-2}}

\newcommand{\wk}{\WS {k}}

\newcommand{\wkm}{\WS {k-1}}

\newcommand{\wkmm}{\WS {k-2}}

\newcommand{\wsn}{ W^{[\leq s]}(\nu)}

\newcommand{\wsmn}{ W^{[\leq s-1]}(\nu)}

\newcommand{\wsk}{ W^{[\leq k]}(\nu)}

\newcommand{\pws}{(\pr_1 W)^{\{\leq s\}}}

\newcommand{\hpws}{(\wh{\pr_1 W})^{\{\leq s\}}}

\newcommand{\hpwsm}{(\wh{\pr_1 W})^{\{\leq {s-1}\}}}

\newcommand{\mws}{(\pr_0 W)^{\{\leq s\}}}

\newcommand{\pwk}{(\pr_1 W)^{\{\leq k\}}}

\newcommand{\mwk}{(\pr_0 W)^{\{\leq k\}}}

\newcommand{\pwkm}{(\pr_1 W)^{\{\leq k-1\}}}

\newcommand{\pwsm}{(\pr_1 W)^{\{\leq s-1\}}}

\newcommand{\pwkmm}{(\pr_1 W)^{\{\leq k-2\}}}

\newcommand{\mwkmm}{(\pr_0 W)^{\{\leq k-2\}}}

\newcommand{\mwkm}{(\pr_0 W)^{\{\leq k-1\}}}

\newcommand{\mwsm}{(\pr_0 W)^{\{\leq s-1\}}}

\newcommand{\mwkp}{(\pr_0 W)^{\{\leq k+1\}}}

\newcommand{\wmok}{W^{\langle k\rangle}}

\newcommand{\wmokp}{W^{\langle k+1\rangle}}

\newcommand{\wmokm}{W^{\langle k-1\rangle}}
\newcommand{\wmokmm}{W^{\langle k-2\rangle}}

\newcommand{\wmos}{W^{\langle s\rangle}}

\newcommand{\wmosm}{W^{\langle s-1\rangle}}

\newcommand{\wmoo}{W^{\langle 0\rangle}}

\newcommand{\wwk}{\( \wmok , \wmokm \)}

\newcommand{\wwkp}{\( \wmokp , \wmok \)}

\newcommand{\wwkm}{\( \wmokm , \wmokmm \)}

\newcommand{\wws}{\bigg( \wmos , \wmosm \bigg)}

\newcommand{\wasn}{W^{\lc s\rc}(\nu)}

\newcommand{\wakn}{W^{\lc k\rc}(\nu)}

\newcommand{\dqr}{\pr_- Q_r}

\newcommand{\ds}{\pr_s}

\newcommand{\dsm}{\pr_{s-1}}

\newcommand{\dow}{\pr_0 W}

\newcommand{\daw}{\pr_1 W}

\newcommand{\hdaw}{\wh{\pr_1 W}}

\newcommand{\hwm}{H_*\( \wmok, \wmokm\)}

\newcommand{\hkwm}{H_k\( \wmok, \wmokm\)}

\newcommand{\hwmp}{H_*\( \wmokp, \wmok\)}

\newcommand{\hkwmm}{H_k\( \wmokm, \wmokmm\)}

\newcommand{\hkwmp}{H_k\( \wmokp, \wmok\)}

\newcommand{\yz}{Y_k(v)\cup Z_k(v)}

\newcommand{\Gama}{{\nazad{ \Gamma}}}
\newcommand{\ug}[1]{\llcorner {#1} \lrcorner}
\newcommand{\npqv}{n(\bar p, \bar q; v)}
\newcommand{\fms}{f:M\to S^1   }

\newcommand{\nkpqv}{n_k(\bar p, \bar q; v)}

\newcommand{\GLT}{\GG lt}\title{$C^0$-generic properties of boundary operators in Novikov Complex}
\date{\today}
\author{A.Pajitnov}
\address{Universit\'e de Nantes, 
Facult\'e des Sciences et Techniques, 
2, rue de la Houssini\`ere, 44072, 
Nantes, Cedex}
\email{pajitnov@math.univ-nantes.fr}
\maketitle
\newcommand{\spa}{\hspace{0,5cm}}

\newcommand{\sspa}{\hspace{1cm}}
\newcommand{\df}{\dotfill}

\centerline{\bf Table of Contents}
\pa
\pa

CHAPTER 1. Introduction\df \pageref{ch:i}
\pa
\spa \ref{s:incneg} Introduction to Novikov Complex and
Novikov Exponential Growth Conjecture
 \dotfill
\pageref{s:incneg}

\sspa \ref{su:mcmf} Morse Complex of a Morse function \dotfill
\pageref{su:mcmf}

\sspa \ref{su:nc} Novikov Complex of a Morse map to the circle
\df\pageref{su:nc}

\sspa \ref{su:mnt} Morse-Novikov theory \df\pageref{su:mnt}

\sspa \ref{su:nc} Novikov Exponential Growth Conjecture
 \df\pageref{su:negc}

\spa \ref{s:cgrcp} $C^0$-generic rationality of \nics.
The contents of the paper
 \dotfill
\pageref{s:cgrcp}

\sspa \ref{su:cpp}
The contents of the paper.
\df\pageref{su:cpp}

\sspa \ref{su:pmm}
Preliminaries on Morse maps $M\to S^1$
\df\pageref{su:pmm}

\sspa \ref{su:fsp}
First step of the proof
\df\pageref{su:fsp}

\sspa \ref{su:opmt}
Outline of the proof of the main theorem
\df\pageref{su:opmt}

\sspa \ref{su:fr}
Further remarks
\df\pageref{su:fr}

\pa
\spa \ref{s:termi}
Terminology
\df\pageref{s:termi}

\sspa \ref{su:tc}
Terminological conventions
\df\pageref{su:tc}

\sspa \ref{su:diffnot}
Some differences between the notation of this paper and that of
\cite{pasur}, \cite{paepr}
\df\pageref{su:diffnot}

CHAPTER 2. Morse functions and their gradients\df \pageref{ch:mftg}
\pa
\spa \ref{s:bdc} Basic definitions and constructions \dotfill
\pageref{s:bdc}

\sspa \ref{su:bd} Basic definitions \df\pageref{su:bd}

\sspa \ref{su:lsmfg} Local structure of Morse function and its
gradient in a \nei~
of a critical point
\df\pageref{su:lsmfg}

\sspa  \ref{su:gst} Global structure of $v$-trajectories
\df\pageref{su:gst}

\sspa   \ref{su:tddcp}  Thickenings of descending discs and
compactness properties
\df\pageref{su:tddcp}

\sspa  \ref{su:tcrl} Transversality conditions and Rearrangement lemma
\df\pageref{su:tcrl}

\sspa  \ref{su:ahc} Adding a horizontal component to a vector field
nearby the boundary
\df\pageref{su:ahc}

\pa\pa

\spa \ref{s:ssags} $s$-submanifolds and $gfn$-systems
\df\pageref{s:ssags}

\sspa \ref{su:ssubm}  $s$-submanifolds
\df\pageref{su:ssubm}

\sspa \ref{su:gfsn} Good fundamental systems of neighborhoods
\df\pageref{su:gfsn}
\pa
\spa \ref{s:hlf} Handle-like filtrations
\df\pageref{s:hlf}

\sspa  \ref{su:sepdef} $\d$-separated \gr s: definition
\df\pageref{su:sepdef}

\sspa  \ref{su:tfo} Definition of handle-like filtrations
\df\pageref{su:tfo}

\sspa  \ref{su:hpof} Homotopical properties of the filtrations
\df\pageref{su:hpof}
\pa
\spa \ref{s:mc}  Morse complex
\df\pageref{s:mc}

\sspa \ref{su:dmc} Definition of the Morse complex
\df\pageref{su:dmc}

\sspa \ref{su:pic} Properties of incidence coefficients
\df\pageref{su:pic}
\pa
\spa  \ref{s:gdm} Gradient descent map
\df\pageref{s:gdm}

\sspa \ref{su:gpt} General properties of tracks
\df\pageref{su:gpt}

\sspa \ref{su:trssb}  Tracks of $s$-submanifolds
\df\pageref{su:trssb}

\sspa \ref{su:totss} Tracks of $ts$-submanifolds
\df\pageref{su:totss}

\pa
\spa \ref{s:pcic} Properties of continuity of integral curves
of vector fields with respect to $C^0$ small perturbations of
of vector fields
\df\pageref{s:pcic}

\sspa  \ref{su:mwb} Manifolds without boundary
\df\pageref{su:mwb}

\sspa  \ref{su:mb} Manifolds with boundary
\df\pageref{su:mb}

\sspa\ref{su:gmf} Gradients of Morse functions
\df\pageref{su:gmf}

\pa\pa CHAPTER 3. Quick flows\df\pageref{ch:qf}

\spa \ref{s:isr} Introduction and statement of results
\df\pageref{s:isr}

\sspa \ref{su:tergf} Theorem on the existence of rapid gradient flows
\df\pageref{su:tergf}

\sspa \ref{su:t} Terminology.
\df\pageref{su:t}

\sspa \ref{su:tlsg} Two lemmas on standard gradients in $\RRR^n$.
\df\pageref{su:tlsg}

\sspa \ref{su:teqf} Theorems on the existence of quick flows.
\df\pageref{su:teqf}
\pa
\spa \ref{s:sc} S-construction
\df\pageref{s:sc}

\sspa  \ref{su:aux} Auxiliary constructions and the choice of $\nu_0$.
\df\pageref{su:aux}

\sspa \ref{su:mfg} Morse function $F$ and its gradient $u$.
\df\pageref{su:mfg}

\sspa  \ref{su:pfg} Properties of $F$ and $u$.
\df\pageref{su:pfg}

\pa
\spa \ref{s:pteqf} Proof of the theorem on the existence of quick flows.
\df\pageref{s:pteqf}

\sspa \ref{su:pteqfint} Introduction
\df\pageref{su:pteqfint}

\sspa  \ref{su:fwcp} Functions without critical points.
\df\pageref{su:fwcp}

\sspa  \ref{su:fwocp} Functions with one critical point.
\df\pageref{su:fwocp}

\pa
\spa\ref{s:fr} Final remark\df\pageref{s:fr}

\pa
 CHAPTER 4. Regularizing the gradient descent map
\df\pageref{ch:rgdm}

\spa \ref{s:cc}  Condition $(\gC)$
\df\pageref{s:cc}

\sspa  \ref{su:cc} Condition $(\gC)$:
the definition and the statement of the density theorem
\df\pageref{su:cc}

\sspa  \ref{su:cogct} $C^0$ openness of $(\gC)$
\df\pageref{su:cogct}

\sspa  \ref{su:cdi} $C^0$ density of $(\gC)$:
 The idea of the proof
\df\pageref{su:cdi}

\sspa  \ref{su:cdt} $C^0$ density of $(\gC)$:
 terminology
\df\pageref{su:cdt}

\sspa  \ref{su:cdpgg} $C^0$ density of $(\gC)$:
 perturbation of a given gradient
\df\pageref{su:cdpgg}

\sspa  \ref{su:cyc} Cyclic cobordisms and condition $(\gC\YY)$
\df\pageref{su:cyc}

\spa \ref{s:hgd} Homological gradient descent
\df\pageref{s:hgd}

\sspa  \ref{su:midp}
 Main idea of the proof and the definition of $\HH_s(-v)$.
\df\pageref{su:midp}

\sspa  \ref{su:sykxk} Sets $T_k(v), t_k(v)$.
\df\pageref{su:sykxk}

\sspa  \ref{su:poh} Properties of $\HH_s$.
\df\pageref{su:poh}

\sspa  \ref{su:couple} A pairing between $C_s(u_1)$ and $C_s(v)$.
\df\pageref{su:couple}

\sspa  \ref{su:eqhgd} Homological gradient descent: equivariant version
\df\pageref{su:eqhgd}

\pa
 CHAPTER 5. $C^0$-generic rationality of the boundary
operators in Novikov Complex
\df\pageref{ch:gpbon}

\spa \ref{s:bt}
Basic terminology
\df\pageref{s:bt}

\sspa  \ref{su:nrbt} Novikov rings
\df\pageref{su:nrbt}

\sspa  \ref{su:mmbt} Morse maps $M\to S^1$
\df\pageref{su:mmbt}

\sspa  \ref{su:ncd}  Novikov Complex: the definition
\df\pageref{su:ncd}

\spa \ref{s:cccgr}
Condition $(\gC\CC)$ and $C^0$ generic rationality
of Novikov incidence coefficients with values in $\ZZZ((t))$.
\df\pageref{s:cccgr}

\spa \ref{s:icvc}
Incidence coefficients with values in completions of group rings
\df\pageref{s:icvc}

\sspa  \ref{su:icdef} Definitions
\df\pageref{su:icdef}

\sspa  \ref{su:ncalg} Some non-commutative algebra
\df\pageref{su:ncalg}

\sspa  \ref{su:eeg} Elements
of exponential growth in Novikov ring
\df\pageref{su:eeg}

\sspa  \ref{su:rnic} Rationality of \nics
\df\pageref{su:rnic}

\pa
FIGURES
\df\pageref{ch:fig}

\pa
 REFERENCES
\df\pageref{refer}

\newpage

\chapter{Introduction}
\lb{ch:i}

\section[Novikov Complex]{Introduction to Novikov Complex and Novikov 
exponential growth conjecture }
\lb{s:incneg}
\pa
\subsection{ \mc of a \Mf}
\lb{su:mcmf}

We  start with a recollection of one of the
basic constructions of the Morse theory -- Morse complex.
Let $M$ be a closed \ma, and $g:M\to\RRR$
be a \Mf.
Let $S(g)$ be the set of critical points of $g$, and $S_k(g)$
be the set of critical points of $g$ of index $k$.
Consider the free abelian group
$C_k(g)$, freely generated
by
$S_k(g)$.
We shall introduce \hos~
$\pr_k: C_k(g)\to C_{k-1}(g)$
which will endow $C_*(g)$
with the structure of a chain complex, called
{\it Morse complex}.
The \hos~depend on some additional choices. Namely one must choose:
\been\item
A \glvf~ $v$ for $g$, satisfying \TA.\footnote{
See \cite{milnhcob}, Def 3.1 for the definition of \glvf.}
\item
Orientations of stable \mas~of zeros of $v$
\enen

Having made these choices, we proceed as follows.
Let $p,q\in S(f), \ind p=\ind q+1$.
Let $L(p,q)$ be the set of orbits of $v$ joining $p$ with $q$.
The \ta~ for $v$ implies that $L(p,q)$ is finite.
The choice of orientation of stable \mas~ allows to associate to each
$\g\in L(p,q)$ a sign
$\ve(\g)\in\{-1, +1\}$.
Set
$n(p,q)=\sum_{\g\in L(p,q)}\ve(\g)$.
The integer $n(p,q)$ is called
{\it incidence coefficient}
corresponding to $p,q$.
Now define the \ho~$\pr_k$
on the generators $p\in S_k(f)$ by the following formula
$$
\pr_k(p)=\sum_{r\in S_{k-1}(f)}n(p,r) r
$$
One can prove that
$\pr_k\circ\pr_{k+1}=0$
and that the \hog~of the resulting complex is isomorphic to
$H_*(M)$.
(See for example \cite{patou}, Appendix, for the details.)

\subsection{Novikov Complex of a Morse map to the circle}
\lb{su:nc}

In the early 80s S.P.Novikov generalized the
construction of Morse complex to the
case of Morse maps $\fms$.
We shall give a construction in details a bit later
(see \sub~ \ref{su:mmbt}), now we just sum up the properties
of the resulting chain complex.

So let $\fms$
be a Morse map, where $M$ is a closed connected \ma.
We assume that $f$ is not homotopic to zero.
As before, let $S(f)$ be the set of critical points of $f$, and
$S_k(f)$ be the set of critical points of $f$
of index $k$.
As the classical Morse complex, the \nc is a free complex generated in
dimension
$k$ by the set
$S_k(f)$, but the base ring here is much larger.
Namely,
consider the ring of all Laurent power series in one variable
 with integral coefficients and
finite negative part, that is
\begin{multline*}
\ZZZ((t))=
\{\l=\sum a_kt^k\mid  \forall k a_k\in\ZZZ 
\\
\mbox{ and } \exists
N=N(\l)\mbox{ \sut~ } a_k=0\mbox{  if  } k<N\}
\end{multline*}

Let $\CC:\bar M\to M$
be the infinite connected cyclic cover,
\sut~ $f\circ \CC$  is homotopic to zero
(note that $\CC$ is uniquely determined by this last property).
To define \nc~one must make some additional choices.
Namely, we must choose:
\been\item
A \glvf~ $v$ for $f$, satisfying \TA
(
the definition of \glvf~for the maps to 
the circle copies that for Morse functions).
\item
Orientation of stable \mas~ of critical points of $f$.
\item
For every $p\in S(f)$ a lifting $\bar p$ of $p$ to $\bar M$
(that is, a point $\bar p\in \bar M$, satisfying $\CC(\bar p)=p$).
\enen

Having made these choices, one applies the procedure,
exposed in \sub~ \ref{su:mcmf}
(with corresponding modifications, see\sub~ \ref{su:ncd} for details)
and gets for every pair $p,q$ of critical points of
$f$, satisfying $\ind p=\ind q+1$
an element
$\npqv\in\ZZZ((t))$, called
{\it \nic}
of $p,q$.
(This coefficient depends on the choice
of liftings $\bar p, \bar q\in \bar M$, see the formula \ref{f:chbase}).
We shall also write
$\npqv=\sum_i\nkpqv t^k$.

One can prove that $\pr_{k+1}\circ\pr_k=0$
for every
$k$. Thus we obtain a chain complex
$C_*(f,v)$ \sut~$C_p(f,v)$
is a free  $\ZZZ((t))$-complex
freely generated by $S_p(f)$.
One can prove that
$$H_*(C_*(f,v))\approx H_*(\bar M)\tens{\ZZZ[t, t^{-1}]}\ZZZ((t))$$
(See \cite{patou} for details.)

\bere\lb{re:glf}
In the present paper we introduce a notion
of $f$-gradient which is a generalization
of the notion of
\glvf~for $f$
(definition \ref{defgrad}).
The gradient of $f$ \wrt~an arbitrary Riemannian metric is an
$f$-gradient,
as well as a \glvf~for $f$.
One defines \nic~
and \nc for $f$-gradients in a way similar to the exposed above
(see \sub~\ref{su:ncd}).
We shall use this in the following sections of the Introduction.
\enre

\subsection{Morse-Novikov theory}
\lb{su:mnt}
The \nc was introduced by S.P.Novikov in \cite{novidok}.
In this paper the analogs of the 
classical Morse inequalities were obtained.
These inequalities (called {\it Novikov inequalities})
were then studied in a series of papers\footnote{ We make
no attempt to give here a complete survey.}.

M.~Farber
\cite{farber} proved that these inequalities are sharp 
in the case where $\pi_1
M=\ZZZ$, $\dim M\geq 6$. J.-C.~Sikorav [25] applied these inequalities
to the theory of Lagrangian intersections. In \cite{padok}, \cite{pasbor}
 the author
obtained some results concerning the sharpness of these inequalities for
$\pi_1 M=\ZZZ^m$, $\dim M\geq 6$.
The study of the Novikov
complex itself advanced slower. In \cite{sikothese} yet another
version of the Novikov complex, defined over a completion of the group
ring of $\pi_1(M)$, was given. In \cite{paodense}, \cite{patou}
 I proved that the 
simple homotopy 
type
of this chain complex equals the simple homotopy type of the
completed chain simplicial complex of the universal covering of $M$. In
\cite{paura}, \cite{pasur}
 I proved a theorem on realization of the abstract chain
complexes as the Novikov complexes of Morse maps $M\to S^1$. This theorem
implies the results of \cite{farber} and of \cite{pasbor}.
 The algebraic result by Ranicki
\cite{rani} shows that the main theorem of \cite{paura}, \cite{pasur}
 implies also the
Farrell fibration theorem \cite{farrell}.
In \cite{latour} 
F.Latour has developed an approach to Morse-Novikov theory, 
based on the Cerf's method \cite{cerf}.

\subs{Novikov exponential growth conjecture}
\lb{su:negc}

A priori the \nics ~
$\npqv$ are elements of $\ZZZ((t))$, that is, power series of the form
$\sum_{i=-N}^{\infty} \nkpqv t^k$, where $\nkpqv $
are integers.
From the outset S.P.Novikov conjectured that these powers 
series had some nice 
analytic properties. Namely, he conjectured that
the positive part of the power series 
$\npqv$
converges in a disc of a non-zero radius with the center $0$
(equivalently: the coefficients
$\nkpqv$ grow at most exponentially with $k$).
During several years the conjecture has been discussed 
and adjusted. Now there are several versions of the conjecture. To 
discuss these versions it is convenient to introduce a definition.

\bede\lb{d:egdef}
\been\item
Let $f:M\to S^1$
be a Morse map, $v$ be an $f$-gradient,
satisfying \TA.
The pair $(f,v)$ will be called {\it Morse-Novikov pair}.
\item
Let $(f,v)$ be a  {\it Morse-Novikov pair}.
We say that the {\it \egp}
holds for the pair $(f,v)$ if for every $p,q\in S(f)$
with $\ind p=\ind q+1$ there are
$C,D>0$
\sut~\fe~$k\in\ZZZ$:
$$
|\nkpqv|\leq Ce^{kD}
$$
\enen
\end{defi}

The initial version of \negc~
said that the \egp~
holds \fe~  Morse-Novikov pair $(f,v)$.

As far as I know, the first published version
of the conjecture
appeared in the paper of V.~I.~Arnold 
``Dynamics of Intersections'', 1989. Arnold wrote
in \cite{Arno}, p.83:
\begin{quote}
``The author is indebted to S.~P.~Novikov
who has communicated the following conjecture,
which was
the starting point of the present paper. Let
$p:\widetilde M\to M$ be a covering of a compact
manifold $M$
with fiber
$\ZZZ^n$, and let $\a$ be a closed 1-form on $M$
such that $p^*\a=df$, where $f:M\to \RRR$ is a Morse
function.
The Novikov conjecture states that, ``generically''
the number of the trajectories of the vector field
$-\grad f$ on $M$ starting at a critical
point
$x$ of the function $f$ of index $k$ and
connecting it
with the critical points $y$ having index $k-1$
and satisfying $f(y)\geq f(x)-n$, grows in $n$
more slowly than some exponential, $e^{an}$.''
\end{quote}

Thus the version of
\negc, cited in \cite{Arno}, says in particular, that
\egp~holds for every generic pair
$(f,v)$, where $v$ is a 
{\it riemannian gradient} of $f$.
Note that this version is stated in somewhat more general
 setting than our present one:
it deals with Morse 1-forms, and not only with Morse maps to the circle.
Note also that Arnold has conjectured something a bit stronger that the 
original \negc.
Namely, this version deals with the {\it absolute number}
of trajectories joining $p$ with $q$.

Here is the version of the conjecture,
published in 1993 by S.~P.~Novikov \cite{noviquasi}, p.229:

\begin{quote}
``CONJECTURE. For any closed quantized
analytic 1-form
$\omega$, the boundary operator $\partial$
in the Morse complex has all coefficients
$a_{pq}^{(i)}\in K$
with positive part convergent in some region
$\vert t\vert \leq r_1^*(\omega)\leq r_1$. (If we replace
$S$ by $(-S)$ and $t$ by $t^{-1}$, we have to
replace $r_1$
by $r_1^{-1}$.) In particular, $r_1^*\not=0$ and
there will be coefficients
$a_{pq}^{(i)}$ with radius of convergence not
greater than
$r_1(M,[\omega])$ if at least one jumping point really exists.''
\footnote{For convenience of the reader, we explain the terminology
in this citation. S.~P.~Novikov denotes by $K$ the ring
denoted by $\ZZZ((t))$ in the present paper; $p,q$ correspond to
critical points with $\ind p = \ind q + 1$, and $a_{pq}^{(i)}$ is
$n(p,q;v)$ in our notation. A closed 1-form $w$ is called {\it
quantized} if its cohomology class belongs to $H^1(M,\ZZZ)$.}
\end{quote}

Thus the version of \cite{noviquasi}
says that \egp~holds for every 
Morse pair $(f,v)$ where $f$ is analytic and $v$ is a 
Riemannian gradient of $f$.

We proceed to versions of \negc~which are in a sense stronger than \negc.

\bede\lb{d:ratdef}
We say that 
the {\it rationality property}
holds for the Morse-Novikov pair
$(f,v)$, if \fe~ $p,q\in S(f)$ with $\ind p=\ind q+1$
the power series $\npqv$
is the Laurent series of a rational function.
\end{defi}

It is obvious that rationality property for a \mnp
$(f,v)$
implies the \egp~for $(f,v)$.

M.Farber has conjectured that the \rp~ holds \fe~ \mnp~
~$(f,v)$. P.Vogel conjectured, that the \rp~ holds \fe~ \mnp
$(f,v)$. He motivated it by the fact, that the \nics~are results of
one and the same procedure
applied to the fundamental cobordism of $f$
(see the definition of the fundamental cobordism of $f$
below in the \sub~\ref{su:pmm}).

\newpage

\section[The contents of the paper]{$C^0$-generic rationality of Novikov
incidence coefficients.
The contents of the paper}
\lb{s:cgrcp}

\subsection{The  contents of the paper}
\lb{su:cpp}

In this subsection we give a brief outline of the 
contents of the present paper; the reader will find more information
 on the principal ideas
of the work  in their relation to the main theorem
in \sub~\ref{su:opmt} and in the comments and introductions to 
different chapters.

About 1995 I had proved that the \nics~are $C^0$-generically rational 
functions.
For the first time the results appeared in the e-print dg-ga 9603006
(14 March 1996).
They were announced in \cite{pamrl}
and published in the paper \cite{pastpet}, which is a revised
for publication version of \cite{paepr}.

The aim of the present paper is to generalize this result and to provide
 a systematic treatment of the subject.
For a Morse map $f:M\to S^1$
we define a new class of vector fields, which contains the \glvf s for $f$
and riemannian gradients for $f$.
These vector fields are called $f$-\gr s (see Definition \ref{defgrad}).
We give (Chapter 2) the systematic 
and self-contained 
presentation
of all the necessary results and constructions from
 Morse theory; so this chapter can be considered as an introduction 
to the corresponding part of Morse
theory.

 \sub s \ref{su:tcrl}, 
\ref{su:ahc},
\ref{su:tfo},
\ref{su:hpof},
\ref{su:dmc},
\ref{su:pic}
contain the standard Morse-theoretical techniques, generalized to 
the case of $f$-gradients.
In \sub s
\ref{su:tddcp}
and  Sections
\ref{s:ssags}, \ref{s:gdm}
we generalize the techniques of
\cite{pastpet}
to the case of $f$-gradients.

The chapter \ref{ch:qf}
contains a construction which is a Morse-theoretic analog
of triangulation of a \ma~ with simplices of arbitrary small
diameter.
This construction makes no use of triangulations of 
the \ma;
we do not leave the framework of Morse theory, \Mf s
and their gradients.
The contents of this chapter does not differ much from the 
corresponding chapter 
of \cite{pastpet},
we have only revised the exposition
on order to make it easier to read.

Chapter \ref{ch:rgdm}
is technically 
the central part of the work.
It deals with gradients of Morse functions
on \cob s.
We construct a small perturbation of a given \gr, so that the \gr
~descent map, associated to the perturbed \gr,
respects a given handle decomposition of the components of the boundary
of the \cob. 

Chapter \ref{ch:gpbon}
contains the proof of the main results of the work.
Let $\GT(f)$
be the set of all $f$-\gr s, satisfying \TA.
We introduce a subset $\GCCT(f)\sbs \GT(f)$
which is open and dense in $C^0$ topology
(see Proposition \ref{p:cgener}).
The next theorem 
asserts that $C^0$ generically
the \nics~are rational functions.
Denote by $\wi L$
the subring of the ring $\QQQ(t)$
formed by all the rational functions of the form
$\frac {P(t)}{t^mQ(t)}$
where $m\in \NNN$ and
$P,Q$ are polynomials with integer coefficients
and $Q(0)=1$.

\begin{theocitbeznom}  {(Sect. 5, Th. \ref{t:cccgr})}
Let $v\in \GCCT(f)$. Let $x,y\in S(f)$,
$\ind x=\ind y+1$. Then $n(x,y;v)\in \wi L$,
\end{theocitbeznom}  

The gradients belonging to $\GCCT(f)$
have remarkable properties.
In particular, the incidence coefficients $n(x,y;v)$
are invariant under $C^0$ small perturbations
of $\GCT(f)$, if this perturbation has support outside $S(f)$.
Here is the precise statement.

\begin{theocitbeznom}
{(Sect. 5, Th. \ref{t:continc})}
Let $v\in\GCCT(f)$. Let $U$ be an open \nei~
of $S(f)$.
Then there is $\d>0$ \sut~ for every
$w\in \GCCT(f)$
 with $\Vert w-v\Vert<\d$ and $v\mid U=w\mid U$,
and for every $x,y\in S(f)$ with $\ind x=\ind y+1$ we have:
$n(x,y;v)=n(x,y;w)$.
\end{theocitbeznom}

There are some naturally arising versions of Novikov Complex
defined over completions of group rings, for example
the completion of the group ring $\ZZZ[\pi_1(M)]$.
(One  may consider these versions as "equivariant analogs" of the 
usual Novikov Complex).
There is the corresponding generalization of the
Theorem
\ref{t:cccgr}
to this situation.
To cite here this theorem, let $f:M\to S^1$ be a Morse map,
which we assume to be primitive, that is
the \ho~$f_*:\pi_1(M)\to \pi_1(S^1)=\ZZZ$
is epimorphic.
To the \ho~  $\xi=f_*$
one associates a certain completion $\Lxi$
 of the group ring $\L=\ZZZ[\pi_1(M)]$
(see \sub s    \ref{su:nrbt} for definition).
There is a non-commutative localization procedure,
which allows to define a certain subring of $\Lxi$
which is a natural analog of the ring $\wi L$.
This subring is denoted by $\Im \ell$
(see \sub~  \ref{su:ncalg} for definition).
The following statement is a part of the theorem \ref{t:eqcccgr}.

\begin{thbeznom}
Let $v\in\GCCT(f)$. Then for every $x,y\in S(f)$
with $\ind x=\ind y+1$
the Novikov incidence coefficient
$\wh n(\wh x,\wh y;v)$ 
is in $\Im\ell$.                                        
\end{thbeznom}

\subsection{Preliminaries on Morse maps $M\to S^1$}
\lb{su:pmm}

In this \sub~we introduce some terminology which will be used also in
Chapter \ref{ch:gpbon}.

Let $M$ be a closed connected \ma~and $f:M\to S^1$
be a Morse map, non homotopic to zero.
We shall assume that the homotopy class 
$[f]\in [M, S^1]=H^1(M,\ZZZ)$
is indivisible
(the general case is easily reduced to this one by considering
the map
$\wi f=f/n$ with some $n\in \NNN$).

The map $f$ lifts to a Morse function
$F:\bar M\to \RRR$.
Let $t$ be the generator of the structure group
($\approx\ZZZ$)
of $\CC$, \sut~for every $x$ we have:
$F(xt)< F(x)$.
(We invite the reader to look at the figure
on the page \pageref{fig:cyccov}
 where one can visualize these constructions.) 
We shall assume that $f$ is primitive, that is
$f_*:H_1(M)\to H_1(S^1)$
is epimorphic
(one can always reduce the situation to this case by considering
the map $F/n$ with some $n\in\NNN$.)
In this case we have:
$F(xt)=F(x)-1$.
To simplify the notation we shall assume that
$1\in S^1$
is a regular value of $f$,
and \th~ every $n\in\NNN$ is a regular value of
$F$.
Denote $f^{-1}(1)$ by $V$.
Set
$V_\a=F^{-1}(\a),\quad W=F^{-1}\([0,1]\), \quad
V^-=F^{-1}\(]-\infty,1]\)$.
Note that $V_\a t=V_{\a-1}$.

The cobordism $W$ can be thought of as the result of
cutting $M$ along $V$.
We shall sometimes call $W$ a {\it fundamental cobordism} of $f$.
Denote $Wt^s$ by $W_s$; then $\bar M$
is the union
$\cup_{s\in\ZZZ} W_s$, the neighbor copies
$W_s$ and $W_{s+1}$
intersecting by
$V_{-s}$.
For any $n\in\ZZZ$ the restriction of $\CC$ to $V_n$
is a diffeomorphism
$V_n\to V$.

The $\CC$-preimage of a subset $A\sbs M$
will be denoted by
$\bar A$.
On the other hand, if
$x\in M$ we shall reserve the symbol
$\bar x$ for liftings of $x$ to
$\bar M$
(that is $\bar x\in \bar M$ and $\CC(\bar x)=x$).
Let $v$ be an $f$-gradient.
Lift $v$ to $\bar M$;
we obtain a complete vector field on $\bar M$,
which is obviously an $F$-gradient.
We shall denote this $F$-gradient by the same symbol $v$, since \noconf.
To restriction of this $F$-gradient to any cobordism
$W_{[\l, \mu]}=F^{-1}\([\l,\mu]\)$
(where $\l,\mu$ are regular values of $F$)
is an
$(F\mid W_{[\l,\mu]})$-gradient.
Note that if $v$ satisfies \TA,
then the $F$-gradient
$v\mid W_{[\l,\mu]}$
satisfies \TA.

The set of all $f$-gradients will be denoted by $\GG (f)$, and the 
set of all $f$-\gr s satisfying \TA will be denoted by
$\GT(f)$.

The \nic~
$n_i(p,q;v)$
equals by definition
to the incidence coefficient
$n(p,qt^i;v)$
of the critical points $p$ and $qt^i$
of the \Mf~ $F$, restricted to the \cob~ $F^{-1}\([a,b]\)$,
where the interval $[a,b]$ is chosen to be
sufficiently
large as to contain
$F(p)$ and $F(qt^i)$.

\subsection{First step of the proof}
\lb{su:fsp}
We adopt here the terminology of the preceding subsection.

Let $p,q\in W$ be critical
points of $F$,
$\ind p=\ind q+1$. Set $n=\dim M$
and $l+1=\ind p$.
For simplicity, assume that $F(p)$ is the lowest
critical level of $F|W$, so that the interval $]0,F(p)[$
is regular.
Similarly, we assume that the interval $]F(q),1[$ is regular and
that
$\bar p\in W_0$, $\bar q\in W_1$.
Let $S(p)$ denote the
intersection of the stable manifold
of $p$ with $V_0$; this is an oriented
$l$-dimensional
embedded sphere in $V_0$. The intersection $S(q)$
of the unstable manifold of $q$ with $V_1$
is a cooriented $(n-l-1)$-dimensional
sphere in $V_1$.

Let $v\in\GT (f)$. 
Let $K_1$ be the union of all unstable discs
of critical points of $F|W$ and $K_0$
be the union of all 
stable discs of critical points of $F|W$.
Let $\stexp {(-v)}$
be the diffeomorphism 
of $V_1\sm K_1$ to $V_0\sm K_0$
which sends each point $x$ to the point
of intersection of the $(-v)$-\tr~
starting at $x$ with $V_0$.
Using the presentation
of $\bar M$ described above,
we obtain the identity
\begin{equation}
n_{k}(p,q;v)=(\stexp {(-v)})^k\big(S(p)\big)\krest S(q)  \lbl{krest}
\end{equation}
Here are some comments on this formula.
1) The symbol $\krest$ stands for the algebraic
intersection index of two manifolds; it can be
checked that this number is well defined in our situation.
2) $\stexp {(-v)}$ is not defined on all of $V_1$; nevertheless, 
we denote
by $\stexp {(-v)}(A)$ the set $\stexp {(-v)}(A\sm K_1)$, where
$A\subset V_1$. 
3) We identify
$V_1$ and $V_0$ with $V$, so that the
iterations of $\stexp {(-v)}$ make sense
and are applicable to $S(p)\subset V_0$.

Should $\stexp {(-v)}$ be an everywhere defined
diffeomorphism, the intersection
index in (\ref{krest}) would coincide with the
intersection index of some homology classes,
and our result would be a consequence of linear
algebra. Indeed, let
$[p]$ be a homology class of the sphere $S(p)$ in
$H_{l}(V)$, let $]q[{}\in H^{l}(V)$
be the cohomology class dual
to $S(q)$, and let $L:H_*(V)\to H_*(V)$
be the homomorphism induced by $\stexp {(-v)}$; then
$((-v)^{\rightsquigarrow})^k(S(p))\krest S(q)=\langle ]q[,
L^k([p])\rangle t^k$ (here and in what follows we denote by
$\langle \cdot , \cdot \rangle$ the canonical pairing between homology
and cohomology). A simple computation shows
that the power series
$\sum_{k=0}^\infty \big\langle ]q[, L^k([p])t^k\big\rangle$
is a rational function with
denominator $1-\det L\cdot t$.
(see Lemma \ref{l:cramer} for details).

These considerations were well known to 
S.P.Novikov and V.I.Arnold in th end of 80s.
They form a motivation for \negc.
I learned the algebraic lemma
\ref{l:cramer} from V.I.Arnold lecture in Moscow in 1990.
The paper
\cite{Arno}
gave rise 
to anew domain of investigation: the theory of topological complexity
of intersection.
The basic object of study in this theory
is the intersection
of the form
$A^n X\cap Y$, where $X$ and $Y$ are \sma s of a \ma~$M$
and $A$ is a diffeomorphism~$M\to M$.
One tries to establish the asymptotic properties 
of this intersection
when $n\to\infty$.
See the paper 
\cite{Arnoprob}
and the bibliography there for the actual state
of this theory.

\subsection{Outline of the proof of the main theorem}
\lb{su:opmt}

The exposition in this section is rather informal,
 nevertheless we think that it will be useful to the
 reader to go through it, since it explains the main ideas of
what follows.

We continue with the arguments of the previous subsection.
If $f$ has critical points, then the
map $\stexp {(-v)}$ is no longer an everywhere defined
diffeomorphism, and the indices $n_i(p,q;v)$ do
not admit any simple homological interpretation.
We shall show
 that for a $C^0$-generic $f$-gradient $v$ we can still
obtain a formula of type $n(p,q;v)=
\sum_{k=0}^\infty \big\langle ]q[, L^k([p])\big\rangle t^k$,
where, however, the meaning of the objects
$[p], ]q[, L$ will be more complicated than in the model situation
discussed in \sub~\ref{su:fsp}.

Namely, we shall
show that for any $f$ and any $p,q\in S(f)$ one can find
an $C^0$-open-and-dense subset $G\subset \GT(f)$, disjoint compacts
$A,B\subset V$, and a homomorphism
$L:H_{l}(V\sm B, A)\to
H_{l}(V\sm B, A)$ such that the following condition is fulfilled
for any $v\in G$:

$S(p)$ does not intersect $B$, $S(q)$ does not intersect $A$,
the group $H_*(V\sm B, A)$ is finitely generated, and

\begin{equation}
n(p,q;v)=\sum_{k=0}^\infty\big\langle  \lb{f:genfun}
]q[, L^k([p])\big\rangle t^k,
\end{equation}

where $[p]$ is the image in
$H_l(V\sm B,A)$ of the
fundamental class
$[S(p)]$, and
$]q[$ is an element of $H^l(V\sm
B,A)$, dual to $S(q)$.

\bere
The homomorphism $L$ is {\it not } induced by any continuous
map of the pair $(V\sm B, A)$ to itself.
\enre

Now we describe what are the sets
$A,B$, the homomorphism $L$, and the
$C^0$-small perturbations of $v$ satisfying (\ref{f:genfun}). We use
the terminology of \sub s \ref{su:tddcp}, \ref{su:tcrl}.
We start with an auxiliary Morse function
$\phi\:V\to \RRR$ and a
$\phi$-gradient $u$ satisfying the almost transversality condition.
For $i\leq n-1$ consider the set
$D(\indl i; u)$
(recall from
Definition \ref{d:dthick} that it is the union of all the stable \ma s
of $v$ of dimension $\leq i$.)
These sets form a sort of cellular decomposition of $V$,
 and we shall call
$D(\indl i;u)$
{\it $i$-skeleton} of $V$.
(This is a certain abuse of notation, but in this subsection
we shall not use the term
{\it skeleton}
in its usual meaning ).
The sets 
$D(\indl j ;-u)$ form the dual decomposition of $V$.
We have also "thickenings" of this skeleton:
the sets
$B_\d(\indl i;u),
D_\d(\indl i;u),
C_\d(\indl i;u)$.
Here 
$B_\d(\indl i;u)$
is an open set, and
$D_\d(\indl i;u),
C_\d(\indl i;u)$
are compacts.
One may think of 
$D_\d(\indl i;u)$
as of union of all handles of $V$ of thickness $\d$ and indices $\leq i$.

For $\d$ sufficiently small the sets
$D_\d(\indl i;u)$ and
$D_\d(\indl {n-i-2};-u)$
are disjoint. 

Put 
$A=
D_\d(\indl {l-1};u),\quad
B=
D_\d(\indl {n-l-2};-u)$.

Returning to the cobordism $W$
with boundary $\partial W=V_0\sqcup V_1$,
we note that for every geometric
object associated with $V$ there are two
copies of it that correspond to $V_0$ and to $V_1$, because both $V_0$
and $V_1$ are diffeomorphic to $V$.
These objects will be denoted by the same letter but with
indices 0 or 1, respectively.
In this way, we obtain Morse functions $\phi_0,\phi_1$,
compacts $A_0,A_1$, etc.

The \ho~ $L$ will be constructed for vector fields $v$,
 satisfying additional restriction, namely condition
$(\gC\YY)$ (see Definition \ref{d:condcy}). The principal part 
of this condition is given by the following two
formulas:

 \begin{gather*}\stv \bigg(C_\d(\indl j;u_1)\bigg)
\cup
\bigg(D_\d\(\indl {j+1};v\)\cap \pr_0 W\bigg)
\sbs
B_\d(\indl j,u_0)
\\
\mbox{ for every } j\tag{B1}
\end{gather*}

\begin{gather*}\st v\bigg(C_\d(\indl j;-u_0)     \bigg)
\cup
\bigg( D_\d(\indl {j+1}; -v)\cap \pr_1 W\bigg)
\sbs
B_\d(\indl j;-u_1)
\\
\mbox{ for every } j\tag{B0}
\end{gather*}

(see Definitions \ref{d:condcy}, \ref{d:cc} for the details).

The conditions (B0), (B1)
 may be regarded as analogs of the cellular approximation
conditions; however,an additional restriction
is imposed: the soles of the $\d$-thickened handles of $v$ of
dimension less than or equal to $j+1$ must belong to the
$\d$-thickened $j$-skeleton of $V$.

If condition $(\gC\YY)$
 is satisfied, then we construct
the homomorphism
$L=\HH_l(-v):H_l(V\sm B, A)\to
H_l(V\sm B,A)$.
To clarify the idea of this construction, we 
make a simplifying assumption: we assume
that the \hog~class 
$x\in H_l(V\sm B, A)$
is represented by an oriented closed \sma~ $N$ of
$V\sm B$
of dimension $l$.
Set $N'=\stexp {(-v)}(N_1)$;
then $N'$ is a \sma~of $V_0$ 
(non compact in general).
The condition (B1) implies 
$N'\sbs V_0\sm B_0$.
The condition (B0) implies that of the discs $D(a, -v)$
the \ma~ $N_1$ can intersect only those for which $\ind a\leq l$.
Therefore 
the set
$N'\cup (D(\indl l;v)\cap V_0)$ is compact,
and the set $N'\sm B_\d(\indl {l-1}; u_0)$
is also compact.
Thus the fundamental class of $N'$ modulo $A_0$
is well defined as an element 
of $H_l(V_0\sm B_0, A_0)$.
Identifying $V_0$ and $V$ we obtain an element 
of $H_*(V\sm B, A)$, which is by definition $\HH_l(-v)(x)$.

The detailed construction of $\HH_l(-v)$ is 
the subject of the section \ref{s:hgd}.

Now we explain how to 
 gain the 
properties (B0), (B1) by $C^0$-small perturbation of the
initial $f$-gradient $v$ (that is, how to prove that the condition
$(\gC\YY)$ is $C^0$-dense.)

As already noted, for a compact submanifold $R$ in $V_1$,
the submanifold $\stexp {(-v)}(R)$ in $V_0$
is not necessarily compact; to make it
compact, we must add the soles of the
descending disks of $v$ (thus losing
the manifold property). Therefore, it should be expected
that the correspondence
$\cdot\mapsto\stexp {(-v)}(\cdot)$ is well defined on an appropriate
category of stratified \ma s. As shown in \ref{s:ssags}, this
is indeed so. The results of \ref{s:ssags} imply that (under some
mild transversality conditions)
the set
\begin{equation}
{\RR}_j=\stexp {(-v)}\big(D(\indl j;u_1)\big)\cup
\big(D(\indl {j+1};v)\cap V_0\big)
\end{equation}
is compact for every $j$. Since this set is a union of submanifolds
of $V_0$ of dimension not exceeding
$j$, we can push ${\RR}_j$ away from
$D(\indl {n-j-2}; -u_0)$; then,
applying a diffeomorphism $\Phi(T,-u_0)$ for $T$
sufficiently large, we
obtain an isotopy of $V_0$ that pushes ${\RR}_j$ into
$B_\d(\indl j;u_0)$. After some work with
thickenings we obtain
an isotopy of $V_0$ such that the set
$\stexp {(-v)}(D_\d(\indl j;u_1))\cup (D_\d(\indl {j+1} ; v)\cap V_0)$
is also pushed  inside $B_\d(\indl j;u_0)$.
This isotopy may be realized by modifying
$v$; the new gradient will satisfy (B0) and (B1). 
This procedure is described
in detail in Section \ref{s:cc}.
The resulting
gradient can be chosen so as to be $C^0$-close to the
initial one. In other words, the pair
$(\phi_0,u_0)$ can be chosen in such a way that
$\Phi(\tau, -u_0)$ pushes
$V_0\sm B_\d(\indl {n-j-2}; -u_0)$ inside
$B_\d(\indl j;u_0)$ already for $\tau$ sufficiently small.
Such flows will be called
{\it quick flows}. The construction of quick
 flows is the subject of the Chapter \ref{ch:qf}.

As we have already mentioned, our construction of quick flows
differs from the usual
constructions of ``small handle decomposition''
in that it uses no triangulations of
the manifold. A description of the basic idea of this construction
can be found at the beginning of Chapter \ref{ch:qf}.

\subsection{Further remarks}
\lb{su:fr}

\been\item
Note that the $C^0$ topology is very important for us during 
all the work. One of the reasons
why the $C^0$ topology is relevant, lies in the following
property of the integral
curves of differential equations:
 the integral
curves (on a given finite interval) of
a {\it $C^0$-small} perturbations
$v'$ of a $C^1$-vector field $v$ are $C^0$-small perturbations
of the integral curves of $v$ (see, e.g., \cite{biro}, Chapter 5, \S4]).

See the section \ref{s:pcic}
for the detailed exposition of the $C^0$ continuity properties.
We think in general that $C^0$ topology is much more natural and efficient
in Morse theory than it seems at present.

\item
The section \ref{s:icvc}
of the paper is devoted to the generalization of
\negc~ for the case of
incidence coefficients
with values in group rings.
These results (in a somewhat weaker form)
were announced in \cite{pamrl}.
Unfortunately I have no time to include in the present paper
the proofs
of the other
results announced in \cite{pamrl};
that will be the subject of the second part of the present work.
\enen

\section{Terminology}
\lb{s:termi}

\subsection{Terminological conventions}
\lb{su:tc}

Here I have gathered some definitions and notational conventions
which I use throughout the paper.

The  closed disc of radius $r$ and center $a$ in the Euclidean
space $\RRR^n$ is denoted by $D^n(a,r)$.

The  open disc of radius $r$ and center $a$ in the Euclidean
space $\RRR^n$ is denoted by $B^n(a,r)$.

Let $v$ be a vector field on a manifold $W$
(possibly, with boundary). We denote by
$\gamma(x,t;v)$
the value at $t$ of the integral curve
of $v$ starting at $x$.
If $A\subset \RRR$ and $\gamma(x,\cdot ;v)$
is defined
on $A$, then $\g(x,A;v)$ denotes the set
$\{\,\gamma(x,t;v)\mid t\in A\,\}$.

A subset $X\subset W$ is said to be {\it $v$-invariant} if
$\g(x,t;v)$ is defined for every $t\geq 0$
and every $x\in X$
and $\g(x,[0,\infty[;v)\subset X$. A subset $X\subset W$
is {\it $(\pm v)$-invariant} if it is $v$-invariant
and $(-v)$-invariant.

A subset $X\subset W$ is said to be {\it weakly
$v$-invariant} if
$\g(x,t;v)\in X$ whenever $x\in X$, $t\geq 0$, and
$\gamma(x,t;v)$ is defined. Defining
the notion of weak $(\pm v)$-invariance is left to the reader.

For any (time dependent) vector field
$w$ on a closed manifold $M$
and any $t\in\RRR$, we denote by $\Phi(w,t)$
the diffeomorphism $x\mapsto\gamma(x,t;w)$
of $M$.

The closure of a set $U$ is denoted by
$\overline U$. The short bar is reserved for the objects
related to the infinite cyclic covering (e.g.,
$\PP \:\bar M\to M$).

Let $R:A\to B$ be a diffeomorphism of riemannian \ma s.
For $x\in A$ set $E(R,x)=\Vert R'(x)\Vert$.
The number
$$\max\big(\sup_{x\in A} E(R,x), \sup_{y\in B} E(R^{-1}, y)\big)$$
is called 
{\it expansion constant} of $R$.

{\it A cobordism} is 
a compact manifold $W$
together with a presentation
$\pr W=\pr_0W\sqcup \pr_1W$
where
$\pr_1W$ and $\pr_0W$
are compact manifolds of dimension $\dim W-1$ 
without
boundary
(one or both of them can be empty). The \ma
~$W\sm\pr W$
will be denoted by $\Wkr$.

{\it A riemannian cobordism}
is a cobordism endowed with  a riemannian metric.

Let $W$ be a riemannian \cob~ of dimension $n$,  $q\in \Wkr$.
Let $\d_0>0$ be so small, that the restriction of the exponential map 
$\exp_q:T_q W\to W$
to the disc 
$B^n(0,\d_0)$
is a diffeomorphism onto its image, and this image belongs to 
$\Wkr$.
 Let $0<\d<\d_0$.
Denote by $B_\d(q)$, resp. $D_\d(q)$ the $\exp_q$-images
of the open euclidean ball $B^n(0,\d)$,
resp. closed euclidean ball $D^n(0,\d)$ with the center in $0$.
Then 
$B_\d(q)$, resp. $D_\d(q)$ are riemannian $\d$-discs (open, and resp.
 closed).
In what follows we shall use the notation 
$B_\d(q),  D_\d(q)$ 
only when $0<\d_0<\d$.

The end of the proof of a theorem, lemma etc. is marked by $\square$;
 the end of a definition or a remark is marked by $\triangle$.

\subsection{Some differences between the notation of this paper
and that of
\cite{pasur}, \cite{paepr}}
\lb{su:diffnot}

The gradients satisfying the almost 
transversality condition
 were called ``almost good'' in \cite{paepr}, \cite{paepri},
 \cite{pasur}. 

The diffeomorphism
denoted by $\stexp {(-v)}$ in the present paper was denoted by $\st v$ in 
\cite{paepr}, \cite{paepri}, \cite{pasur}.Our 
present notation is the same as that of \cite{pastpet};
it seems more natural because this
diffeomorphism is a shift in the {\it positive} direction along
the integral curves of the field
$-v$. In  \cite{pasur} the sets $B_{\delta}(\indl s;v)$,
$D_{\delta}(\indl s;v)$, $D(\indl s;v)$ were denoted by
$B_{\delta}(\leq s;v)$, $D_{\delta}(\leq s;v)$, $D(\leq s;v)$. In 
\cite{pasur}
the set $D(\indl s, v)$ was denoted by $K(\indl s;v)$.

\chapter{Morse functions and their gradients}
\lb{ch:mftg}

\section{ Basic definitions and constructions }
\lb{s:bdc}
\pa
\subsection{ Basic definitions}
\lb{su:bd}

Recall that  {\it cobordism} for us is 
a compact manifold $W$
together with a presentation
$\pr W=\pr_0W\sqcup \pr_1W$
where
$\pr_1W$ and $\pr_0W$
are compact manifolds of dimension $\dim W-1$ 
without
boundary
(one or both of them can be empty). The \ma
~$W\sm\pr W$
is denoted by $\Wkr$.

\begin{defi}
A {\it Morse function}
$f:W\to[a,b]$
is a $\smo$ map $f:W\to\RRR$, \sut ~
$f(W)\subset[a,b], f^{-1}(b)=\pr_1W,
f^{-1}(a)=\pr_0W$,
all the critical points of $f$ are 
non-degenerate and
belong to
$\Wkr$. 
\end{defi}

We shall denote $f^{-1}(\l)$ by $V_\l$ and 
$f^{-1}([\a,\b])$ by
$W_{[\a,\b]}$ if \noconf.

Recall that for a critical point $a$ of $f$
 there is
a well defined bilinear form
on $T_aM$, called {\it  Hessian}
of $f$, and denoted by $d^2f(a)$.
For a non-degenerate critical point $a$    
 the dimension
of the
negative eigenspace of $d^2f(a)$
is called {\it  index} of $a$. The set of 
critical
point of $f$
is denoted by $S(f)$, the set of critical 
points of
$f$ of
index $k$ is denoted by $S_k(f)$.

To introduce the notion of gradient vector
field for
$f$ recall first some notions and
results from the general theory of vector fields.
Our reference for that theory is the book by R.Abraham and J.Marsden
\cite{abrob}, 
but we modify a bit their terminology. Let 
$M$ be a \ma~ 
without boundary. We say that a zero $p$ 
of a vector field $v$ on $M$
is {\it elementary}
if the matrix of the derivative $v'(p)$ has 
no eigenvalues with zero 
real part.
Recall the following Hadamard-Perron local 
stable manifold theorem (\cite{abrob}, p.85)

\begin{theo}\label{locstma}
 Let $p$ be an elementary zero 
of a vector field $v$. Then there is a \nei~ 
$U(p)$ and two submanifolds 
$W^{un}_{loc}(p,v)\sbs U(p),
W^{st}_{loc}(p,v)\sbs U(p)$,
having the following properties:
\begin{multline*}
W^{un}_{loc}(p,v)=
\{ x\in U(p)\mid \g(t,x;v) 
\mbox{ is defined and belongs to }   
 U(p) \mbox{ for all } 
t\leq 0, \\
\mbox{ and }
\g(t,x;v)\rad{t\to -\infty} p
\end{multline*}

\begin{multline*}
W^{st}_{loc}(p,v)=
\{ x\in U(p)\mid \g(t,x;v) 
\mbox{ is defined and belongs to }\\  
  U(p) \mbox{ for all } 
t\geq 0, \mbox{ and }
\g(t,x;v)\rad{t\to  \infty} 
p
\end{multline*}

$W^{un}_{loc}(p,v)\ptf
W^{st}_{loc}(p,v)$,

$W^{un}_{loc}(p,v)\cap
W^{st}_{loc}(p,v)=\{p\}$.
\end{theo}

The tangent spaces to
$W^{st}_{loc}(p,v)$, resp.
$W^{un}_{loc}(p,v)$ will be denoted by
$T_-(p;v)$, resp. $T_+(p;v)$.

Let $W$ be a \ma~~ with boundary.
 Let $v$ be a vector field on $W$, 
\sut~~ $\forall x\in\pr W: v(x)\not=0$.
We can apply the theorem above to
every elementary zero $p$ of $v\mid \Wkr$.
Thus we obtain:
a \nei~ $U(p)$, local unstable manifold
$W^{un}_{loc}(p,v)$,  local stable 
manifold
$W^{st}_{loc}(p,v)$ and the corresponding
  subspaces
$T_+(p;v),  T_-(p;v)$ of $T_pW$.

\begin{defi}    \label{defgrad}
Let $\fcob$ be a Morse function on a 
cobordism $W$.
A vector field $v$ on $W$ is called 
{\it gradient for $f$},
or simply {\it $f$-gradient}, if
\begin{enumerate}
\item For every $x\in W\sm S(f)$
we have:
$df(x)\((v(x)\)>0$.
\item
The set of zeros of $v$ equals to $S(f)$
and each zero of $v$ is elementary.
\item For every $p\in S(f)$
the vector subspaces
$T_-(p,v)$ and $T_+(p,v)$ are orthogonal
\wrt~ the bilinear form
$d^2f(a)$. $\triangle$
\end{enumerate}
\end{defi}

\begin{rema}\label{r:dgop}
The same definition is valid of course for
Morse functions on manifolds without boundary
(the case of closed manifolds is already contained in the definition).
\end{rema}

The basic example of an $f$-gradient is 
riemannian
gradient. Namely, let $G$ be a riemannian 
metric
on $W$; recall that {\it riemannian gradient} 
$\grad \! f$ is the vector 
field defined
by
the formula

$$
\langle\grad f(x),h\rangle=df(x)(h) 
\for x\in M
$$

\begin{lemm}
Riemannian gradient of a Morse function $f$ 
is an
$f$-gradient in the sense of the preceding
definition.
\end{lemm}

\Prf
The property 1) and the first half of 2) are 
obvious.
Proceeding to the local properties in a \nei~~
of a critical point $p\in S(f)$,
consider the geodesic coordinate system
 centered in
$p$.
We can assume therefore, that $f$ and 
$v=\grad\! f$
are defined in a \nei~
~$U$ of $0\in\RRR^n$, and that for every $x\in U$ the matrix 
$G(x)$
of our riemannian
metric in $U$ satisfies
$G(x)=I+O(\vert x\vert^2)$.
Write the Taylor decomposition
of $f$ as
$f(x)=f(0)+1/2 \langle Ax,x\rangle+ 
O(\vert x\vert^3)
$
where $\langle\cdot,\cdot\rangle$
stands for the Euclidean scalar product,
$\vert\cdot\vert$ stands for the Euclidean norm,
and
$A$ is a non degenerate symmetric linear 
operator.
Then we have:
$f'(x)h=\langle Ax, h\rangle+
O(\vert x\vert^2) (h)$,
and
$v(x)=Ax+O(\vert x\vert^2)$.
Thus $v'(0)=A$ and all the eigenvalues of 
$v'(0)$ are real and nonzero.
$\qs$

{\bf Remark on the terminology}
\lb{remterm}

Our present terminology
is different from that of
\cite{paepr} and from that of \cite{pasur}.
For the convenience of the reader we cite here 
the
definitions from
\cite{milnhcob},
\cite{pasur} and \cite{paepr}, and establish relations between them.
\pa

Here is the definition of \glvf~ (\cite{milnhcob}, Def. 3.1):

\begin{quote}

Let $f$ be a Morse function for the triad $(W^n; V, V')$.
A vector field $\xi$ on $W^n$
is a \underline{\glvf ~ for $f$} if
\been\item
$\xi(f)>0$ throughout the complement of the set of critical 
points of $f$, and
\item
given any critical point $p$ of $f$
there are coordinates
$(\vec x, \vec y)=(x_1,..., x_\l, x_{\l+1}, ..., x_n)$
in a \nei~  $U$ of $p$ so that
$f=f(p)-|\vec x|^2
+
|\vec y|^2$
and $\xi$ has coordinates
$(-x_1,..., -x_\l, x_{\l+1}, ..., x_n)$
throughout $U$.
\enen
\end{quote}

Here is the definition of $f$-gradient and 
related objects from
 \cite{pasur},\S 2.2:

\begin{quote}

     Let $f:W \to [a,b]$ be a Morse function 
on $W$
(for
the definition see
\cite{milnhcob},
 def 2.3).
Denote by $S(f)$ the
set of critical points of $f$ and by $S_k (f)$ 
the
set of
critical points of $f$ of index $k$. For each 

$p \in
S(f)$
choose a chart $\Phi(p) : U(p) \to  B^n(0,r)$ 
(where
$U(p)$ is
 a neighborhood of $p$ in $W$, and $n = $ dim $W$),
such that
$f \circ \Phi(p)^{-1} (x_1,...,x_n) = 
f(p)-(x_1^2 +...+x_k^2)
+
(x_{k+1}^2 +...+x_n^{2}) $ where $k = $ ind $p$ . We
call this
chart {\it standard chart} and we consider it as a
part of the
data of the Morse function $f$. The neighborhood 
$U(p)$
is called
 {\it standard coordinate neighborhood}. We say
(cf. \cite{milnhcob},
  \S 3),
 that a vector field $v$ on $W$ is a {\it gradient}
for $f$, or an
{\it $f$-gradient} if 1) $df(v) > 0$, except for 
the
critical points
of $f$, and
2) for any $p\in S(f)$ the vector fields $\Phi(p)_*(v)$
 and $(-x_1,...,-x_k, x_{k+1},...,x_n)$
are equal in some neighborhood of $0$ (where 
$k=\ind
p$).

\end{quote}

Here is the definition 
of $f$-gradient and related objects 
from \cite{paepr}, def 1.1

\begin{quote}

 Let $f:M\to \RRR$ be
 a Morse function on a closed manifold
$M$. Denote $\dim M $ by $ n$.
The set of critical points of $f$ will be
denoted by $S(f)$.
A chart $\Phi_p:U_p\to B^n(0,r_p)$
(where $p\in S(f), U_p$ is a neighborhood
of $p$, $r_p >0$) is called \it standard
chart for $f$ around $p$ of radius $r_p$
\rm (or simply $f$-\it chart\rm ) if
there is an extension of $\Phi_p$ to a chart
$\widetilde\Phi_p:V_p\to B^n(0,r'_p)$,~
(where $\overline {U_p}\subset V_p$ and  $r'_p>r_p$),
such that
$(f\circ \widetilde\Phi_p^{-1})(x_1,...,x_n)=f(p)+
\sum_{i=1}^n\alpha_i x_i^2$, where
 $\alpha_i<0$ for $i\leq \ind_fp$ and
$\alpha_i>0$ for $i>\ind_fp$.
 The domain $U_p$
is called \it standard coordinate neighborhood.
\rm
Any such extension $\widetilde\Phi_p$ of $\Phi_p$
will be called
\it standard extension of $\Phi_p$
\rm

...

A family
${\UU}= \{\Phi_p:U_p\to B^n(0,r_p)\}_{p\in S(f)}$
 of $f$-charts is called $f$-\it chart-system\rm ,
 if the family $\{\overline{U_p}\}$ is disjoint.

...

Given an $f$-chart system
${\UU}= \{\Phi_p:U_p\to B^n(0,r_p)\}_{p\in S(f)}$
we say, that a vector field $v$ on $M$
 is an
\it $f$-gradient with respect to $\UU$,
\rm
if

1) $\forall x\in M \setminus S(f)$ we
have $df(v)(x)>0$;

2) $\forall p\in S(f)$ we have
$(\widetilde\Phi_p)_*(v)=(-x_1,...,-x_k,x_{k+1},...,x_n)$,
where $k=\ind_f p$, and $\widetilde\Phi_p$ is some
standard
extension of $\Phi_p$.

We say that a vector field $v$ is an
\it $f$-gradient \rm
if there is an $f$-chart system
 $\UU$ ,
such that $v$ is an $f$-gradient
with respect to $\UU$.

\end{quote}

The notion of \glvf~from \cite{milnhcob}
coincides essentially with that of \cite{pasur},
except that in \cite{pasur}
we considered the standard coordinate system around the critical point
as a part of the data of a Morse function.

The following simple lemma shows that the notion of
 $f$-gradient in the sense of
\cite{paepr}
and the notion of
\glvf~ are the same.

\bele\lb{l:grads}
Let $f:W\to\RRR$ be a \Mf
хn a \cob~ $W$, $v$ be a vector field on $W$.
Then $v$ is a \glvf, \ifff~
$v$ is an $f$-\gr~ in the sense of \cite{paepr}
\enle

\Prf To prove that every $f$-\gr~in the sense of \cite{paepr}
 is a \glvf~for $f$
(the other sense is obvious)
consider the standard coordinate system
$\Phi:U\to V$
for $f$ and $v$ around a critical point
$p\in S_k(f)$. We have:
$(f\circ \widetilde\Phi_p^{-1})(x_1,...,x_n)=f(p)+
\sum_{i=1}^n\alpha_i x_i^2$, where
 $\alpha_i<0$ for $i\leq k$ and
$\alpha_i>0$ for $i>k$
and
$(\widetilde\Phi_p)_*(v)=(-x_1,...,-x_k,x_{k+1},...,x_n)$.
Define a \dfm~
$\Psi:\RRR^n\to\RRR^n$
by
$\Psi(x_1,..., x_n)=
( \sqrt{|\a_1|} x_1,..., 
\sqrt{|\a_n|} x_n)$.
Then the chart $\Psi\circ \Phi$
is a standard chart for $f$ around $p$, and it is obvious that 
$(\Psi\circ\Phi)_*(v)(x)=
(-x_1,...,-x_k,x_{k+1},...,x_n)$.
$\qs$

Of course the notion of $f$-gradient suggested in the present paper
is strictly wider than that of \glvf.
In the present paper we shall sometimes use \glvfs.
For brevity we shall sometimes use the term
"$f$-Gradient" instead of "\glvf".
The reason for using in our paper the term "$f$-gradient" already used 
in the previous paper in the different sense 
(we apologize to the reader for every inconvenience
 caused by this)
is that we try to find the widest possible class of 
"gradients " for which the
 usual
Morse-theoretic
constructions
like "cellular-like decomposition, formed by the descending
separatrix discs"
are still valid.
The extension which we propose in this paper (Definition
\ref{defgrad}
)
is motivated by the fact that riemannian gradient of
a Morse function is
now an $f$-gradient (this was not true in general with
the definitions
of \cite{pasur}, \cite{paepr}).

\subsection{Local structure of 
Morse function and its gradient in a \nei~ of
a critical point.}\label{su:lsmfg}

In this subsection $W$ is a cobordism, $\fcob$ 
is a Morse function on $W$, $p$ is a critical 
point of $f$ of index $k$, $v$ is an $f$-gradient.
Recall the
famous Morse Lemma, asserting that $f$
 is equivalent to its Hessian
$d^2f(p)$
locally in a \nei~ of $p$.
The next lemma is a bit sharper version of the Morse
Lemma.

\begin{lemm}\label{morselemma}
Let $L_-, L_+$ be the subspaces of
$T_pW$, \sut~
$d^2f(p)\mid L_-<0,
d^2f(p)\mid L_+>0$,
and $L_-$ is orthogonal to $L_+$ \wrt~ the
bilinear form $d^2f(a)$.

Then there is a chart
$\vphi:U\to V$ at the point $p$, \sut~ $\vphi(p)=0$,
and

\begin{equation}f\circ \vphi^{-1}(x_1,..., x_n)=
f(p)-
\sum_{i=1}^k x_i^2
+
\sum_{i=k+1}^n x_i^2   \lbl{standchart}
\end{equation}
      and
$\vphi'(p)(L_-)=\RRR^k,
\vphi'(p)(L_+)=\RRR^{n-k}$,
where $\RRR^{n-k}$ stands for
the subspace
$\{x_j=0\mid
1\leq j\leq k
\}$, and
$\RRR^k$
stands for
the subspace
$\{x_j=0\mid 
 k+1\leq j\leq n\}$.
\end{lemm}
\Prf  

The Morse lemma provides a chart 
$\Phi$ which satisfies (\ref{standchart})
and $\Phi(p)=0$.
Thus $(f\circ\Phi^{-1})(x)=f(p)+A(x,x)$
where
$A$ is the bilinear form on $\RRR^n$
with the matrix
\begin{equation}
A_{ij}=
\begin{cases}
0, & if~i\not= j\\
1, & if~ i=j\leq k\\
-1, & if~ i=j>k
\end{cases}
\end{equation}

Let $K_+=\Phi'(p)(L_+),
K_-=\Phi'(p)(L_-)$.
It suffices to find a linear isomorphism 
$\xi:\RRR^n\to\RRR^n$
\sut~
$\xi(K_+)=\RRR^{n-k},
\xi(K_-)=\RRR^{k}$.
Note that 
$A|K_+>0, A|K_-<0, A(K_+, K_-)=0$.
From these properties one deduces easily the existence of the isomorphism 
$\xi$ satisfying the properties above. $\qs$

A chart $\vphi:U\to V$, satisfying 
 the conclusions of \ref{morselemma}
is called {\it 
 Morse chart for $f$ at $p$ (with respect 
to $L_+,L_-$)
}.  Morse charts provide a natural tool   for 
the study of 
 the local behavior
of $f$ and $v$ nearby the critical point 
$p$.

Note that the linear subspaces $T_+(p,v), T_-(p,v)\sbs T_pW$
satisfy:
\begin{gather*}
d^2f(p)|T_+(p, v)>0,\quad
d^2f(p)|T_+(p, v)<0,\\
d^2f(p)(T_+(p, v), T_-(p,v))=0
\end{gather*}
(Indeed, the third property is already in the definition of
$f$-\gr. Since
$f|W_{loc}^{st}(p,v)$
has a local maximum in $p$, we have
$d^2f(p)|T_+(p, v)\geq 0$, and using the third property together with
non-degeneracy of $d^2f$
it is not  difficult to obtain the strict positivity of
$d^2f|T_+$.)

Set 
$L_-=T_-(p;v); 
L_+=
T_+(p;v)$, and let
 $\vphi:U\to V\sbs\RRR^n$ be a Morse chart 
for $f$
at $p$
with respect to
$L_+,
L_-$.
In the following definition we develop the notation 
for subsets, maps, etc. in 
$\RRR^n$, related to $\vphi_*(v)$ and to 
$f\circ\vphi^{-1}$ (see the figure on the page
\pageref{fig:morse}).

\begin{defi}\label{morseterm}
Denote $\vphi_*v$ by $\wi v$, $f\circ\vphi^{-1}$ 
by $\wi f$.

We write a point $x\in\RRR^n$ as $(x_+,x_-)$, 
where $x_+\in\RRR^{n-k}$,
and $x_-\in\RRR^k$. 
The orthogonal projection of $\RRR^n $ onto $\RRR^k$
will be denoted by $p_-$.
The orthogonal projection of $\RRR^n $ onto $\RRR^{n-k}$
will be denoted by $p_+$.

Set
\begin{gather*}
Q_r=D^k(0,r)\times D^{n-k}(0,r)  \\
Q_r^\circ=B^k(0,r)\times B^{n-k}(0,r) \\
Q_r^-=Q_r\cap\{(x_+,x_-)\mid \vert x_-\vert\geq 
\vert x_+\vert\} \\
Q_r^+=
Q_r\cap\{(x_+,x_-)\mid \vert x_-\vert\leq 
\vert x_+\vert\}  \\
C=\{(x_+,x_-)\mid \vert x_+\vert=
\vert x_-\vert \} \\
\end{gather*}
The boundary $\pr Q_r$ is the union of two 
subsets:

\begin{gather*}
\pr_-Q_r=\{x\in \pr Q_r\mid \vert x_-\vert =r\} \\
\pr_+ Q_r=\{ x\in \pr Q_r \mid \vert x_+\vert = r\}
\end{gather*}

Set 
\begin{gather*}
\pr_-^\circ Q_r=\pr_- Q_r\sm\pr_+ Q_r\\
\pr_+^\circ Q_r=\pr_+ Q_r\sm\pr_- Q_r
\end{gather*}

Let $W^-$ be the local stable \ma, and
$W^+$ be the local unstable 
\ma~\wrt~a \nei~$U(p)$ of $p$.
We can assume $U(p)\sbs U$.
Denote $\vphi(W^{\pm}(p,v))$ by $N_\pm$.
Since $\wi v$ is an $\wi f$-gradient, we have
$Q_r\cap N_\pm\sbs Q_r^\pm$.
Let $\wi v'(0)=A$. By definition,
$A$ preserves $\RRR^{n-k}$ and $\RRR^k$, and
there is $\nu>0$ \sut~ 
$\langle A x_+,x_+\rangle\geq \nu\langle
x_+,x_+\rangle,
\langle A x_-,x_-\rangle\leq -\nu\langle
 x_-,x_-\rangle$.
We have
$\wi v(x)= Ax+\d(x)$,
where
$\d(x)=o(|x|)$. 
\end{defi}

In what follows we shall  assume 
that the size of the set
$Q_r$ is sufficiently small. The following condition
(R) is satisfied for every $r$ sufficiently small.

$$\lb{R}
\mbox{ (R) } \quad\begin{aligned}
 \bullet~  &  \wi v \mbox{ points inward }    Q_r \mbox{ in the 
points of } \pr_-^\circ Q_r
\mbox{ and outward  } Q_r \mbox{ in  } \pr_+^\circ Q_r\\
\bullet~   &  p_-\mid (N_-\cap Q_r)
\mbox{ is a diffeomorphism onto } 
D^k(0,r)\\
\bullet~   &  p_-\mid (N_+\cap Q_r)
\mbox{ is a diffeomorphism onto } 
D^{n-k}(0,r)\\
 \bullet~  &  \mbox{ for  } x\in Q_r \mbox{ we have } 
 \vert\d(x)\vert\leq 
\frac {\nu}{10}   x\\
 \bullet~  &  Q_r\sbs \vphi(U(p))
\end{aligned}$$

Up to the end of this subsection we assume that $r$ 
satisfies the condition 
(R).
Let $N_+(r), N_-(r)$
be the local stable, resp. local unstable \ma s
of $\wi v|Q_r^\circ$.
We have obviously
$N_\pm(r)\sbs N_\pm\cap Q_r^\circ$.
Actually the inverse inclusion is also true.
Indeed, let
$x\in N_-\cap Q_r^\circ$.
If the \tr~ 
$\g(x,\cdot; \wi v)$
quits 
$Q_r^\circ$,
it must intersect 
$\pr_+^\circ Q_r^\circ$
and then at the moment of the intersection
the value
$\wi f(\g(x,t;\wi v))$
will be strictly greater than  $\wi f(0)$
which contradicts $x\in N_-$.
Thus
$\g(x,t;\wi v)\in Q_r^\circ$
for all $t\geq 0$,
and $x\in N_\pm (r)$.

The next lemma will be useful in the sequel.

\begin{lemm} \label{staysconv}
Assume that a $\wi v$-trajectory
$\g(x,t;\wi v)$ stays in $Q_r$
for all $t\geq 0$. Then
$\g(x,t;\wi v)\rad{t\to\infty} p$.
\end{lemm}

\Prf 
Using the fact that $\wi f(\g(x,t;\wi v))$
is an increasing function of $t$, it is easy to
prove that
$\g(x,t;\wi v)\in Q_r^-$
for every $t\geq 0$.
Write
$$\g'(t)=(\g_+'(t),\g_-'(t))=(A\g_+(t),
A\g_-(t))+ \d(\g(t))$$
So 
$\langle\g_-',\g_-\rangle=
\langle A\g_-, \g_-\rangle + \langle
\d(\g),\g_-\rangle$.
Since $\vert\g\vert\leq 2\vert\g_-\vert$
in $Q_r^-$ we have $\langle\g_-',\g_-\rangle
\leq -\frac {\nu}{2} 
\langle\g_-,\g_-\rangle$. Therefore
$\g_-(t)\to 0$ as $t\to \infty$ and
the same is true for
$\g_+(t)$. $\qs$

\subsection{ Global structure of $v$-trajectories}
\lb{su:gst}

In this subsection $\fcob$ is a Morse function
 on $W$, and $v$ is an $f$-gradient.

We consider here the behavior of
$v$-trajectories on the whole of
$W$, not only in a \nei~ of a critical point.
Since $v\mid \pr W$ is transversal to $\pr W$
and points inward $W$ on $\pr_0 W$ and outward
on $\pr_1 W$, there are 4 possible types of behavior
for an  integral curve
of $v$ ($a,b$ in the formulas below are
finite real numbers):
\begin{enumerate}
\item
$\g$ is defined on $[a,\infty[,~~
\g(a)\in\pr_0 W, ~~ \g(]a,\infty[)\sbs \Wkr$

\item
$\g$ is defined on $]-\infty,a],~~
\g(a)\in\pr_1 W, ~~\g(]-\infty, a[)\sbs \Wkr$

\item
$\g$ is defined on $[a,b],~~
\g(a)\in\pr_0 W, ~~ \g(b)\in\pr_1 W, ~~\g(]a,b[)
\sbs
\Wkr$

\item
$\g$ is defined on $]-\infty,\infty[,~~
 \g(]-\infty,\infty[)\sbs \Wkr$
\end{enumerate}

(see, e.g., \cite{pasur}, Prop. 8.9.) Since
$v$ is a
gradient of a Morse function we can say
more about the behavior of $v$-trajectories.
First we
study an integral curve in a \nei~ of a
critical point.
In the next lemma $p$ is a critical point of $f$, 
and
$r$ is sufficiently small so that (R)
 holds. Set $I_r=\{t\in\RRR\mid\g(t)\in\vphi^{-1}(Q_r)\}$.

\begin{lemm} 
One of the
following possibilities holds:
\begin{enumerate}
\item $I_r=[a,\infty[, \quad \g(t)\rad{t\to\infty} p$,
\item $I_r=]-\infty,a], \quad \g(t)\rad{t\to-\infty} p$,
\item $I_r=[a,b]$ \quad (where $a=b$  not excluded)
\item $I_r=]-\infty,\infty[,\quad \g(t)=p$ for every $t$.
\item $I_r=\emptyset$
\end{enumerate}
\end{lemm}

\Prf Assume that $I_r\not=\emptyset$, and
let $\tau_0\in I_r$. Consider the connected
component
$J$ of $\tau_0$ in $I_r$; then $J$ is a
closed
interval, finite or
infinite from one or both sides.
The proof goes by considering subsequently
these
cases. We do three of them and leave the
rest to the reader.
\begin{enumerate}\item
$J=\RRR$. Then
$\Imm \g\sbs U(p)\cap (W_{loc}^{st}\cap
W_{loc}^{un})=\{p\}$
(by Lemma \ref{staysconv})
and we obtain the case $(4)$.

\item $J=]-\infty,\a]$. 
Denote $\vphi(\g(t))$ by $\wi\g(t)$.
Then $\wi\g(t)\in Q_r$
for every
$t\in]-\infty,\a]$, and by
Lemma \ref{staysconv} we have $\g(t)\in W_{loc}^{un}$.
Note that
$\wi\g(\a)\in\pr Q_r$,
and since $\wi v$ points inward $Q_r$ on
$\pr_-^\circ Q_r$, we have
$\wi\g(\a)\in \pr_+ Q_r$. Therefore
$f(\g(\a))\geq f(p)$, and
for every $t>\a$ we have $f(\g(t))>f(p)$.
Moreover for $t>\a$ sufficiently close to $\a$,
we have
$\wi\g(t)\notin Q_r$. Assume now that
$I_r\not= J$, let
$t_0$ be
the least number $>\a$, \sut~~  $\wi\g(t)\in Q_r$.
Then $\wi\g(t_0)\notin \pr_+^\circ Q_r$,
therefore
$\wi\g(t_0)\in \pr_- Q_r$, which is impossible since
$\wi f\mid \pr_- Q_r \leq f(p)$. Thus $J=I_r$.

\item $J=[\a,\b]$. The same reasoning as in the
preceding point,
proves that
$\wi\g(\b)\in \pr^+ Q_r$,
and that
$\wi\g(\tau)\notin Q_r$ for $\tau >\b$.
 By symmetry
 we have:
$\wi\g(\a)\in \pr_- Q_r$ and
$\g(\tau)\notin Q_r$
for $\tau<\a$.

\end{enumerate}
$\qs$

From this lemma it follows easily that for a
non-constant
$v$-trajectory $\g$
the four possibilities can occur:
\begin{enumerate}
\item
$\g$ is defined on $[a,\infty[,
\g(a)\in\pr_0 W, \g(t)\rad{t\to\infty} p$, where
$p\in S(f)$.

\item
$\g$ is defined on $]-\infty, a[, ~~ 
\g(a)\in\pr_1 W, ~~ \g(t)\rad{t\to-\infty} q$, where
$q\in S(f)$.

\item
$\g$ is defined on $[a,b], ~~ 
\g(a)\in\pr_0 W, ~~ \g(b)\in\pr_1 W, ~~ \g(]a,b[)
\sbs \Wkr$

\item
$\g$ is defined on $]-\infty,\infty[,~~ 
 \g(t)\rad{t\to\infty} p$, where $p\in S(f)$, and
$\g(t)\rad{t\to-\infty} q$, where $q\in S(f)$.
\end{enumerate}

\begin{defi} Let $p$ be a critical point of $f$.
We denote by $D(p,v)$ the set of all $x\in W$, 
\sut~~ 
$\lim_{t\to\infty}\g(x,t;v)=p$, and we call
it {\it descending disc}.

We denote by $D(p,-v)$ the set of all $x\in W$, 
\sut~
\break
$\lim_{t\to\infty}\g(x,t;-v)=p$, and we call
it {\it ascending disc}.

\end{defi}

{\bf Remarks.} 
\been\item
$D(p,v)$ is an analog of the global stable
manifold
in the dynamical system theory. Note, however, that
the flow
generated by $v$ is not defined on the whole of $\RRR$
in the case
$\pr W\not=\emptyset$ ( this is one of the
 reasons to prefer this notation to the
notation $W^{st}(p)$, used in the literature on
dynamical systems).
\item
By the Grobman-Hartman theorem the behavior of the
$v$-trajectories belonging to
$D(p,v)$ nearby $p$
is {\it topologically} equivalent to that of
the solutions of the standard system
$x'=-x$ in $\RRR^n$.
This is not the case in the differentiable category,
and the behavior
of trajectories nearby $p$
can be complicated.
The following  lemmas show nevertheless that globally
$D(p,v)$
is a good geometrical object. $\qt$
\enen

In the next three lemmas $p$ is  a 
critical point of $f$.
 In these lemmas we use the terminology 
from Definition
\ref{morseterm} and we assume that $r$ 
satisfies (R).

\begin{lemm}\label{locstabglob}
$D(p,v)\cap \vphi^{-1}(Q_r)=
W_{loc}^{st}(p,v)\cap\vphi^{-1}(Q_r)$
\end{lemm}
\Prf     Let $\wi x=\vphi(x)\in Q_r$. If
$\g(\wi x,\cdot;\wi v)$
stays forever in $Q_r$, then
$x\in W_{loc}^{st}(p,v)$ 
by Lemma \ref{staysconv}. 
If
$\g(\wi x,\cdot;\wi v)$
leaves $Q_r$, then it must intersect
$\pr_+ Q_r$
(since in $\pr_-^\circ Q_r\sm C$
the vector field $\vphi_*(v)$ points 
inward $Q_r$)
say at a moment $t_0$.
Then for $t>t_0$ we have:
$f(\g(x,t;v))>f(p)$
and $x\notin D(p,v)$. $\qs$

Denote $f^{-1}(\l)$
by $V_\l$ and $f^{-1}([\l,\m])$ by $W_{[\l,\m]}$.

\begin{lemm}\label{stabdiscmanif}
We have:
\begin{enumerate}\item
$D(p,v)\cap\Wkr$ is a submanifold of $\Wkr$ 
of dimension
$\ind p$.
\item
$D(p,v)\cap V_\l$ is a submanifold of
 $V_\l$
of dimension $\ind p-1$.
\end{enumerate}
\end{lemm}

\Prf      The second part  follows from the first, since
$D(p,v)$ is transversal to 
$V_\l$.
To prove 1) note that it suffices to prove 
that
$D(p,v)\cap \Wkr$
is locally a submanifold of $\Wkr$ in a 
\nei~ of a
point
$x\in\vphi^{-1}(Q_\e)$, where $\vphi$ is a 
Morse
chart. But this follows
immediately from the preceding Lemma.
$\qs$

The following Lemma justifies the using
the term
{\it descending disc}
for $D(p,v)$.

\begin{lemm}\label{discisdisc}
Let $\l\in[a,f(p)]$ be a regular value of $f$.
Assume that for every $x\in D(p,v)\cap W_{[\l,b]}$
the $(-v)$-trajectory reaches $V_\l$.
Then the pair
$(D(p,v), D(p,v)\cap V_\l)$
is diffeomorphic to
$(D^k,\pr D^k)$, where $k=\ind p$.
\end{lemm}
\Prf For 
$0\leq \rho\leq r$ set 

$$\D_\r=N_-\cap\{(x_+,x_-)\mid \vert x_-
\vert\leq \r\},
\quad$$
$$\Sigma_\r=\pr\D_\r= N_-\cap\{(x_+,x_-)\mid 
\vert
x_-\vert= \r\} $$
$$A(\r,r)= N_-\cap\{(x_+,x_-)\mid \r\leq\vert
 x_-\vert\leq
r\}$$

For $0<\r<r$ we have:
$\D_\r$ is diffeomorphic to a $k$-disc, $A(\r,r)$
is diffeomorphic to the annulus
$S^{k-1}\times [0,1]$.
Each $(-v)$-trajectory, starting from a point
$x\in A(\r,r)$
reaches $\Sigma_r$ at a moment $t(x)$ and reaches
$V_\l$ at a moment $T(x)>t(x)$, both $t(x)$ and $T(x)$
being $\smo$ functions of $x$.
It is easy now to construct a diffeomorphism
$\Phi:\D_r\to D(p,v)\cap V_{[\l, f(p)]}$. Namely,
set $\Phi\mid {\D_\r}=\id$ and to define $\Phi$ outside
$\D_\r$ accelerate the movement along each $(-v)$-trajectory
$\g(x,t;-v)$ for $x\in\Sigma_\r$. We leave the details
to the reader.
$\qs$

\subsection{Thickenings of descending discs and compactness
properties}
\label{su:tddcp}

Let
$\fcob$ be a Morse function on a riemannian cobordism $W$ of dimension
$n$, and
$v$ be an $f$-gradient.

\begin{defi}\label{d:dthick}
For $p\in S(f)$ set

\begin{gather*}
B_\delta(p,v)=\{x\in M~\vert~ \exists
t\geq 0 :
 \gamma (x,t;v)\subset B_\delta (p)\}\\
D_\delta(p,v)=\{x\in M~\vert ~\exists
t\geq 0 :
\gamma (x,t;v)\subset D_\delta (p)\}, 
\end{gather*}

Recall from \sub~ \ref{su:tc}
that in the notation $B_\d(p), D_\d(p)$ it is implicitly
understood
that $\d$ is sufficiently small.

We denote by $D_\d(v)$, resp. by $B_\d(v)$
the union of all the $D_\d(p,v)$, resp. $B_\d(p,v)$.
We denote by $D(v)$ the union of all the descending
discs.

For $s\in \NNN, 0\leq s\leq n$
we denote by $D (\indl s;v)$
the union of   $
D(p,v)$ where $p$ ranges over critical
points of $f$ of index $\leq s$.
We denote by $B_\d (\indl s;v)$, resp.
by $D_\d(\indl s;v)$
the union of
$B_\d(p,v)$, resp. of
$D_\d(p,v)$, where $p$ ranges over critical
points of $f$ of index $\leq s$. 
\end{defi}

\bere
We shall use
similar
notation like $D_\d(\inde s;v)$ or
$B_\d(\indg s;v)$, which are now clear
without special definition. 
\enre

\begin{lemm}\label{bopen}
$B_\d(p,v)$ is open.
\end{lemm}
\Prf Follows immediately from the continuous dependence of the integral 
trajectories of vector fields on their initial values.
$\qs$
\pa

The things are not that simple for $D_\d(p,v)$ or for
$D(p,v)$.
 These sets are not necessarily closed.
The general reason for that is the existence of $v$-trajectories,
joining critical points:
if there is a $(-v)$-trajectory joining two critical 
points
$p$ and $q$,
then clearly
$q\in\overline{D(p,v)}\sm D(p,v)$.
But for the union of  all the $D(p,v)$ the compactness
property is valid:

\begin{lemm}\label{dcomp}
\begin{enumerate}
\item $D(v)$ is compact.
\item $D_\d(v)$ is compact
\item $D_\d(v)=\overline{B_\d(v)}$.
\end{enumerate}
\end{lemm}
\Prf 1) The complement $W\sm D(v)$ consists of all points
$x\in W$, \sut~~ $\g(x,\cdot;v)$ reaches $\pr_1 W$.
This set is open by Corollary \ref{co:reach}\footnote
{We shall use in this chapter the results of Section \ref{s:pcic}.
The results of the section \ref{s:pcic}
makes no use of the results of the preceding sections; we
placed it in the end in order not to interrupt the exposition.}.

2) similar to 1).

3) It follows from the  inclusions
$B_\d(v)\sbs D_\d(v)\sbs\overline{B_\d(v)}$
$\qs$

As $\d$ goes to $0$, these $\d$-thickenings, for example,
$B_\d(p,v)$,
 should collapse to their
"skeleton"
$D(p,v)$, that is, $B_\d(p,v)$ should form a fundamental
system of
neighborhoods for $D(p,v)$.

Again in general it is not true, 
but the analog of this property is true for
the union of all
$D(p,v)$;
that is the subject of the following proposition.

\begin{prop}\label{fundsys}
\begin{enumerate}
\item
For every $\theta>0$ we have:
$D_\theta(v)=\cap_{\d>\theta}B_\d(v)$.
\item
$D(v)=\cap_{\d>0} B_\d(v)$.
\item $\{B_\d(v)\}_{\d>\theta}$ form a fundamental system of \nei s
for
$D_\theta(v)$.
\item $\{B_\d(v)\}_{\d>0}$ form a fundamental system of \nei s
for
$D(v)$.

\end{enumerate}
\end{prop}
\Prf 1) A point $x\in W$ {\it does not}
belong to
$D_\theta (v)$
if and only if
the trajectory $\g(x,\cdot;v)$ reaches $\pr_1 W$ without
 intersecting the discs $D_\theta(p), p\in S(f)$.
Then it does not intersect the discs $D_\d(p)$ for
some $\d>\theta$. Similar argument proves 2).

3) Denote $B_\d(v)$ by $R_\d$, and $D_\theta(v)$ by
$K$
(here $\d$ is in some interval $I_=]\theta,\d_0[$).
Then:
\begin{enumerate}
\item $\overline{R_\d}\sbs R_{\d'}$, if $\d<\d'$.
\item $\cap_{\d>\theta} R_\d=K$.
\item $K$ and $W$ are compacts.
\end{enumerate}

These three properties imply in general, that
$\{R_\d\}_{\d>\theta}$
form a fundamental system
of neighborhoods
of $K$. This is an exercise from an elementary general
topology, which we leave to
the reader ( see also Lemma \ref{gfimplf}, where a similar argument 
is exposed).
$\qs$

\subsection{Transversality conditions and Rearrangement
Lemma}
\lb{su:tcrl}

In this section $\fcob$ is a Morse function on a cobordism
$W$ and $v$ is
an $f$-gradient.
\begin{defi} We say that
 $v$ satisfies
           {\it Transversality Condition }, if 
$$\big(
x,y \in S(f) \big)
\Rightarrow \big( D(x,v)\cap \Wkr \pitchfork D(y,-v)\cap
\Wkr \big)$$

\vskip0.1in

 We say that 
 $v$ satisfies {\it Almost Transversality Condition},
if
   $$\big( x,y \in S(f) ~\&~
\ind x \leq \ind y \big) \Rightarrow \big(
D(x,v)\cap\Wkr \pitchfork D(y,-v)\cap\Wkr\big)$$
\vskip0.1in
\end{defi}

One of the basic instruments in the Smale's theory
of handlebodies is the construction
from any Morse function on a cobordism $W$
a so-called self-indexing Morse function, that is a
function, \sut
~$f(p)=\ind p$ for every $p\in S(f)$.
The details of this construction can be found in
\cite{milnhcob}, \S 4. 
The construction actually gives a bit more:
if we start with an arbitrary \Mf~ $g$ and a \glvf~ $v$
for $g$, satisfying  \TA, then it is possible to find a 
self-indexing \Mf~ $f$ on the cobordism, so that $v$ is also
a $f$-gradient.
The same is true for $f$-gradients; the corresponding
Rearrangement Lemma is the main subject of the present
subsection.

\begin{defi} 
\label{d:adjust}
A Morse function $\phi:W\to[\a,\beta]$
is
called {\it adjusted to $(f,v)$}, if:

 1) $S(\phi)=S(f)$,

 2) the function $f-\phi$ is constant in a \nei~
of $V_0$, in a \nei~ of $V_1$, and in a \nei~ of each
point of
$S(f)$,

3) $v$ is also a $\phi$-\gr.
\end{defi}

\begin{prop}\label{rearrlemm}
Assume that $S(f)=\{p,q\}$, and that
$D(p,\pm v)\cap D(q,\mp v)=\ems$.
Let $\a,\b\in]a,b[$.
Then there is a Morse function
$g:W\to[a,b]$, adjusted to $(f,v)$, \sut~  
~$g(p)=\a, g(q)=\b$.
\end{prop}
\Prf Note first, that we can assume that $f(p)\not=
f(q)$
(this can be always achieved by a small perturbation
of $f$ with support
in a \nei~of $S(f)$).
Assume, for example, that $f(p)<f(q)$.
Denote
$D_\d(p,v)\cup D_\d(p,-v)$ by $\DD_\d(p)$
and $B_\d(p,v)\cup B_\d(p,-v)$ by $\BB_\d(p)$.
Similarly we define
$ \DD_\d(q)$ and
$ \BB_\d(q)$.
Note, that $\BB_\d(q)$ and $\BB_\d(p)$ are open.

\begin{lemm}
For $\d>0$ \sufsm~we have:
$\BB_\d(p)\cap \BB_\d(q)=\ems$.
\end{lemm}
\Prf Choose any $\l$ \sut~      $f(p)<\l<f(q)$.
Denote the cobordism
$f^{-1}([a,\l])$ by $W_0$, and the cobordism
$f^{-1}([\l,b])$ by $W_1$. Applying to $W_0$ and $W_1$
Proposition \ref{fundsys}
we deduce from $D(p,-v)\cap D(q,v)=\ems$
the following:

For $\d>0$ \sufsm~ the sets
$D_\d(p,-v)\cap f^{-1}(\l)$ and
$D_\d(q,v)\cap f^{-1}(\l)$ are disjoint. This implies
the Lemma.
$\qs$

It follows from the Lemma  that for $\d>0$ \sufsm~
the sets $\DD_\d(p), \DD_\d(q)$ are disjoint and compact.
Denote their intersections with
$\pr_0 W$ by
$K_p$, resp. $K_q$. Choose a $\smo$ function
$\a:\pr_0 W\to [0,1]$, \sut~
$\a\mid K_p =1,
 \a\mid K_q = 0$.
Note, that $W'=W\sm (D(v)\cup D(-v))$
is the union of all the $v$-trajectories, starting
at a point of
$\pr_0 W\sm D(v)$, and reaching $\pr_1 W$. Therefore
there is a unique extension $A$ of $\a$ to
$W'$, \sut~$A$ is constant on each
such trajectory. This function is smooth on $W'$. Remark
that
$W=W'\cup \BB_\d(p)\cup \BB_\d(q)$ and it follows
from the definition that
$A\mid W'\cap \BB_\d(p)=1,
A\mid W'\cap \BB_\d(q)=0$.
Therefore $A$ extends to $\smo$ function on $W$, which
is constant on every $v$-trajectory.
The rest of the proof is standard (see \cite{milnhcob},
\S 4). We reproduce it briefly for the sake of completeness.

One constructs first an auxiliary real $\smo$ function
of two real arguments
$G:[a,b]\times [0,1]\to[a,b]$, which will be considered
as a real function
$[a,b]\to[a,b]$ depending on the real parameter $s\in
[0,1]$.
This function must satisfy the following conditions:
\begin{enumerate}
\item
$\frac{\pr G_s(x)}{\pr x}
>0$ for every $s,x$.
\item
$G_s(x)=x$ for $x$ close to $a$ or $b$ and for every $s$.
\item $G_0(f(q))=\b$ and $\frac {\pr G_0(x)}{\pr x}=1$
in a \nei~of $x=f(q)$.
\item $G_1(f(p))=\a$ and $\frac {\pr G_1(x)}{\pr x}=1$
in a \nei~of $x=f(p)$.
\end{enumerate}

Then the function $x\mapsto G_{A(x)}(f(x))$
satisfy the conclusion of our Proposition.
$\qs$

\begin{coro}\label{order}
Assume that $v$ satisfies the almost transversality condition.
Let $\prec$ be any partial order on the set $S(f)$, \sut
~for every $p,q\in S(f)$
with $\ind p<\ind q$ we have:
$p\prec q$.
Then there is a Morse function $q: W\to[a,b]$,
 adjusted to $(f,v)$, \sut
~$f(p)<f(q)\Leftrightarrow p\prec q$. $\qs$
\end{coro}
\begin{coro}
Assume that $v$ satisfies the almost transversality condition.
Then for every $s$ the set
$D(\indl s,v)$ is compact.
\end{coro}
\Prf By the preceding Corollary there is a function 
$\phi:W\to[a,b]$ and a regular value
$\l$ of $\phi$, \sut
~$\ind p\leq s\Leftrightarrow \phi(p)<\l$.
Now apply  Lemma \ref{dcomp}
to the cobordism $\phi^{-1}([a,\l])$. $\qs$

The preceding Corollary has an analog for $\d$-thickenings of
descending discs.
To formulate it we need two more definitions.

\begin{defi}
We say that $f$ is {\it ordered}
with an {\it ordering sequence} $a_0,...,a_{n+1}$, if
$a=a_0<a_1<...<a_{n+1}=b$ are regular values of $f$ 
\sut~ $S_i(f)\subset
 f^{-1}(]a_i,a_{i+1}[)$. 
\end{defi}

\begin{defi}
We say that $v$ is
{\it
$\delta$-ordered
} (where $\d>0$),
if there is an ordered Morse function
$\phi:W\to[a,b]$, adjusted to $(f,v)$,
with an ordering sequence
$a_0,...,a_{n+1}$, such that
\begin{equation}\forall p\in S_k(f):\quad
D_\d(p)\sbs\phi^{-1}(]a_k,a_{k+1}[)  \lbl{ordergrad}
\end{equation}
\end{defi}

\begin{rema}\label{orderedplus}
 If $f$-gradient $v$ is $\d$-ordered,
then it is $\d'$-ordered for some $\d'>\d$. Every 
$f$-gradient, satisfying
\ata, is $\d$-ordered for some $\d$ (it
 follows from Corollary \ref{order}).
\end{rema}

{\bf Warning:} Being $\d$-ordered does not imply \ata.

\begin{prop}\label{fundsysord}
If $v$ is  $\e$-ordered then
$\forall\d : 
0< \d\leq \e$ and
$\forall
s: 0\leq s\leq n$
we have:
\begin{enumerate}
\item
$
D_\delta (\indl s;v)~ \text{~~is compact}~~;$
\item $
D (\indl s;v) ~ \text{~~is compact}~~;$
\item $
{\bigcap} _     {\theta>0} B_\theta (\indl s;v) =
D (\indl s;v);$
\item $
{\bigcap} _{\theta>\delta} B_\theta (\indl s;v)
= D_\delta (\indl s;v);$
\item $
\overline{B_\delta (\indl s;v)} = D_\delta (\indl s;v).$
\end{enumerate}
\end{prop}
\Prf Let $\phi:W\to[a,b]$ be the ordered \Mf~, 
satisfying
(\ref{ordergrad}). Apply Lemmas \ref{fundsys} and \ref{dcomp}
to the cobordism $\phi^{-1}([a_0,a_{s+1}])$. $\qs$

\subsection{Adding a horizontal component to a vector field 
nearby the boundary}
\label{su:ahc}

In this subsection we describe a (fairly standard) construction:
modifying a vector field on a cobordism by adding a horizontal component
near the boundary. See for example \cite{milnhcob}, page 62.
The aim of this subsection is to give explicit estimates and
 to derive some 
properties  
for the resulting vector fields.

Let $W$ be a cobordism, $v$ be a vector field on $W$, 
transversal to $\pr_0 W$ (that is, for $x\in\pr_0 W$
 we have $v(x)\notin T_x(W)$).
We shall expose here the procedure of changing $v$ nearby $\pr_0 W$. 
( If $v$ is transversal to $\pr_1 W$,  the similar
 procedure can be performed nearby 
$\pr_1 W$; we leave the statement of results to the reader.)

Let $a>0$; assume that for every 
$x\in\pr_0 W$ the $v$-trajectory $\g(x,\cdot;v)$ is defined on 
$[0,a]$ and
$\g(x,[0,a];v)\sbs W\sm\pr_1 W$.
The set $\g(\pr_0 W, [0,a];v)$ will be called
{\it closed collar} of $\pr_0 W$
of height $a$ \wrt~ $v$, and denoted by
$C_0$.
Consider a map
$\Psi:\pr_0 W\times[0,a]\to W$
defined by
$\Psi(x,t)=\g(x,t;v)$.
Then $\Psi$ is a diffeomorphism on its image $C_0$.
The tangent space $T_{(x,t)}\(\pr_0 W\times[0,a]\)$
being naturally isomorphic to
$T_x(\pr_0 W)\oplus \RRR$,
any vector field $\eta$
on $\pr_0 W\times [0,a]$
can be written  as follows:
$\eta=\(\eta_1(x,t), \eta_2(x,t)\)$
where $\eta_1$ is a vector field on $\pr_0 W$
(it will be referred to as {\it the first component} of $\eta$)
and $\eta_2$ is a function $\pr_0 W\to\RRR$
(it will be referred to as {\it the second component} of $\eta$).
In this notation we have:
$\(\Psi^{-1}\)_*(v)=(0,1)$.

Let $\xi$ be a vector field on $\pr_0 W$, and let
$h:[0,a]\to\RRR$ be a positive 
$\smo$ function \sut~ 
$h(t)=0$ for  $t$ in a \nei~ of  $0$ and in a \nei~
of  $a$.
Define a vector field $v'$ on $\pr_0 W\times [0,a]$
by
$v'(x,t)=\(h(t)\xi(x), 1\)$
and define a vector field $u$ on $W$
by the following formula:
\begin{gather*}
u=v\mx{ in } W\sm C_0;\\
u=\Psi_*(v')\mx{ in } C_0
\end{gather*}

We shall call this vector field $u$
{\it the result of adding to $v$ a  horizontal component 
proportional to
$\xi$
nearby $\pr_0 W$} (or simply the result of AHC-construction).
The triple
$(a,\xi,h)$
will be called {\it input data} for AHC-construction.
We shall also write $u=AHC(a,\xi,h;v)$.
We shall establish three basic properties of this construction.
They are all easily checked, and we leave the proofs to the reader .

\pa
{\it AHC1. Estimate of $u-v$}\label{ahc1}
\pa
Let $E$ be the expansion constant of $\Psi$.
We have:
$$\Vert u-v \Vert\leq E\Vert \xi\Vert\Vert h\Vert$$
(where $\Vert\cdot\Vert$ is  the $C^0$ norm ).
Moreover, the map $(\xi,h)\mapsto u$ is continuous \wrt~
$\smo$ topology.
\pa
{\it AHC2. The behavior of $(-u)$-trajectories}\label{ahc2}
\pa
Let $\g$ be a $(-v)$-trajectory starting at 
$x_0\in W\sm C_0$, and such that 
$y=\g(x,t_0;-u)\in \pr_0 W$.
Then the $(-u)$-trajectory starting at $x_0$
reaches $\pr_0 W$ at the moment $t_0$
and 
$\g(x_0,t_0;-u)=\Phi(-\xi, T_0)(y)$,
where
$T_0=\int_0^a h(\tau)d\tau$.
Thus

$$
\stind {(-u)}ba = \Phi(-\xi, T_0)\circ\stind {(-v)}ba 
$$
\pa
{\it AHC3. Gradient property. } 
\pa
Assume that $v$ is an $f$-gradient, where $\fcob$
is a Morse function on $W$.
Set 
$m=\min_{x\in C_0} df(v)(x), 
 \quad M=\max_{x\in C_0}\Vert df(x)\Vert$.
If $\Vert \xi\Vert\Vert h\Vert<\frac m{EM}$
then $u$ is also an $f$-gradient.

\section{$s$-submanifolds and $gfn$-systems
}\label{s:ssags}
\pa
This section is an extended version 
of \S 2.1, 2.2 of \cite{pastpet}.

\subsection{ $s$-submanifolds }\label{su:ssubm}

Here we define and study $s$-submanifolds.
The basic example of an $s$-\sma~is the set of all the 
descending discs of an $f$-gradient $v$,
satisfying the almost transversality condition 
(where $f:M\to\RRR$ is a Morse function
on a closed \ma.)
The notion of $s$-submanifold of a \ma ~$M$ is close to the notion
of stratified \ma ~in the sense of Whitney-Thom.
We have chosen to introduce the notion of $s$-\sma ~and 
to deal with these objects, 
rather than to stick to the classical definition, 
since $s$-\sma s carry less structure, and are more flexible
 than the classical stratified manifolds.
It seems that if we impose additional restriction on $v$, these 
$s$-manifolds
are actually Whitney-stratified manifolds. We hope to return to 
these questions in further publications.

Proceeding to the precise definitions, let
$\aa=
\{A_0,...,A_k\}$
be a finite sequence of subsets of
a topological space $X$.
We denote  $A_s$ also by $\aa_{(s)}$, and
we denote  $A_0\cup...\cup A_s$
by $A_{\leq s}$ and
also by $\aa_{(\leq s)}$.
We say that $\aa$ is a
{\it compact family }
if  $\aa_{(\leq s)}$
is compact for every $s$.

\begin{defi}\label{ssubman}
Let $M$ be a manifold without boundary.
 A finite sequence $\xx =  \{X_0,...,X_k\} $
 of subsets of $M$ is called  {\it
$s$-submanifold of $M$ }
($s$ for stratified) if
\begin{enumerate}
\item  $\xx$ is
disjoint and each
 $X_i$ is a submanifold of $M$ of dimension $i$
with the trivial normal bundle.
\footnote{ We impose the triviality of the normal bundle since 
it makes easier
some of the proofs, see for example the proof of the Proposition 
\ref{transvssub}.}
\item $\xx$ is a compact family. \quad$\triangle$
\end{enumerate}
\end{defi}
For an $s$-submanifold $\xx=\{X_1,...,X_s\}$
the largest $k$ \sut $X_k\not=\emptyset$
is called {\it dimension } of $\xx$,
and denoted  by $\dim \xx$.
 For a diffeomorphism $\Phi :M\to N$ and
 an $s$-submanifold
 $\xx$ of $M$ we denote by $\Phi (\xx)$
the $s$-submanifold
of $N$ defined by $\Phi (\xx)_{(i)} =
\Phi(\xx_{(i)})$.
\vskip0.1in
If $V$ is a submanifold of $M$ and $\xx$
 is an $s$-submanifold of $M$,
then we say that $V$ {\it is transversal
to $\xx$ }
(notation: $V\pitchfork \xx$) if
 $V\pitchfork \xx_{(i)}$ for each $i$.
If $V$ is a compact submanifold of $M$,
 transversal to an
 $s$-submanifold $\xx$, then the family
$\{\xx_{(i)}\cap V\}$
is an $s$-submanifold of $V$
 which will be denoted
 by $\xx \cap V$.
\vskip0.1in
If $\xx , \yy$ are two $s$-submanifolds
of $M$, we say, that $\xx$
 is {\it transversal to $\yy$ }
 (notation: $\xx \pitchfork \yy$) if
$\xx_{(i)}\pitchfork \yy_{(j)}$ for every
$i,j$; we say that $\xx$
is {\it almost transversal to }
   $\yy$ (notation: $\xx\nmid\yy$) if
$\xx_{(i)}\pitchfork \yy_{(j)}$ for
$i+j<\dim M$. 

\begin{rema}\label{transvinters}
Note, that
$\xx\nmid\yy$ if and only if $ \xx _{(\leq i)}
\cap  \yy _{(\leq j)} =
 \emptyset$
whenever $i+j<\dim M$.
\end{rema}

Let $\fcob$ be a Morse function
on a \cob~$W$, $v$ be an $f$-\gr. 
We denote by $\dd(v)$
the family $\{D(\inde i;v)\}_{0\leq i \leq \dim W}$
of subsets of $W$. For $\l\in\RRR$ denote by 
 $\dd_\l(v)$ the family
$\{D(\inde {i+1};v)\cap f^{-1}(\l)\}_{0\leq i\leq \dim W-1}$
of subsets of $V_\l=f^{-1}(\l)$.

\begin{lemm}
Assume that $v$ satisfies \ata. Then:
\begin{enumerate}
\item $\dd(v)$ is a compact family.
\item If $M$ is a closed manifold, then
$\dd(v)$ is an $s$-submanifold, transversal to $f^{-1}(\l)$
for every regular value $\l$ of $f$.
\item If $M$ is any cobordism, then for every regular value
$\l$ of $f$ the family $\dd_\l(v)$ is an $s$-submanifold 
of $f^{-1}(\l)$.
\end{enumerate}
\end{lemm}
{\it Proof.\quad}
1) follows immediately from 
Proposition \ref{dcomp}. Note that the normal bundle
to $D(p,v)$ is always trivial, which implies 2) and 3). \quad$\square$

We proceed to an analog for $s$-submanifolds of the Thom theorem
on density of transversality property.
In the proposition \ref{transvssub}
we shall prove that
 if $\xx,\yy$ are two $s$-submanifolds, it is always possible to
 perturb a little one of them
say, $\xx$, so that the result is almost transversal to $\yy$.
First we 
need  
a definition.

Let $M$ be a \ma~without boundary. Denote by $\vem$ the subset of 
$\Vect(M\times [0,1])$, formed by
the vector fields
$v\in \Vect(M\times [0,1])$, satisfying
\begin{enumerate}
\item $v$ vanishes together with all partial derivatives
in 
$M\times\{0\} \cup M\times\{1\}$.

\item The coordinate of $v$, corresponding to the factor $[0,1]$ 
is equal to zero.
\end{enumerate}

One can consider the elements of $\vem$
as the $C^\infty$ vector fields on $M$, depending smoothly
on the parameter $t\in[0,1]$. There is
a natural topology on 
$\Vect(M\times[0,1]$, and $\vem$
is a closed subspace of
$\Vect(M\times [0,1])$
(see, e.g., \cite{pasur},\S 8).
From now on and to the end of the subsection we assume that $M$
is compact. In this case the space 
$\Vect(M\times [0,1])$ as well as the space
$\vem$ has a $sup$-norm
(with respect to a riemannian metric),
which induces the {\it  $C^0$ topology } on this space.

For every $v\in \vem$ the integral curve $\g(x,t;v)$
is defined for all $t\in\RRR$. Recall that 
$\g(x,t;v)$ is also denoted by 
$\Phi(v,t)(x)$, and that
for a fixed $t$ the map
$x\mapsto \Phi(v,t)(x)$
is a diffeomorphism of $M$.

\begin{prop}\label{transvssub}
Let $\xx,\yy$ be $s$-submanifolds of a closed \ma~$M$.
Then the set of $v\in\vem$, \sut~
$\Phi(v,1)(\xx)\nmid\yy$
is open and dense in $\vem$.
\end{prop}

\Prf We shall give a sketch of the proof. A reader familiar 
with the differential topology will easily reconstruct the details.
Denote by $\TT(\xx,\yy)$ the subset of all $v\in\vem$,
\sut~
$\Phi(v,1)(\xx)\nmid\yy$.
Note first that 
$\TT(\xx,\yy)$ 
is open  in $C^0$-topology,
which follows easily from the results of \S 5 and 
the Remark
\ref{transvinters}.
To prove its $C^\infty$-density we proceed by induction.
 Consider two
auxiliary assertions:

$A(n)$:\quad $\TT(\xx,\yy)$ is dense for every $\xx$ and $\yy$ with
$\dim\xx+\dim\yy\leq n$.

$A_0(n)$:\quad For every $\xx,\yy$ with $\dim\xx+\dim\yy\leq n$
the closure of 
$\TT(\xx,\yy)$ 
contains $0$.

It is easy to show that 
$A_0(n)\Rightarrow A(n)$.
Therefore it is sufficient to show that 
$A(n-1)\Rightarrow A_0(n)$.
For this, let 
$\xx=\{X_0,...,X_k\},
\yy=\{Y_0,...,Y_l\}$
be $s$-submanifolds of $M$,
$\dim\xx=k,\dim\yy=l, k+l=n$.
Denote by $\xx',\yy'$ the $s$-submanifolds 
$\{X_0,...,X_{k-1}\}$,
resp.
$\{Y_0,...,Y_{l-1}\}$.
By the induction assumption there is a small
vector field 
$w\in\vem$, \sut~
$\Phi(w,1)(\xx')\nmid\yy$
and
$\Phi(w,1)(\xx)\nmid\yy'$.
The manifolds which are not yet in general position
are
$X_k'=\Phi(w,1)(X_k)$ and $Y_l$. It suffices to consider 
the case $k+l<\dim M$. 
Consider three disjoint compacts:
$\xx'_{(\leq k-1)}, X_k'\cap Y_l, \yy_{(\leq l-1)}$.The compact 
$X_k'\cap Y_l$ is in $X_k$.
Recall that the normal bundle to $X_k'$ is trivial. Let 
$U_0$ be a \nei~ of $X_k'\cap Y_l$ in $X_k'$.
Choose a tubular \nei~
$U$ of $U_0$, fiberwise diffeomorphic to the product
$U_0\times\RRR^{n-k}$, \sut
~ $\overline{U}\cap 
(\xx'_{(k-1)}
\cup \yy_{(l-1)})=
\emptyset$.
A standard argument using Sard Lemma 
(see, e.g., \cite{milnhcob}, p.46)
implies that there is
a small isotopy $\Psi_t$ with support in $U$, \sut~
$\Psi_1(X_k')\pitchfork Y_l$.
$\qs$

\subsection{ Good fundamental system of neighborhoods}
     \label{su:gfsn}

The aim of this subsection is to introduce the notions of
$gfn$-systems and $ts$-submanifolds.
The $ts$-submanifolds are "thickenings of $s$-submanifolds" in 
the same
 sense as 
$B_\d(v)$ is a thickening 
of $D(v)$. The notion of $gfn$-system is a generalization of the 
notion
of $ts$-submanifold to the situation when we 
replace $s$-submanifolds by compact families of subsets of 
topological spaces.
They appear in our setting when we deal with
Morse functions on cobordisms with non-empty boundary, so that 
the descending discs are manifolds with boundary, non compact 
in general.

Proceeding to the definitions we start with good fundamental systems.
Recall that a {\it fundamental system of neighborhoods}
of a subset of $A$ in a topological space $X$ is a
family $\{ U_i\}_{i\in I}$
of open subsets
of $X$, \sut~ for every open set $U\supset A$ there is
$i\in I$ with $U_i\sbs U$.

\begin{defi} 
Let $X$ be a topological space, $A\subset X$
 a closed subset, $I$ an open interval
$]\d_0,\d_1[$. A
{\it
good fundamental system of neighborhoods of $A$}
(abbreviation: $gfn$-system for $A$)
is a family
$\{ A(\delta )\}_{\delta\in I}$
of open subsets of $X$,
 satisfying the following conditions:
\begin{enumerate}
\item for each $\delta\in I$ we have $A\subset A(\delta)$. Moreover, 
$\delta_1<\delta_2\Rightarrow A(\delta_1)\subset A(\delta_2)$,
\item for each $\delta\in I$ we have
 $\overline{A(\delta)}=\bigcap\limits_{\theta>\delta}A(\theta)$,
\item $A=\bigcap\limits_{\theta>\d_0}A(\theta)$.
\end{enumerate}
\end{defi}
\begin{lemm}\label{gfimplf}
 Assume that $X$ is compact. Let $I=]\d_0,\d_1[$. Let
$\{ A(\delta )\}_{\delta\in I}$
be a $gfn$-system for $A$.
Then: 
\begin{enumerate}
\item      The family
$\{ A(\theta )\}_{ \theta>\d_0}$
is a fundamental system of neighborhoods of $A$.
\item For every $\delta\in I$ the family
$\{ A(\theta )\}_{ \theta>\delta}$
is a fundamental system of neighborhoods of $\overline{A(\delta)}$.
\end{enumerate}
\end{lemm}
\Prf
(1) Let $U$ be an open neighborhood of $A$.
The sets $\{X\setminus \overline{A(\theta)}\}_{\theta >\d_0}$
form an open covering of the compact $X\setminus U$.
There is a finite subcovering, therefore there is $\theta >\d_0$,
 such that
$\overline{A(\theta)}\subset U$.
(2) is proved similarly.$\qs$

Next we define good fundamental systems of families of subsets.

\begin{defi} 
Let $X$ be a topological space,
$\aa = \{A_0,...,A_k\}$ be a  compact  family  of  subsets  of
$X$~,~   $I$  an  open  interval  $]\e_0,\e_1[$.  A  {\it  good
fundamental system of neighborhoods of $\aa$ }  (abbreviation:
$gfn$-system  for  $\aa$)  is  a  family  $\AAA=\{  A_s(\delta
)\}_{\delta\in  I,0\leq  s\leq  k}$  of  open  subsets  of  $X$,
satisfying the following conditions:

      (FS1) For every $s : 0\leq s\leq k$ and
 every $\delta\in I$ we have $A_s\subset A_s(\delta)$

      (FS2) For every $s : 0\leq s\leq k$
 we have $\delta_1<\delta_2\Rightarrow
A_s(\delta_1)\subset A_s(\delta_2)$

      (FS3)  For every $\delta\in I$  and every $j$
with $0\leq j\leq k$
we have
$\overline{A_{\leq j}(\delta)} =
{\bigcap}_{\theta >\delta}\big( A_{\leq j}(\theta)\big)$

      (FS4)  For every   $j$
with $0\leq j\leq k$ we have
$ A_{\leq j} =
{\bigcap}_{\theta >0}\big( A_{\leq j}(\theta)\big)$.

$I$ is called {\it interval of definition
}
of the system. $\aa$ is called the {\it core } of $\AAA$.

We shall denote $A_s(\d)$ also by $\AAA_{(s)}(\d)$ and
$A_{\leq i}(\d)$ also by $\AAA_{(\leq i)}(\d)$.
We shall denote $A_s$ also by $A_s(\e_0)$.
\quad$\triangle$
\end{defi}

\begin{rema}
Assume that $X$ is compact. 
Let $I=]\e_0,\e_1[$.
Let $\AAA=\{  A_s(\delta
)\}_{\delta\in  I,0\leq  s\leq  k}$  
be a $gfn$-system for a compact family 
$\aa=\{A_0,..., A_k\}$ of subsets of $X$. 
Then for every $\d\in [\e_0, \e_1[$ and every $0\leq s\leq k$
the family
$\{A_{\leq s}(\theta)\}_{\theta\in]\d,\e_1[}$
is a $gfn$-system for $\ove{A_{\leq s}(\d)}$.
\end{rema}

   If $M$ is a manifold without
boundary, $\xx$ is an $s$-submanifold of $M$,
and $\XXX$ is a $gfn$-system for $\xx$, then
we say that $\XXX$ is a
{\it $ts$-submanifold  with the core $\xx$. }
For a $ts$-submanifold
$\XXX =\{ X_s(\delta )\}_{\delta\in I,0\leq s\leq k}$
we shall denote $X_i(\delta)$ by $\XXX _{(i)}(\delta)$,
and $X_{\leq i}(\delta)$
  by
$\XXX _{(\leq k)}(\delta)$.

The next lemma follows
immediately from  Proposition \ref{fundsysord}. 
It provides the basic example of $gfn$-system and of $ts$-manifold.

\begin{lemm}\label{l:gfndisc}
 Let $W$ be a
  riemannian 
cobordism, $n=\dim M$.
 Let $f:W\to[a,b]$ be a Morse function
and $v$ be an $\e$-ordered $f$-gradient.
Then
 the family
$\{ B_\delta (\inde s;v)\}_{\delta\in ]0,\epsilon[,0\leq s\leq n}$
is a $gfn$-system for $\dd (v)$.
If $W $ is a closed manifold, this family
is a $ts$-submanifold with the core $\dd(v)$.
\quad$\square  $
\end{lemm}

 The $gfn$-system,
introduced in the lemma,
 will be denoted
by $\DDD (v)$. Note that there is no canonical choice
of the interval  of definition for this system.

\section{Handle-like filtrations}\label{s:hlf}

It is well known since Smale (see e.g. \cite{smapoi}, \cite{smale})
that to any Morse function on a closed manifold $M$
one can associate a so-called
{\it handle decomposition} of $M$.
This section contains a version of this construction, which 
suits better to our purposes, than 
the original one. The object in our setting, which replaces the 
classical "handle corresponding to a critical point
$p$" is the set $D_\d(p,v)$
(see Definition \ref{d:dthick}).
With our definition a systematic treatment of
"$\d$-thin handles with $\d\to 0$" is possible.
On the other hand these objects are not very well suited for 
studying \dfm ~type of \ma s
(which was the original aim of the seminal papers of S.Smale, 
cited above).
Indeed, we do not even introduce a smooth 
structure
on $D_\d(\indl s;v)$.
We only identify the \hot~type of
$D_\d(\indl s;v)$
(under some mild restrictions on $v$).

\subsection{$\d$-separated \gr s: definition}
\label{su:sepdef}

Recall from \cite{dold}, Ch. 5, that a {\it filtration}
of a topological space $X$ is a family
$\{ X_i\}_{i\in\ZZZ}$ of subsets of $X$ \sut~$i\leq j 
\Rightarrow X_i\sbs X_j$.
We shall say that a filtration $\{ X_i\}$ is {\it complete}
if $\cup_{i\in\ZZZ} X_i=X$.

Let $f:W\to[a,b]$ be a Morse function on 
a riemannian \cob~ $W$ of dimension $n$,
 and $v$ be an $f$-gradient.
The very first filtration of $W$  to introduce is of course the
 following:
$W_s=D_\d(\indl s;v)$
(where $\d>0$ is some fixed number). In general
these  $W_s$ do not form a good filtration;
$W_s$ can well be non-compact.
To obtain a  filtration with "handle-like" properties
we shall impose more restrictions on $f$ and on $v$.

Let $\phi:W\to[a,b]$ be an ordered Morse function with an
ordering sequence
$(a=a_0<a_1... <a_{n+1}=b)$.
 Let $\d>0$ and let $v$ be a $\phi$-gradient.
Denote $\phi^{-1}\([a_i,a_{i+1}])$ by $W_i$.

\begin{defi}
We say that $v$ is {\it $\d$-separated \wrt~
 the ordering sequence $(a_0,...,a_{n+1})$}, if

i) for every $i$ and every $p\in S_i(f)$ we have
$D_\d(p)\sbs \Wkr_i$.

ii) for every $i$ and every $p\in S_i(f)$
there is a Morse function
$\psi:W_i\to[a_i,a_{i+1}]$, adjusted to
$(\phi\mid W_i, v)$
and a regular value $\l$ of $\psi$ \sut~
$$D_\d(p)\sbs\psi^{-1}(]a_i,\l[)$$
and for every $q\not= p, q\in S_i(f)$ we have
$$D_\d(q)\sbs \psi^{-1}(]\l,a_{i+1}[)$$

If \noconf~we shall say that $v$ is a $\d$-separated $\phi$-\gr.

\end{defi}

\begin{defi}
Let $\fcob$ be an arbitrary Morse function , $v$ be an $f$-gradient.
We say that $v$ is {\it $\d$-separated }if 
it is $\d$-separated \wrt~some ordered Morse
function $\phi:W\to[a,b]$, adjusted to $(f,v)$.
\end{defi}

It is obvious that each $f$-gradient satisfying \ATA~ is $\d$-separated
for
some $\d>0$.

\subsection{Definition of handle-like filtrations}
\lb{su:tfo}

Now let $\phi:W\to[a,b]$ be an ordered Morse function with an ordering
 sequence
$(a_0,...,a_{n+1})$, where $n=\dim W$,
let $v$ be a $\d$-separated $\phi$-gradient; let $\nu\in]0,\d]$. Let 
$s$ be an integer between $0$ and $n+1$. Set
\begin{gather}
W^{\{\leq s\}}
=\phi^{-1}([a_0, a_{s+1}]); \quad W^{\{\geq s\}}=
\phi^{-1}([a_s,a_{n+1}])\\
W^{[\leq s]} (\nu)
= D_\nu(\indl s ;v)\cup \pr_0 W;\\
W^{\lc s\rc}(\nu)=
 \fii 0s
\cup D_\nu(\inde s;v)
\end{gather}

Note that
$W^{\{\leq 0\}}=\fii 01, W^{\{\geq 0\}}=W$,
$W^{[\leq 0]}(\nu)=
\(\cup_{p\in S_0(f)} D_\d(p)\)
\cup\pr_0 W=W^{\lc 0 \rc}$.
Extend  the definition to the negative values of
$s$ as follows. For $s\leq -2$
define all these sets (except $W^{\{\geq -s\}}$)
to be empty.
For $s=-1$ the definitions of
$\ws$,$\wsn$ are good as they are; we have
$W^{\{\leq -1\}}=\pr_0W= W^{[\leq -1]}(\nu)$.
For $s\leq 0$ set
$W^{\{\geq 0\}}=W$. Set  $W^{\lc -1\rc}(\nu)=
W^{[\leq -1]}(\nu)=
\dow$.

There is an ambiguity in this notation, since
$W^{\{\leq s\}}$ depend on $\phi$, and
$W^{[\leq s]} (\nu),
W^{\lc s\rc}(\nu)$
depend on $v$.
Each time we use this notation it will be clear what are
$\phi$ and $v$.
For the rest of the section \ref{s:hlf}
 we fix both $\phi$ and $v$, so \noconf.
 We have obviously:
$W^{\{\leq s\}} \sps
W^{\lc s\rc}(\nu)\sps
W^{[\leq s]} (\nu) $.
All the three filtrations will 
be referred to as {\it handle-like filtrations}.

We introduce one more family of subsets setting 

$$U_s= \fii 0s
\cup \bigg(\bigcup_{\ind p=s} D(p,v)\bigg)$$
for $s\geq 0$
and setting
$U_{-1}=\pr_0 W$, $U_s=\emp$ for $s\leq -2$.
The space $U_s$ is obtained from
$\fii 0s
$ by attaching the following 
family of $s$-dimensional cells:
$\{D(p,v)\cap\phi^{-1}\([a_s,a_{s+1}]\)\}_{p\in S_s(\phi)}$

\subsection{Homotopical properties of the filtrations}
\label{su:hpof}
\begin{prop}\label{p:fil}
Let $s\in\ZZZ$ and $\nu\in ]0,\d]$.
\begin{enumerate}
\item The  inclusions
$W^{[\leq s]} (\nu) \sbs
\wasn\sbs
W^{\{\leq s\}} $
are homotopy equivalences.
\item
The inclusion $U_s\sbs \wasn$ is a homotopy
equivalence.
\item
The inclusions of pairs
$$\bigg(W^{[\leq s]} (\nu),  W^{[\leq s-1]} (\nu)\bigg)
\sbs
\bigg(W^{\lc s\rc},W^{\{\leq  s-1\}}\bigg)
\sbs
\bigg(
W^{\{\leq s\}} , W^{\{\leq s-1\}}\bigg)
 $$
are homotopy equivalences.
\item The  inclusion
$$\bigg( U_s, W^{\{\leq s-1\}} \bigg) \sbs
\bigg( \wasn, W^{\{\leq s-1\}}\bigg)$$
is a homotopy equivalence.
\end{enumerate}
\end{prop}
\Prf
We consider only the case $s\geq 0$, the rest is trivial.
The assertions 1 and 2 
follow immediately from 3 and 4. The proofs of 3 and 4 are based
on the following lemma.
\bele
\label{l:defret}
Let $(Y,B)\sbs (X,A)$
be pairs of topological spaces.
Assume that there is a homotopy
$\Phi_t:X\to X, t\in[0,a]$, verifying the following
three properties:
\begin{enumerate}
\item[(i)] 
$\Phi_0=\id$,
\item[(ii)]
 the subsets $ Y,A$ and $B$
are invariant under $\Phi_t$ for all $t$,
\item[(iii)]
$\Phi_a(X)\sbs Y, \Phi_a(A)\sbs B$.
\enen

Then the inclusion
$i:(Y,B)\rInto (X,A)$
is a homotopy
equivalence.
\enle
\Prf
The map
$\Phi_a:(X,A)\to (Y,B)$
is the homotopy inverse for $i$.
Indeed,
$\Phi_a\circ i :(Y,B)\to(Y,B)$
equals to
$\Phi_a\mid (Y,B)$ which
is homotopic to the identity map
$(Y,B)\to (Y,B)$ via the homotopy
$\Phi_t\mid (Y,B)$. 
Further, the map
$i\circ \Phi_a: (X,A)\to (X,A)$
is homotopic to the identity via the homotopy $\Phi_t$. $\qs$

To deduce the assertion (3)  from Lemma \ref{l:defret}
define a homotopy $\Phi_t:W\to W$ by $\Phi_t(x)=\G(x,t;-v)$
(see Definition \ref{d:ascend} and Corollary \ref{co:ascend}
 for the definition 
and properties 
of $\G(\cdot, \cdot ; \cdot)$).
Here $t\geq 0$; note that for $T>0$ sufficiently large we have for every
$s$:
 $\Phi_T(\ws)\sbs \wsn$.
(Indeed, there is $\l>0$ \sut~ $df(v)(x)\geq\l$ 
for $x\in\ws\sm\cup_{p\in S(f)} B_\nu(p)$.
Therefore if $T$ is sufficiently large, there exists 
$t_0\leq T$ \sut~
$\G(x,t_0;-v)\in\dow\cup B_\nu(\indl s; v)$.)

The part 4 of our assertion 
 is a bit more delicate:
we can not do only with deformations
of the type
"shift along the $(-v)$-trajectories".
Note that the compacts
$D_\d(p,v)\cap \wa s{s+1}$
where $\ind p=s$ are disjoint, therefore
it suffices to prove that for each $p\in S_s(f)$
the space
$D(p,v)\cup \waa s$
is a strong deformation retract of
$D_\d(p,v)\cup \waa s$.
\footnote{Recall that $A\sbs X$ is {\it a strong deformation retract
of $X$ } if there is a homotopy
$f_t:X\to X$, \sut~$f_0=\id, f_t\mid A=\id, f_1(X)\sbs A$.}
We shall proceed as if there were only one critical point $p$ of
index $s$; the general case differs from this one only by complications
in the notation.

Let $\vphi$ be a Morse chart for $f$ at $p$.
We shall use the terminology from
Definition
\ref{morseterm}.
We assume that the condition (R) from the definition
\ref{morseterm}
is true, and also that 
$\vphi^{-1}(Q_r)\sbs D_\nu(p)$.
Recall that the vector field $(-\wi v)$ points
outward $Q_r$
in the points of $\pr_-^\circ Q_r$
and inward
$Q_r$ in the points of $\pr_+^\circ Q_r$.
Set 
$$R=\{z\in\wa s{s+1}\mid 
z=\g(x,t;-v)\mx{ where } x\in \vphi^{-1}(Q_r), t\geq 0\};$$
$$R_-=\{z\in\wa s{s+1}\mid 
z=\g(x,t;-v)\mx{ where } x\in \vphi^{-1}(\pr_-Q_r), t\geq 0\}$$

The sets $R$ and $R_-$ are weakly $(-v)$-invariant.
Set $W_0=\waa s$.
We have inclusions:
$$
(U_s,W_0)\sbs(R_-\cup U_s, W_0)
\sbs
(R\cup W_0, W_0)
\sbs \(\wasn\cup W_0, W_0\)
$$
The argument similar to the one used above shows
that the first and the third inclusions are homotopy equivalences.
To prove that the second inclusion is a \hot~ equivalence,
it suffices to 
prove that $R_-\cup U_s$
is a strong deformation retract of $R\cup W_0$.
For this it suffices in turn to prove that
$\dqr\cup(D(p,v)\cap Q_r)$ is a strong deformation retract 
of $Q_r$. The last assertion is
 obvious since the pair
$\(Q_r, \dqr\cup (D(p,v)\cap Q_r)\)$
is homeomorphic to
$(Q_r, \dqr\cap \RRR^k)$.$\qs$
\section{Morse complex}\label{s:mc}
\pa
\subsection{Definition of the Morse complex}\label{su:dmc}

In this section  $W$ is a riemannian cobordism of dimension $n$, 
$\fcob$
is a Morse function on  $W$, $v$
is an $f$-gradient.
Let $s$ be an integer between $0$ and $n$.
Consider the free abelian group $C_s$
generated by all critical points of $f$ of index $s$.
Assuming that $v$ satisfies \TA~
we shall introduce for every $s$
a \ho~
$\ds: C_s\to C_{s-1}$
satisfying $\pr_{s-1}\circ\pr_s=0$,
so that we shall obtain a chain complex over $\ZZZ$, called
{\it Morse complex}.
The construction of $\ds$
makes use mainly of the orbits of
$v$ joining critical points of $f$, but it depends
also on an additional arbitrary choice
(namely on the  choice of orientations
of descending discs).
The homology of the complex is isomorphic to
$H_*(W,\pr_0 W)$.

Proceeding to precise definitions, choose an ordered Morse function
$\phi:W\to[a,b]$
with an ordering sequence
$(a_0=a, a_1,...,a_{n+1}=b)$, adjusted to $(f,v)$.
In the preceding subsection we have constructed a filtration
$\ws$. The filtration $\ws$ is complete;
we obtain thus a complete filtration in the singular chain complex of
$(W,\pr_0 W)$.
From the proposition \ref{p:fil}
we know that the pair
$\(\ws,\wsm\)$ is \hot~equivalent to
$\(U_s,\wsm\)$, and $U_s$ is obtained from
$\wsm$ by attaching
the following family
of $s$-dimensional cells:
$\{D(p,v)\cap\fii s{s+1}\}_{p\in S_s(\phi)}$.
Therefore we have:
\begin{equation}
H_*\(\ws, \wsm\)=
\left\{
\begin{aligned}
0 &\mx{ if } *\not= s\\
(\ZZZ)^{\#S_s(\phi)}&\mx{ if } *=s
\end{aligned}
\right.\lbl{f:cell}
\end{equation}
We shall denote
$H_s\(\ws,\wsm\)$ by $C_s$.
The property (\ref{f:cell})
implies that the filtration
$\{\ws\}$
is a {\it cellular filtration}
in the sense of \cite{dold}, Ch 5, \S 1.
(Note that in the subsection \ref{su:tfo}
we have defined $\ws$
for all $s\in\ZZZ$; in particular
$W^{\{-1\}}=\dow$).
Therefore, by \cite{dold}, Ch. 5, Prop. 1.3
 there is an isomorphism
$$
H_*(W,\dow)\approx H_*(C_*,\pr_*).
$$
Here $\ds:C_s=H_s\(\ws,\wsm\)\to C_{s-1}=H_{s-1}\(\wsm,\wsmm\)$
is the composition of the boundary homomorphism
of the exact sequence of the pair
$\(\wsm,\wsmm\)$
and the projection
$$H_{s-1}\(\wsm\)\to H_{s-1}\(\wsm,\wsmm\)$$
The complex $(C_*,\pr_*)$ is called {\it Morse complex}, and
denoted by $C_*(\phi,v)$.

By Theorem \ref{t:iceq} and Remark \ref{r:indind} below
this complex  depends only on $v$, but 
not on the choice of ordered Morse function $\phi$,
of which $v$ is a gradient.
So we shall often denote it by $C_*(v)$.

It occurs that the boundary operator $\pr_*$
is computable from the geometric data
(see Theorem \ref{t:iceq} below).

{\bf Historical remark. \quad}
The essential in this construction was discovered by S.Smale
and used in his proof of the Poincare conjecture
in dimensions $\geq 5$ (see \cite{smapoi}).
In the book \cite{milnhcob}, \S 7
one finds almost the statement of the theorem \ref{t:iceq}.
Much later the interest to the subject was stimulated by the paper
\cite{witt}, where one finds the theorem \ref{t:iceq}
in its present form
(for the case of riemannian gradients).
The proof suggested in \cite{witt}
is, however, based on different ideas, than those of
Smale. A proof of
Theorem \ref{t:iceq} is contained 
in the Appendix to my paper \cite{patou}
(only the case of pseudo-gradients
 is treated there, but we shall see a bit later that this does not
diminish the generality).
$\qt$

Choose 
for each $p\in S_s(f)$ an orientation of the tangent
space $T_-(p,v)$ to the local stable
manifold. Then the manifold
$D(p,v)\cap \Wkr$
obtains an orientation.
Further, the manifold $D(p,-v)$
obtains a co-orientation
(recall that a {\it co-orientation} of a submanifold
$N$ of a manifold $M$
is an orientation of the normal bundle to $N$ in $M$).
The fundamental class of $D(p,v)$ modulo its boundary defines
an element  $d(p)\in C_s$.
These elements form a free base in $C_s$.
Write
$$
\ds d(p)=\sum_{q\in S_{s-1}(f)} n(p,q)\cdot q
$$
where $n(p,q)\in\ZZZ$.
We call $n(p,q)$
{\it homological incidence coefficient}
corresponding to $p$ and $q$.
We shall also write
$n(p,q;v)$ to stress that a priori this index depends on $v$
(we shall see later that often it does not ).

Assume now that  $v$  satisfies
\TA. Denote $D(p,v)\cap\phi^{-1}(a_s)$
by
$S(p,v)$ and $D(q,-v)\cap \phi^{-1}(a_s)$
by $S(q,-v)$.
For each $p\in S_s(f), q\in S_{s-1}(f)$ the sets $S(p,v), S(q,-v)$
are compact 
submanifolds of complementary dimension in $V_{a_s}$,
 and they are transversal to each other.
Therefore their intersection $N(p,q;v)$ is a finite set.
Its points are in 1-1 correspondence with
the orbits of $v$ joining $q$ to $p$.
Due to our choice of orientations each point $x\in N(p,q;v)$ obtains
a sign $+$ or $-$, denoted by $\ve(x)$.
(Namely, the manifold
$S(p,v)$ is oriented as the boundary of $D(p,v)$. Further,
$S(q,-v)$ is a cooriented submanifold
of $V_{a_s}$, since $D(q,-v)$ is cooriented.
Then in each point of $x\in S(p,v)\cap S(q,-v)$
we obtain an orientation of
$T_x(S(p,v))$
and a coorientation of
$T_x(S(q,-v))$.
Since these two subspaces of $T_xV_{a_s}$ are
complementary, we obtain a sign).

 \label{finite}

Set 
$$
\nu(p,q;v)=\sum_{x\in N(p,q;v)}\ve(x)
$$
\bere\label{r:indind}
The number $\nu(p,q;v)$
depends only on $v$, but not on the choice
of ordered Morse function $\phi$,
adjusted to $(f,v)$.
Indeed, the above definition is easily reformulated
in terms of orbits of $v$: the number
$\nu(p,q;v)$
is the algebraic number of $v$-orbits joining $p$ with $q$.)
\enre

\subsection{Properties of incidence coefficients}
\lb{su:pic}

\beth\label{t:iceq}
For all $p,q$ with $\ind p=\ind q+1$ we have:
$n(p,q;v)=\nu(p,q;v)$.
\enth
\Prf

By the remark \ref{r:indind}
we can assume that all the critical points of $\phi$ of
index $s$ lie on the same critical level of $\phi$, say
$c_s\in]a_s,a_{s+1}[$.

We shall say that two $\phi$-gradients $v$ and $w$ are {\it tangent}
if for every critical point $p\in S(\phi)$ we have:
$T_-(p,v)=T_-(p,w), T_+(p,v)=T_+(p,v)$.
Using partitions of unity, it is not difficult to 
show that every $f$-gradient is 
tangent to an $f$-Gradient.

\bepr
\label{p:indinc}
Let $w$ be another $\phi$-gradient, satisfying \TA.
Assume that $v$ and $w$ are tangent. 
 Then
for every $p,q\in S(f)$
with $\ind p=\ind q+1$
we have
\been\item
$\nu(p,q;v)=\nu(p,q;w)$
\item
$n(p,q;v)=n(p,q;w)$.
\enen
\enpr
\Prf 1) 
The number $\nu(p,q;v)$
is the intersection index of submanifolds
$S(p,v), S(q,-v)$ in $V_{a_s}$
(the first being oriented, the second cooriented).
To prove $\nu(p,q;v)=\nu(p,q;w)$
it suffices to show that $S(p,v)$ is
homotopic to
$S(p,w)$
and $S(q,-v)$ is homotopic to $S(q,-w)$.

The proof of
$S(p,v)\sim S(p,w)$
makes use of
Morse  chart $\varphi$ for
$\phi$ around $p$.
Since $v$ and $w$ are tangent, we can 
use the same Morse chart $\vphi$ for both $v$ and $w$.
We assume here therefore the terminology
of Definition \ref{morseterm}.
Let $r$ be so small that the condition (R)
from Definition
\ref{morseterm} is valid for $v$ and for $w$. Let
$t\in [0,1]$ and
consider the vector field $v_t=tw+(1-t)v$.
In general it is not a $\phi$-gradient, but it satisfies
$d\phi(v_t)(x)>0$ for every $x\notin S(\phi)$.
This implies that every $(-v_t)$-trajectory, starting
 in $\dqr$, reaches $V_{a_s}$.

Consider a continuous map
$\psi_t:\dqr\to \dow$
which associates to each point $y$ of $\dqr$ the point of intersection
of
$\g(y,\cdot;-v_t)$ with $V_{a_s}$.
Set $\s(p,v)=\dqr\cap D(p,v)$,
$\s(p,w)=\dqr\cap D(p,w)$.
These are two oriented embedded spheres in $\dqr$,
which are homotopic to each other; denote the corresponding homotopy by
$H_t$.
Then $\psi_t\circ H_t$ is the \hot~
which we were seeking for.
The proof of $S(q,-v)\sim S(q,-w)$
is similar. 

2) The proof is similar to 1): deforming $w$ to $v$

one shows that the discs $D(p,v)$ and $D(p,w)$ define the same homology
classes in
$H_*\(\ws,\wsm\)$.
We leave the details to the reader.
$\qs$
\pa
In view of this proposition we can assume that the 
image of $v$ in the Morse chart corresponding to $p$ is
the standard vector field $\Vv_s$ (where $s=\ind p$).
Denote $\phi^{-1}\([a,a_{s-2}]\)$
by $W_0$, 
$\phi^{-1}\([a,a_{s-1}]\)$
by $W_1$, 
$\phi^{-1}\([a,c_{s-1}[\)$
by $W'$.

Fix some $p\in S_s(f)$
and denote
$S(p,v)$ by $M$.
Consider its fundamental class $[M]$ as an element
of
$H_{s-1}(W_1, W_0)$.
Our theorem describes the decomposition
of $[M]$ \wrt~to the base $\{d(q)\}_{q\in S_{s-1}(f)}$ of
$H_{s-1}(W_1, W_0)$
in terms of the algebraic intersection index of $M$ with
the
 spheres $S(q,-v)$.
In order to find this decomposition we shall perform a homotopy of
$M$ in $W_1$.
The inclusion
$(W_1, W_0)\rInto (W,W')$
is a \hot~ equivalence;
therefore it suffices to calculate the image of $[M]$
 in the group
$H_{s-1}(W_1,W')$.
Let $q\in S_{s-1}(\phi)$ and let 
$x_1,...,x_l\in M$ be the points of intersection
of $M$ with $D(q,-v)$. 
Denote the closed $(s-1)$-disc $D^{s-1}(0,1)$ by $\D$ 
and $D^{s-1}(0,\l)$
by $\D_\l$. The interiors of these discs will be denoted by
$\D^\circ$, resp. $\D_\l^\circ$.
For each $i:1\leq i\leq l$
choose
an embedding
$E_i:\D\to M$
\sut~$E_i(0)=x_i$ and 
$E_i$
conserves the orientation.
We assume also that the images of $E_i$ are disjoint.
The set
$M'=M\sm\bigcup_i E_i(\D^\circ)$
is a compact subset of the domain of
definition of
$\stind v{a_s}{a_{s-1}}$
therefore the gradient descent along $(-v)$ will push
$M'$ below $c_{s-1}$. (In precise terms,
denote $\G(x,t;-v)$ by $H_t(x)$. Then 
for $T>0$ sufficiently large we shall have:
$H_t(M')\sbs\phi^{-1}\([a,c_{s-1}[\)$.)
Therefore the image of $[M]$ in $H_*(W_1, W')$
equals to the sum
(over all $i$) of the fundamental 
classes modulo boundary of singular discs
$H_t\circ E_i:\D\to (W_1, W')$.
The same is of course true if instead of the embeddings $E_i$ we
consider their restrictions to $\D_\l, 0<\l<1$.
Consider for example the point $x_1\in M$, assume that
$x_1\in D(q,-v)\cap D(p,v)$,
where $q\in S_{s-1}(\phi)$.
Let $\ve$ be the  intersection sign of $S(p,v)$ and $S(q,-v)$ in $x_1$.
Let $\vphi:U\to V$
be a Morse chart for $\phi$ around $q$.
We refer again to Definition \ref{morseterm}
for the notation and definitions.
Let $T$ be so large that
$\g(x_1,T; -v)$ belong
to $\Int\vphi^{-1}(Q_r)$.
The map $\rho:y\mapsto\vphi\(\g(E_1(y), T;-v)\)$
restricted to $y\in \D_\l$, where $r$ is sufficiently small,
is a diffeomorphism onto its image, which is in
$\Int Q_r$.

Consider now the local situation in $Q_r$
The map
$\rho'(0): \RRR^{s-1}\to\RRR^n=\RRR^{s-1}\oplus\RRR^{n-s-1}$
is an injection, moreover, its first coordinate is an isomorphism
$H:\RRR^{s-1}\to\RRR^{s-1}$.
Note that the sign of $\det H$ equals to $\ve$.
Consider the linear map $\a:\RRR^{s-1}\to\RRR^{s-1}\oplus\RRR^{n-s+1}$
defined by
$x\mapsto\bigg(\frac r{2\Vert H\Vert} Hx, 0\bigg)$.
The restriction of $\a$ to the unit disc $D^{s-1}(0,1)$ 
is a diffeomorphism onto its image
$K$ which is contained in $\Int Q_r \cap L_-$.

Apply now to $\rho$ the diffeomorphism
$\Phi(T',-\Vv_{s-1})$, where $T'$ is a large enough real number.
It is clear that
if $\l>0$ is small, the embedded disc
$\Phi(T', -\Vv_{s-1})\(\rho(\D_\l)\)$
will become larger and more and more flat.
More precisely, for every $\e>0$ there is $\l>0$ and $T'>0$
\sut~the embedded disc
$K'=\Phi(T',\Vv_{s-1})(\rho(\D_\l))$
is $\e$-close to
the embedded disc $K$.
(This is a matter of an easy computation
which will be left to the reader.)
Therefore the homology class of $(K', \pr K')$
in $(Q_r, Q_r^-)$
is the same as that of the negative disc modulo
its boundary, multiplied by $\ve$. $\qs$

\section{Gradient descent map}
\label{s:gdm}

In this section we deal with the following situation.
Let $W$ be a cobordism, $\fcob$ be a Morse function
on $W$, and $v$ be an $f$-gradient.
Let $U_1\sbs\pr_1 W$
be the set of all $x$, \sut~
the $(-v)$-trajectory $\g=\g(x,t;-v)$
starting at $x$
 reaches $\pr_0 W$.
For $x\in U_1$ the point of intersection of $\g$ with
$\pr_0 W$ will be denoted by
$\stexp {(-v)} (x)$.
We obtain a $\smo$ diffeomorphism
$\stexp {(-v)}: U_1\to U_0$, where $U_0$ is an open
subset of
$\pr_0 W$.
In general of course $U_i\not=\pr_i W$.
For $A\sbs\pr_1 W$
denote
$\stexp {(-v)}(A\cap U_1)$ by
$\stexp {(-v)}(A)$. Note that if $L\sbs\pr_1 W$ is a
submanifold of
$\pr_1 W$, then
$\stexp {(-v)}(L)$
is a submanifold of $\pr_0 W$, but
 compactness of $L$ does
not imply compactness of
$\stexp {(-v)}(L)$. To make
$\stexp {(-v)}(L)$
compact, we should add the soles of descending
discs, and the resulting object will not be a
manifold. So one can think that the gradient descent map
is well defined on the category of stratified
manifolds. It is indeed so, and
the corresponding construction is given in subsection \ref{su:trssb}
(subsection \ref{su:gpt} contains some preliminaries).
To explain it briefly, let
$\aa=\{ A_0,...,A_k\}$
be an $s$-submanifold of $V_1$
(see the definition in Subsection \ref{su:ssubm}).
Let $A'_i$ denote the union of
all
$(-v)$-trajectories, starting in $A_i$.
Set $B_i=A'_i\cup D(\inde i, v)$. The main
observation is that
$\{B_i\}_{0\leq i\leq k}$
is a compact family, provided that
$\aa\nmid \dd(-v)$
(see Lemma \ref{l:trcomp}).
It is natural therefore to define an $s$-submanifold
$\stexp {(-v)}(\aa)$
of $\pr_0 W$ setting
$\stexp {(-v)}(\aa)_{(i)}=B_{i+1}\cap\pr_0 W$.
The most technical part of the section is Subsection \ref{su:totss}.
Here we deal with thickenings. We construct a $gfn$-system
 for
$\{B_i\}$
starting with a $gfn$-system for $\aa$ and the $gfn$-system
$\DDD(v)$.

\subsection{General properties of tracks}
\label{su:gpt}
In this subsection $W$ is a cobordism, $f$ is a Morse
function
on $W$, $v$ is an $f$-gradient.

\begin{defi}
Let $X\sbs W$. The set
$\{\g(x,t;-v)\nmid x\in X, t\geq 0\}$
is denoted by
$T(X,-v)$
and called
{\it track of $X$ \wrt~ $(-v)$}. 
\end{defi}
The aim of this subsection is to establish the basic
 properties of $T(X,-v)$.
Next lemma is obvious.

\begin{lemm}\label{tractriv}
\begin{enumerate}
\item For every $A\sbs W$ we have
$T(\ove{A},-v)\sbs\ove{T(A,-v)}$.
\item
If $\{ X_i\}_{i\in I}$ is a family of subsets of
$\pr_1 W$, then
$\cap_i T(X_i,-v)=T(\cap_i X_i,-v)$.
$\qs$
\end{enumerate}
\end{lemm}

In the following lemma we list the compactness properties of sets
$T(X,-v)$. It can be considered as an analog
of
\ref{dcomp}.

\begin{lemm}\label{tracomp}
\begin{enumerate}
\item
 If $X$ is compact,
 then  $T(X,-v)\cup D(v)$ is compact.

\item
 If $X$ is compact, and
every $(-v)$-trajectory starting from a point of $X$
reaches $\pr_0 W$ then $T(X,-v)$ is compact.

\item
 For any $X$ we have $\overline{T(X,-v)\cup D(v)} =
T(\overline X,-v)\cup D(v)$.

\item
 For any $X$
and $\delta>0$
 we have $\overline{T(X,-v)\cup B_\delta (v)} =
T(\overline X,-v)\cup D_\delta (v)$.
\end{enumerate}
\end{lemm}

\Prf
\begin{enumerate}
\item A point $y\in W$ {\it does not} belong to
$T(X,-v)\cup D(v)$, \ifff the trajectory
$\g(x,\cdot; -v)$ reaches $\pr_1 W$ without
 intersecting $X$. Therefore
the complement of
$T(X,-v)\cup D(v)$ in $W$
is open.

\item
For $x\in X$ let
$\tau(x)\geq 0$ be the
moment when
$\g(x,t;-v)$ intersects $\pr_0 W$. Then
$\tau(x)$ is continuous in $x$, so it is bounded by
 a constant, say, $C>0$. Now let
$y_i=\g(x_i,t_i;-v)\in T(X,-v)$
where $x_i\in X, 0\leq t_i\leq C$.
Choose a subsequence $i_k$, \sut~ $\{x_{i_k}\}$ and
$\{t_{i_k}\}$ converge, and we are over.

\item Since $\ove{X}$ is compact,
$T(\ove{X},-v)\cup D(v)$
is also compact, and the inclusion
\break
$\ove{T(X,-v)\cup D(v)}
\sbs
T(\ove{X},-v)\cup D(v)$
holds. The inverse inclusion follows from the preceding lemma.

\item Similar to (3).
\end{enumerate}
$\qs$

\subsection{Tracks of $s$-submanifolds}
\label{su:trssb}

In this and in the following
subsection $\fcob$ is a Morse function on a cobordism
$W$, $v$ is an $f$-gradient, satisfying  the  \ATA,
$\aa=\{A_0,..., A_k\}$
is an $s$-submanifold of $\pr_1 W$, satisfying $\aa\nmid\dd_b(-v)$.

\begin{defi}
For $0\leq i\leq k+1$ set $TA_i(-v)=T(A_{i-1} ,-v)
\cup D(\inde i ;v)$ ( we set
$A_{-1}=\emptyset$ by definition).
Introduce a family $\ttt(\aa,-v)$ of subsets of $W$ by
 $\ttt (\aa ,-v)=\{ TA_i(-v) \}_{0\leq i\leq k+1}$
 For $\l\in[a,b]$
introduce a family $\stind {(-v)}b\lambda (\aa)$ of subsets of $V_\l$
  as follows:
$\stind {(-v)}b\lambda (\aa)=
\{ TA_{i+1}
(-v)
\cap f^{-1}(\lambda)\}_{0\leq i\leq k}$
of subsets of $f^{-1}(\lambda)$.
If the values $b,\lambda$ are clear from the context
 we shall abbreviate
$\stind {(-v)}b\lambda (\aa)$ to $\st {(-v)} (\aa)$.
The family $\ttt (\aa ,-v)$ will be called {\it track }
of $\aa$, and the family $\stind vb\lambda (\aa)$
will be called
{\it $\st v$-image } of $\aa$.
\quad$\triangle$
\end{defi}
\begin{lemm}\label{l:trcomp}
\begin{enumerate}\item    $\ttt (\aa ,-v)$
and
$\stind {(-v)}b\lambda (\aa)$ are compact families .

\item    If $\lambda$ is a regular value of $f$, then
$\stind {(-v)}b\lambda (\aa)$
is an $s$-submanifold of  $f^{-1}(\lambda)$.
\end{enumerate}
               \end{lemm}
\Prf
(1) Let
$0\leq s\leq k+1$. Denote
$ T(A_{\leq s-1} ,-v)\cup D(\indl s;v)$
by $Q(s)$.
Let $\phi :W\to[a,b]$ be an ordered Morse
function, adjusted to $(f,v)$, and
$(a_0,...,a_{n+1})$ be the ordering sequence for
$\phi$.
The all the critical points of $\phi$ of indices $>s$
are above $a_s$ and
those of indices $\leq s$ are below $a_s$.
Since
$\aa\nmid (\dd (-v)\cap V_1)$, every
 $(-v)$-trajectory, starting at a point
of $A_{\leq s-1}$ reaches
$\phi^{-1}(a_{s+1})$. Then
$Q(s)\cap\phi^{-1}([a_{s+1} ,a_{n+1}])$
is compact by \ref{tracomp} (2)
and
$Q(s)\cap\phi^{-1}([a_0, a_{s+1}])$
is compact by \ref{tracomp} (1), therefore
$Q(s)$ is compact.
(2) is obvious.$\qs$

The assumption $\aa\nmid\dd_b(-v)$,
adopted in the beginning of the present subsection is
not very restrictive: given $\aa$, we can always perturb 
$v$ so that this assumption will hold. A more general 
statement is the subject of the following Proposition.

\begin{prop}\label{trper}
Let $\bb$ be an $s$-submanifold of $\pr_1 W$, $\cc$
be an $s$-submanifold of $\pr_0 W$.
Let $U$ be a \nei~  of $\pr W$.
 Then there is an $f$-gradient $w$,
satisfying \ATA, \sut~ 
\begin{enumerate}\item
$\dd_b(-w)\nmid\bb$;
\item $\stexp {(-w)} (\bb)\nmid \cc;\quad \st w(\cc)\nmid \bb$
\end{enumerate}
Moreover, $w$ can be chosen arbitrarily $\smo$ close to $v$,
and \sut~ $\supp (u-w)\sbs U$, and 
$\supp (u-v)\cap \pr W=\emp$.
\end{prop}
\Prf Choose $a>0$ \sut~
the collars
$C_0=\g(\pr_0 W,[0,a];v)$ and
$C_1=\g(\pr_1 W, [0,a]; -v)$
are disjoint. Our new $f$-gradient $w$ will satisfy
$\supp(w-v)\sbs \Int (C_0\cup C_1)$,
which imply immediately \ATA.
Apply to $v$ the AHC-construction nearby the upper boundary 
$\pr_1 W$ of $W$.
We choose the input data $(a,\xi,h)$ as follows.
Fix first a function $h$ as to satisfy
$\int_0^ah(\tau)d\tau =1$. Choose then a vector field
$\xi$ on $\pr_1 W$ so that 
$\Phi(\xi, 1)\(\dd_b(-v)\)\nmid\bb$.
(This is possible by \ref{transvssub}),
and, moreover, we can choose $\xi$
so $\smo$ small
that
$v_0=AHC(a,\xi,h;v)$
is $\smo$ close to
$v$, and that 
$v_0$ is still an $f$-gradient. 
Now AHC2 from Subsection \ref{su:ahc}
 implies that
$\Phi(\xi,1)\(\dd_b(-v)\)=\dd_b(-v_0)$
and the point 1)
of our conclusions is true with $v_0$ instead of $w$.
Therefore $\ttt(\bb, -v_0)$
is defined and the $s$-submanifold $\ll=\stexp {(-v_0)}(\bb)$ of
$\pr_0 W$ is defined. Applying
to $\ll$ and to $\cc$ a similar construction we obtain 
an $f$-gradient $w$, satisfying the 1) and the first half 
of 2) of our conclusion. Note that the first half of 2) 
implies the second, and the proof is over. $\qs$

\subsection{ Tracks of $ts$-submanifolds}
\label{su:totss}

In this subsection we assume that our $f$-gradient $v$
is $\d_1$-ordered.

\begin{defi}
 Let
 $\AAA =\{ A_s(\delta )\}_{\delta\in ]0,\delta_0[ ,0\leq s\leq k}$
be a $ts$-submanifold of $V_1$ with the core
$\aa$.
 For $0<\delta<\min (\delta_0,\delta_1)$
and $0\leq s\leq k+1$
set
$TA_s(\delta , -v) = T(A_{s-1} (\delta),-v)
\cup
 B_\delta (\inde s; v)$ (where we set by definition
$A_{-1}(\delta)=\emptyset$).
\quad$\triangle$
\end{defi}
We shall prove  that for \sufsm~  $\e>0$ the system
$\{TA_s(\delta,-v)\}_{\d\in]0,\e[, 0\leq s\leq k}$
 is a $gfn$-system for
$\ttt (\aa ,-v)$.
Up to the end of this section
$Q(s,\delta)$ stands for
$T(\overline{A_{\leq s-1}(\delta)},-v)
\cup  D_\delta (\indl s;v)$.

\begin{prop}\label{tracgfn}
There is $\epsilon \in ]0,\min (\delta_0,\delta_1)[$
such that
\break
$\{TA_s(\delta , -v) \}_{\delta\in]0,\epsilon[,~0\leq s\leq k}$
is a gfn-system for $\ttt (\aa ,-v)$.
\end{prop}
\Prf
We shall first prove two lemmas, establishing the properties of
the introduced family $\{TA_s\}$.

\begin{lemm}
For every $\d$ with $0<\d\leq \min(\d_0,\d_1)$
and for every $s$
we have
  \begin{gather}
TA_{\leq s}(-v)={\cap}_{\theta>0} TA_{\leq s}(\theta,-v)
\\
Q(s,\delta)={\cap}_{\theta>\delta} TA_{\leq s}(\theta,-v),
\\
TA_{\leq s}(\delta,-v)\subset Q(s,\delta)
\subset \overline{TA_{\leq s}(\delta,-v)}
\end{gather}
\end{lemm}

\Prf To prove (1) note first that
$TA_{\leq s}(-v)={\cap}_{\theta>0} TA_{\leq s}(\theta,-v)$
(by \ref{tractriv})
\break
${\bigcap} _     {\theta>0} B_\theta (\indl s;v) =
D (\indl s;v);$
(by \ref{fundsysord})
Since both
$TA_{\leq s}(\theta,-v)$
and
$B_\theta (\indl s;v)$
are decreasing in $\theta$ families of subsets,
the (1) holds.
(2) is proved similarly.
To prove (3), note that the first inclusion is
straightforward, the second follows
from
\ref{tractriv} $\qs$

\begin{lemm}
Let $0<\epsilon\leq\min(\delta_0,\delta_1)$.
Then (i)$\Leftrightarrow$(ii).
\begin{itemize}\item[(i)] For every $0<\delta<\epsilon$ and every
$s:0\leq s\leq k+1$ the set
$Q(s,\delta)$ is compact
\item[(ii)]
$TA_s(\delta , -v)_{0<\delta<\epsilon,
0\leq s\leq k+1}$
is a $gfn$-system for $\ttt (\aa,-v)$
\end{itemize}\end{lemm}
{\it Proof.\quad}
The properties (FS1), (FS2), (FS4) hold for the system
\break
$TA_s(\delta , -v)_{0<\delta<\epsilon,
0\leq s\leq k+1}$
for every $\e$. So we have to check that (i) is equivalent to (FS3).
And this follows immediately from the preceding lemma.
$\qs$

Now the question is reduced to the compactness of the subsets
$Q(s,\d)$.
Since $v$ is $\d_1$-ordered, there is a Morse function
$\phi: W\to[a,b]$, adjusted to
$(f,v)$ with an ordering sequence
$(a_0,a_1,...,a_{n+1})$,
\sut~ for every
$p\in S_s(f)$
we have:
$D_{\d_1}(p)\sbs \phi^{-1}(]a_s,a_{s+1}[)$.

\begin{lemm}
Let $0<\d\leq\min(\d_0,\d_1)$. Then
$(i) \Leftrightarrow (ii)$.

(i) $Q(s,\d)$ is compact.

(ii) every $(-v)$-trajectory starting at
a point of $\ove{A_{\leq  s-1}(\d)}$
reaches $\phi^{-1}(a_{s+1})$.
\end{lemm}

\Prf Assume that (ii) holds. Then
$Q(s,\d)\cap\phi^{-1}([a_{s+1},a_{n+1}])$
is compact by \ref{tracomp} (2)
and
$Q(s,\d)\cap\phi^{-1}([a_0,a_{s+1}])$
is compact by \ref{tracomp} (1).

Assume now that (ii) does not hold. Then there is
$x\in\ove{A_{\leq s-1}(\d)}$
\sut~
$\g(x,t;-v)$ converges to
$p\in S(f)$ as $t\to\infty$,
and $p\in\phi^{-1}(]a_{s+1}, a_{n+1}[)$.
But in this case
$p\in \ove{Q(s,\d)}\sm Q(s,\d)$.
$\qs$

This lemma implies in particular that if
$0<\d'\leq\d\leq\min (\d_0,\d_1)$
and $Q(s,\d)$ is compact,
then $Q(s,\d')$ is also compact.

\begin{defi}
Let $\l>0$. Let
$\AAA$ be an $ts$-\sma~ of $\pr_1 W$
with interval of definition
$I=]0,\d_0[$.
We say that the pair $(\AAA, v)$ is $\l$-ordered, if
\begin{enumerate}
\item
 $v$ is $\l$-ordered
\item $\l$ is in the interval of definition of $\AAA$
\item
$Q(s,\l)$ is compact for every $s$. 
\end{enumerate}
\end{defi}

Three preceding lemmas imply that if $(\AAA,v)$ is $\l$-ordered, then
\break
$\{TA_s(\delta,-v)\}_{\d\in]0,\l[, 0\leq s\leq k}$
is a $gfn$-system.
Therefore to prove our proposition it suffices to show that there is
$\l$ \sut~ $(\AAA,v)$ is $\l$-ordered.
But this is obvious. Indeed, every
$(-v)$-trajectory starting  in $A_{\leq s-1}$
reaches $\phi^{-1}(a_{s+1})$.
Recall that    $\{ A_{\leq s-1}(\theta)\}_{\theta>0}$
is a fundamental system of neighborhoods of the compact
$A_{\leq s-1}$. An easy compactness argument
shows that for some $\theta >0$ every
$(-v)$-trajectory, starting in
$A_{\leq s-1}(\theta)$
reaches
$\phi^{-1}(a_{s+1})$ as well. $\qs$

 We shall denote the $gfn$-system,
introduced in \ref{tracgfn} by $\TTT(\AAA,-v)$
and call it {\it track } of $\AAA$ (\wrt~ $-v$).
Note that there is no canonical choice of
the interval of definition
for this system.

\section[$C^0$ continuity of integral curves]{Properties of
 continuity of integral curves
of vector fields with respect to $C^0$ small perturbations of
of vector fields}
\lb{s:pcic}

Let $v$ be a $C^1$ \vf~ on a \ma~$M$,
$\g$ be an integral curve of $v$.
It is well known that given $t\geq 0$
the value $\g(t)$
depends continuously on the initial value
$\g(0)$
and on the \vf~ itself (where the topology
on the space of \vfs~ is the {\it $C^1$- topology}.
Actually this continuity property is still true if we replace
$C^1$ topology
by {\it $C^0$-topology}.
\footnote{
By the time when the preprint  \cite{paepr}
was written I had no references for this property at hand.
I am grateful to O.Viro
for the indication that
this "$C^0$ continuity property"
is known and must be exposed somewhere in the literature.}
It is due to this fact that the $C^0$ topology
plays an important role in our work.
The present section
contains the statements of the continuity theorems we shall need
in the sequel. We start with citing two propositions from \cite{pastpet}.
Afterwards we deduce from them some properties (see
Corollary \ref{c:cont3}, Proposition \ref{p:coco}, Proposition
\ref{p:co})
in order to obtain the results in the form
ready for the applications in the further sections.
The results which are used in the sequel are gathered
in \sub~ \ref{su:mb}.
Essentially all the contents of of the present section is already
present in \cite{pastpet}; we have only rearranged the material.

\subs{Manifolds without boundary}
\lb{su:mwb}

For a manifold $M$ (without boundary)
 we denote by
$\hrrr$ (resp. $\vemm$) the vector space
of $C^1$ vector
fields on $M$ (resp. the vector space of
 $C^1$ vector
fields on $M$ with compact support).
 In this subsection $M$ is a
riemannian manifold
without boundary, $v\in\vemm , n=\dim M$.

\begin{propcit}{(\cite{pastpet}, Prop. 5.1)}
\lb{p:5.1}
Let $a,b\in M, t_0\geq 0$
and $\gamma (a,t_0;v)=b$.
Then for every open
 neighborhood $U$ of
$\gamma(a,[0,t_0];v)$ and
every open
 neighborhood $R$ of $b$
there exist $\delta>0$ and an
open neighborhood
 $S\subset U$ of $a$ such that
$\forall x\in S$ and
$\forall w\in\vemm$ with
$\nr w-v\nr <\delta$ we have:
  $\gamma(x,t_0;w)\in R$ and
$\gamma(x,[0,t_0];w)\subset U$.$\qs$
\end{propcit}

\begin{propcit}{(\cite{pastpet}, Prop. 5.1)}
\lb{p:5.2}
Let $a,b\in M,t_0\geq 0$
and $\gamma (a,t_0;v)=b$.
Let $E$ be a submanifold
without boundary
of $M$
of codimension 1, such that
$b\in E$
and $v(b)\notin T_bE$.

Let $U$ be an open neighborhood
 of
$\gamma(a, [0,t_0];v)$ and $R$ be an
 open neighborhood of $b$.
 Then
there is $\theta_0>0$ \sut~ for every $\theta\in ]0,     \theta_0[$
there exist $\delta>0$ and an open neighborhood
 $S\subset U$ of $a$, such that $\forall x\in S$
and $\forall w\in\vemm$ with $\nr w-v\nr<\delta$
we have:
\begin{enumerate}
\item $\gamma(x,[-\theta, t_0+\theta];w)\subset U$,
and $\gamma(x, [t_0-\theta,t_0+\theta];w)\subset R$.
\item There is a unique
$\tau_0=\tau_0(w,x)\in [t_0-\theta,t_0+\theta]$,
 such that
$\gamma(x,\tau_0;w)\in E$.
\item If $E$ is compact and $t_0$
is the unique $t$ from $[0,t_0]$
 such that $\gamma(x,t_0;v)\in E$,
 then $\forall y\in S$
the number $\tau_0(w,y)$ is the
unique $\tau$ from
$[-\theta,t_0+\theta]$ such that
$\gamma(y,\tau;w)\in E$.$\qs$
\end{enumerate}
\end{propcit}

For the further use we shall need a version of Proposition \ref{p:5.1};
this version 
(Proposition \ref{p:co} below)
is proved by the same method, so we shall expose
only briefly the proof.

\begin{defi}
We denote by $\Vert\cdot\Vert$ the $sup$-norm on $\vemm$.
Let $\d>0, v\in\vemm$.
We denote by
$B(v,\d)$ the set of all
$w\in\vemm$, \sut~ $\Vert w-v\Vert<\d$.
The topology, induced in $\vemm$
by the $sup$-norm will be called $C^0$-\top.
\end{defi}

\bepr\lb{p:co}
Denote by
$\EE:
M\times\RRR \times\vemm
\to
M$
the map defined by
$$
(x,t,v)\mapsto \g(x,t;v)
$$
Then $\EE$ is continuous \wrt~the $C^0$ \top~in $\vemm$.
\enpr

\Prf
Let
$(a,t_0,v_0)
\in
M\times\RRR \times\vemm
$; 
let $b=\g(a,t_0;v_0)$.
The continuity of $\EE$ at the point
$(a,t_0, v_0)$
is equivalent to the following assertion:

for every \nei~ $U$ of $b$ there are:
a \nei~ $S$ of $a$, a number $\d>0$, a number $\theta>0$, \sut

~$$
\(x\in S, w\in B(v,\d), \vert t-t_0\vert<\theta\)
\Rightarrow
\(\g(x,t;w)\in U\)$$

Subdividing the interval $[0,t_0]$ into small intervals
it is easy to see
that we can restrict our assertion
to the case when
the curve
$\g(a,[0,t_0];v_0)$
belongs to the domain of a single chart
$\Phi:V\to V'\sbs\RRR^n$.
This case is reduced to the case of \vfs~ with compact support
in $\RRR^n$.
This last case is settled via the next lemma, which in turn is a direct
consequence of
\cite{biro}, \S 4, Th. 3.

For  $v\in \verr$
we denote by $\Vert v\Vert_1$
the $sup$-norm of the derivative
$dv:\RRR^n\to L(\RRR^n ,\RRR^n)$.

\begin{lemm}
Let $u,w\in\verr ,
\quad\nr u -w\nr<
\alpha,\quad\nr u\nr_1\leq D$,
where $D>0$.
Let $\gamma,\eta$
be trajectories of, respectively, $u,w$,
and assume
that $\vert \gamma(0)-\eta(0)\vert
\leq\epsilon$.
Then for every $t\geq 0$ we have:
 $\vert\gamma(t)-\eta(t)\mid\leq
\epsilon e^{Dt} +\frac \alpha D (e^{Dt} -1)$. \quad$\square$
\end{lemm}

\begin{coro}
\lb{c:cont3}
Let $a,b\in M, t_0\geq 0, \g(a,t_0;v)=b$.
Let $E$ be a submanifold
without boundary
of $M$
of codimension 1, such that
$b\in E$
and $v(b)\notin T_bE$.
Let $U$ be an open \nei~ of $b$.
There is $\theta_0>0$, \sut~ for every
$\theta\in]0,\theta_0[$
there is
$\d>0$
and a \nei~ $S$ of $a$
and a function
$\tau_0:B(v,\d)\times S\to \RRR$, \sut
~\been\item
$\forall (w,x)\in \bv\times S$
we have
$y\in\gam xw \in E\cap U$
and
$w(y)\notin T_yE$.
\item
$\tau_0(w,x)\in ]t_0-\theta, t_0+\theta[$
and $\tau_0(w,x)$
is the unique number $\tau$
of this interval \sut~ $\g(x,\tau;w)\in E$.
\enen
\end{coro}

\Prf
Choose a \nei~ $R\sbs U$
of $b$ and $\d_0>0$, \sut~for every $w\in B(v,\d_0)$ and every
$y\in R\cap E$ we have $w(y)\notin T_yE$.
Then apply  Proposition \ref{p:5.2}
with this $R$.
There is $\theta_0>0$, \sut
~for every
$\theta\in]0,\theta_0[$
there is
$\d>0$ and a \nei~
$S$ of $a$ and a function
$\tau_0: \bv\times S\to\RRR$, satisfying the conclusion
of \ref{p:5.2}.
For $(w,x)\in\bv\times S$
set
$y(w,x)=\gam xw$.
It follows from
our choice of $R$, that
$  w(y(w,x))   \notin    T_{y(w,x)}E   $.
It remains only to diminish $\d$ and $S$ in order to
force $\tau_0(w,x)$ to belong to
$]t_0-\theta, t_0+\theta[$
instead of the interval
$[t_0-\theta, t_0+\theta]$
from (2).
For this just take $\d'$ and $S'$ corresponding to $\theta/2$.$\qs$

\bepr\lb{p:coco}
Let $a,b\in M, t_0\geq 0, \g(a,t_0;v)=b$.
Let $E$ be a submanifold
without boundary
of $M$
of codimension 1, such that
$b\in E$
and $v(b)\notin T_bE$.
Let $\d>0, \theta>0$, let $S$ be a \nei~of $a$, and assume that
$\tau_0:\bv\times S\to\RRR$
is a function
satisfying (1) and (2) from the Corollary \ref{c:cont3}.
Then $\tau_0$ is continuous on its domain.
\enpr
\Prf Let $x_0\in S, w_0\in\bv$, and consider the $w_0$-trajectory
$\g(x_0, \cdot; w_0)$.
Set $t_0'=\tau_0(w_0, x_0)$ and
$y_0=\g(x_0, t_0'; w_0)$.
Then $w(y_0)\notin T_{y_0} E$.
Thus we can apply the preceding Corollary to the trajectory
$\g(x_0, \cdot; w_0)$.
Let $\e>0$.
We apply Corollary \ref{c:cont3}
with $\theta'_0$, satisfying
$\theta'_0<\e, ~
[t_0'-\theta',
t_0'+\theta']\sbs
]t_0-\theta,
t_0+\theta[$.
We obtain  a number $\d'>0$
a \nei~ $S'$ of $x_0$ and a unction
$\tau':B(w,\d')\times S'\to\RRR$, \sut~
\been\item
For every $(w,x)\in B(w_0, \d')\times S'$ we have:
$y=\g(x,\tau'(w,x);w)\in E$ and
$w(y)\notin T_y E$
\item
$\tau'(w,x)\in
]t_0'-\theta',
t_0'+\theta'[$
and
$\tau'(w,x)$ is the unique number $\tau$ of this interval
\sut~ $\g(x,t;w)\in E$.
\enen

We can diminish $\d'$ and $S'$ so that
$B(w_0, \d')\sbs \bv$ and
$S'\sbs S$.
Then the domain of definition of $\tau'$
is in the domain of definition
of $\tau_0$ and our choice of $\theta'$ implies that
$\tau'(w,x)=\tau_0(w,x)$.
The condition (2) implies then that
$\tau_0(w,x)\in
]t_0-\e,
t_0+\e[$
for
$(w,x)\in B(w_0, \d')\times S'$.
Since $\e>0$ was arbitrary this proves continuity.
$\qs$

In the proposition \ref{p:comap}
we assume that the \sma~ $E$
and the \vf~ $v$ satisfy the following assumptions:
\been\item[$(\EE 1)$]
There is $\a>0$ \sut~ for every $w\in B(v,\a)$
every $w$-trajectory intersects $E$ at most once.
\item[$(\EE 2)$]
for every $a\in E$
we have:
$v(a)\notin T_a E$.
\enen

\bepr\lb{p:comap}
Let $K\sbs M$ be a compact.
Let $U$ be an open subset of $E$.
 Assume that for every $x\in K$
the \tr~$\g(x,\cdot; v)$ reaches $E$ at the moment
$\tau(x)$
and that $\g(x,\tau(x);v)\in U$.
Then there is $\d>0$ and a continuous function
$\tau:\bv\times K\to\RRR_+$, \sut~
\been\item
For every $(w,x)\in\bv\times K$ we have:
$\ga xw \in U$
\item
Denote by $A:\bv\times K\to U$
the map defined by
$A(w,x)=\ga wx$.
For $w\in \bv$
denote by $A_w$ the map $K\to U$
defined by
$A_w(x)=A(w,x)$.
Then $A$ is continuous, and $A_w$ is continuous \fe~$w$ and $A_w$
is homotopic to $A_v$ \fe~ $W$.
\enen
\enpr
\Prf
Apply Corollary \ref{c:cont3}
and choose
\fe~ $a\in K$
a \nei~ $S(a)$ of $a$, a number $\d(a)$ and a function
$\tau_{(a)} :B(v,\d(a))\ti S(a)\to \RRR_+$
\sut~ \fe~ $x\in S(a)$
and
$w\in B(v, \d(a))$ we have:
$\gama wx  \in U$.
Moreover the function $\tau_{(a)}$
is continuous at $(v,a)$. We assume that $\d(a)<\a$ \fe~
$a$.
Note that if a pair $(w,x)$ is in the intersection
of
$B(v, \d(a))\ti S(a)$
and
$B(v, \d(b))\ti S(b)$
(where $a,b\in K$)
then
$\tau_a(w,x)=\tau_b(w,x)$
(it follows from $(\EE 1)$).
Choose a finite covering of $K$ by sets
$S(a_i), 1\leq i\leq N$,
and set
$\d=\min_i \d(a_i), S=\cup_i S(a_i)$.
We obtain a continuous function
$\tau:\bv\ti S\to \RRR_+$.
This function is continuous.
Now the continuity of the map $A$, defined
 in the point 2) of our Proposition is obvious:
just apply Proposition \ref{p:co}
and get that $A$ is a composition of continuous maps.
To obtain a homotopy
between $A_w$ and $A_v$
set
$H_t=A_{v+t(w-v)}$.
$\qs$

\subs{Manifolds with boundary}
\lb{su:mb}

Let $W$ be a compact \ma~with boundary $\pr W$.
Let
$\phi:\pr W\ti[0, \infty[\to W$
be a collar of $\pr W$
(see \cite{hirsch}, p. 113 for the definition of
collar and the theorem of existence).
Consider the smooth \ma~ $W'$ without boundary
obtained by gluing $W$ to
$\pr W\ti ]-1, 1[$
via the diffeomorphism
$\phi | \pr W\ti[0,1[$.
The manifolds
$W, \pr W\ti]-1,1[$
are embedded in the obvious way to $W'$, so we
shall consider
$W, \pr W\ti]-1,1[$
as submanifolds of $W'$;
Let $\xi:\pr W\ti]-1,1[\to ]-1,1[ $
be the projection onto the second factor.
Denote by $\vew$ the set of
all $C^1$ vector fields $v$ on $W$ \sut~ that \fe~
$x\in \pr W$
we have:
$v(x)\notin T_x(\pr W)$.
Let
$v\in \vew$.
Denote by $\pr^+_v W$, resp. by
$\pr^-_v W$
the set of $x\in \pr W$, \sut~ $v(x)$ points
inward $W$, resp. outward $W$.
(If the value of $v$ is clear from the context we shall write simply
$\pr^- W, \pr^+ W$.)
Since $W$ is compact, the set
$\vew$
is open in
$\mbox{ Vect } ^1(W)$
\wrt~  $C^0$ topology.
Moreover for given $v\in\vew$
there is $\d>0$
\sut~
for $\Vert v-w \Vert<\d$
we have:
$\pr^+_v W=\pr^+_w W, \quad
\pr^-_v W=\pr^-_w W$.

A \vf~ $u$ on $W'$ will be called {\it good extension}
of $v$, if
\been\item[(E1)]
$u | W=v$;
\item[(E2)]
$\supp u $ is compact;
\item[(E3)]
the restriction of $u$ to
$\pr_0 W\ti [0, 1[$ satisfies
$\xi'(x,t)(u(x,t))\geq 0$
and
\break
$\xi'(x,t)(u(x,t))> 0$
\fe~ $x$ and $t$ in some interval $[0,\b]$.
\item[(E4)]
the restriction of $u$ to
$\pr_1 W\ti [0, 1[$ satisfies
$\xi'(x,t)(u(x,t))\leq 0$
and
\break
$\xi'(x,t)(u(x,t))< 0$
\fe~ $x$ and $t$ in some interval $[0,\b]$.
\enen

\bere\lb{r:ext}
It is not difficult to show using the partitions of unity
 that for any
$v\in \vew$ good extensions of
$v$ exist.
Moreover, one can show that if
$v,w\in \vew$
and
$\pr_v^{\pm} W=
\pr_w^{\pm} W$,
then one can choose good extensions
$\wi v, \wi w$ of $v$, resp. of $w$,
\sut~
$\Vert \wi v - \wi w\Vert_{W'} \leq 2\Vert v-w \Vert_{W}$.
\enre

\bele\lb{l:e1e2}
Let $v\in\vew$
and let $\wi v$
be a good extension of $v$.
Then the vector field
$\wi v$ and the submanifold $\pr^+ W\sbs W'$ satisfy
$(\EE 1)$ and $(\EE 2)$
from Subsection \ref{su:mwb}.
\enle
\Prf
The condition $(\EE 2)$ goes by definition.
To prove
$(\EE 1)$
note first that every $\wi v$-\tr~ $\g$, which intersects
$\pr^+ W$
at a moment say $t_0$,
leaves necessarily
the domain $\pr^+ W\ti[0,\b]$, passing by a point of
$\pr^+W\ti \{\b\}$.
Afterwards it stays forever in
$\pr^+W \ti[\b, 1[$.
So $(\EE 1)$
holds for $v$.
If $w$ is sufficiently $C^0$ close to
$\wi v$
the property $(E 3)$ above is still true for $w$
and thus $(\EE 1)$
holds for $w$.

\begin{coro}\lb{co:reach}
Let $v\in \vew$. let $x\in W$ and assume that
the \tr~$\g(x,\cdot;v)$ reaches $\pr^+ W$.
Then there is an open \nei~ $S$ of $x$ \sut~  
$\forall y\in S$
the \tr~  $\g(y,\cdot; v)$
reaches $\pr^+ W$.$\qs$
\end{coro}

\bede\lb{d:ascend}
For $x\in W, t\geq 0$ define the element
$\G(x,t;v)\in W$ as follows:
\been
\item
if $\g(x,\cdot; v)$
is defined on $[0,t]$
then $\G(x,t;v)=\g(x,t;v)$
\item
if $\g(x,t_0; v)\in \pr^+ W$
for some $t_0\leq t$, then
$\G(x,t;v)=\g(x,t_0;v)$
\enen
\end{defi}

\begin{coro}\lb{co:ascend}
The map
$(x,t)\mapsto \G(x,t;v)$
is continuous.
\end{coro}
 
\Prf
Let $x\in W$.
If 
$\g(x,\cdot;v)$
does not reach $\pr^+ W$
then the continuity of $\G$ in $(x,t)$ for every $t$
is obvious.
If $\g(x, \cdot;v)$
reaches $\pr^+ W$, then
the same is true for any $y$ in a \nei~  
of $x$;
denote the moment of intersection
of 
$\g(x,\cdot;v)$
with $\pr^+ W$
by $\tau_0(x)$. Then
$\G(x,\cdot;v)=
\g(x,\min(t, \tau_0(x));\wi v)$
and the continuity is proved. $\qs$

\begin{coro}
\lb{co:homot}
Let $v\in \vew$.
Let $U$ be an open subset of $\pr^+ W$, and
$K\sbs W$ be a compact.
Assume that \fe~ $x\in K$ the $v$-\tr~
$\g(x, \cdot; v)$ starting at $K$ reaches $\pr^+ W$
 and intersects it at $\a_v(x)\in U$.

Then the map $x\mapsto \a_v(x)$ is a continuous map $K\to U$.
There is is $\l>0$ \sut~
\fe~ $w$ with $\Vert w-v\Vert <\l$ every $w$-trajectory starting at $K$
reaches $\pr^+ W$, intersects it at a point $\a_w(x)\in U$, and the map
$x\mapsto \a_w(x)$ is a continuous map $K\to U$ homotopic to $\a_v$.
\end{coro}

\Prf
Let $\wi v$ be a good extension of $v$.
The \vfs~$\wi v$ and the \sma~ $\pr^+ W$ of $W4$
satisfy $(\EE 1)$ and $(\EE 2)$, as it follows from
Lemma \ref{l:e1e2}.
Therefore Proposition \ref{p:comap}
applies, and we obtain a number
$\d>0$
and a continuous function
$A: B(\wi v,\d)\ti K\to U$, \sut
~$A(u,x)=\g(x,\tau(x,u);u)$.
Set $u=\wi v$; we obtain a continuous function
$A_{\wi v} :K\to U$
and a \cf $x\mapsto \tau(x,\wi v)$, which we shall
denote by $\tau_{\wi v}$.
Since $K\sbs W$ for every $x\in K$ we have
$\tau_{\wi v} (x)\geq 0$ and
$\g(x,[0, \tau_{\wi v} (x)]; \wi v)
\sbs W$.
Therefore
$\a_v(x)=A_{\wi v}(x)$
and the map
$\a_v$ is continuous.

Assume now that $\d$ is so small that \fe
~$\xi\in \veww$
with
$\Vert \xi-v\Vert_W<\d$
we have:
$\xi\in \vew$
and $\pr^\pm_\xi W=
\pr^\pm_v W$.
Let $\l<\d/2$; let $w\in \veww$
and $\Vert v-w\Vert_W<\l$.
Let $\wi w$ be a good extension of $v$, satisfying
$\Vert\wi v-\wi w\Vert_{W'}<\d$.
We have:
$\ga x{\wi w}\in U\sbs E$ for $x\in K$.
Note that $\tau(x,\wi w)\geq 0$
(indeed, if
$\tau(x,\wi w)<0$
then for some
$\theta\in]\tau(x,\wi w), 0[$
we have:
$\g(x,\theta; \wi w)\in \pr^+W\ti]0,1[$,
and the relation
$\g(x,0;\wi w)\in W$ is no more possible).
Further,
$\g(x,[0,\tau(x,\wi w)]; \wi w)\sbs W$, \th
~~
$\ga x{\wi w} = A(\wi w, x)=\a_w(x)\in U$
and $\a_w$ is continuous in $x$.

To prove that $\a_w$ is homotopic to
$\a_v$
for
$\Vert v-w\Vert_W<\l$, consider a \vf~
$\xi(t)=\wi v+t(\wi w - \wi v)$
where
$t\in [0,1]$.
Note that
$\Vert\xi(t) -\wi v\Vert_{W'}<\d$, \th~
\fe~$t\in[0,1]$
the \vf
~$\xi(t)$
is in
$B(\wi v, \d)$.
The map
$[0,1]\ti K\to U$
given by the following
formula:
$(t,x)\mapsto A(\xi(t), x)$
is then a \hot~between
$\a_w$ and $\a_v$.
$\qs$

For some further applications 
we shall need a generalization of Corollary
\ref{co:homot}
to the case
of pairs of compacts $(K_1, K_2)$. For the 
further references we shall formulate it as a separate statement.
The proof is immediate generalization of
the proof of \ref{co:homot}.

\begin{coro}
\lb{co:hompair}
Assume the hypotheses and terminology of Corollary \ref{co:homot}.
Let $K'$ be a compact subset in $K$, and $U'$ be an open subset of $U$.
Assume that $\a_v(K')\sbs U'$.

There is is $\l>0$ \sut~
\fe~ $w$ with $\Vert w-v\Vert <\l$ every $w$-trajectory starting at $K'$,
resp. at $K$ 
reaches $\pr^+ W$, intersects it at a point $\a_w(x)\in U'$,
resp $\a_w(x)\in U$,
 and the map
$x\mapsto \a_w(x)$ defines  a continuous map 
of pairs 
$(K, K')\to (U, U')$, which is homotopic to 
 $\a_v$, considered as a map of pairs
$(K, K')\to (U, U')$.$\qs$
\end{coro}

\subsection{Gradients of Morse functions}
\lb{su:gmf}

In this subsection $\fcob$
is a \Mf~on a riemannian cobordism
$W$, $v$ is an $f$-gradient.
We start with an elementary lemma.

\bele\lb{l:adjpert}
Let $\phi:W\to[\a,\b]$ be a \Mf~ adjusted to
$(f,v)$.
Then there is $\d>0$, \sut~ every
$f$-gradient
$w$, satisfying
$\Vert w-v\Vert<\d$ is also a $\phi$-\gr.
\enle
\Prf
To check the condition 1 from Definition \ref{defgrad}
for the pair $(\phi,w)$
(the other two conditions hold obviously), choose
for each $x\in S(f)$
a neighborhood $U(x)$
such that
$(f-\phi)\vert U=const$.
Set $U=\cup_xU(x)$
and
$K=W\sm U$.
We have: $d\phi(x)(v(x))>0$ for every $x\in K$.
Since
$K$ is compact, this implies
that for $w$ \sclv~we have:
$d\phi(x)(w(x))>0$
\fe~$x\in K$.
For every
$x\in U\sm S(f)$
we have:
$d\phi(x)(wx))=
df(x)(wx))>0$, and the condition 1 is proved. $\qs$

\begin{coro}
\lb{c:seppert}
\been\item
Assume that $v$ is a $\d$-ordered $f$-gradient. then
 every $f$-gradient $w$, \sclv
~is $\d$-ordered.
\item
Assume that $v$ is a $\d$-separated $f$-gradient.
Then every $f$-gradient $w$ \sclv~is $\d$-separated.$\qs$
\enen
\end{coro}

In the following two propositions
$\l\in[a,b]$ is a regular value of $f$,
$p$ is a critical point of $f$, satisfying $f(p)\geq \l$, $U$ is an open
subset of $f^{-1}(\l)$.
We denote $f(p)$ by $c$, and $f^{-1}(\l)$ by $V'$.

\bepr\lb{p:discpert}
Assume that \fe~
$x\in D(p,v)\cap f^{-1}\([\l, c]\)$
the $(-v)$-trajectory
$\g(x,\cdot; -v)$
reaches $V'$.
Assume that $D(p,v)\cap V'\sbs U$.

Then there is $\d>0$, \sut~\fe~ $f$-gradient
$w$, satisfying
$\Vert v-w\Vert<\d$
we have:

For every
$x\in D(p,v)\cap f^{-1}\([\l, c]\)$
the $(-w)$-trajectory
$\g(x,\cdot; -w)$
reaches $V'$,
and  $D( p,v)\cap V'\sbs U$.
\enpr
\Prf
Let $W'=f^{-1}\([\l,b]\)$.
For every $q\in S(f)\cap W'$
we have
$D(q,-v)\cap D(p,v)=\emp$,
\th~applying Rearrangement Lemma \ref{rearrlemm} several times
we can find a \Mf~
$\phi:W'\to[\l, b]$, adjusted to $(f,v)$
and a regular value $\mu$
of $\phi$, \sut~
$p$ is the unique critical point
of $\phi$
in the domain $\phi^{-1}\([\l,\mu]\)$.
Consider the cobordism
$W''=\phi^{-1}\([\l,\mu]\)$.
Recall from Lemma \ref{l:adjpert}
that any $f$-gradient $w$, \sclv~
is also a $\phi$-\gr.
Therefore for
$w$ \sclv~
a $(-w)$-\tr~starting at $p$ necessarily reaches $V'$, since
there is only one critical point of $\phi$
in $W''$.
Further, let $K=V'-U$.
Every $v$-\tr~starting at $K$ reaches
$\phi^{-1}(\mu)=\pr_1 W''$.
Therefore the same is true \fe~ \vf~$w$, \sclv
(by Corollary \ref{co:homot}).
Thus $D(p,w)\cap V'\sbs V'\sm K=U$                $\qs$

Assume  the hypotheses of the preceding proposition.
Choose an orientation of $D(p,v)$.
Set
$\Sigma(v)=D(p,v)\cap V'$.
Then $\Sigma(v)$
is an embedded $s$-dimensional oriented sphere in $V'$,
where
$s+1=\ind p$
(see Lemma \ref{discisdisc}).
Therefore its fundamental class defines an element
$\s(v)\in H_s(U)$.
By Proposition \ref{p:discpert}
the same is true for any $f$-gradient $w$ if only $w$ is
\sclv.

\bepr
\lb{p:dischom}
Assume the hypotheses of Proposition \ref{p:discpert}.
Let $U$ be 
a \nei~of $p$.

Then there is  $\d>0$, \sut~ \fe~$f$-\gr~ $w$, satisfying 
$v|U=w|U$ and 
$\Vert v-w\Vert<\d$
we have:
$\s(v)=\s(w)$.
\enpr
\Prf
We assume the terminology from the proof of
\ref{p:discpert}.
Set $\nu=\phi(p)$
and let
$\e>0$ be so small that
$D(p,v)\cap\phi^{-1}\([\nu-\e, \nu]\)\sbs U$.
Set
$\Sigma_0(v)=D(p,v)\cap \phi^{-1}(\nu-\e), \quad
\Sigma_0(v)=D(p,v)\cap \phi^{-1}(\nu-\e), \quad$.
Then
$$
\Sigma(v)=
\stind {(-v)}{\nu-\e}\l (\Sigma_0(v)),\quad
\Sigma(w)=
\stind {(-w)}{\nu-\e}\l (\Sigma_0(w))
$$

Note that the condition $v|U=w|U$
implies 
$\Sigma_0(v)=
\Sigma_0(w)$, \th~ to obtain the conclusion of our Proposition
just apply the Corollary
\ref{co:homot}
to the 
\cob~ $\phi^{-1}\([\l, \nu-\e]\)=W''$
and the vector field
$(-v)| W''$.
$\qs$

\chapter{Quick Flows}\label{ch:qf}

Let $M$ be a closed \ma~ of dimension $n$,
 $f$ be a Morse function on $M$,
$v$ be an $f$-gradient
satisfying the Transversality Condition
Consider the union of all the descending discs of
indices $\leq k$, denote it by $M_k$.
These descending discs form a filtration of $M$,
which is a close analog of
cellular decomposition. (In many cases
the
descending discs 
 {\it do} give a real cell decomposition, but
we  neither study this subject
in our paper, nor use the
corresponding results.)
The ascending discs form a "dual decomposition", denoted by
$\{\wh M_k\}_{0\leq k\leq n}$.
The transversality condition implies that
$M_k\cap \wh M_{n-k-1} =\emp$.
Let $U$ be an open \nei~  of $M_k$,
and $V$ be an open \nei~
of
$\wh M_{n-k-1}$. It is clear that the gradient descent map
$\Phi(T,-v)$
pushes $M\sm V$ inside $U$, if only
$T$ is large enough.
The main aim of this section
is to find for any $M$ such a function $f$ and
its gradient $v$, and the  neighborhoods $U,V$,
that
the number $T$ from  above and the $C^0$ norm of $v$ are small.
We demand also that $V$ and $U$ be sufficiently small neighborhoods of
$M_k$, respectively of
$\wh M_{n-k-1}$. 
The flows generated by such gradients are called {\it rapid flows},
since they push $M\sm V$ to $U$ rapidly.

The precise statement of the result is  Theorem \ref{t:quickpush}
which is the main result 
 of the present chapter. 
This theorem will be used in the proof of $C^0$-density 
of condition $(\gC)$
(Theorem \ref{t:cc});
this last theorem is essential for the proof of Theorem \ref{t:cccgr}.
Later on we  shall use only
the statement of Theorem \ref {t:quickpush}, but
not the proof, so  the reader interested only in the applications
of the theorem can read  the present introduction and 
Subsection \ref{su:tergf} and skip the rest of the chapter.

One may try to give a construction of rapid flows  
using triangulations of $M$. Indeed,
choose a smooth triangulation of $M$.
Pick a Morse function $f$ on $M$, \sut~  the critical
 points of $f$ of index
$p$
are the barycenters of the simplices of the
triangulation.
(The existence of such a function seems to
be known in folklore since sixties; there is a recent
work \cite{pozniak} of M.Pozniak on this subject.)
It seems that if the triangulation of $M$ is sufficiently
fine (that is, the simplices are of small diameter)
the function    $f$ and
some $f$-gradient would give a solution of our problem.

I have chosen another approach, in which it is easier to
obtain explicit estimates for times and norms. I do not make any use
of triangulations of the manifold.
Instead this approach suggests an alternative construction
of "small handle decomposition" of  a   manifold, which may prove 
useful in other domains
of topology.

Now I shall give a brief informal description of the construction.
Start with an arbitrary Morse function $f:M\to\RRR$
\sut the values of $f$ at all the critical points are different. Then
choose a finite sequence of regular values
$\a_1,...,\a_N$ of $f$, such that between $\a_i$ and $\a_{i+1}$
there is at most one critical value of $f$. We then think of $M$ as 
of the union
of slices $W_i=f^{-1}([\a_i,\a_{i+1}])$, which we choose to be 
very thin.
Consider a slice $W_i$ which has no critical points inside it.
It looks like a cylinder $V_i\times [0,1]$, where $V_i=f^{-1}(\a_i)$
(see the figure A on the page  \pageref{fig:skl} ). 
Now modify $f$ inside $W_i$ so as
to make a fold (see the figure B on the page \pageref{fig:skl}
 )\footnote{ "fold" is  "skladka"
in Russian,
 which explains the term "S-construction"}.
Glue  together the functions $\widetilde f_i$ 
and obtain a smooth function on $W$. 
Note  that
the trajectories of its  gradient have become much shorter:
they can  intersect no  more than three successive  slices $W_i$.
However each  function $\widetilde f_i$ is not a Morse function:
there are two level surfaces,
consisting only of critical points. This must be repaired by
 perturbing $\widetilde f_i$
nearby these level surfaces with the help of an auxiliary
Morse function on $V_i$ (see the figure C on the page 
\pageref{fig:skl} ).
 Gluing together the perturbed  functions 
we obtain the Morse function sought.

What is going on in one slice is the subject of section \ref{s:cc}
 of this chapter,
"gluing back" is done in section \ref{s:pteqf}, which contains
also the end of the proof of 
the Theorems \ref{t:sq} and \ref{t:q}.
Section \ref{s:isr}  contains all the new terminology we need.
Unfortunately we do need a lot, since the precise
estimates of times which a trajectory of a gradient
of a Morse function
can spend in different domains is not a standard item.
In particular we shall  distinguish between 
{\it rapid flows} (Definition \ref{d:rapid}) and {\it quick flows} 
(Definition \ref{d:quickflo}). In both definitions 
the main item is
 the time which 
a trajectory of a flow can spend outside a \nei~
$U$ of $S(f)$; but in the definition of
rapid flows $U$ is the union of some riemannian balls
around the critical points of $f$, and
in the definition of quick flows $U$ is the union of
coordinate neighborhoods corresponding to Morse charts.
In the section \ref{s:isr} we also state the main result of this
 chapter
(Theorem \ref{t:quickpush}) and do  first steps in its proof.

\section{Introduction and statement of results}\label{s:isr}

\subsection{Theorem on the existence of rapid gradient flows}
\label{su:tergf}

In this subsection $f:M\to\RRR$ is a Morse function
on a closed manifold of dimension $m$.

\begin{defi}\label{d:rapid}
Let $v$ be an $f$-gradient, satisfying \ATA.
Let $\e>0, t\geq 0$.
We say that the flow, generated by $v$ is {\it $(t,\e)$-rapid},
if
for every $s:0\leq s\leq m$ we have:
\begin{equation}
\Phi(v,t)\(M\sm B_\e(\indl s;v)\)\sbs
B_\e\(\indl {m-1-s}; -v\)\lbl{qpu1}
\end{equation}
\begin{equation}
\Phi(-v,t)\(M\sm B_\e(\indl s;-v)\)\sbs
B_\e\(\indl {m-1-s}; v\)\lbl{qpu2}
\end{equation}
\end{defi}

Of course the definition would make sense for
arbitrary $f$-gradient, not necessarily satisfying
\ATA.
The reason to restrict ourselves to
the gradients, satisfying \ATA,
is the following lemma.

\begin{lemm}\lb{l:rap}
For every $\e>0$ there is
$t\geq 0$
\sut~ the flow generated by
$v$ is $(t,\e)$-rapid.
\end{lemm}
\Prf
If the lemma is true for some
$\nu$, then it is true for every $\nu'\geq \nu$,
so it is sufficient to prove it for
such $\nu$, that $v$ is $\nu$-ordered. For such $\nu$
consider the compact subset
$$
K=M\sm\bigg( B_\e(\indl {m-s-1}, -v)\cup
\bigcup_{\ind p\leq s} B_\e(p)\bigg)
$$
There are no critical points in $K$, therefore
$df(v)(x)>C>0$
for
$x\in K$. This implies the existence of $D>0$, \sut~
for every $x\in K$
there is $t_0\leq D$, \sut~
we have:
$\g(x,t_0;-v)\in M\sm K$.
Since $M\sm B_\e(\indl {m-s-1}, -v)$
is $(-v)$-invariant
we have:
$\g(x,t_0,-v)\in\cup_{\ind p\leq s} B_\e(p)$,
and since
$B_\e(\indl s,v)$ is
also $(-v)$-invariant, we have
\break
$\g(x,D;-v)\sbs B_\e(\indl s, v)$.$\qs$

Similarly, given $t$ and $\e$ one easily constructs an
 $f$-gradient which is $(t,\e)$-rapid, but
has maybe a very large $C^0$ norm.
(Indeed, multiply a given $f$-gradient
satisfying \ATA~
by a large constant outside $B=\cup_{\ind p\leq s} B_\e(p)$
and smooth the result nearby the boundary $\pr B$.)
The aim of the construction of rapid flows
is to produce flows $v$ which have small $C^0$ norm and are
$(t,\e)$-rapid for a small $t$ and $\e$.
This is achieved by the theorem \ref{t:quickpush}
which is the main result of the present chapter.
We need a definition.

\begin{defi}\label{d:equiv}
Let $\fcob$ be a Morse function on
a compact cobordism, and $u,v$ be $f$-gradients.
We say that $u$ and $v$ are
{\it equivalent},
if $u(x)=\phi(x)v(x)$,
where $\phi:W\to\RRR$ is a $\smo$ function, \sut~$\phi(x)>0$
for every $x$
and $\phi(x)=1$
for $x$ in a \nei~ of $S(f)$.
\end{defi}

\begin{theo}\label{t:quickpush}

Let $M$ be a closed riemannian manifold. Let $C>0, t>0$.

There is a Morse function
$f:M\to\RRR$ and an $f$-gradient $v$, satisfying \ATA~
and for every $\e>0$ there is
a $(t,\e)$-rapid $f$-gradient $u$ with $\Vert u\Vert\leq C$,
equivalent to $v$.
\end{theo}

We shall deduce Theorem \ref{t:quickpush} from Theorem \ref{t:sq},
see page \pageref{pftrtsq}.

\subsection{ Terminology}\label{su:t}
The rapid $f$-gradients which we shall construct will be \glvf s.
(The time which a $v$-\tr~ spends in some domain nearby a critical point
of $f$ is much easier computed 
or estimated for a \glvf~than for a general $f$-\gr.)

Let $n, s$ be positive integers and $0\leq s\leq n$; set
\begin{equation}
\QQ_s(x)=
-\sum_{i=1}^s x_i^2 -
\sum_{i=s+1}^n x_i^2
\end{equation}
\begin{equation}
\gotn_s(x)=
(-x_1,...,-x_s,x_{s+1},..., x_n)
\end{equation}

\begin{defi}\label{d:defmorse}

{\bf 1. Quasi-Morse charts ($qM$-charts) for Morse functions.}

Let $p$ be a critical point of $f$ of index $s$.

A chart $\Phi:U\to B^n(0,r)$, where
$U$ is a neighborhood
of $p$, and $r >0$ is called {\it quasi-Morse
chart for $f$ around $p$ of radius $r$}
 (or simply $qM$-{\it chart for $f$} ) if $\Phi(p)=0$ and
there is an extension of $\Phi$ to a chart
$\widetilde\Phi:V\to B^n(0,r')$,~
where $\overline {U}\subset V\aand r'>r$,
such that
$$(f\circ \widetilde\Phi^{-1})(x_1,...,x_n)=f(p)+
\sum_{i=1}^n\alpha_i x_i^2$$
 with
 $\alpha_i<0$ for $i\leq s$ and
$\alpha_i>0$ for $i> s$.
 The domain $U$
is called {\it quasi-Morse coordinate neighborhood.
}
Any such extension $\widetilde\Phi$ of $\Phi$
will be called
{\it quasi-Morse  extension of $\Phi$.}
We shall often abbreviate "quasi-Morse" to $qM$.
The set $\Phi^{-1}(\RRR^s\times\{0\})$, resp.
$\Phi^{-1}(\{0\}\times\RRR^{n-s})$,
is called {\it negative disc, }
resp. {\it positive disc. }

\pa
{\bf 2. Quasi-Morse chart systems.}

Assume that for every critical point $p\in S(f)$ there is given
a quasi-Morse chart
$\{\Phi_p:U_p\to B^n(0,r_p)\}_{p\in S(f)}$.
The family
${\UU}= \{\Phi_p:U_p\to B^n(0,r_p)\}_{p\in S(f)}$
 is called {\it quasi-Morse  chart system} ,
 if the family $\{\overline{U_p}\}$ is disjoint
and no one of $\ove{U_p}$ intersects $\pr W$.
We denote $\min_p r_p$ by $d({\UU})$, and
 $\max_p r_p$ by $D({\UU})$.
If all the $r_p$ are equal to $r$, we shall say
that $\UU$
{\it is of radius $r$. }
The set $\Phi_p^{-1}(B^n(0,\lambda))$, where
 $\lambda\leq r_p$ will be denoted by $U_p(\lambda)$.
For $\lambda\leq d(\UU)$ we denote
$\cup_{p\in S(f)}U_p(\lambda)$ by $\UU(\lambda)$.
In the case when $W$ is riemannian, we       denote
       $\max_p\GG(U_p,\Phi_p)$ by
$\GG(\UU)$.
Let $\UU=\{\Phi_p:U_p\to B^n(0,r_p)\}_{p\in S(f)}$,
$\UU'=\{\Phi'_p:U'_p\to B^n(0,r'_p)\}_{p\in S(f)}$
be two $qM$-chart-systems for $f$. We say, that $\UU'$
is
{\it
a restriction of
}
 $\UU$, if for every $p\in S(f)$
we have:
 $r'_p\leq r_p, ~U'_p\subset U_p, \Phi'_p=\Phi_p\mid U'_p$.

\pa
{\bf 3. $f$-Gradients}
\label{d:Grad}

Given an $qM$-chart system
${\UU}= \{\Phi_p\}$
we say, that a vector field $v$ on $W$
 is an
{\it $f$-Gradient with respect to $\UU$,
}
if

1) $\forall x\in W \setminus S(f)$ we
have $df(v)(x)>0$;

2) $\forall p\in S(f)$ we have
$(\widetilde\Phi_p)_*(v)(x)=(-x_1,...,-x_k,x_{k+1},...,x_n)=
\gotn_k(x)
$,
where $k $ is the index of $p$, and $\widetilde\Phi_p$ is some $qM$
extension of $\Phi_p$.

Note that a vector field $v$ is a \glvf~for $f$ \ifff~it is
an $f$-Gradient \wrt~  some qM-chart system.
\pa

{\bf  4. $\MM_0$-flows,  $\MM$-flows, $\AA\MM_0$-flows, $\AA\MM$-flows}

If $\UU=\atlas$ is a $qM$-chart system for $f$ and $v$ is
an $f$-Gradient of $f$ \wrt~ to $\UU$, we shall say
that the triple
$\VV=(f,v,\UU)$ is an {\it $\MM_0$-flow}.
 We       say that $\VV$
{\it is
of radius $r$, }
if~  $\UU$ is of radius $r$. We set $d(\VV)=d(\UU),~D(\VV)=D(\UU)$.

If $W$ is riemannian,
we shall say that the triple
$\VV=(f,v,\UU)$ is an {\it $\MM$-flow} if moreover
 for each chart $\chart$ from $\UU$ the coordinate
 frame in $p$ is orthonormal with respect to
 the riemannian metric.

An $\MM_0$-flow  $(f,v,\UU)$ (resp. an $\MM$-flow)
will be called {\it $\AA\MM_0$-flow}
(respectively, $\AA\MM$-flow), if $v$ satisfies the almost 
transversality
condition.

Let $\flow$ be an $\MM$-flow on $W$,
where
$\UU=           \atlas$.
We       denote
       $\max_p\GG(U_p,\Phi_p)$ by
$\GG(\VV)$.
\pa
{\bf  5. Subordinate flows}

Let  $\VV_1=(f,v_1,\UU_1)$ and
$\VV_2=(f,v_2,\UU_2)$
be two $\MM_0$-flows on $W$
(the Morse function is the same for both), where
$\UU_1= \{\Phi_p:U_p\to
B^n(0,r_p)\}_{p\in S(f)}$, and
$\UU_2= \{\Psi_p:U_p\to
B^n(0,r'_p)\}_{p\in S(f)}$
are the corresponding $qM$-chart systems;
We say that $\VV_2$ is {\it subordinate to $\VV_1$},
if
\begin{enumerate}\item
$\UU_2$ is a restriction of $\UU_1$
\item$v_2=\phi\cdot v_1$, where
$\phi:W\to\RRR$
is a strictly positive
$\smo$ function
equal to $1$ in a \nei~
of the closure of $\cup_p U_q$.
\end{enumerate}

\end{defi}

(Note that if $\VV_2$ is subordinate to $\VV_1$ then $v_1$ and $v_2$
are equivalent in the sense
of
Definition \ref{d:equiv}.)

\bere

Here are two observations, which may serve as some justification
for having introduced so many definitions.
\begin{enumerate}
\item
As we have already noted (Remark \ref{remterm})
each $f$-gradient is actually a \glvf.
Still we keep the qM-charts in our 
terminology, since we shall often consider
different \Mf s with the same $f$-\gr~
(for example, we shall multiply \Mf s by a constant).
\item
The condition or orthonormality in the definition
of $\MM$-flow
seems artificial, but it
 makes easier some constructions.
Namely, it implies that the riemannian $\d$-balls around critical 
points and the $\d$-balls \wrt~ standard charts
are very close to each other.
 The basic invariants of
such flows , like $N(\VV), T(\VV)$
(see Definition \ref{d:quickflo})
can be defined of course
for every $\MM_0$-flow.$\qt$
\end{enumerate}

\end{rema}

\begin{lemm}\label{l:qmtom}
Let $f:\RRR^n\to\RRR$ be a quadratic form
$f(x)=
\sum_{i=1}^n\alpha_i x_i^2$, where
 $\alpha_i<0$ for $i\leq s$ and
$\alpha_i>0$ for $i> s$.
(Note that $\gotn_s$ is an $f$-gradient.)

Let $R>R'>r'>r>0$. Then  there is a Morse function $g:\RRR^n\to\RRR$
\sut~
\begin{enumerate}\item $g(x)=f(x)$ outside $B^n(0,R')$
\item $g(x)=C\QQ_s(x)$ in $B^n(0,r)$
where $C>0$.
\item $S(g)=\{0\}$, and $g'(x)\(v(x)\)>0$ for every $x\not= 0$.
\end{enumerate}
\end{lemm}
\Prf Choose a $\smo$ function $h:\RRR^n\to[0,1]$
\sut~ $\supp h\sbs B^n(0, R')$ and $h(x)=1$ for $x\in B^n(0,r')$.
It is obvious that the function
$x\mapsto C h(x) \QQ_s(x) + (1-h(x)) f(x)$ satisfies
the conclusion of the lemma if only
$C$ is \sufsm~.
$\qs$

It follows from this lemma that if $v$ is an $f$-Gradient \wrt~
a $qM$-chart system
${\UU}= \{\Phi_p:U_p\to B^n(0,r_p)\}_{p\in S(f)}$
 and $\g$
is a $v$-trajectory, then for every $p\in S(f)$ and every $\l\leq r_p$
 the set $I=
\{t\in\RRR \mid \g(t)\in U_p(\l)\}$ is an open interval, empty, or
finite, or infinite from one side.

\begin{defi}\label{d:quickflo}
{\bf ( Quick flows) }

Let $\flow$ be an $\MM_0$-flow on $W$.
 Let $\gamma$ be a $v$-trajectory.

A) The number
of sets $U_p=U_p(r_p)$,
intersected by $\gamma$, will be denoted by
$N(\gamma)$. The number $\max_\gamma N(\gamma)$
will be denoted by $N(\VV)$.
The set
$\{t\in\RRR \vert \gamma(t)\notin\cup_{p\in S(f)}
U_p\}$ is a finite union of disjoint
closed intervals; the sum of their lengths
 will be denoted by $T(\gamma)$.
The number $\max_\g T(\g)$ is denoted by $T(\VV)$.

B)
Let $\beta>0,~C>0$. We say, that $\VV$ is
{\it
$(C,\beta)$-quick,
}
if
$\nrv \leq C$
and  $T(\VV)\leq\beta$.

C) Let
$\delta\leq d(\UU)$.
The number
of sets $U_p(\delta)$,
intersected by $\gamma$, will be denoted by
$N(\gamma,\delta)$. The number
$\max_\gamma N(\gamma,\delta)$
will be denoted by $N(\VV,\delta)$.
The set
$\{t\in\RRR \vert \gamma(t)\notin\UU(\delta)\}$
 is a finite union of
disjoint
closed intervals
and the sum of their lengths
 will be denoted by $T(\gamma,\delta)$.
The number $\max_\g T(\g,\d)$ will be denoted
by $T(\VV,\d)$.

D)
Let $\beta>0,~C>0, \d\leq d(\UU)$. We say, that $\VV$ is
{\it
$(C,\beta,\delta)$-quick,
}
if
$\nrv \leq C$
and  $T(\VV,\delta)\leq\beta$.

\end{defi}

\subsection{ Two lemmas on standard gradients in $\RRR^n$.}
\label{su:tlsg}

In  this subsection
$n,k$ are positive integers, and $0\leq k\leq n$.
We refer to ${\bold R}^n$ as to
the product $\RRR ^k\times \RRR ^{n-k}$;
a point $z\in \RRR ^k$ is therefore written down as
 $(x,y)$, where $x\in \RRR^k,y\in \RRR^{n-k}$.
The Euclidean norm in $\RRR^n$ is denoted by $\vert\cdot\vert$.

We have:
$\gotn_k(x,y)=(-x,y);\quad \QQ_k(x,y)= 
-\vert x\vert^2+\vert y\vert^2$.
The $\gotn_s$-trajectories are of the form
$(x_0e^{-t}, y_0 e^t)$. using this fact
it is not difficult to prove the following
two lemmas .

\begin{lemm}\label{l:standtime}
Let $R>r>0$ and $\gamma$ be a
$\gotn_s$-trajectory. Then the time,
which $\gamma$ spends in $B^n(0,R) \setminus
B^n(0,r)$ is not more than
$\ln \big( (\frac Rr)^2
 +\sqrt{(\frac Rr)^4 -1}~\big) $
 and the length
 of the corresponding part of $\gamma$ is
 not more than $2R$.
\end{lemm}

\Prf
Let $\g(t)=(x(t),y(t))$ be a $\gotn_s$-trajectory.
We can assume that both $x$ and $y$ are never equal to $0$.
(Indeed, if $x(t_0)=0$, then $x(t)=0$ for all $t$
and the time, respectively length, cited in the statement of 
the lemma, are equal to
$\ln \frac Rr$, respectively, to $R-r$.)
Set $\LL(x)=\ln(x+\sqrt{x^2-1})$.
There is a unique $t_0$, \sut~
$\vert x(t_0)\vert=\vert y(t_0)\vert$.
This point is the minimum of $\vert\g(t)\vert$.
Reparametrize $\g$ so as to have
$t_0=0$ and denote
$\vert x(0)\vert=\vert y(0)\vert$ by $\a$.
Assuming $\sqrt 2 \a\leq R$, we calculate
easily  that
$\g$ quits
$B(0,R)$
at the moment
$\frac 12 \LL(\frac {R^2}{2\a^2})$.
If $r<\sqrt 2\a$, we are over.
If not, the time spent between two spheres equals to
$\LL\(\frac {R^2}{2\a^2}\) -\LL\(\frac {r^2}{2\a^2}\)$
and it suffices to prove that
$\LL(x) -\LL(y)\leq \LL(\frac xy)$
for $x>y>1$,
but this is obvious (take the derivative in $x$).
To obtain the estimate of the length, just note that the 
length can not be more, than
the length of the intersection of
the 1-dimensional hyperbole
$H=\{(x,y)\in\RRR^2\mid xy=\a\}$
with $B^2(0,R)$.
$\qs$

\begin{lemm}\label{l:standtimetwo}
Let $r>0$ and $\gamma$ be a $\gotn_k$-trajectory.
 Then the time, which $\gamma$ spends in the set
$\QQ_k^{-1}([-r^2,r^2]) \setminus B^n(0,r)$ is
not more than $2$.
\end{lemm}
\Prf
The function
$\QQ_k'(x)(v(x))$
is bounded below in the domain by $r^2$
and the variation of $\QQ_k$ in this domain is not more than 2.
$\qs$

\begin{coro}\label{co:radrestr}
Let $\flow$ be an $\MM$-flow on
                        a cobordism.
 Assume, that
$\VV$ is $(C,\beta,\alpha)$-quick.
Then
 $\VV$ is $(C,\beta + 8N(\VV),\alpha /2)$-quick.
\end{coro}
{\it Proof.\quad}
By Lemma \ref{l:standtime}
the time, which a $v$-trajectory can
 spend in $\UU(\a)\sm\UU(\a/2)$
is not more than
$8N(\VV )$.
$\qs$

\subsection{ Theorems on the existence
of quick flows
}\label{su:teqf}

In this subsection $W$ is a riemannian cobordism of dimension $n$
($\pr_0 W$ and $\pr_1 W$ can both be empty).

\begin{theo}\label{t:sq}
Let $C>0,
\beta>0, A>1$.
There is an $\AA\MM$-flow $\flow$ on $W$ of radius $r$,
  and for every
$\mu\in ]0, r[$
there is an
$\AA\MM$-flow $\WW = (f,w,\UU ')$,
subordinate to $\VV$, and such that

\text{\rm (1)}\quad $\WW$ is of radius $\mu$ and
 is $(C,\beta )$-quick;\quad

\text{\rm (2)}\quad $N(\WW  ) \leq n2^n$;\quad

\text{\rm (3)}\quad $\GG (\WW )\leq A$.

\end{theo}

\begin{theo}\label{t:q}
Let $B>0,\mu_0>0$.
Then there is an  $\AA\MM$-flow $\flow$
on $W$ with the following properties:

\text{\rm (1)}\quad $\GG (\VV )\leq 2\aand D(\VV ) <\mu_0 $;
\qquad\qquad

\text{\rm (2)}\quad$N(\VV )\leq n2^n$;

\text{\rm (3)}\quad$\VV$ is $(B,S(n))$-quick,
where $S(n) = 30+n2^{n+5}$.
\end{theo}

The scheme of the proof of \ref{t:sq} and \ref{t:q} is
as follows. 
We shall use the notation \ref{t:sq}(n) which means 
"Theorem \ref{t:sq} for manifolds of dimension $n$.
The proof of \ref{t:q} $\Rightarrow$
\ref{t:sq} is given below. The proof of
 \ref{t:sq}(n) $\Rightarrow$  \ref{t:q}(n+1)
 is done in \S 2 and \S 3 with the help of
(S)-construction. The same argument proves also 
\ref{t:q} for 1-dimensional
manifolds.

{\it Proof of Theorem \ref{t:quickpush} from \ref{t:sq} }\label{pftrtsq}

Let $(f,v,\UU)$ be the $\AA\MM$-flow of radius $r$, satisfying
the conclusion of \ref{t:sq} \wrt~ the manifold $M$, the numbers
$C,t$ and $A=2$.
Let $\k>0$ be any number, so small that
\begin{enumerate} \item
$\k<r$ and $\k<\e$
\item $v$ is $\k$-ordered
\item for every $p\in S(f)$ the disc $D_\e(p)$ is in
$U_p$.
\end{enumerate}
Let $\WW=(f,w,\UU')$ be the $\MM$-flow, subordinate to
$\VV$ and satisfying \ref{t:sq} \wrt~ the number
$\k/2$.
I claim that it satisfies our conclusions.
Indeed, $\Vert w\Vert\leq C$ by definition.
To show that the flow generated by $u$ is
$(t,\e)$-rapid, it suffices to check
(\ref{qpu1}, \ref{qpu2})
with $\k$ instead of $\e$.
The set
$M\sm B_\k(\indl {n-s-1}; -v)$
is
$(-v)$-invariant and $(-w)$-invariant,
so we must only check that for every
$M\sm B_\k(\indl {n-s-1}; -v)$
there is $p\in S(f)$ and $t_0\leq t$
\sut~
$\g(x,t_0;-w)\in B_\k(p)$.
This is possible, since
$U_p(\k/2)\sbs B_\k(p)$
and $\WW$ is
$(C,t,\k/2)$-quick.
The proof of \ref{qpu2} is similar.
$\qs$

Before proceeding to the proof of
\ref{t:q}$\RA$\ref{t:sq}
we explain informally the main idea
of the proof.
The principal difference between the statements
of these two theorems
is that
\ref{t:q}
asserts the existence of {\it one} 
 flow with prescribed quickness properties, and
\ref{t:sq}
asserts the existence of a whole family
of such flows, subordinate to a given one.
All these flows must have the prescribed quickness properties,
and in this family
there are  flows of arbitrary small radius.
Assume for example that we have a flow $\flow$
of radius $r$, satisfying the conclusions of
\ref{t:q}.
The problem is of course that
any $v$-trajectory which passes very close to
a critical point $p$ spends a lot of time
nearby $p$, in particular in the spherical slice
$U_p(r)\sm U_p(\mu)$, if $\mu$ is small.
But if we seek the flow of radius $\mu$
with good quickness properties, we can multiply
$v$ in
$U_p(r)\sm U_p(\mu)$
by a large positive real function, so that the (Euclidean) norm
of the resulting field becomes constant in this domain, which
will diminish dramatically the time spent in
$U_p(r)\sm U_p(\mu)$.

These considerations are formalized in the following proposition.
In this proposition we present a construction, producing
from a given
$\MM$-flow $\VV$ another one $\VV'$, for which $D(\VV')$
is
as small as we like, but $T(\VV')$ is still under control.

\begin{prop}\label{p:accel}
Let $\VV=(f,v,\UU)$
be an $\MM$-flow on $W$ with
 $\GG(\VV)\leq 2, \nrv\leq B$.
Write its atlas as $\UU=\atlas$.
Let $C\geq B$. Assume that to each $p\in S(f)$
               a number $\mu_p\in ]0,r_p[$ is assigned.
Then there is an $\MM$-flow $\VV'=(f,w,\UU')$, subordinate
to $\UU$, with
$\UU'=\{\Phi'_p:U_p'\to B^n(0,\mu_p)\}_{p\in S(f)}$
and $\Phi_p'=\Phi_p\mid U'_p$ and satisfying:
$\nrw\leq C$ and
\begin{equation}
C\cdot T(\VV')\leq B\cdot T(\VV)+5N(\VV)\cdot D(\VV)\lbl{f:accel}
\end{equation}  \end{prop}
{\it Proof.\quad}
Outside $\cup_pU_p$ we set $w=\frac CB v$.
We shall construct $w$ in each $U_p$ separately
with the help of the standard chart.
Namely, $(\Phi_p)_*(v)$ is the standard vector
field $\gotn_k$ on $\RRR^k\times \RRR^{n-k}$,
where $k=\ind p$. Let $\l:\RRR_+\to\RRR$ be a $\smo$
function \sut                          $$
\left\{\begin{gathered}
\l(t)=1\ffor 0\leq t\leq \mu_p\\
\min (1, \frac C{2t})\leq\l(t) \leq \max (1, \frac C{2t})\mbox{ for }
  t\in [\mu_p,\mu_p+\d]\\
\l(t)=\frac C{2t}\ffor \mu_p+\d\leq t\leq r_p-\d\\
\min (\frac CB, \frac C{2t})\leq\l(t) \leq \max (\frac CB, \frac C{2t})
\mbox{ for }
  t\in [r_p-\d,r_p]\\
\l(t)=\frac{C}{B}\ffor r_p\leq t
\end{gathered}
\right.                                  $$
(see the figure on the page  
\pageref{fig:lambda}). The number $\d\in]0,\frac{r_p-\mu_p}{2}[$
will be chosen later.

Define a vector field $w_0$ in $\RRR^n$ by
$w_0(x)=\l(\vert x\vert)\gotn_k(x)$. Note that
$(\Phi_p^{-1})_*(w_0)$ can be glued to $w$
previously defined on $M\sm\cup_pU_p$
(with the help of a standard extension);
performing this construction for each $p$
we obtain an $\MM$-flow $\VV'=(f,w,\UU')$,
subordinate to $\VV$.

The Euclidean norm of $w_0(x)$
 equals $\l(\vert x\vert)\vert x\vert$
and using $\GG(\VV)\leq 2$ it is easy to check that
$\nrw\leq C$. To obtain the estimate for
$T(\VV)$ note that the time which a
$w$-trajectory spends outside $\cup_pU_p$
is not more than $T(\VV)\frac BC$.
Therefore to obtain (\ref{f:accel}) it suffices
to prove that for each $p\in S(f)$ the
time which a $w$-trajectory
 can spend in $U_p\sm U_p(\mu_p)$
is not more than ${\frac 5C}r_p$.
This is a consequence of two following observations
(where $\g$ stands for a $w_0$-trajectory).
\begin{enumerate}
\item The time spent by $\g$ in each of the annuli
$\AA_1=\{\mu_p\leq\vert x\vert \leq\mu_p+\d\},\quad
\AA_2=\{r_p-\d\leq\vert x\vert \leq r_p\}$
goes to zero when $\d\to 0$ (it follows from
\ref{l:standtime}      )
\item The Euclidean length of the part of $\g$
in the annulus
 $\AA_3=\{\mu_p+\d\leq\vert x\vert\leq r_p-\d\}$
is not more than $2r_p$. The Euclidean
norm of $\g'$ in $\AA_3$ equals $\frac C2$.
$\qs$ \end{enumerate}

{\it Proof of Theorem  \ref{t:sq} from  Theorem \ref{t:q}.\quad}
We are given $C,\beta,A$. Choose $B>0,\mu_0>0$
as to satisfy
$$
B\leq C,\quad B\leq C\beta\frac{1}{3S(n)},\quad
\mu_0\leq C\beta\frac{1}{15n2^n}
$$
and let $\VV_0=(f,v,\UU_0)$ be an $\MM$-flow
satisfying the conclusions of \ref{t:q}
\wrt~ $B,\mu_0$.
Write $\UU_0=\{\Psi_p:U_p\to B^n(0,r_p)\}_{p\in S(f)}$.

Choose $r<d(\VV_0)$ so small that
$\GG(U_p(r),\Psi_p)<A$ for every $p\in S(f)$.
(This is possible since the coordinate frame
of $\Psi_p$ at $p$ is orthonormal \wrt~ the
 riemannian metric of $W$).
For each chart
$\{\Psi_p:U_p\to B^n(0,r_p)\}_{p\in S(f)}$
of $\UU_0$ consider its restriction
$\Phi_p=\Psi_p\vert U_p(r)$, form the corresponding
$qM$-chart system $\UU=\{\Phi_p\}$ and denote
by $\VV$ the $\MM$-flow $\VV=(f,v,\UU)$.

I claim that $\VV$ satisfy the conclusions of
\ref{t:q}. Indeed, let $0<\mu<r$. Apply the
proposition \ref{f:accel} to the flow $\VV_0$
and to the set of numbers
 $\mu_p=\mu$ for each $p\in S(f)$. We obtain
an $\MM$-flow $\WW=(f,w,\UU')$, subordinate to
$\VV_0$, where $\UU=\{\Phi_p'\}_{p\in S(f)}$
and $\Phi'_p=\Phi_p\vert U_p(\mu)$. This $\MM$-flow
is $(C,\beta)$-quick, since $\nrw\leq C$ and
(by (9))
$T(\VV')\leq \left(\frac {\beta T(\VV_0)}{3S(n)}+
{\frac 5C }n2^n\mu_0\right)\leq {\frac 23}\beta$.
The properties (2) and (3) go by the construction.
$\qs$

\section{ S-Construction
 }\label{s:sc}
The aim of this section is to present the S-construction.
Given a Morse function without critical points on a
cobordism $W$, the S-construction produces another
Morse function on $W$. This new function behaves like
 "fold" (see the figure on the page 
 \pageref{fig:skl} and the informal introduction
in the beginning of the present chapter).

The corresponding theorem
is stated below (theorem \ref{t:sc}), and the
 proof occupy the rest of the section.
In Subsection \ref{su:aux} we make some auxiliary constructions and
give some definitions. The output of the construction, namely,
Morse function $F$ and its gradient $u$ are introduced in 
Subsection \ref{su:mfg}.
The proof of the properties of $F$ and $u$, listed in the theorem
\ref{t:sc},
occupies the end of Subsection \ref{su:mfg} and Subsection
\ref{su:pfg}.

Proceeding to the statement of the theorem \ref{t:sc},
let $W$ be a  riemannian cobordism
of dimension $n+1$, where $n\geq 0$.
In the present section we denote by $\vert\cdot\vert$
the norm induced
by the metric 
on the tangent spaces, 
and by $\Vert\cdot\Vert$
the corresponding $C^0$ norm
in the space of vector fields.

\begin{theo}\label{t:sc}
Let $g:W\to [a,b]$ be a Morse function
without critical points,
$g^{-1}(a) =\pr_0 W, g^{-1}(b) =\pr_1W$.
Let $C>0$ and let $w$ be a $g$-gradient,
such that $\nrw \leq C$.
Denote $g^{-1}(\frac {2a+b}3)$ by $V_{1/3}$,~
$g^{-1}(\frac {2b+a}3)$ by $V_{2/3}$,~
$g^{-1}(\frac {a+b}2)$ by $V_{1/2}$.
Denote
$g^{-1}([a,\frac {2a+b}3])$ by $W_0$,
$g^{-1}([\frac {2b+a}3,b])$ by $W_1$,
$g^{-1}([\frac {2a+b}3,\frac {2b+a}3])$
by $W_{\frac 12}$.
Denote by $\text{\rm grad}(g)$ the riemannian
gradient of $g$ and by $\text{\rm grd}(g)$ the
 vector field
$\text{\rm grad}(g)/
\vert \text{\rm grad}(g) \vert$.
Denote by $T$ the maximal length
of the interval of definition
of a $\text{\rm   grd}(g)$-trajectory.
Then there is $\nu_0>0$ such that:

For every $s$-submanifold $\xx$ of $\pr_0 W$,
$s$-submanifold $\yy$ of $\pr_1W$,
every             $\AA \MM$-flow
$\VV_1=
(F_1, u_1, \UU_1)$
on $\VODIN$ and
$\AA \MM$-flow
 $\VV_2=(F_2, u_2, \UU_2)$
on $\VDVA$, and every \break
$\mu\leq \min (\nu_0,d(\VV_1),d(\VV_2))$
there is an
$\AA\MM$-flow $\VV=(F,u,\UU)$ on $W$,
having the following properties:
\begin{enumerate}
\item
$F:W\to[a,b], S(F)=S(F_1)\cup S(F_2)$;
in a neighborhood of $\partial W$
we have: $F=g,u=w$.

\item    $\VV $ is of radius $\mu$;\quad
      $\GG (\VV)\leq {\frac {3}{2}}
\max (\GG (\VV _1), \GG (\VV _2))$.

\item $V_\tret, V_\dvet, W_{\polo}$ are
$(\pm u)$-invariant,
$W_0$
 is $u$-invariant and weakly $(-u)$-invariant,
 $W_1$ is $(-u)$-invariant and weakly
$u$-invariant.
\item
Assume that $\VV _1$ is $(C_1,\beta_1,\mu/2)$-quick
and  $\VV _2$ is $(C_2,\beta_2,\mu/2)$-quick.
Then $\VV$ is
\break
$(3/2(C+C_1 +C_2),~\beta_1 +
\beta_2+5+4\frac {T}{C}~, \mu)$-quick.
\item
Let $\gamma$ be an $u$-trajectory.
If $     \Im \gamma\subset W_\polo$, then
$N(\gamma,\mu)\leq N(\VV _1,\mu)+N(\VV _2,\mu)$.
If $     \Im \gamma\subset W_0$, then
 $N(\gamma,\mu)\leq N(\VV _1,\mu)$.
If $     \Im \gamma\subset W_1$, then
$N(\gamma,\mu)\leq N(\VV _2,\mu)$.
\item
$\dd_a(u)\nmid \xx;\quad \dd_b(-u)\nmid\yy.$
\end{enumerate}
\end{theo}

{\it \quad Comments.}
To visualize the constructions we suggest three figures
(see the page \pageref{fig:skl})
which we have already 
mentioned.
What is drawn exactly is the case 
$\pr_0 W=\pr_1 W=S^1$.
On the figure A
we see the cobordism 
$W\approx S^1\times [0,1]$, embedded in $\RRR^3$,
the function $g$ being the vertical projection onto $[0,1]$.
The figure B
presents the manifold $W'$, diffeomorphic to $W$ but embedded 
differently into $\RRR^3$.
The projection $\pi$ onto the vertical, restricted to $W'$,
has the same behavior as the functions 
$\chi\circ g, \chi\circ f$ (see Subsection \ref{su:aux}).
The function $\pi\mid W'$ has two critical
levels: $\frac 23$ and $\frac 13$, and
the level surfaces correspond to $V_{\frac 13}$
and to
$V_{\frac 23}$.
Finally, on the figure C
we find a further  perturbed embedding 
$W''$ of $S^1\times [0,1]$ to $\RRR^3$, and our resulting function
$F$ is the vertical projection, restricted to $W''$.

The submanifolds $V_{\frac 13}, V_{\frac 23}$ correspond to
the critical level surfaces of $\widetilde f_i$ on the
figure B. The flows $\VV_1$ and $\VV_2$ correspond 
to the perturbations
of $\widetilde f_i$ which are necessary to turn $\widetilde f_i$
to a Morse function. \quad$\triangle$

                            \vskip0.1in

In the proof we  shall assume that $a=0, b=1$, since the
general case is easily reduced to this one
by an affine transformation of $\RRR$.

\subsection{ Auxiliary constructions and the choice of $\nu_0$}
\label{su:aux}
${}$
\vskip0.1in
{\it\quad   1. Function $f$ and its gradient $v$}
           \vskip0.1in
\begin{lemm}
There is a Morse function $f:W\to [0,1]$
 without critical points and an $f$-gradient
$v$, such that:
\begin{enumerate}
\item  $\nrv \leq C$;
\item in a neighborhood of
$\pr W$ we have: $ f=g$
and $ v=w$;
\item $f^{-1}(1/3)=V_{1/3},~ f^{-1}(2/3)=V_{2/3},~
f^{-1}(1/2)=V_{1/2}$;
\item for $x$ in a neighborhood of
 $V_{1/3}\cup  V_{1/2}\cup V_{2/3}$
we have:
\break
$ \vert v(x)\vert =C\aand df(v)(x)=C$;
\item for
$\lambda=1/3,~1/2,~2/3\aand x\in g^{-1}(\lambda)$
we have: $v(x)\bot g^{-1} (\lambda)$;
\item The maximal length of the domain of a
 $v$-trajectory is not more than $2T/C$.
\end{enumerate}
\end{lemm}
{\it\quad  Proof.}
Let $U$ be an open neighborhood of
 $\pr W$, such that
$U\cap(\VODIN\cup\VDVA\cup\VM)=\emptyset$
and let $h:W\to [0,1]$ be a $\smo$ function such
that $\supp ~h\subset U$ and for $x$ in a
neighborhood of $\partial W$ we have $h(x)=1$. Set
$v(x)=h(x)w(x) + (1-h(x))C\text{grd}(g)$. It is
 obvious that $v$ is a $g$-gradient, satisfying
 (1) and (5).
We have also $v(x)=w(x)$ nearby $\partial W$,
as well as $\vert v(x)   \vert=C      $ nearby
$\VODIN\cup\VDVA\cup\VM$.
(6) holds also, if only $U$ was chosen sufficiently small.
To construct $f$ define it first of all in a \nei~ of
$\VODIN\cup\VDVA\cup\VM$ directly by (4), and in 
a \nei~ of $\partial W$
as equal to $g$. Now extend $f$ to each of the cobordisms
$W_0, W_1, g^{-1}([\frac 13, \polo]), g^{-1}([\polo, \frac 23])$
using the usual procedure of gluing two functions
having the same gradient
(this procedure is described in details
in [21, Corollary 8.14]).$\qs$

For $\lambda\in[0,1]$ we denote $f^{-1}(\lambda)$ by
 $V_\lambda$. Fix some
$\epsilon \in ]0,\frac {1}{12} [$.
For $\nu>0$ sufficiently small
the map
$(x,\tau)\mapsto\gamma(x,\tau;v/C)$ is defined on
$V_i\times[-\nu,\nu]$, where $i=1/3, 1/2, 2/3$, on
$V_0\times[0, \nu]$
and on $V_1\times[-\nu, 0]$. The corresponding
embeddings will be denoted by
\begin{gather}
\Psi_0(\nu) : V_0 \times [0,\nu ] \to W,\quad\\
\Psi_1(\nu) : V_1 \times [-\nu,0 ] \to W,\quad\\
\Psi_i(\nu) : V_i \times [-\nu,\nu ] \to W,\quad\ffor i=\frac 13,
\polo, \frac 23
\end{gather}
 The image of $\Psi_i(\nu)$ will be denoted by
$\Tb_i(\nu)$ (for $i=0,\frac 13, \polo, \frac 23, 1$).

\pa
{\it 2. Choice of $\nu_0$.}
\pa

Let $\nu_0$ satisfy the following condition ($\frak R$):
$$
(\gR)
\left\{
\begin{aligned}
 (1)\quad  &2\nu_0<C \\
(2)\quad &\text{\rm ~For~} i=0, 1/3, 1/2, 2/3, 1  \text{\rm~we have:~}\\
\qquad  &(2A)\quad
 \Tb_i(2\nu_0)\subset ]i-\e, i+\e[\\
\qquad   &(2B)\quad\text{\rm
   ~The riemannian metric induced by~} 
\Psi_i(2\nu_0)\text{\rm ~from~}    W\\
   &\text{\rm~
and the product metric on the domain of~}
 \Psi_i(2\nu_0)
\text{\rm~
  are~} \\
&
 \polt-
\text{\rm~
equivalent.}\\
(3)\quad
&\text{\rm~
For~}
i=\frac 13, \polo, \frac 23
\text{\rm~ we have:~}
\vert v(x)\vert=df(v(x))=C
\text{\rm~
for~}
x
\text{\rm~
 in a}\\
&\text{\rm~    neighborhood of~}
\Tb_i(\nu_0)
\end{aligned}
\right.
$$
                   \vskip0.1in
We shall prove that $\nu_0$ satisfy the conclusions of our
theorem. So let $\xx,\yy, \VV_1=(F_1,u_1,\UU_1),
\VV_2=(F_2, u_2, \UU_2)$
be as in the statement of the theorem, where
$\UU_1=\{\Phi_{1p}:U_{1p}\to B^n(0,r_{1p})\}_{p\in S(F_1)}     $ 
and
$\UU_2=\{\Phi_{2q}:U_{2q}\to B^n(0,r_{2q})\}_{q\in S(F_2)}      $.
Let $\mu\leq\min(\nu_0, d(\VV_1), d(\VV_2))$.
We shall denote $\Psi_{i}(\mu)$ by
$\Psi_{i}$.
       Using $(\frak R)$ it is not difficult to check the following
property:
\begin{multline}
\ffor \vert\l\vert\leq 2\mu \aand i=\frac 13, \polo, \frac 23
\text{\rm ~we have~} \\
f^{-1}(i+\l)=\Psi_i(V_i\times\{\l\})\lbl{f:psi}
\end{multline}

\pa
{\it 3. Auxiliary functions, vector fields and $s$-submanifolds}
\pa

 For $\lambda\in [0,1/3]$ we denote by
 $\ll_\lambda (u_1)$ the
$s$-submanifold
$\stind v{1/3}\lambda (\dd (u_1))$ of $V_\lambda$.
For $\lambda\in [1/3,1]$ we denote
 by $\ll_\lambda (u_1)$
 the $s$-submanifold
 $\stind {(-v)}{1/3}\lambda \dd (u_1)$ of
$V_\lambda$.
For $\lambda\in [2/3,1]$ we denote
 by $\ll_\lambda (-u_2)$
 the $s$-submanifold $\stind {(-v)}{2/3}\lambda (\dd (-u_2))$
and
for $\lambda\in [0,2/3]$ we denote by
$\ll_\lambda (-u_2)$ the
 $s$-submanifold $\stind v{2/3}\lambda (\dd (-u_2))$.
Note that since $v$ have no zeros,
$\stind v\alpha\beta$ is a diffeomorphism
of $V_\alpha$ onto $V_\beta$
for any $\beta<\alpha$.
Note also
that $\ll_{2/3} (-u_2)$ equals to
$\dd (-u_2)$ and
$\ll_{1/3} (u_1)$ equals to
$\dd (u_1)$.

Let $\chi:[0,1]\to [0,1]$ be a $\smo$
 function, with the following properties:
\begin{gather*}
\chi (x)=x\ffor x\in [0,\epsilon]\cup [1-\epsilon , 1];\\
\chi (x) =1-x\ffor x\in [1/2-\epsilon, 1/2+\epsilon];\\
\chi (x) = 2/3-(x-1/3)^2
\for
 x\in [1/3-\epsilon,1/3+\epsilon];\\
\chi'(x)>0
\for
 x\in[0,1/3[~\cup ~]2/3,1];
\\
\chi (x) = 1/3+(x-2/3)^2\for
 x\in [2/3-\epsilon,2/3+\epsilon];\\
\chi'(x)<0 \for x\in ]1/3 , 2/3[.
\end{gather*}
(see the graph of $\chi$ on the  page \pageref{fig:chi}).

Denote $\chi\circ f$ by $\psi$. The following property
follows immediately from $(\frak R)$:
     \begin{multline}
  (\psi\circ\Psi_\tret)(x,\tau)=2/3-\tau^2,   \lbl{f:psi1}    \\
(\psi\circ\Psi_\dvet)(x,\tau)=1/3+\tau^2,\\
       (\psi\circ\Psi_\polo)(x,\tau)=1-\tau  
                                  \end{multline}

Let $B:\RRR\to [0,1]$ be a $\smo$ function
such that
$B(t)=0\ffor \vert t\vert\geq 5\mu/3\aand
B(t)=1\ffor  \vert t\vert\leq 4\mu/3.$
Let $B_1:\RRR\to\RRR^+$ be a $\smo$ function
such that $\supp B_1\subset ]\mu, 2\mu[$
and that $\int_0^{\infty} B_1(t)dt=C$.

Let $z_0$ be a $\smo$ vector field on $V_0$ such
that $\Phi(z_0,1)(\ll _0 (u_1))\nmid\xx$.
Let $z_1$ be a $\smo$ vector field on $V_1$ such
that $\Phi(z_1,1)(\ll _1 (-u_2))\nmid\yy$.
Let $z_{1/2}$ be a $\smo$ vector field on $V_{1/2}$
such that $\Phi(z_{1/2},1)(\ll _{1/2} (u_1))
\nmid\ll _{1/2} (-u_2)$.
We shall assume that $z_i$ (where $i=0,1/2,1$)
are chosen so small, that
$\sup_\tau\vert B_1(\tau)\vert
\cdot\Vert z_i\Vert<C/9$.

\subsection{ Morse function $F$ and its gradient $u$}
\label{su:mfg}
${}$
\pa
{\it\quad 1.   Morse function $F$.}
\begin{gather*}
\text{\rm Let } a_1,a_2>0. \text{\rm ~~
 Set:}\quad
 F(y)
 =\psi (y) \ffor y\in W
\setminus (\Tb _{\frac 13}(2\mu)\cup\Tb _{\frac 23}   (2\mu));\\
(F\circ\Psi _\tret)(x,\tau)
  =
a_1 B(\tau)F_1(x) +2/3-\tau^2 \ffor (x,\tau)\in
 \VODIN\times [-2\mu,2\mu]; \\
(F\circ\Psi _\dvet)(x,\tau)
  = a_2 B(\tau)F_2(x) +1/3+
\tau^2\ffor (x,\tau)\in \VDVA\times [-2\mu,2\mu]
\end{gather*}
It follows from (3.2) that these formulas define correctly
                          a smooth function
$F:W\to\RRR$, which equals to
$f$ nearby $\partial W$.
To find critical points of $F$ note that $S(F)$ is
contained obviously in
$\Tb _\tret(2\mu )\cup\Tb _\dvet(2\mu)$. For
$(x,\tau)\in\VODIN\times[-2\mu,2\mu]$ we have
$$
d(F\circ\Psi_1)(x,\tau)=
(a_1 B'(\tau) F_1(x) -2\tau)d\tau +
a_1 B(\tau) dF_1(x).
$$
For $a_1$ small enough this can
vanish
only for $\tau =0$. We conclude therefore
(applying the same reasoning to $\Tb _\dvet(2\mu)$),
 that $S(F)=S(F_1)\cup S(F_2)$ if only $a_i$ are
small enough, and we make this assumption from now on.
To prove that $F$ is a Morse function we shall
explicit the standard charts for $F$.
 Let $p\in S(F_1)$ and write
$(F_1\circ\Phi_{1p}^{-1})(x)=
F_1(p)+
\sum_i\alpha_ix_i^2$.
 Consider the chart
$\Phi_{1p}\times\id :U_{1p}\times ]-\mu,\mu[\to B^n(0,r_{1p})\times
 ]-\mu,\mu[$ of the manifold $\VODIN\times
]-\mu,\mu[$ around the point $(p,0)$.
We have
$$
F\circ\Psi_\tret\circ (\Phi_{1p}\times\id)^{-1}(x,\tau) =
a_1 F_1(p)+2/3+a_1 \sum_i \alpha_ix_i^2 -\tau^2
$$
 therefore the chart
$$
\big( (\Phi_{1p}\times\id)
\mid (\Phi_{1p}\times\id)^{-1} (B^{n+1}(0,\mu)\big)\circ
(\Psi_\tret\mid \Tb_\tret(\mu))^{-1}
$$
is a standard chart of radius $\mu$
for $F$ at $p\in S(F_1)$.
 These charts
together with the similar ones for $q\in S(F_2)$
give an $F$-chart-system of radius $\mu$.
We shall denote this system
by $\UU$.
Note      that $\ind _F p=\ind _{F_1}p+1$ for
$p\in S(F_1)$ and that $\ind _F q=\ind _{F_2}q$
 for $q\in S(F_2)$.
 Note also     that if $a_1$ and $a_2$ are chosen
sufficiently small, then
$F(\Tb _\tret(2\mu))\subset
[2/3 -\epsilon, 2/3 +\epsilon]$, and
$F(\Tb _\dvet(2\mu))\subset [1/3 -\epsilon, 1/3 +\epsilon]$.
In this case also $F^{-1}(1/2) =V_{1/2} \cup
 V_{\theta} \cup V_{\theta'}$, where
$\epsilon <\theta <1/3 -\epsilon,\quad
2/3 +\epsilon <\theta' < 1-\epsilon$.
\vskip0.1in
{\it\quad  2.  $F$-gradient $u$   }
\vskip0.1in

Set $\text{Tub}=
\Tb _0 (2\mu) \cup
\Tb _{\frac 13} (2\mu) \cup
\Tb _{\frac 12} (2\mu) \cup
\Tb _{\frac 23} (2\mu) \cup
\Tb _1 (2\mu)
$

Define a vector field $u$ on $W$ as follows:

\begin{equation*}
u(y)=v(y)\mbox{  for  } y\in (W_0\sm\Tub)\cup(W_2\sm\Tub)
\end{equation*}
\begin{equation*}
u(y)=-v(y)\mbox{  for  } y\in (W_{\frac 12}\sm\Tub)
\end{equation*}
\begin{equation*}
\big( (\Psi _0^{-1})_* (u) \big) (x,\tau)=
(-B_1(\tau)z_0(x),~ C);
\end{equation*}
\begin{equation*}
\big( (\Psi _1^{-1})_* (u) \big) (x,\tau)=
(B_1(-\tau)z_1(x),~C);
\end{equation*}
\begin{equation*}
\big( (\Psi _\tret^{-1})_* (u) \big) (x,\tau)=
\big(B(\tau)u_1(x), ~-B(\tau)\tau -
(1-B(\tau))\cdot C\cdot\text{sgn} \tau\big);
\end{equation*}
\begin{equation*}
\big( (\Psi _\polo^{-1})_* (u) \big) (x,\tau)=
(-B_1(-\tau)z_{1/2}(x),~
 -C);
\end{equation*}
\begin{equation*}
\big( (\Psi _\dvet^{-1})_* (u) \big) (x,\tau)=
 \big(B(\tau)u_2(x), ~B(\tau)\tau +
(1-B(\tau))\cdot C\cdot\text{sgn} \tau\big).
\end{equation*}

An easy computation
 using the definition of $F$ and $u$
shows that  these formulas
define correctly a
$\smo$ vector field on $W$, and that
$(F,u,\UU)$ is an $\MM$-flow on $W$
of radius $\mu$, if only $z_0, z_1, z_{1/2}$
are small enough (which assumption we
make from now on). It is also easy to see
that
$u$ is an $f$-gradient in
 $W_0\setminus V_{1/3}$
and $W_1\setminus V_{2/3}$; ~
$(-u)$ is an $f$-gradient in
$W_\polo\setminus (V_{1/3}     \cup V_{2/3}  )   $.

We claim that $\VV=    (F, u, \UU)$
satisfies all the conclusions of 3.1.
(1) and (2) follow immediately
from the construction. To prove
(3) note that  $V_{1/3}$ is a closed submanifold
of $\kr {W}$
 and $u\mid \VODIN$ is tangent to $\VODIN$;
therefore $\VODIN$ is $(\pm u)$-invariant
 (same for $\VDVA$).
 To prove that $W_\polo$ is $u$-invariant,
let $x\in\kr {W}_\polo$.
The trajectory $\g(x,\cdot; u)$ cannot intersect
$V_\tret$ or $V_\dvet$, since they are
$(\pm u)$-invariant. Therefore $\g(x,\cdot;u)$
has to stay in $W_\polo^\circ$ forever.

\subsection{ Properties of $F$ and $u$.}\label{su:pfg}
${}$
\vskip0.1in
{\it\quad    1. Estimate of the quickness of $\VV$}
                                \vskip0.1in
To obtain the estimate of      $\Vert u\Vert$
note that the inequality
$\vert u(x)\vert\leq
3/2(C+C_1+C_2)$ is to be checked
only for $x\in\cup_i\Tb_i(2\mu)$,
where it is a matter of a simple computation;
we leave it to the reader.
 To estimate the
time, which an $u$-trajectory spends outside
 $\UU (\mu)$, note first that for a trajectory, starting at
 a point of $\VODIN$ (resp. $\VDVA$),
this time is not more than $\beta _1$ (resp.
$\beta_2$),
since it is actually a $u_1$- (resp. $u_2$)-trajectory.
Now let  $\gamma$
be an
$u$-trajectory, passing by
a point of $\kr {W}_\polo$. In the end of
the preceding subsection
  we proved, that $\g$
 stays in $\kr {W}_\polo$ forever.
 Since $u$ is a $(-f)$-gradient in $\kr {W}_\polo $,
the function $t\mapsto f(\gamma (t))$
is strictly decreasing and
$\lim_{t\to -\infty} f(\gamma (t)) = 2/3,\quad
\lim_{t\to \infty} f(\gamma (t)) =1/3$.
 This implies
 $\lim_{t\to -\infty} \gamma (t) =q\in S(F_2)$ and
$\lim_{t\to \infty} \gamma (t) =p\in S(F_1)$.
We shall estimate the time which $\gamma$
spends outside $\UU (\mu)$
between the various level surfaces of $f$.

1) $f(\gamma (t))\in[2/3 -\mu/2,2/3]$. \quad
By (\ref{f:psi}) and (\ref{f:psi1}) this condition is equivalent to:
$\gamma (t)\in\Psi_\dvet(V_{2/3}\times[-\mu/2,0[)$.
The curve $\Psi_\dvet^{-1}(\gamma(t))$
is a product of
an $u_2$-trajectory and the curve
$\tau\mapsto \alpha e^{\tau}$
with $\alpha<0$.
 Since
 $\Psi_\dvet^{-1}(\UU(\mu))$
contains
$U_p(\mu/2)\times ]-\mu/2,\mu/2[$
for every $p\in S(F_2)$
the time which $\gamma (t)$ spends in
$\Psi_\dvet(V_{2/3}\times [-\mu/2,0[)$
outside $\Psi_\dvet^{-1}(\UU(\mu))$
is not more than $\beta_2$.

2) $f(\gamma (t))\in[2/3 -2\mu,2/3-\mu/2]$, in
other words,
$\gamma(t)\in \Psi_\dvet(V_{2/3}\times [-2\mu,-\mu/2])$.
The vector field $(\Psi_\dvet^{-1})_*(u)$
equals to
$(B(\tau) u_2(x),\varkappa (\tau))$, where
$\varkappa (\tau)\leq\tau$. Therefore the total time
which $\gamma$ can spend in this domain
is not more than
$\ln \big(2\mu/(\mu/2)\big) < 2$.

3) $f(\gamma (t))\in[1/2, 2/3-2\mu]$.\quad
In the domain
 $f^{-1}([1/2,2/3-2\mu])$
we have $u=v$. Therefore the time is
$\leq 2T/C$.

4) $f(\gamma (t))\in[1/2-2\mu, 1/2]$.\quad
Here $\gamma (t)\in\Psi_\polo(V_{1/2}\times
[-2\mu, 0])$. The second coordinate
of $(\Psi_\polo^{-1})_*(u)$ is equal to $(-C)$,
therefore the total time
spent here is not more than $2\mu/C\leq 1$.

Similarly to the cases 1) - 3) above
one shows that the time which
$\gamma$ spends in
$f^{-1}(]1/3,1/2 -2\mu])\setminus \UU (\mu)$
is not more than
$2T/C+2+\beta_1$. Summing
 up, we obtain that the time
which $\gamma$
can spend in $\kr {W}_\polo\setminus\UU (\mu)$
is not more than
$\beta_1+\beta_2+4T/C+5$.

Similar analysis of behavior of
 $u$-trajectories in $W_0$ and $W_1$
shows that the same estimate
holds in these cases also.
                                           \vskip0.1in
{\it\quad    2.    Estimate of $N(\gamma,\mu)$}
                                                      \vskip0.1in

Let $\gamma(\cdot)$ be an $u$-trajectory in $W_\polo$.
If $\gamma$ is in $V_{1/3}$, (resp. in $V_{2/3}$),
then obviously
$N(\gamma,\mu)\leq N(\VV_1,\mu)$
(resp.
$N(\gamma,\mu)\leq N(\VV_2,\mu)$).
Assume that $\Im\gamma\subset\kr {W}_\polo$.
 Let $\alpha\in\RRR$
(resp. $\beta\in\RRR$) be the
unique number, such that
$f(\gamma(\alpha)) =2/3-\mu$
(resp.
$f(\gamma(\beta)) =1/3+\mu$).
We have
$\cup_{p\in S(F_i)} U_p(\mu)
\subset
\Tb_i(\mu)$
 (for $i=1,2$). Therefore
$\gamma(t)\in\UU (\mu)$
can occur
only if
 $t\geq\beta$ or $t\leq \alpha$.
For $t\leq\alpha $ (resp. $t\geq\beta$)
the curve $\Psi_\dvet^{-1}(\gamma(t))$
(resp. $\Psi_\tret ^{-1}(\gamma(t))$)
 is an integral curve of the vector
field
$(u_2(x),\tau)$
(resp. $(u_1(x),-\tau)$).
Since an integral curve of $u_2$
(resp. $u_1$)
 can intersect no more than
$N(\VV _2,\mu)$
(resp. $N(\VV _1,\mu)$)
standard coordinate neighborhoods of
radius $\mu$, the first part
of (5) follows. The case of curves
in $W_0$ and $W_1     $
is considered similarly.
                                    \vskip0.1in
{\it\quad 3.   Transversality properties}
                                      \vskip0.1in
 The next lemma implies that
$\VV$ satisfies  Almost Transversality Condition.
It implies also the point (6)
 of our conclusions.
\begin{lemm}
\begin{enumerate} \item
 The family
$\{     \cup_{p\in S_i(F_1)} D_p(-u)\}_
{0\leq i\leq n}$
 equals to $\dd (-u_1)$.

              The family
$\{\cup_{q\in S_i(F_2)} D_q(u)\}_
{0\leq i\leq n}$
 equals to $\dd (u_2)$.
\item For each $\lambda\in [0,1]$ the family
$\{(\cup_{q\in S_i(F_2)} D_q(-u))
\cap V_\lambda\}_{0\leq i\leq n}$
is an $s$-submanifold of $V_\lambda$,
which is equal to:
\begin{enumerate}
\item
$\Phi (z_1 ,1)(\ll _0(-u_2))\iiif
\lambda =1$,~~
\item
       $\ll _\lambda (-u_2)\iiif
\lambda\in
[1/2, 1 - \epsilon ]$,~~
\item
$\emptyset\iiif
\lambda\leq 1/3$.
\end{enumerate}

\item For each $\lambda\in [0,1]$ the family
$\{(\cup_{p\in S_i(F_1)} D_p(u))
\cap V_\lambda\}_{0\leq i\leq n}$
is an $s$-submanifold of $V_\lambda$,
which is equal to:
\begin{enumerate}
\item $
\Phi (z_0,1)(\ll _0(u_1))\ffor \lambda =0~~;$
\item  $
\ll _\lambda (u_1)\ffor \lambda\in
[\epsilon,1/2-2\mu ]~~;           $
\item   $
\Phi (z_{1/2},1) (\ll _{1/2}(u_1))\ffor
 \lambda =1/2;                   $
\item    $
~\emptyset\ffor \lambda\geq 2/3;$
\end{enumerate}

 \item $u$ satisfies  Almost Transversality Condition,
and we have: $\dd(u)\cap V_0\nmid \xx, ~\dd(-u)\cap V_1\nmid \yy$
\end{enumerate}
\end{lemm}
{\it \quad Proof.}
 (1) For every $p\in S_i(F_1)$ the positive disc
of $F$ in $p$  belongs to $V_{1/3}$, which implies
immediately the first
assertion;
the second is proved similarly.

(2)
It is not difficult to see that
it suffices to prove the assertion
for $\lambda\in[-2\mu+2/3,2\mu+2/3]$.
For these values of $\lambda$ it follows from
the analysis of the behavior of $u$-trajectories
in $\Tb_\dvet (2\mu)$, carried out
in           Subsection \ref{su:pfg} 1) above.
 (3) is proved similarly.
(4): To prove that $u$ satisfies  Almost Transversality Condition
let
$p,q\in S(F),~ \ind p\leq \ind q$. Let
$\gamma$ be an $(-u)$-trajectory joining
$p$ with $q$. The case $p,q\in S(F_1)$
or $p,q\in S(F_2)$
follows easily from (1). Let
  $p\in S(F_1), q\in S(F_2)$.
Denote by $z$ the (unique) point
of intersection of
$\gamma$ with
$V_{1/2},~ \ind_Fp$ by $k$,
~ $\ind_Fq$ by $r$. Then
$z\in \big(\Phi(z_{1/2},1)(\ll _{1/2}(u_1))\big)_{k-1}$
and $z\in \ll _{\frac 12}(-u_2)_{n-r}$,
which is impossible by the choice of $z_{1/2}$.
 The last point is already proved.
$\qs$

\section{Proof of the theorem on the existence of quick flows}
\label{s:pteqf}

\subsection{Introduction}\label{su:pteqfint}

Let $W$ be a riemannian cobordism of dimension $n+1$.
We shall deduce from \ref{t:sq} the theorem
\ref{t:iq} below,
which is somehow stronger than \ref{t:q}.
Indeed, we have included in the conclusions of the
theorem \ref{t:iq}
not only the restrictions (1) -- (3) from \ref{t:q}
on the resulting $\MM$-flow $\VV$, but also the
 restrictions
$(\pr 1) - (\pr 3)$ which concern the behavior
of $\VV$ nearby $\pr W$
(the conditions $(\pr 1) - (\pr 3)$ are empty when
$\pr W=\emp$).

\begin{theo}\label{t:iq}
Let $g:W\to[a,b]$ be a Morse function on $W$.
Let $B>0,\mu_0>0$.
Let $u$ be a $g$-gradient with $\Vert u\Vert\leq B/2$.
Let $\xx$ be an $s$-submanifold of $\pr_0 W$, $\yy$ be
an $s$-submanifold
of $\pr_1 W$.
There is an $\AA\MM$-flow
$\VV=(f,v,\UU)$ on $W$, such that
\begin{enumerate}
\item $\GG(\VV)\leq 2,~D(\VV)\leq\mu_0$;
\item $N(\VV)\leq (n+1) 2^{n+1}$;
\item $\VV$ is $(B,S(n))$-quick, where
$S(n)=30+n2^{n+5}$
\end{enumerate}

Moreover:
\begin{enumerate}\item[($\pr 1$)]
$f(W)=[a,b]$ and in a \nei~ of $\pr W$ we have:
$f=g,u=v$.
\item[($\pr 2$)]
$\dd_a(v)\nmid\xx; \quad\dd_b(-v)\nmid\yy$
\item[($\pr 3$)]
If a $v$-trajectory $\g$ starts at a point of $\pr_0 W$
 then $\g(t)$ converges to a critical point of $f$
as $t\to\infty$.
If a $(-v)$-trajectory $\theta$ starts at a point of
$\pr_1 W$ then $\theta(t)$ converges to a critical point
of $f$ as $t\to\infty$.
Moreover for such trajectories $\g$ we have:
$N(\g)\leq n2^n, \quad T(\g)\leq 8+n2^{n+4}$.
\end{enumerate}
\end{theo}
We shall prove the theorem by induction on the number of critical
 points of $g$.
The lemma \ref{l:indst} below
shows how to do it: if \ref{t:iq}
is true for two "halves" of a cobordism $W$, then it is true also
for the whole $W$.
To make our notation shorter we shall
introduce some terminological conventions.
In the next definition and in the lemma \ref{l:indst}
$g:W\to[a,b]$ is a Morse function.

\begin{defi}
\begin{enumerate}\item
An {\it input data} is
a quintuple
$(B,\mu_0,u,\xx,\yy)$
where $B>0,\mu_0>0$ are real numbers, $u$ is
 a $g$-gradient with $\Vert u\Vert\leq B/2$,
$\xx,\yy$ are $s$-submanifolds of $\pr_0 W$, resp. of
$\pr_1 W$.
\item The pair $(g,W)$ will be called
{\it good}, if the theorem \ref{t:iq}
is true for $(g,W)$ and any input data
$(B,\mu_0,u, \xx,\yy)$.
\item A $v$-trajectory starting at a point
of $\pr_0 W$ or finishing at a point of $\pr_1 W$ will
be called
{\it boundary trajectory}.     \end{enumerate}           \end{defi}

Let $\l$ be a regular value of $g$. Set
$W_0=g^{-1}([a,\l]),\quad W_1=g^{-1}([\l,b])$.
Set $g_0=g\mid W_0, g_1=g\mid W_1$.

\begin{lemm}\label{l:indst}
If $(g_0,W_0)$ and $(g_1,W_1)$ are good, then
$(g,W)$ is good.
\end{lemm}
\Prf Let
$(B,\mu_0,u,\xx,\yy)$
be an input data.
Set $u_0=u\mid W_0, u_1 = u\mid W_1$. Apply \ref{t:iq}
to
$(g_0,W_0)$ and to the following input data:
$(B,\mu_0,u_0,\xx,\emp)$.
We obtain an $\AA\MM$-flow
$\VV_0=(f_0,v_0,\UU_0)$
on $W_0$.
Next apply Theorem \ref{t:iq}
to the function $g_1:W_1\to[\l,b]$
and to the following input data:
$(B,\mu_0,u_1,\dd_\l(-v_0),\yy)$.
We obtain an $\AA\MM$-flow
$\VV_1=(f_1,v_1,\UU_1)$ on $W_1$. Note that the condition
$(\pr 1)$ from \ref{t:iq}
allows to glue together the
$\MM$-flows $\VV_0$ on $W_0$ and $\VV_1$ on $W_1$
to an $\MM$-flow $\VV=(f,v,\UU)$ on $W$.
I claim that $\VV$ satisfies the conclusions of
\ref{t:iq}.
We must check first that $\VV$ is an $\MM$-flow. This
follows, since $\VV_0$ and $\VV_1$ are $\AA\MM$-flows and
$\dd_\l(v_1)\nmid\dd_\l(-v_0)$.
Next we check that $\VV$ verifies the conclusions
(1) -- (3)
of \ref{t:q}.
The (1) follows from the properties of $\VV_0,\VV_1$.
To prove (2) let $\g$ be a $v$-trajectory.
If it stays forever in one of $W_0,W_1$, then we
are over since the corresponding condition
is true for $\VV_0,\VV_1$.
If $\g$ intersects $V_\l=f^{-1}(\l)$, say,
 $\g(t_0)\in V_\l$, then
$\g(t)\in W_1$ for
$t\geq t_0$
and
$\g(t)\in W_0$
for $t\leq t_0$.
Now the inequality $N(\g)\leq (n+1)2^{n+1}$
follows from
the condition
$(\pr 3)$ which holds for $\VV_0$ and for $\VV_1$.
The proof of (3) is similar.

Now we proceed to check $(\pr 1) - (\pr 3)$ for
$\VV$.
The properties $(\pr 1)$ and $(\pr 3)$ follow
 immediately from
the corresponding properties of $\VV_0, \VV_1$.
To check $(\pr 3)$ note that
no $v$-trajectory starting at a point of $W_0$
 reaches $\pr_1 W$ (this follows since the property
$(\pr 3)$ holds for $\VV_1$).
Therefore the condition
$\dd(-v)\nmid\yy$
is equivalent to
$\dd(-v\mid W_1)\nmid\yy$, which holds by the
assumption. $\qs$

Therefore it suffices to deduce
\ref{t:iq}(n+1)
from \ref{t:sq}(n)
for
Morse functions $g:W\to[a,b]$ without critical points
and for those with one critical point.
This is the aim of the subsections \ref{su:fwcp}, \ref{su:fwocp}.
We shall not use the full force of \ref{t:sq}(n),
but only the following corollary of
\ref{t:sq}(n).

\begin{lemm}\label{l:wsq}
Let $N$ be a closed riemannian manifold of dimension $n$.
Let $D>0,\b>0,A>1$.
There is $\nu_0>0$ \sut~ for every $\mu<\nu_0$ there
 is an $\AA\MM$-flow $\WW$, \sut~
\begin{enumerate}
\item
$\WW$ is of radius $\mu$ and is
$(D,\beta)$-quick;
\item
 $N(\WW)\leq n2^n$;
\item
$\GG (\WW )\leq A$. $\qs$
\end{enumerate}
        \end{lemm}

\subsection{Functions without critical points}\label{su:fwcp}
Here we prove
\ref{t:sq}(n+1) for Morse functions without critical points.
The theorem \ref{t:sc} applies in this situation and we
shall first of all state a corollary
of (\ref{t:sc} + \ref{l:wsq}).

\begin{coro}\label{co:score}
Let $g:W\to[a,b]$ be a Morse function without critical points.
Let
$(B,\mu_0,w,\ll,\kk)$ be an
input data for $g$. Denote by $T$ the maximal length of the
interval of definition of
$\grd(g)$-trajectory.

Then there is an $\MM$-flow
$\VV=(F,u,\UU)$ of radius $\mu$ on $W$,
\sut~
\begin{enumerate}\item
In a \nei~ of $\pr W$ we have: $F=g,u=w$.
\item $\GG(\VV)\leq 2, \mu\leq \mu_0$
\item $N(\VV)\leq n2^{n+1}$,
\item $\VV$ is
$(B,7+n2^{n+4}+8\frac TB$)-quick
\item Every boundary trajectory
 $\g$ of
$v$ converges to a critical point of $F$.
 Moreover we have for such trajectories:
$N(\g)\leq n2^n$.
\item $\dd_a(u)\nmid\ll,\quad \dd_b(-u)\nmid\kk$.
\end{enumerate}     \end{coro}
\Prf
As in the proof of \ref{t:sc}
set
$g^{-1}(\frac {2a+b}3)=V_{\frac 13},
g^{-1}(\frac {2b+a}3)=V_{\frac 23}$.
Apply the lemma \ref{l:wsq} to
$\VODIN$ and to $\VDVA$
and obtain two $\AA\MM$-flows $\VV_1, \VV_2$ of some
 radius $\mu<\mu_0$
on $\VODIN, \VDVA$,
\sut~
$\GG(\VV_i)\leq 2, N(\VV_i)\leq n2^n$
and that $\VV_i$ are
$(\frac B{12}, 1)$-quick.
Note that $\VV_i$ are
$(\frac B{12}, 1+n2^{n+3}, \mu/2)$-quick
(by \ref{co:radrestr}).
Now just apply \ref{t:sc}
and produce an $\AA\MM$-flow
$\VV=(F,u,\UU)$ on $W$ of radius $\mu$, which is
$(B,7+n2^{n+4}+8\frac TB )$-quick.
We have also $N(\VV)\leq n2^{n+1}$.
The other points of our conclusion also follow from the
properties of
$\VV$ listed in the \ref{t:sc}. $\qs$

Of course this corollary does not yet give
an $\AA\MM$-flow satisfying the conclusions of \ref{t:iq}
since $T$ can be arbitrarily large.
So what we shall do is to cut
$W$ to thin slices so that for each slice $W_s$
the corresponding time $T_s$ is small,
then apply the \ref{co:score} to each of $W_s$ and then
glue the
resulting flows together.

Proceeding to the precise constructions,
recall that we are given the input data
$(B,\mu_0,u,\xx,\yy)$.
Choose a sequence $a=a_0<a_1<a_2<...<a_N=b$
of regular values of $g$, so close to each other, that
 for each $s$
the time which a $\grd(g)$-trajectory spends in
$W_s=g^{-1}([a_s,a_{s+1}])$
is not more than $B/8$.
We shall define by induction in $s$ a sequence
of $\AA\MM$-flows
$\VV_i=(F_i,u_i,\UU_i)$ on $W_s$.
Assume that $\VV_i$ is already constructed for $i\leq s-1$.
Apply \ref{co:score} to
$W_s$, choosing the following input data:
$(B,\mu_0,u,\ll_s,\kk_s)$ where
$\ll_s$ and $\kk_s$ are defined as follows:

\begin{align*}
s &=0: &\ll_s&=\xx,\quad\kk_s=\emp \\
0 &<s<N-1:  &\ll_s&=\dd_{a_s}(-u_{s-1}), \quad\kk_s=\emp \\
s &=N-1: &\ll_s&=\dd_{a_s}(-u_{s-1}),\quad \kk_s=\yy
\end{align*}

We obtain thus a sequence of $\AA\MM$-flows $\VV_i$ of
radius $\mu_i\leq\mu_0$, each of $\VV_i$
is
\break
$(B,8+n2^{n+4})$-quick,
and satisfies $N(\VV_i)\leq n2^{n+1}$.
Now glue all the $\VV_i$ together (this is possible,
since nearby every $g^{-1}(a_i)$
the functions $f_i,f_{i+1}$ coincide with $g$ and the
vector fields
$u_i, u_{i+1}$ coincide with $u$). We obtain an
$\MM$-flow $\VV=(F,u,\UU)$ on $W$,
and I claim that it satisfies
the conclusions of \ref{t:iq}.
Indeed, note first that $\VV$ is an
$\AA\MM$-flow (this follows from the restrictions
$\dd_{a_s}(u_s)\nmid\dd_{a_s}(-u_{s-1})$,
imposed on $\VV_s$ in the process of construction).
Further, (1) and (2) follow immediately
(actually, we obtain here a bit stronger
estimate $N(\VV)\leq n2^{n+1}$).
As to the estimate for the quickness, note that no
$v$-trajectory
can intersect more than two of $W_s$. Thus
$\VV$ is
$(B, 2(8+n2^{n+4}))$-quick. The other conditions follow
immediately. $\qs$

\subsection{Functions with one critical point}\label{su:fwocp}
Here we prove \ref{t:iq}(n+1)
for functions
$g:W\to[a,b]$ with one critical point. Let $p$ be this
critical point, let
$k$ be its index, and let
$(B,\mu_0,u,\xx,\yy)$ be our input data for construction
of the flow $\VV$.
Changing if necessary the $g$-gradient
$u$ in a \nei~ of $p$, we can assume that there is
a Morse chart $\Phi:U\to B^n(0,r)\sbs\RRR^{n+1}$ for $g$
\sut~ $\Phi(p)=0$ and
$\Phi_*(v)=\gotn_k$.

\begin{lemm}\label{l:six}
There is $\d_0,  0<\d_0<r$, \sut~ for every $\d<\d_0$ the time which
 a $v$-trajectory can spend in
$g^{-1}([g(p)-\d^2,g(p)+\d^2])\sm
U_p(\d/2)$
is not more than $6$.
\end{lemm}

\Prf
 Denote
$ g^{-1} ([g(p)-\delta^2,g(p)+\delta^2])$ by $S_\delta $.
Let $\gamma$ be a $w$-trajectory.
The function $dg(w)$ is bounded from below in
$S_\delta
\setminus U_p$, therefore for
$\delta >0$ sufficiently small the time which
$\gamma$ spends in $S_\delta\setminus U_p$
 is less
than 1. The time, which $\gamma$ spends in
$U_p\cap \big(S_\delta
\setminus U_p(\delta)\big)$
is not more than 2 (by 1.10), and
 the time which $\gamma$ spends
 in $U_p(\delta)\setminus U_p(\delta/2)$
is not more than
$\ln 8\leq 3$ (by 1.9).
$\qs$

Now let $\d>0$ be so small that
\begin{enumerate}
\item $\d<\mu_0,\d<\d_0$
\item $\GG(U(\d),\Phi\mid U(\d))\leq 2$
\end{enumerate}

Set
\begin{align*}
W_+ &=g^{-1}\([g(p)+\d^2, b]\),  &g_+ = g\mid W_+ \\
W_- &=g^{-1}\([a, g(p)-\d^2]\),  &g_- = g\mid W_- \\
W_- &=g^{-1}\([g(p)-\d^2, g(p)+\d^2]\),
 &g_0 = g\mid W_0 \\
V_- &= g^{-1}(g(p)-\d^2) & \\
V_+ &= g^{-1}(g(p)+\d^2) & \\
S_+ &= D(p,-v)\cap g^{-1}(g(p)+\d^2)    & .
\end{align*}

By abuse of notation we keep the symbol $u$ to
denote all the restrictions
of $u$ to: $W_+, W_-,...$.
Apply to $g_+$ Theorem \ref{t:iq}
(proved in \ref{su:fwcp}).
Choose as the input data the following string:
$(B,\mu_0, u, S_+,\yy)$.

Denote the resulting $\AA\MM$-flow by
$\VV_+=(f_+,v_+,\UU_+)$.
Glue it to the $\AA\MM$-flow $\VV_0=
(g\mid W_0, u, \{\Phi\mid U(\d/2)\})$.
We obtain an $\AA\MM$-flow
on $W_+\cup W_0$, which will be denoted by
$\VV_1=(F_1,v_1, \UU_1)$.
Apply now Theorem \ref{t:iq} to
$g_-$, choosing as the input data the following string:
$(B,\mu_0,u,\xx, \dd_{f(p)-\d^2}(v_1))$.
We obtain an $\AA\MM$-flow. Denote it by
$\VV_-=(f_-,v_-,\UU_-)$
and glue it to $\VV_1$.

Denote the resulting $\MM$-flow by
$\VV=(f,v,\UU)$.
I claim that it satisfies the conclusions of
\ref{t:iq}.
Indeed, the condition (1) goes by construction,
as well as the condition
$\Vert v\Vert\leq B$.
To check (2) note that it is already proved for
 trajectories of $v$, which stay forever in $W_-$ or
$W_+$. Assume that a $v$-trajectory $\g$ intersects $V_+$.
If $\lim_{t\to-\infty}\g(t)=p$, then
$N(\g)\leq 1+n2^n$.
If $\g$ intersects $V_-$ then $N(\g)\leq 1+n2^{n+1}$.
The time which a $v$-trajectory $\g$ can spend outside
$\cup_p U(p)$ is estimated from above using the lemma
\ref{l:six} and the estimates of Corollary \ref{co:score}.

The properties
$(\pr 1) - (\pr 3)$
of Theorem \ref{t:iq} are verified without difficulties.
$\qs$

\section{Final remark}
\lb{s:fr}

The results of this chapter
are essentially contained in \cite{pastpet}.
In particular, the theorem \ref{t:sq}
is a generalization of the Theorem 1.12 from \cite{pastpet}
to the case of
\cob s (in \cite{pastpet} it was stated only for
closed \ma s).
The theorem \ref{t:q}
is a generalization of the theorem 1.13 from \cite{pastpet}
to the case of \cob s
(also we diminished the corresponding constant
from 100+$n 2^{n+7}$
in \cite{pastpet} to
30+$n 2^{n+5}$;
we think that one can diminish it still more).
The theorems 1.12 and 1.13 of \cite{pastpet}
are in turn exactly the theorems 1.17, resp. 1.18
of the preprint \cite{paepr}.

\chapter{Regularizing the gradient descent map}
\label{ch:rgdm}

Let $W$ be a riemannian \cob~ of dimension $n$,
$\fcob$ be a \Mf~and $v$ be an $f$-\gr.
Consider the open set $U_1\sbs \pr_1 W$ of all the points $x$,
\sut~ the \tr~ 
$\g(x,\cdot; -v)$
reaches $\dow$
and intersects it it at some point $y$.
The correspondence $x\mapsto y$
is a \dfm~of $U_1$
onto an open subset
$U_0\sbs \dow$, we denote
it by
$\stexp {(-v)}$.
The aim of this chapter is to prove 
an analog of cellular approximation theorem for the map $\stexp {(-v)}$.
The problem is that 
$\stexp {(-v)}$
in general can not be extended to a continuous map on the whole of $\daw$.

Recall however from \ref{s:gdm}
that 
$\stexp {(-v)}$
is well defined on the category
of $s$-\sma s
(under some mild transversality condition, see 
\sub~ \ref{su:trssb}). The collection of descending discs
 of a \gr~of a \Mf~ form
an 
$s$-\sma~(if the \gr~ satisfies \ATA); we pick a \Mf~
$\phi_1$ on $\daw$
and a $\phi_1$-\gr~ $u_1$, satisfying \ATA.
Consider the corresponding $s$-\sma~ 
$\dd(u_1)$,
formed by the descending discs
of
$u_1$. Denote $\dd(u_1)$ by $\aa$ and 
$\cup_{i\leq j} D(\inde i ; u_1)$
by $A_{\leq j}$.
Proceed similarly with $\dow$ and obtain
the $s$-\sma~$\bb=\dd(u_0)$, where $u_0$ is a gradient of a \Mf on $\dow$,
and $u_0$ satisfies \ATA.
The decomposition of a \ma~ to the union of the 
descending discs of a gradient
resembles much of a cellular decomposition, so one can try 
to reproduce the cellular approximation theorem in this situation.
We can not, of course, expect, that, perturbing $v$, we obtain for the
new \gr~ $u$ the direct analog of the cellular approximation condition:
$\stexp {(-u)} (A_{\leq j}) \sbs B_{\leq j}$.
We can, though, deform $v$ so that
 for the  new \gr~$u$ the following will hold:

\begin{equation}
\stexp {(-u)}(A_{\leq j} )\sbs 
D_\d(\indl j; u_0)  \lb{f:opyat}
\end{equation}

Indeed, first perturb 
$v$ so as to make the \sma~ 
$\stexp {(-u)} (\aa)$
transversal to the $s$-\sma~ 
$\wh{\bb}=\dd(-u_0)$. Then
$\stexp {(-u)} (A_{\leq j})$
does not intersect 
$\dd_{\leq n-2-j}(-u_0)$
and for
$T$ sufficiently large
the \dfm~ 
$\Phi(T, -u_0)$ pushes
$\stexp {(-u)} (A_{\leq j})$
into a small \nei~ 
of
$\dd_{\leq j} (-u_0)$.
We can modify $v$
in a small \nei~of $\dow$, so that for the new vector field $\wi v$
the following holds:
$
\stexp{(-\wi v)}=
\Phi(T,-u_0)\circ \stexp{(-v)}$.
This is done by adding to
$v$ a horizontal component proportional to $u_0$
(see \sub~ \ref{su:ahc} for definition of the construction 
of adding a horizontal component).
Thus for $\wi v$
the 
(\ref{f:opyat})
holds.

It is essential for us that the property
(\ref{f:opyat})
(and actually a bit stronger property)
can be achieved by a $C^0$-small perturbation
$\wi v$ of $v$.
The proof of this result
(Theorem \ref{t:cc})
is based on the main theorem of the previous chapter.

\section{Condition $(\gC)$}
\label{s:cc}

\subsection{Condition $(\gC)$: the definition
and the statement of the density theorem
}
\label{su:cc}
${}$
\pa
{\it Terminology}
\lb{ssu:term}

Let $g:W\to[a,b]$ be a Morse function on a cobordism $W$
of dimension $n$, $v$ be a $g$-gradient.

Let $\l\in[a,b]$ be a regular value of $f$. Let $A\sbs\daw$.
We say, that {\it $v$ descends $A$ to $f^{-1}(\l)$},
if every
$(-v)$-trajectory starting in $A$ reaches $f^{-1}(\l)$.
Let $X\sbs W$. We say that
{\it $v$ descends $A$ to $f^{-1}(\l)$
without intersecting $X$}, if $v$ descends $A$ to
$f^{-1}(\l)$
and
$T(A, -v)\cap X=\emp$.

We leave to the reader to  define the meaning of
{\it "$v$ lifts $B$ to $f^{-1}(\l)$"}
and
{\it "$v$ lifts $B$ to $f^{-1}(\l)$
without intersecting $Y$"}, where
$B\sbs \dow, Y\sbs f^{-1}\([a,\l]\). \qt$

Recall the sets
$B_\d(\indl s;v), D_\d(\indl s;v)$
from Definition \ref{d:dthick}
and introduce one more set
with similar properties:
set
$C_\d(\indl s; v)=W\sm B_\d(\indg {n-s-1};-v)$.
The following lemma
 is proved by standard Morse-theoretical arguments.

\bele\lb{l:cdelta}
Assume that $v$ is $\d$-separated
\wrt~$\phi$ and the ordering sequence
$(a_0,...,a_{n+1})$.
Then for every $s$
we have:
$$
\dow\cup B_\d(\indl s;v)\sbs
\phi^{-1}\([a_0,a_{s+1}]\)
\sbs
C_\d(\indl s;v)\cup \dow
$$
and these inclusions are homotopy equivalences.
$\qs$\enle

Now we can formulate the condition $(\gC)$.
\bede\lb{d:cc}

We say, that
{\it $v$ satisfies condition $(\gC)$}
if there are objects 1) - 4), listed below,
with the properties (A), (B0), (B1) below.
\pa
{\it Objects:}
\pa
\begin{enumerate}
\item[1)] An ordered Morse function $\phi_1$
on $\pr_1 W$ with ordering sequence
$(\a_0,...,\a_n)$, and a $\phi_1$-gradient $u_1$.
\item[2)]  An ordered Morse function $\phi_0$
on $\pr_0 W$ with ordering sequence
$(\b_0,...,\b_n)$, and a $\phi_0$-gradient $u_0$.
\item[3)]  An ordered Morse function $\phi$
on $ W$ with ordering sequence
$(a_0,...,a_{n+1})$,
adjusted to $(f,v)$.
\item[4)]  A number $\d>0$.
\end{enumerate}
\pa
{\it Properties:}
\pa
(A)\quad  $u_0$ is $\d$-separated \wrt~ $\phi_0$,
$u_1$ is $\d$-separated \wrt~ $\phi_1$, $v$ is
$\d$-separated \wrt~ $\phi$.
\pa

 \begin{equation*}\stv \bigg(C_\d(\indl j;u_1)\bigg)
\cup
\bigg(D_\d\(\indl {j+1};v\)\cap \pr_0 W\bigg)
\sbs
B_\d(\indl j,u_0)
\mbox{ for every } j\tag{B1}
\end{equation*}

\begin{equation*}\st v\bigg(C_\d(\indl j;-u_0)     \bigg)
\cup
\bigg( D_\d(\indl {j+1}; -v)\cap \pr_1 W\bigg)
\sbs
B_\d(\indl j;-u_1)
\mbox{ for every } j\tag{B0}
\end{equation*}

\end{defi}

The set of all $f$ gradients satisfying $(\gC)$ will be denoted by
$\GgC(f)$.
\pa
{\it Comments:}
\begin{enumerate}\item
As we have already mentioned in the introduction,
the conditions (B0), (B1) together form
 an analog of cellular approximation condition.
Indeed, the sets $C_\d(\indl j;u_1)$ and $B_\d(\indl j;u_0)$
are "thickenings" of
$D(\indl j;u_1)$, resp. of $D(\indl j; u_0)$ (the first one
being "thicker" than the second). Thus the condition
(B1)
requires that $\stv$ sends a certain thickening of a
$j$-skeleton of $\daw$
to
a certain thickening of a $j$-skeleton of $\dow$.
Warning: (B1) actually requires more than that:
for every $j$
the soles of $\d$-handles of $W$ of indices $\leq j+1$
must belong to the
set
$B_\d(\indl j;u_0)$.

\item
The condition $(\gC)$ is stronger than the condition
$(\CC)$ from \cite{pator}, \S 1.2. This is obvious, since
$(\gC)$
is formulated similarly to $(\CC)$, except that in $(\CC)$ we
have demanded
$\stv \(\phi^{-1}\([a_0, a_{j+1}]\)\)
\sbs
B_\d(\indl j;u_0)$. The condition $(\CC)$ is in turn stronger, than
the condition (RP) from \cite{pastpet}, \S 4, this is a bit less obvious,
and we shall not prove it here.

\item
If $v$ satisfies condition $(\gC)$,
then $v$ obviously satisfies \ata. Note that if
 we are given a $\phi_0$-gradient $u_0$, a $\phi_1$-gradient $u_1$,
a $\phi$-gradient $v$, which satisfy Almost
Transversality Condition,
then the condition (A) is always true if
$\d$ is sufficiently small.

\end{enumerate}

\begin{theo}\lb{t:cc}
$\GgC(f)$ is open and dense in $\GG(f)$ \wrt~  $C^0$ topology.
Moreover, if $v_0$ is any $f$-gradient then one can
 choose a $C^0$-small perturbation $v$
of $v_0$ \sut~ $v\in\GgC(f)$, 
$v=v_0$ in a \nei~ of $\pr W$, and
$v=v_0$
in a \nei~ of every critical point of $f$.
\end{theo}

The proof  occupies the subsections
\ref{su:cogct} -- \ref{su:cdpgg}.

{\it Convention:}
We shall say that $(\gC)$ is {\it open}, resp. {\it dense }
if the set of all
$f$-gradients satisfying $(\gC)$ is open, resp. dense in $\GG(f)$.
$\qt$

\pa

\subsection{ $C^0$-openness of $\gC$}
\quad\lb{su:cogct}
\pa
Let $v$ be an $f$-gradient satisfying $(\gC)$.
Choose the corresponding functions
$\phi_0, \phi_1, \phi$, their gradients
$u_0, u_1$, and a number $\d>0$ satisfying (A) and (B0), (B1).
Fix now the string
$\SS=(\d,u_0, u_1,\phi_0, \phi_1,\phi)$.
For a $\phi$-\gr~ $w$ we shall denote the condition (A) with
respect to $w$ and $\SS$
by (A)($w$). Similarly, we shall denote
 the conditions  (B1), (B0) \wrt~ $w$ and $\SS$
by  (B1)($w$), resp. (B0)($w$). We know that
(A)($v$)$\&$ (B1)($v$) $\&$ (B0)($v$)
 holds; we shall prove that
(A)($w$)$\&$ (B1)($w$) $\&$ (B0)($w$)             holds
for every $w$ sufficiently $C^0$ close to $v$.

 Let $w$ be another $f$-gradient. Corollary \ref{c:seppert}
 implies that
if $\Vert w-v\Vert$ is sufficiently small,
then $w$ is
a $\phi$-gradient and is
 also $\d$-separated, so the condition A($w$) is
satisfied also for $w$.

It is convenient to reformulate the conditions (B1)($w$)$\&$(B0)($w$).
Introduce three new conditions:
\been
\item[$\b_1(w)$:]
For every $j$ we have:
$w$ descends $C_\d(\indl j;u_1)$ to $\phi^{-1}(a_{j+1})$
without intersecting $\cup_{\ind p\geq j+1} D_\d(p)$.

\item[$\b_2(w)$:]
For every $j$ we have:
$w$ lifts $C_\d(\indl {n-j-2};-u_0)$ to $\phi^{-1}(a_{j+2})$
without intersecting $\cup_{\ind p\leq j+1} D_\d(p)$.

\item[$\b_3(w)$:]
For every $j$ the sets:

$$R_j^+(w)=\stind {(-v)}b{a_{j+2}} (C_\d(\indl j;u_1)),$$

$$R_j^-(w)=\stind {v}a{a_{j+2}} (C_\d(\indl {n-j-2};-u_0))$$

are disjoint.
\enen

Note that $\b_1(w)$ and $\b_2(w)$ are in a sense dual to each other.

\bele
Assume that A($w$) holds. Then
(B1)($w$)$\&$(B0)($w$)$\Leftrightarrow
\b_1(w)\&\b_1(w)\&\b_3(w)$
\enle
\Prf $\Rightarrow$

The condition $\b_1$ follows from (B1), $\b_2$ follows from (B0).
To prove $\b_3$, note that if
$R_j^+(w)\cap R_j^-(w)\not=\emp$, then there is
a $(-v)$-trajectory joining
$x\in C_\d(\indl j;u_1)$ with
$y\in C_\d(\indl {n-j-2}; -u_0)=
\dow\sm B_\d(\indl j; u_0)$.
And this contradicts (B1).

$\Leftarrow$

 $\b_1$ says that for every $j$ we have:
$$
\bigg(\bigcup_{\ind p\geq j+1} D_\d(p,-v)\bigg)
\cap\daw\sbs
\daw\sm C_\d(\indl j;u_1)=B_\d(\indl {n-2-j}; -u_1).
$$
That is,
$D_\d(\indl {n-1-j}; -v)\cap \daw\sbs B_\d(\indl {n-2-j};-u_1)$.
Further, $\b_3$ implies that for every $j$
every $v$-trajectory starting in
$C_\d(\indl {n-j-2}; -u_0)$
and reaching $\daw$, cannot intersect $\daw$ at a point of
$C_\d(\indl j;u_1)$. Therefore it intersects
$\daw$ at a point of
$B_\d(\indl {n-2-j};-u_1)$.
Therefore (B0)($w$) holds. The proof of (B1)($w$) is similar. $\qs$

Returning to the proof of $C^0$-openness of $(\gC)$, note that the
conditions
$\b_1(w)$ and $\b_2(w)$ are obviously open.
We know that $\b_1(v)$ and $\b_2(v)$ hold, therefore
for every $j$ the sets
$R_j^+(v), R_j^-(v)$
are compact. Since they are disjoint (by $\b_3(v)$),
we can choose two disjoint open subsets
$U_j, V_j$ of
$\phi^{-1}(a_{j+2})$ \sut
~$R_j^+(v)\sbs U_j, R_j^-(v)\sbs V_j$. Then for every
$w$
sufficiently close to $v$ and for every
$j$
we have (by Corollary \ref{co:homot}):
$$
\stind {(-w)}b{a_{j+2}} \(C_\d(\indl j;u_1)\)\sbs U_j, \quad
\stind wa{a_{j+2}}\( C_\d(\indl {n-2-j}; -u_0)\)\sbs V_j
$$
and $\b_3(w)$ holds. Therefore
$\b_1(w)\&
\b_2(w)\&
\b_3(w)$
is $C^0$ open and we have proved $C^0$ openness of $(\gC)$.
\pa

\subsection{ $C^0$-density of $(\gC)$: The idea of the
 proof }
\lb{su:cdi}

In this subsection I present the idea of the
proof on the intuitive level.

Choose some pair $(\phi_0, u_0)$, where $\phi_0$ is
an ordered Morse function on $\dow$, and
$u_0$ is an $\phi_0$-gradient.
Let us try to perturb a given $f$-gradient
$v_0$ in order to obtain (B1).
The set
$$\gN_j(v,\d)=
\stv \bigg(C_\d(\indl j;u_1)\bigg)
\cup
\bigg(D_\d\(\indl {j+1};v\)\cap \pr_0 W\bigg)
$$
contains in particular the set
$$
\gN_j(v)=
\stv \(D(\indl j;u_1)\)
\cup
\(D\(\indl {j+1};v\)\cap \pr_0 W\)
$$

We know from Proposition
\ref{trper} and Lemma \ref{l:trcomp}
that perturbing $v$ if necessary we can assume
that
$\{\gN_j(v)\}$
form an $s$-\sma, transversal to to $\dd(-u_0)$.
In particular
$\gN_j(v)$ is disjoint from
$D(\indl i, -u_0)$ for $i+j<n-1$.
Using Proposition \ref{tracgfn} we conclude that for
$\e>0$ small enough the sets
$\gN_j(v,\e)$
and
$D_\e(\indl i; -u_0)$
are disjoint if $i+j<n-1$.

Now  we shall modify $v$ by adding a horizontal component nearby $\dow$.
Let $\xi$ be a vector field equivalent to $u_0$
(see Definition \ref{d:equiv}).
Then for some $T>0$ we have
$$
\Phi(T,-\xi)(\dow\sm B_\e(\indl i, -u_0))
\sbs
B_\e(\indl {n-2-i}; u_0
$$
Thus if we apply to the AHC-construction \wrt~
$\xi$
and to the function $h$ satisfying
$\int_0^l h(\tau) d\tau=T$,
we shall have for the resulting vector field $w$:
$$
\gN_j(w)\sbs B_\e(\indl {n-2-i}, u_0)
$$
(see AHC2, page \pageref{ahc2};
actually we should also take some care in order that
$w$ be an $f$-gradient).

Of course if the vector field $u_0$ and $\xi$
were chosen arbitrarily
then $T$ may be large and therefore $\Vert v-w\Vert$
may also be large.

Now I shall explain
why it is possible to choose $T$ and $\Vert v-w\Vert$
small. This step uses in essential way the Theorem \ref{t:quickpush}.
We shall chose the pair $(\phi_0, u_0)$ in a special way.
Let $C>0, t>0$ be real numbers and $(\phi_0, u8а)$ be the
pair, satisfying the conclusions of \ref{t:quickpush}
\wrt~ $C,t$.
Choose $\wi u_0$, equivalent $u_0$ with
$\Vert \wi u_0\Vert\leq C$, and having
the property that the flow generated
by $(-\wi u_0)$ pushes
$\dow\sm B_\e(\indl i; -u_0)$
inside
$B_\e(\indl {n-2-i}; u_0)$,
and the time necessary for this procedure is $\leq t$.
Now set $\xi=\wi u_0$ and to $v$ a horizontal component proportional
to $\xi$ nearby $\dow$.
The input data $(a,\xi, h)$ for AHC-construction
will be chosen so that
$\int_0^a h(\tau) d\tau=t$.
The estimates from
AHC1 (page \pageref{ahc1})
shows that $\Vert v-w\Vert$
will be small enough if only $C$ and $t$ are chosen
sufficiently small.

\subs{$C^0$-density of $(\gC)$: terminology }
\label{su:cdt}
\pa
Let $v_0$ be an $f$-gradient.
Choose regular values of $\a,\b\in]a,b[$
of $f$, \sut~ $a<\a<\b<b$ and that the intervals
$[a,\a], [\b,b]$ be regular.
Denote $f^{-1}(\l)$ by $V_\l$, so that $\pr_1 W=V_b, \pr_0 W=V_a$.
Denote $f^{-1}([\a,\b])$ by $W'$.
The map $\stind  {(-v_0)}b\b :V_b\to V_\b$
and the map
$\stind  {(v_0)}a\a :V_a\to V_\a$
are diffeomorphisms.
We shall denote these maps by $\Psi_1$, resp. $\Psi_0$.
Set
$$m=\inf_{x\in W\sm W'}\Vert df(v)(x)\Vert,\quad
M=\sup_{x\in W\sm W'}\Vert df(x)\Vert$$

Let $E_0, E_1$ be the expansion constants for
$\Psi_0$, resp. $\Psi_1$.
Set $E=\max (E_0, E_1)$.

Now  we shall define  the horizontal
components which will be added to $v$.
Let $\e>0$. Our aim is to  construct a vector field $v$,
\sut~(B1) and (B0) hold and
$\Vert v-v_0\Vert<\e$.
Choose some closed collars $C_0=\g(\pr_0 W,[0,l];v)$
and
$C_1=\g(\pr_1 W, [0,l]; -v)$
so that
$C_0\sbs \phi^{-1}([a,\a[), C_1\sbs \phi^{-1}(]\b,b])$.
Choose any positive $\smo$ function $h$ on $[0,l]$
\sut~$h(x)=0$ for $x$ in a \nei~of $0$ and in a \nei~of $l$.
Denote $\sup_{\tau\in[0,l]} h(\tau)$ by $H$.
Set $t_0=\int_0^lh(\tau)d\tau$.
Choose a real number $A>0$ \sut~
\begin{equation}
A<\frac{m}{EMH}, \quad
A<\frac\e{2EH} \lbl{smnr}
\end{equation}

\bere\label{r:choice}
It is easy to deduce from AHC1, AHC3,
that with this choice of $A, t_0, h, H$ we have:
Let $\xi$ be a vector field $\dow$ with $\Vert\xi\Vert\leq A$.
Let $v$ be an $f$-gradient.
Let $w$ be the result of AHC-construction with the input data
$(l,\xi,h)$, applied to $v$.
Then $w$ is also an $f$-gradient annd
$\Vert w-v\Vert\leq \e/2$.
(Similar assertion holds for adding a horizontal component
nearby $\daw$.)
\enre

Let $\phi_0:\pr_0 W\to\RRR$ be a Morse function and
$u_0$ be a $\phi_0$-gradient, satisfying the conclusions of
Theorem \ref{t:quickpush}
with respect to
$C=A, t=t_0$. We can obviously assume that $\phi_0$ is ordered;
let
$(\a_0,... \a_{n})$ be the corresponding ordering sequence.

Similarly,
let $\phi_1:\pr_1 W\to\RRR$ be a Morse function and
$u_1$ be a $\phi_1$-gradient, satisfying the conclusions of
Theorem \ref{t:quickpush}
with respect to
$C=A, t=t_0$. We  assume that $\phi_1$ is ordered;
let
$(\b_0,... \b_{n})$ be the corresponding ordering sequence.

Assume that $u_0,u_1$ are $\d$-separated. We obtain the corresponding
$gfn$-systems
$\DDD(u_1), \DDD(u_0), \DDD(-u_1),
\DDD(-u_0)$ with the cores $\dd(u_0), \dd(u_1),...$.
We fix $]0,\d[$ as the interval of definition for these systems.

The diffeomorphism $\Psi_1$ identifies $V_b$
with
$V_\b$, therefore we can transfer to $V_\b$
the function $\phi_1$ and  its gradient $u_1$. The resulting
function and gradient
will be denoted by
$\phi_1'$, resp. $u'_1$.
We obtain also $gfn$-systems:
$\AAA^-=\Psi_1\(\DDD(u_1)\)$ and
$\AAA^+=\Psi_1\(\DDD(-u_1)\)$
in $V_\b$ with the cores $\dd(u_1')$, resp. $\dd(-u_1')$.
({\it Warning:}
$\Psi_1$ is not in general an isometry, so $\AAA^-$ need not
be the same as $\DDD(u_1')$.)
A similar construction gives $gfn$-systems
$\BBB^-, \BBB^+$ on $V_\a$
with the cores $\dd(u_0)$, resp. $\dd(-u_0)$.

\subsection{ $C^0$-density of $\GgC(f)$: perturbation of a given
gradient}\label{su:cdpgg}

In this subsection we  finish the proof of Theorem \ref{t:cc}.
 Namely, we show
how to perturb $v_0$ in order to obtain a gradient $v$, satisfying
the condition $(\gC)$.
\pa
{\it A. First modification of the gradient}
\pa
Applying Proposition \ref{trper}
we obtain an $(f\mid W')$-gradient $v_1$,
having the following properties:
\been
\item[(i)] $v_1$ satisfy \ATA;
\item[(ii)] $v_1=v$ in a \nei~of $V_\a\cup V_\b$;
\item[(iii)] the $s$-submanifolds
$\ttt(\dd(u_1'),-v_1)$ and  $\ttt(\dd(-u_0'), v_1)$
are defined
\item[(iv)]
$\stind {(-v_1)}\b\a\(\dd(u_1'))\nmid \dd(-u_0')$  and
$\stind {(v_1)}\a\b \(\dd(-u_0')\)\nmid \dd(u_1')$.
\enen
Moreover, $v_1$ can be chosen arbitrarily $\smo$ close to
$v_0$, thus we can assume
$\Vert v_1-v_0\Vert<\e/2$.
Recall from
Subsection \ref{su:totss}
that there is $\l>0$ \sut~
the pairs
$(\AAA, v_1)$ and $(\BBB, -v_1)$ are $\l$-ordered.
Therefore we can consider $]0,\l[$ as the interval of definition
of
the $gfn$-systems
$\TTT(\AAA^-, -v_1), \TTT(\BBB^+, v_1)$.

Applying Proposition \ref{tracgfn}
we deduce from (iv) above
that there $\nu<\l$, \sut~
for every $s$:
\begin{multline}
\stind {(-v_1)}\b\a \bigg(\AAA^-_{(\leq s)}(\nu)\bigg)
\cup \bigg(D_\nu(\indl {s+1}; v_1)\cap V_\a\bigg)\\
\mbox{ does not intersect }
\BBB^+_{(\leq n-2-s)}(\nu) \lbl{intmp1}
\end{multline}
and
\begin{multline}
\stind {(v_1)}\a\b \bigg(\BBB^+_{(\leq s)}(\nu)\bigg)
\cup \bigg(D_\nu(\indl {s+1}; -v_1)\cap V_\b\bigg)\\
\mbox{ does not intersect }
\AAA^-_{(\leq n-2-s)}(\nu) \lbl{intmp0}
\end{multline}

We shall choose $\nu$ so small that $v_1$ is $\nu$-separated.
Extend the vector field $v_1$ to the whole of $W$ setting
$v_1(x)=v_0(x)$ for $x\in W\sm W'$.
Then $v_1$ is a $\nu$-separated $f$-gradient. Choose
an ordered Morse function $\phi$ on $W'$, adjusted to
$(f,v_1)$ and extend $\phi$ to the whole of
$W$, setting $\phi(x)=f(x)$ for $x\in W\sm W'$.

\pa
{\it B. Second modification of the gradient}
\pa
Here we start with the gradient $v_1$
constructed
in  previous subsection.
The second modification consists in adding to $v_1$
horizontal components
 nearby $\pr_1 W$ and $\pr_0 W$.
We have chosen a number $\nu>0$.
Choose now (by Theorem \ref{t:quickpush})
a $\phi_0$-gradient $\wi u_0$ and a
$\phi_1$-gradient $\wi u_1$, equivalent to
$u_0$, resp. to $u_1$
and \sut~
$\wi u_0, \wi u_1$
are
$(t_0,\nu)$-rapid
and satisfy
$\Vert \wi u_0\Vert\leq A,
\Vert \wi u_1\Vert\leq A$.
 That is, for $i=0,1$
and for every $s$ we have
$$\Phi(\wi u_i, t_0)(C_\nu(\indl s; \wi u_i))
\sbs
B_\nu(\indl s; \wi u_i)
$$
$$\Phi(-\wi u_i,t_0)(C_\nu(\indl s; -\wi u_i)
\sbs
B_\nu(\indl s; -\wi u_1)$$

Apply to $v_1$ the AHC-construction nearby $\pr_0 W$
with the input data $(l,\wi u_0, h)$
and nearby $\pr_1 W$ with
the input data
$(l,\wi u_1, h)$
(recall that $l$ and $h$ were chosen in the beginning of the
subsection \ref{su:cdt}).

Denote the resulting vector field by $v$. Then
 $\Vert v-v_0\Vert\leq \e$,
and that $v$ is an $f$-gradient
(see Remark \ref{r:choice}).
Since $\supp (v_1-v)\sbs W\sm W'$ and
$f\mid W\sm W' =\phi\mid W\sm W'$
the vector field $v$ is also a $\phi$-gradient.
To compute $\st v$ and $\stv$, denote the diffeomorphism
$\stind {(-v)}b\b$ by $\wi \Psi_1$ and
$\stind va\a$ by $\wi \Psi_0$.

From the construction  of $v$ it follows
(see AHC2), page \pageref{ahc2})
that for every $s$ we have:

\begin{equation}
\wi\Psi_1\bigg(C_\nu(\indl s; u_1  )  \bigg)
\sbs \AAA_{(\leq s)}^-(\nu)\lbl{qp1-}
\end{equation}
\begin{equation}\wi\Psi_1^{-1}
\bigg(V_\b\sm\AAA^+_{(\leq s)}(\nu)\bigg)\sbs
B_\nu(\indl {n-2-s}; -u_1) \lbl{qp1+}
\end{equation}

Similarly,
\begin{equation}
\wi\Psi_0\bigg(C_\nu(\indl s; -u_0    \bigg)
\sbs \BBB^+_{(\leq s)}(\nu)
\lbl{qp0+}
\end{equation}
\begin{equation}
\wi\Psi_0^{-1}\bigg(V_\a\sm\BBB^+_{(\leq n-s-2)}(\nu)\bigg)
\sbs
B_\nu(\indl s; u_0)\lbl{qp0-}
\end{equation}

Write now
$$\stind {(-v)}ba = \wi \Psi_0^{-1}\circ
\stind {(-v_1)}\b\a \circ \wi\Psi_1$$

$$\stind {v}ba = \wi \Psi_1^{-1}\circ
\stind {(v_1)}\a\b \circ \wi\Psi_0$$

and combine
(\ref{qp1-} - \ref {qp0-})
and (\ref{intmp1}, \ref{intmp0})
to note that (B1) and (B0) hold
 for $v$ with $\d=\nu$ and for every $s$.
The proof of Theorem \ref{t:cc} is
finished. $\qs$

\subs{Cyclic cobordisms and condition $(\gC\YY)$ }
\lb{su:cyc}

The condition $(\gC)$
and the related techniques, developed in this chapter
are the basic ingredients in the proof of the main results of our work
(Theorems \ref{t:cccgr}, \ref{t:eqcccgr}).
In this theorems we apply the results of the present chapter in
the special situation when the
~parts $\daw, \dow$ of the
\co~~  are diffeomorphic ( this \co~ is obtained as the result of cutting
a \ma~ $M$ along a regular value of a Morse map $M\to S^1$).
 The present subsection contains
a version of Theorem \ref{t:cc}, adapted to this situation.

\bede\lb{d:cyc}
{\it A cyclic cobordism} $W$ is a riemannian cobordism
together with an isometry
$\Phi:\dow\to \daw$.
\end{defi}

Let $f:W\to[a,b]$ be a Morse function
on a cyclic \co~ $W$ of dimension $n$,
$v$ be an $f$-gradient.

\begin{defi}
\lb{d:condcy}
 We say that $v$ satisfies
$(\gC\YY)$ if $v$ satisfies condition $(\gC)$
(see \sub~ \ref{su:cc}),
and, moreover, the Morse functions $\phi_1, \phi_0$ and their gradients
$u_1, u_0$ from the definition \ref{d:cc}
can be chosen so that
$\phi_1\circ\Phi=\phi_0, u_1=\Phi_*(u_0)$.
The set of all $f$-gradients satisfying
$(\gC\YY)$
will be denoted by
$\GgCY(f)$.
\end{defi}

\beth\lb{t:cyc}
$\GgCY(f)$ is open and dense in $\GG(f)$ \wrt~  $C^0$ topology.
Moreover, if $v_0$ is any $f$-gradient then one can choose a $C^0$
 small perturbation $v$
of $v_0$ \sut~ $v\in\GgCY(f)$ and
$v=v_0$ in a \nei~ of $\pr W$.
\end{theo}
\Prf It goes exactly as the proof of Theorem \ref{t:cc}.
To prove $C^0$-openness of the set
$\GgCY(f)$ we do not need to change whatever it be in the first
part of the proof of Theorem \ref {t:cc}.

For the proof of $C^0$ density just note that
in the proof of $C^0$ density of $(\gC)$
(\sub~ \ref{su:cdt})
we can choose $\phi_0=\phi_1\circ\Phi, u_1=\Phi_*(u_0)$.
$\qs$

\section{Homological gradient descent}\lb{s:hgd}

Now we shall deduce consequences from the condition $(\gC)$.
In this section $\fcob$ is a Morse function on a riemannian \cob~ $W$
of dimension $n$, and $v$ is an $f$-gradient,
satisfying condition $(\gC)$.
Recall that the condition $(\gC)$ requires the existence of Morse
functions
$\phi_0:\pr_0 W\to\RRR, \phi_1:\pr_1 W\to\RRR$
and their gradients $u_0, u_1$, satisfying the conditions (B1) and
(B0) from 
Subsection \ref{su:cc}
There must be  also an ordered Morse function $\phi:W\to[a,b]$,
adjusted to $(f,v)$ with an ordering sequence
$a=a_0<a_1<...<a_{n+1}=b$.
We shall denote $\phi^{-1}(\l)$ by $V_\l$, so that
$\pr_1 W=V_b, \pr_0 W=V_a$.
From the gradients $u_0, u_1$ we obtain the filtrations of
$V_b$ and of $V_a$, namely
$V_b^{\{\leq s\}}, V_a^{\{\leq  s\}}$
($s$ ranges over integers between $0$ and $n-1$).
We have proved  that every $f$-gradient $w$, sufficiently
$C^0$ close to $v$
satisfies the condition $(\gC)$ \wrt~the same objects
$\phi_0, \phi_1, u_0, u_1, \phi$ (see Subsection \ref{su:cogct}).

\begin{theo}\label{t:consh}
For every integer $s:0\leq s\leq n$ there is a homomorphism
$$\HH_s(-v): H_*(V_b^{\{\leq s\}}, V_b^{\{\leq s-1\}} )
\to
H_*(V_a^{\{\leq s\}}, V_a^{\{\leq s-1\}} )
$$
with the following properties:
\begin{enumerate}\item[1)]
Let $N$ be an oriented \sma~ of $V_b$, contained in $V_b^{\{\leq s\}}$
\sut~$N\sm\Int V_b^{\{\leq s-1\}}$ is compact.
Then the manifold $N'=\stind {(-v)}ba (N)$ is in $V_a^{\{\leq s\}}$ and
$N'\sm \Int V_a^{\{\leq s-1\}}$
is compact and the fundamental class of $N'$ modulo
$V_a^{\{\leq s-1\}}$ equals to
$\HH_s(-v)\([N]\)$.
\item[2)] 
There is an $\e>0$ \sut~for every $f$-gradient $w$ with
$\Vert w-v\Vert\leq \e$ and every $s$
we have:
$\HH_s(-v)=\HH_s(-w)$.
\end{enumerate}
\end{theo}
The proof of the theorem occupies Subsections
\ref{su:midp},
\ref{su:sykxk},
\ref{su:poh}.

\pa
\subsection{ Main idea of the proof and the definition of $\HH_s(-v)$}
\label{su:midp}
We shall first give an informal description of $\HH_s(-v)$. Let 
$\a\in H_m(V_b^{\{\leq s\}},
V_b^{\{\leq s-1\}})$.
Assume for simplicity that $\a$ is represented
 by a closed oriented submanifold
$A$ of $V_b$, $\dim A=m$,
contained in $V_b^{\{\leq s\}}$.
Denote $V_{a_{s+1}}$ 
by $V'$ for brevity.
We descend $A$ along the $(-v)$-trajectories from $V_b$ to $V'$.
The condition $(\gC)$ implies
that the whole of $V_b^{\{\leq s\}}$ descends to $V'$
without hanging on critical points of $f$,
that is, every
$(-v)$-trajectory starting in $\vbs$ reaches $\vass$.
Therefore $A'=\stind {(-v)}b{a_{s+1}}(A)$
is a closed oriented submanifold of $V'$.
Try now to apply the same procedure in order to descend $A'$ to $V_a$.
In general $A'$  intersects the upper soles of the ascending
 discs $D(p,-v)$ with $\ind p\leq s$,
and so
$\stind {(-v)}{a_{s+1}}a (A')$
is a submanifold of $V_a$, which is not closed.
But note that the intersection of the upper soles of these ascending
discs with $A'$ is contained in the \nei~
$\cup_{\ind p\leq s} B_\d(p,-v)$
of
$\cup_{\ind p\leq s} D(p,-v)$,
and this \nei~
is carried to $\vasm$ by $\stind {(-v)}{a_{s+1}}a$.
So the fundamental class of $\stexp {(-v)} (A')$ can still be defined
modulo
$\vasm$, and this was our aim.

We proceed now to the precise definitions.
The homomorphism $\HH_s(-v)$ will be defined as the composition of two
homomorphisms: $\HH 1_s(-v)$ and $\HH 0_s(-v)$.
(Informally, $\HH 1$ corresponds to the descent from the level 
$b$ to the
level $a_{s+1}$,
and $\HH 0$ corresponds to the descent from
$a_{s+1}$ to $a$.)
We shall denote the maps
$\stind {(-v)}b{a_{s+1}}$
and
$\stind {(-v)}{a_{s+1}}a$
by $\stexp {(-v1)}$, and respectively $\stexp {(-v0)}$. Every $(-v)$-
trajectory
starting in $\vbs$ reaches $\vass$.
Therefore $\vbs$ is in the domain of
$\stexp {(-v1)}$
and
$\stexp {(-v1)}$
defines a homeomorphism of the pair
$(\vbs, \vbsm)$ to its image, which is a pair of compact
 subsets of $V'$.

Now let $A,B,C$ be subsets of $V'$, satisfying
the following conditions (R1) - (R4):
\begin{equation*}
A\supset B\supset C;\quad \ove{C}\sbs \Int B
 \tag{R1}
\end{equation*}
\begin{equation*}
\stexp{(-v1)}(\vbs)\sbs A, \quad
\stexp{(-v1)}(\vbsm)\sbs B,
\tag{R2}
\end{equation*}
\begin{equation*}
\bigg(\bigcup_{\ind p\leq s} D(p,-v)\bigg)\cap V'\sbs C
\tag{R3}
\end{equation*}
\begin{equation*}
\stexp {(-v0)}(A)\sbs \vas, \quad \stexp{(-v0)} (B)\sbs \vasm
\tag{R4}
\end{equation*}

A bit later we shall check that such subsets exist.
Now we shall
construct $\HH_s(-v)$ assuming that such subsets exist, and
 we shall show that
 $\HH_s(-v)$ does
not depend on the choice of $A,B,C$ satisfying  (R1) - (R4).

The composition of $\stexp {(-v1)}$
with the inclusion defines a continuous map
$$
\phi: (\vbs, \vbsm)\to (A,B)$$
Set $\HH 1_s(-v)=\phi_*$.
Consider the excision isomorphism
\begin{equation}
Exc: H_*(A\sm C, B\sm C)\overset{\approx}{\rTo} H_*(A,B)\lbl{f:exc}
\end{equation}
The conditions (R3) and (R4)
imply that the gradient descent from
$a_{s+1}$ to $a$ defines a continuous map of pairs
\begin{equation}
\rho: (A\sm C, B\sm C)\to (\vas, \vasm)\lbl{f:ro}
\end{equation}
Set $\HH 0_s(-v)=\rho_*\circ (Exc)^{-1}$, and
  $\HH_s(-v)=\HH 0_s (-v)\circ \HH 1_s(-v)$.

Now we  check that the homomorphism $\HH_s(-v)$
does not depend
on the choice
of $A,B,C$.
Let $(A,B,C), (A', B', C')$ be two triplets satisfying (R1) - (R4);
denote the corresponding homomorphisms by
$\HH, \HH'$.
It is easy to see that $A\cap A', B\cap B', C\cap C'$
satisfy again the conditions (R1) - (R4).
Therefore it suffices to verify  $\HH=\HH'$ assuming
$A\supset A', B\supset B', C\supset C'$. This follows from an
easy functoriality
argument, which we  leave to the reader.

In the following subsection we give two examples of triples $(A,B,C)$,
 satisfying
the conditions (R1) - (R4).

\subsection{Sets $T_k(v), t_k(v)$}
\label{su:sykxk}

Let $k$ be an integer, $0\leq k\leq n$ .
\begin{defi}
Denote by $T_k(v)$ the set of all
$y\in\phi^{-1}([a,a_{k+1}])$
\sut~either

i) $\g(y, \cdot; -v)$ reaches $V_a$ and intersects it at a point
$z\in \vakm$ or

ii) $\lim_{t\to\infty}\g(y,t;-v)=p$, where $p\in S(f)$.
\end{defi}

One can reformulate this definition slightly. Consider
the cobordism
$W'=\phi^{-1}([a,a_{k+1}])$.
We have the restriction
$\phi\mid W':W'\to[a,a_{k+1}]$
of the Morse function $\phi$
and the restriction $v\mid W'$ of its gradient $v$.
Then $T_k(v)=
T\bigg(\vakm, v\mid W'      \bigg)
\cup
D(-v\mid W')$
(see the definition of $T(\cdot,\cdot)$ in Subsection (\ref{su:gpt})).
This implies immediately that $T_k(v)$ is compact.
Before proceeding to the following definition recall that we
have fixed $\d>0$ satisfying $(\gC)$.
\begin{defi}
Denote by $t_k(v)$ the set of all
$y\in\phi^{-1}([a,a_{k+1}])$
\sut~either

i) $\g(y, \cdot; -v)$ reaches $V_a$ and intersects it at a point
$z\in \Vakm$ or

ii) $\lim_{t\to\infty}\g(y,t;-v)=p$, where $p\in S(f)$.
\end{defi}

We have obviously $t_k(v)\sbs T_k(v)$. The same argument as used for
$T_k(v)$ shows that $t_k(v)$ is compact.

\begin{lemm}
\label{l:depri}
There is $\d'>\d$ \sut~$v$ is $\d'$-separated and
\begin{equation}
\bigcup_{\ind p\leq k} D_{\d'}(p,-v)\cap V_{a_{k+1}}\sbs t_k(v)
\lbl{f:depri}
\end{equation}
\end{lemm}
\Prf
The condition (B1) implies that
$$
D_\d(\indl k;v)\cap\pr_0 W\sbs B_\d(\indl {k-1}; u_0)
$$
Since $v$ is $\d$-separated, it is $\d''$-separated for some
 $\d''>\d$ and
the interval of definition of the $gfn$-system
$\DDD(v)$ contains $]0,\d''[$.
Therefore
$\{B_\theta(\indl k;v)\cap\pr_0 W\}_{\theta>\d}$
is a fundamental system of neighborhoods for
$D_\d(\indl k;v)\cap\pr_0 W$.
Thus for some
$\d':\d<\d'<\d''$
we have:
$D_{\d'}(\indl k;v)\cap\pr_0 W\sbs 
B_\d(\indl {k-1}; u_0)$. The lemma follows.
$\qs$
\begin{rema}\label{r:depri}
The previous lemma implies that
$$
t_k(v)\cap V_{a_{k+1}}\sbs \Int\big( T_k(v)\cap V_{a_{k+1}}\)$$

Indeed, the inclusion $t_k(v)\sbs T_k(v)$ is obvious.
Further, let $z\in t_k(v)\cap V_{a_{k+1}}$.
If $z\in \cup_{\ind p\leq k} D(p,-v)$, then
$z$ is contained in
$$
\bigg(\bigcup_{\ind p\leq k} B_\d(p,-v)\bigg)\cap V_{a_{k+1}}$$
which is in $t_k(v)$.
If $\g(z,\cdot;-v)$ reaches $V_a$, then it intersects
$V_a$ at a point
$z'\in \Vakm$,
which is in the interior of $\vakm$, and this implies our assertion.
A similar argument proves that  $t_k(v)\cap\phi^{-1}\([a,a_{k+1}[)$
is in
$\Int\(T_k(v)\cap\phi^{-1}\([a,a_{k+1}[\)\)$. 

Returning to the definition of sets $A,B,C$ from 
Subsection \ref{su:midp},
set
\begin{equation}
\left\{
\beal
A &=T_{s+1}(v)\cap V',\\
B &=T_s(v)\cap V',\\
C &=\bigg(\bigcup_{\ind p\leq s} D(p,-v)\bigg)\cap V'
\enal\quad \lbl{f:triodin}
\right.
\end{equation}
It is now obvious that this triple satisfies (R1) - (R4). Another
triple satisfying the conditions is
\begin{equation}
\left\{
\beal 
A &=t_{s+1}(v)\cap V',\\
B &=t_s(v)\cap V';\\
C &=\bigg(\bigcup_{\ind p\leq s} D(p,-v)\bigg)\cap V'
\enal\lbl{f:tridva}
\right.
\end{equation}\end{rema}

\subsection{Properties of $\HH_s$}\label{su:poh}
\pa
{\it Proof of 1)}
\pa
Consider the submanifold $N''=  \stexp {(-v1)}(N)$ 
of $V'=V_{a_{s+1}}$.
Since
$\stexp {(-v1)}$ is a diffeomorphism of
 the open subset $\Int \vbsm$ of $V_b$ onto its
image
in $V'$,
the manifold
$N''\sm \stexp {(-v1)}\(\Int \vbsm\)$
is compact.

Thus the compactness of
$N'\sm\Int V_a^{\{\leq s-1\}}$
follows from the next lemma.

\bele\lb{l:compa}
Let $L$ be a \sma~ of $V'$, \sut~
$L\sm \stexp {(-v1)}(\Int V_b^{\{\leq s-1\}} )$
is compact.
Set $L'=\stexp {(-v0)}(L)$.
Then 
$L'\sm \Int V_a^{\{\leq s-1}\}$
is compact.
\enle
\Prf
The set $L\sm\Int T_s$
is compact. Since $\stexp {(-v0)}$ is defined 
on $V'\sm \Int T_s$, the set
$$
\stexp {(-v0)} (L'\sm\Int T_s)
=
\stexp {(-v0)} (L')\sm \stexp {(-v0)}(\Int T_s)
$$
is  a compact subset in $V_a$.
Note that $\Int \vasm$ is an open set containing
$\vovo(\Int T_s)$, and this implies finally that
$L'\sm\Int\vasm$ is compact.
$\qs$

To verify that
$[N']=\HH_s(-v)\([N]\)$, use
the definition of
$\HH_s(-v)$ with the help of
$A,B,C$ as in (\ref{f:triodin})
and check through the definitions of $N'$ and of $\HH_s(-v)$.
We leave the details to the reader.
\pa
{\it Proof of 2)}
\pa
Here it is convenient to use for the definition of $\HH_s(-v)$ the
triple
$(A,B,C)$, where

\begin{equation}
\left\{
\beal
A &=T_{s+1}(v)\cap V',\\
B &=T_s(v)\cap V';\qquad\lbl{f:tritri}\\
C &=\bigg(\bigcup_{\ind p\leq s} B_\d(p,-v)\bigg)\cap V'
\enal
\right.
\end{equation}
In the rest of this subsection we shall abbreviate $\vbs$ to $R_b$, 
$\vbsm$ to $Q_b$,
$\vas$ to $R_a$, $\vasm$ to $Q_a$.
First we shall show that the triple $(A,B,C)$ defined in (\ref{f:tritri})
satisfies the conditions (R1) - (R4) not only for $v$ but also for 
every $\phi$-gradient $w$ sufficiently $C^0$ close to $v$.
Of course it is true for (R1), since it does not depend on the gradient.

Proceeding to (R2) recall that (by Remark \ref{r:depri})
$t_{s+1}(v)\cap V'\sbs \Int A$, and
$t_{s}(v)\cap V'\sbs \Int B$.
Therefore the pair of
compacts
$$\bigg(\vov (R_b),\vov (Q_b)\bigg)$$
is in $(\Int A, \Int B)$.
Now it follows from  Corollary \ref{co:homot}
that there is $\e>0$ \sut~
for every vector field $w$ on
$\phi^{-1}\([a_{s+1}, b]\)$ with $\Vert w-v\Vert<\e$
each $(-w)$-trajectory
starting in $\vbs$ reaches $V'$ and 
\begin{equation}
\bigg(\stexp {(-w)} (R_b), \stexp {(-w)} (Q_b)\bigg)
\sbs
(\Int A, \Int B)\label{f:incl}
\end{equation}
Therefore for every $\phi$-gradient $w$ with $\Vert w-v\Vert<\e$
the condition (R2) is true \wrt~$w$.

Since $C$ is open, the
 property (R3) \wrt~$w$ is true if $w$ is sufficiently close
to $v$ by Proposition \ref{p:discpert}.
Proceeding to (R4) 
note that
$$
\vovo(A)\sbs \Vas\sbs\Int \vas
$$
So each $v$-trajectory starting in
$V_a^{\{\geq s\}}$
reaches $V'$
and intersects it at a point outside $A$. Then the same is
true for every $w$-trajectory
starting in $V_a^{\{\geq s\}}$
if only $w$ is sufficiently $C^0$ close to $v$.
This implies the first half of the property (R4) \wrt~$w$.
Similarly we verify the second half.

Thus we can use the triple $(A,B,C)$ defined in (\ref{f:tritri})
for the definition of $\HH_s(-w)$.
It is left to show that 
$\HH_s(-w)=\HH_s(-v)$
for all $w$ sufficiently $C^0$ close to $v$.

\bele
\lb{l:h1}
For every $\phi$-gradient $w$
sufficiently $C^0$-close to $v$ we have:
$\HH 1_s(-w)=\HH 1_s(-v)$.
\enle
\Prf
Apply the Corollary \ref{co:hompair}
to the \cob~ 
$\phi^{-1}\([a_{s+1}, b]\)$, and the vector field $(-v)$.
$\qs$

\bele
\lbl{l:h0}
For every $\phi$-gradient $w$
sufficiently $C^0$-close to $v$ we have:
$\HH 0_s(-w)=\HH 0_s(-v)$.
\enle
\Prf Recall  that
$A\sm C$ is in the domain of definition of
$\stind {(-v)}{a_{s+1}}a$.
Moreover, the condition $(\gC)$ implies that
$\stind {(-v)}{a_{s+1}}a$
sends the pair of compacts
$(A\sm C, B\sm C)$
to the pair
$
\(\Vas, \Vasm\)
$
which is in the pair of open sets
$(\Int R_a, \Int Q_a)$.
Therefore (by Corollary \ref{co:hompair})
for every $w$ sufficiently $C^0$ close to
$v$ the map
$\stind {(-w)}{a_{s+1}}a$
is defined on $A\sm C$, sends
the pair
$(A\sm C, B\sm C)$ to
$(\Int R_a, \Int Q_a)$,
and is homotopic to 
$\stind {(-v)}{a_{s+1}}a$
as a map of these pairs. $\qs$

\subsection{ A pairing between $C_s(u_1)$ and $C_s(v)$}
\label{su:couple}

Here we present a construction related to
the homological gradient descent \ho.
This construction will be used in Chapter 4
for computation of Novikov incidence coefficients.

Let $0\leq s\leq n-1$.
As in the previous \sub~ we denote
 $\phi^{-1}(a_{s+1})$ by $V'$.
Every $(-v)$-trajectory starting 
at
$C_\d(\indl s; u_1)$ reaches $V'$.
Therefore the map
$\stind {(-v)}b{a_{s+1}}$
defines a homeomorphism
of the pair
$$\(\Int C_\d(\indl s; u_1),\quad
\Int C_\d(\indl {s-1}; u_1)\)$$
to the pair
$(N_s, N_{s-1})$
of open subsets of $V'$.
Let $\tau$ be the restriction of
$\stind {(-v)}b{a_{s+1}}$
to the pair
$(\pws, \pwsm)$.
Then $\Im\tau$
is a compact subpair $(n_s, n_{s-1})
\sbs (N_s, N_{s-1})$.
Let $q\in S_s(f)$. Consider the submanifold
$S(q)=D(q,-v)\cap V'$ of $V'$;
it is a $(n-s-1)$-dimensional
cooriented embedded sphere.
Denote by $]S_q[$ the cohomology class in
$H^{s}(V', V'\sm S(q))$,
dual to
$S(q)$.
We have:
$S(q)\cap n_{s-1}=\emp$
(this follows from $(\gC)$).
Therefore the pair 
$(N_s, N_{s-1})$
is a subpair of
$(V', V'\sm S(q))$
and we can induce the class $]S_q[$ to the cohomology
of $(N_s, N_{s-1})$.
The resulting cohomology class will be denoted again by $]S_q[$.
For $x\in C_s(u_1)=H_s(\pws, \pwsm)$ set
$$
\langle q, x\rangle_v=\langle ]S_q[, \tau_*(x)\rangle
$$
(The index $v$ reminds that the pairing depends on the vector field $v$.)
This pairing has some natural properties, gathered
in the  proposition \ref {p:couple}.

\bepr\lb{p:couple}
\been\item
Let $L$ be an $s$-dimensional oriented submanifold of $\daw$, \sut~
$L\sbs \pws$
and $L\sm \Int \pwsm$
is compact. Then
$L'=
\stind {(-v)}b{a_{s+1}}(L)$ is an $s$-dimensional oriented submanifold
of $V'$, \sut~$L'\sbs n_s$
and
$L'\sm \Int n_{s-1}$
is compact.
If $L'\pitchfork S(q)$
 then
$$
L'\krest S(q)=\langle q, [L]\rangle_v
$$
\item
Let $U$ be a \nei~ of $S(f)$.
There is $\d>0$, \sut~ for every $\phi$-gradient
$w$ with $v\vert U=w\vert U$
and
$\Vert v-w \Vert<\d$ we have:
$\langle q, x\rangle_v =
\langle q, x\rangle_w$
\enen
\enpr

\Prf 1) is obvious. Proceeding to (2), 
note first that the objects 
\break
$N_s, N_{s-1}, n_s, n_{s-1},
 \tau, S_q$ etc.,
introduced in the beginning of the present subsection 
depend on $v$, and we shall now show it 
in the notation, writing, respectively
\break
$N_s(v), N_{s-1}(v), n_s(v), n_{s-1}(v),
 \tau(v), S_q(v)$ etc.  in order to make explicit the
 dependence on the gradient.

Recall that if $w$ is sufficiently $C^0$ close to $v$
then $w$ still satisfies $(\gC)$ \wrt~ 
the same
string 
$\SS=(\d,u_0, u_1,\phi_0, \phi_1,\phi)$.
In particular, every $(-w)$-trajectory starting in
$\pws$
reaches $V'$.
It follows easily that if only $w$ is sufficiently close to
$v$, we have:
$n_s(w)\sbs N_s(v), n_{s-1}(w)\sbs N_{s-1}(v)$.

The compacts $n_{s-1}(v), S(q,v)$
are disjoint subsets of $V'$.
Therefore there are disjoint open subsets
$\L,\G$ of $V'$ with $n_{s-1}(v)\sbs\L, S(q,v)\sbs\G$.
We have
the following inclusions:
\begin{equation}
\(n_s(v), n_{s-1}(v)\)
\rInto^{i(v)}
(V', \L)
\rInto^{j}
(V', V'\sm \G)
\rInto^{k(v)}
(V', V'\sm S(q,v))\lbl{f:i3}
\end{equation}

Note that, obviously
\begin{equation}
\langle p, x\rangle_v
=
\big\langle  (k(v))^*\(]S_q(v)[\), 
(j\circ i(v)\)_*(\tau_*(v)(x)\)\big\rangle
\lbl{f:inter}
\end{equation}

If $w$ is sufficiently $C^0$ close to $v$, then
still
$n_{s-1}(w)\sbs\L;, S(q,w)\sbs\G$
and the formula
(\ref{f:inter})
still holds if we replace $v$ by $w$.
Therefore it suffices to show that for $w$ sufficiently $C^0$ close to
$v$ we have
\begin{equation}
i(v)\circ\tau(v)
\sim i(w)\circ\tau(w), \quad
\(k(v)\)^* \(]S_q(v)[\)
=
\(k(w)\)^* \(]S_q(w)[\)
\end{equation}

The first  follows from Corollary \ref{co:hompair}.
To prove the second,
note that the class 
$$
\(k(v)\)^* \(]S_q(v)[\)
\in H^s(V', V'\sm \G)
$$
is the Lefschetz-dual to the fundamental class
$[S_q(v)]\in H_{n-s-1}(\G)$.
Therefore
it suffices to prove
that if $v|U=w|U$
and $w$ is \sclv
then the fundamental classes of $S_q(v)$ and of $S_q(w)$
in $H_*(\G)$ are equal.
Now just apply Proposition \ref{p:dischom}. $\qs$

\subs{Homological gradient descent: equivariant version}
\label{su:eqhgd}

Let $p:\wh W\to W$
be a regular covering of $W$ with structure group $H$.
For a subset $A\sbs W$ we denote by $\wh A$ the set
$p^{-1}(A)$.
We obtain thus an $H$-invariant filtration
$\wh V_b^{\{\leq s\}}$
of $\wh V_b$ and 
an $H$-invariant filtration
$\wh V_a^{\{\leq s\}}$ of $\wh V_a$.
Recall from Section \ref{s:bdc} the notation $K^+=
\cup_{x\in S(f)} D(x, -v),
 K^-=
\cup_{x\in S(f)} D(x, v)$.
For $x\in\wh{\daw}\sm \wh K^+$
consider the $(-v)$-trajectory 
$\g(q(x),t;-v)$.
It reaches $\dow$, and the lifting 
$\wh\g$ of $\g$
to $\wh W$ is an integral curve of $-p^*v$.
This integral curve
 reaches $\wh\dow$, and the
 intersection point of $\wh\g$ with $\wh\dow$
will be denoted by
$\stind {(-\wh v)}ba (x)$.
We obtain thus a diffeomorphism
$$
\stind {(-\wh v)}ba (x):
\wh{\daw}\sm \wh K^+
\rTo
\wh{\dow}\sm \wh K^-
$$
which commutes with the action of $H$.

The \hog~ 
$H_*\(\hpws, \hpwsm\)$
is a right $\ZZZ H$-module, and Proposition \ref{p:fil}
implies that it is a free $\ZZZ H$-module with the base formed
by liftings of descending discs 
of critical points of $\phi_1$ of index $s$.
The following theorem is an "equivariant version"
of theorem \ref{t:consh}.
The proof is a direct generalization
of the proof of Theorem
\ref{t:consh}
and we leave it to the reader.

\beth\label{t:eqhgd}
For every integer $s:0\leq s\leq n$ there is a homomorphism
$$\hHH_s(-v): H_*\(\hV_b^{\{\leq s\}}, \hV_b^{\{\leq s-1\}} \)
\to
H_*\(\hV_a^{\{\leq s\}}, \hV_a^{\{\leq s-1\}} \)
$$
with the following properties:
\begin{enumerate}\item[1)]
Let $N$ be an oriented \sma~ of $\hV_b$, contained in 
$\hV_b^{\{\leq s\}}$ and 
\sut~$N\sm\Int \hV_b^{\{\leq s-1\}}$ is compact.
Then the manifold $N'=\stind {(-\wh v)}ba (N)$ is in
 $\hV_a^{\{\leq s\}}$ 
and
$N'\sm \Int \hV_a^{\{\leq s-1\}}$
is compact and the fundamental class of $N'$ modulo
$\hV_a^{\{\leq s-1\}}$ equals to
$\hHH_s(-v)\([N]\)$.
\item[2)] 
There is an $\e>0$ \sut~for every $f$-gradient $w$ with
$\Vert w-v\Vert\leq \e$ and every $s$
we have:
$\hHH_s(-v)=\hHH_s(-w)$. $\qs$
\end{enumerate}
\end{theo}

We proceed to the equivariant analogs of the results of \sub~
\ref{su:couple}. Let $0\leq s\leq n-1$.
We assume the terminology of
\sub~ \ref{su:couple}.
We denote
$H_s\(\hpws, \hpwsm\)$
by $\wh C_s(u_1)$. It is a free right $\ZZZ H$-module.
There is an $H$-equivariant homeomorphism
$$
\wh\tau(v):
\(\hpws, \hpwsm\)\to
(\wh n_s(v), \wh n_{s-1}(v))
$$
Let $q\in S_s(f)$ and choose a lifting $\wh q$ of $q$ to $\wh W$.
Then we obtain  a lifting 
$ S(\wh q,v)$ of $S(q, v)$
to $\wh W$, and a \cog~ class
$] S_{\wh q}  [\in H^s(\wh N_s(v), \wh N_{s-1}(v))$.

For $x\in H_s(\hpws, \hpwsm)$
set
$$
\langle\wh q, x\rangle_v=
\langle]S_{\wh q}  [, (\wh\tau(v))_*(x)\rangle
$$
Thus 
$\langle\wh q, x\rangle_v$
is an integer. If we fix $\wh q$
and consider $x$ as the variable,
 we obtain a \ho~$\wh C_s(u_1)\to\ZZZ$.

Next we give an  equivariant version
of Proposition
\ref{p:couple}.
The first two points of this proposition
are direct generalization of the corresponding points 
of Proposition \ref{p:couple}.
The point (3) is obvious, and (4) follows from the fact 
that for a given singular chain 
$y\in H_*(\wh n_s(v), \wh n_{s-1}(v))$
there is only finite number of $h\in H$, \sut~ 
$\supp y\cap S(\wh q,v)\not= 0$.

\bepr\lb{p:eqcoup}
\been\item
Let $L$ be an $s$-dimensional oriented submanifold of $\hdaw$, \sut~
$L\sbs \hpws$
and $L\sm \Int\hpws$
is compact. Then
$L'=\stind {(-\wh v)}b{a_{s+1}}(L)$
 is an $s$-dimensional oriented submanifold
of $\wh V'$, \sut~$L'\sbs \wh n_s(v)$
and
$L'\sm \Int \wh n_{s-1}(v)$
is compact.
If $L'\pitchfork  S(\wh q,v)$
 then
$$
L'\krest  S(\wh q,v)=\langle \wh q, [L]\rangle_v.
$$
\item
Let $U$ be a \nei~ of $S(f)$.
There is $\d>0$, \sut~ for every $\phi$-gradient
$w$ with $v\vert U=w\vert U$
and
$\Vert v-w \Vert<\d$ we have:
$\langle \wh q, x\rangle_v =
\langle \wh q, x\rangle_w$.

\item Let $h\in H$. Then 
$\langle \wh q h, xh\rangle_v =
\langle \wh q, x\rangle_v$.

\item
For every $x$ the set of $h\in H$ \sut~
$\langle \wh qh , x\rangle_v \not= 0$
is finite.
\enen $\qs$
\enpr
\chapter[$C^0$-generic rationality of  the boundary
 operators]{$C^0$-generic rationality of the boundary
operators in Novikov Complex}
\label{ch:gpbon}

In this chapter we prove the theorems which are the aim of our work.
We refer to the corresponding parts of the Introduction (\ref{su:fsp} and
\ref{su:opmt})
for the exposition of the main ideas of the proof 
and we proceed directly to the precise definitions and statements
of results.

\section{Basic terminology}
\label{s:bt}

In the Subsections
\ref{su:nrbt},
\ref{su:mmbt}
we recall the basic definitions concerning
 algebraic (Subsection \ref{su:nrbt})
and geometric (Subsection \ref{su:mmbt})
framework of Novikov theory, so that our
exposition is self-contained at least in what concerns the terminology.

\subsection{Novikov rings}
\label{su:nrbt}

We shall denote the integer Laurent polynomial ring $\ZZZ[t, t^{-1}]$
by $L$.
The very first example of Novikov ring is the ring
$\ZZZ((t))$, which consists by definition
of all the Laurent
series
$\l=\sum_{-\infty}^\infty a_nt^n$
with integer coefficients and a finite negative part
(so that
$\l\in\ZZZ((t))$ if
$\forall n : a_n\in\ZZZ$ and there is
$N=N(\l)$ \sut~ $a_n=0$ for
$n<N(\l)$).
We shall often denote
$\ZZZ((t))$ by $\wh L$.
The localization of $L$
\wrt~to the set of all the polynomials of the form
$1+tP(t)$, where $P$ is a polynomial, is denoted by
$\wi L$, we have:
$L\sbs \wi L\sbs \wh L$.

To consider "equivariant versions" of the Novikov Complex
one needs to introduce  some more complicated rings.
Let $G$ be a group,
$\xi:G\to\RRR$
be a group \ho.
Set $\L=\ZZZ G$, and denote by
$\wh{\wh\L}$
the abelian group of all the formal linear combinations of the form
$\sum_{g\in G} n_g g$ (infinite in general).
For $\l\in
\wh{\wh\L}, \l=
\sum_{g\in G} n_g g$
set
as usual
$\supp\l=\{g\in G\mid n_g\not=0\}$.
Consider the subset $\Lxi\sbs
\wh{\wh\L}$
defined by
$$
\Lxi=
\{\l\in\lL
\mid \forall c\in\RRR \mbox{  the set  }\supp\l\cap\xi^{-1}([c,\infty[)
\mbox{ is finite}\}
$$
It is not difficult to see that $\Lxi$
has a natural structure of a ring so that the embedding
$\L\rInto\Lxi$
is a \ho.\footnote{
In all my earlier papers on this subject I used
the notation $\L^-_\xi$ for this ring in order to stress
that the
power series $\l\in\L_\xi^-$ are allowed to have the terms
with $\xi(g)$ an arbitrary large {\it negative}
number. Now it seems to me more important to stress that
the definition of the ring uses the completion
procedure.
The symbol $~\wh{}~~$
is included into the notation of the ring in the original paper
\cite{novidok}, and in many other papers, see e.g. \cite{sikoens}.
 }
The most important for us is the case of a group
epimorphism
$\xi:G\to\ZZZ$.
In this case set $H= \Ker\xi$.
 For $n\in\ZZZ$ denote $\xi^{-1}(n)$ by $G_{(n)}$
and
$\{x\in\ZZZ G\mid\supp x\subset G_{(n)}\}$
by $\ZZZ G_{(n)}$.
Denote $\xi^{-1}(]-\infty,-1])$ by $G_-$ and
$\{x\in\ZZZ G\mid\supp x\subset G_-\}$
by $\ZZZ G_-$. Choose $\t\in G_{(-1)}$.
It is easy to see that $\wh{\ZZZ G}_\xi$ is identified with the
ring of power series of the form
$\{a_{-n}\t^n+...+a_1\t+... \mid a_i\in\ZZZ H\}.$

\subsection{Morse maps $M\to S^1$}
\lb{su:mmbt}

The circle $S^1$ is for us the quotient $\RRR/\ZZZ$.
Let $M$ be a closed connected \ma.
A {\it Morse map $f:M\to S^1$}
is a $\smo$ map, \sut~ that all critical points of $f$
 are non-degenerate.
As in

\sub~\ref{su:bd}
the set of all critical points of
$f$ is denoted by $S(f)$, and the set of all critical points
of index $k$ is denoted by
$S_k(f)$.

We leave to the reader to formulate the definition of
$f$-gradient of a Morse map
(following the definition \ref{defgrad}
and Remark \ref{r:dgop}).
The set of all gradients of $f$ is denoted by $\GG(f)$.
Since $M$ is  closed,
 the integral curves of
vector fields on $M$ are defined on the whole of $\RRR$.
For $p\in S(f)$ the stable manifold $D(p,v)$ and the unstable \ma~
$D(p, -v)$ are defined. These are immersed smooth
manifolds in $M$.
We have the corresponding notions of gradients
satisfying \TA~ and, respectively, \ATA.
The set of all $f$-gradients
satisfying \TA, resp. \ATA,
will be denoted by
$\GG\TT(f)$, resp. by
$\GG\AA(f)$.
(Warning:
In our present setting there are no analogs of
Lemma \ref{dcomp}
Proposition \ref{fundsysord}, mainly by the
 reason that the $v$-trajectories
can pass infinitely many times more and more
 closely to the initial position.)

The next lemma implies that $\GT(f)$
is $C^0$ dense in $\GG(f)$.
This is a version of Kupka-Smale theorem.
(It would be more precise to say: "a version of a result, weaker 
than Kupka-Smale theorem", since we make no assertions
 concerning the closed orbits of the considered vector field.)
The proof of Kupka-Smale theorem given in \cite{peix}
works in our case with only minor modifications.
We leave the details to the reader.

\bele\lb{l:kupsm}
Let $v\in \GG(f)$. Let $U$
be a \nei~ of $S(f)$, $\VV$ be a \nei~ of $v$ in $\GG(f)$
\wrt~$\smo$ topology.
Then there is an $f$-\gr~ $v_0$, \sut~ 
$v_0\in\GT(f), v_0\in\VV, \supp(v-v_0)\sbs U$,
and there is a \nei~ $U_0$ of
$S(f)$, such that  
$v_0|U_0=v|U_0$. $\qs$
\enle

We shall assume that $f$ is not homotopic to 0
(otherwise $f$ lifts to a Morse function $M\to\RRR$
and the situation
is the same as considered previously
in Chapter 1).
There is a unique connected infinite cyclic covering
$\CC:\bar M\to M$
\sut~ $f\circ\CC$ is \hot~ to zero.
The map $f$ lifts to a Morse function
$F:\bar M\to \RRR$.
Let $t$ be the generator of the structure group
($\approx\ZZZ$)
of $\CC$, \sut~for every $x$ we have:
$F(xt)< F(x)$.
(We invite the reader to look at the figure
on the page \pageref{fig:cyccov}
 where one can visualize the majority of the constructions of this
subsection.)
We shall assume that $f$ is primitive, that is
$f_*:H_1(M)\to H_1(S^1)$
is epimorphic
(one can always reduce the situation to this case by considering
the map $F/n$ with some $n\in\NNN$.)
In this case we have:
$F(xt)=F(x)-1$.
To simplify the notation we shall assume that
$1\in S^1$
is a regular value of $f$,
and \th~ every $n\in\NNN$ is a regular value of
$F$.
Denote $f^{-1}(1)$ by $V$.
Set
$V_\a=F^{-1}(\a),\quad W=F^{-1}\([0,1]\), \quad
V^-=F^{-1}\(]-\infty,1]\)$.
The cobordism $W$ can be thought of as the result of
cutting $M$ along $V$.
Note that $V_\a t=V_{\a-1}$.
Denote $Wt^s$ by $W_s$; then $\bar M$
is the union
$\cup_{s\in\ZZZ} W_s$, the neighbor copies
$W_s$ and $W_{s+1}$
intersecting by
$V_{-s}$.
For any $n\in\ZZZ$ the restriction of $\CC$ to $V_n$
is a diffeomorphism
$V_n\to V$.
Endow $M$ with an arbitrary riemannian metric and lift it to
a $t$-invariant riemannian metric on $\bar M$.
Now $W$ is a riemannian cobordism, and actually it is
a cyclic cobordism \wrt~the isometry
$t^{-1}:\dow=V_0\to\daw=V_1$.
The $\CC$-preimage of a subset $A\sbs M$
will be denoted by
$\bar A$.
On the other hand, if
$x\in M$ we shall reserve the symbol
$\bar x$ for liftings of $x$ to
$\bar M$
(that is $\bar x\in \bar M$ and $\CC(\bar x)=x$).
Let $v\in\GG(f)$.
Lift $v$ to $\bar M$;
we obtain a complete vector field on $\bar M$,
which is obviously an $F$-gradient.
To restriction of this $F$-gradient to any cobordism
$W_{[\l, \mu]}=F^{-1}\([\l,\mu]\)$
(where $\l,\mu$ are regular values of $F$)
is an
$(F\mid W_{[\l,\mu]})$-gradient.
We shall denote this $F$-gradient by the same symbol $v$, since \noconf.
Note that if $v$ satisfies \TA,
then the $F$-gradient
$v\mid W_{[\l,\mu]}$
satisfies \TA.

We shall consider "equivariant" versions of Novikov Complex;
that requires considering one more covering
(which can be, for example, universal covering).
Let $\PP:\wh M\to M$
be a regular covering with structure group $G$.
Assume that $f\circ \PP$ is homotopic to zero;
then the covering $\PP$ is factored through
$\CC$ as follows:

\begin{diagram}[LaTeXeqno]
\wh M & \rTo^{p}       &   \bar M         \\
      &  \rdTo_{\PP}    &   \dTo_\CC      \quad\lb{f:covs} \\
       &                 &   M             \\
\end{diagram}

Therefore there is a lifting of $f:M\to S^1$ to
a Morse function $\wh F=\bar F\circ p: \wh M\to\RRR$.
Further, the \ho~$f_*:\pi_1M\to\ZZZ$ can be factored through a \ho
~$\xi:G\to\ZZZ$.
Recall that we assume that $\xi$ is epimorphic. So,
applying the terminology of
\sub~ \ref{su:nrbt},
$H=\Ker\xi, G_-=\xi^{-1}\(]-\infty, 0]\)$.

\subsection{Novikov Complex: the definition}
\label{su:ncd}

We assume the terminology
of
\sub~ \ref{su:mmbt}.
Let $v\in\GT(f)$.
For each point $x\in S(f)$ choose a lifting $\bar x$ of $x$ to
$\bar M$.
Moreover, for each point $x\in S(f)$
choose an orientation of the local stable \ma~ $D(x,v)$,
and lift it to an orientation of the corresponding
local stable \ma~ of $\bar x\in S(F)$.

These choices made, we can define the Novikov Complex.
The $r$th term of this complex is defined to be the free
$\ZZZ((t))$-module
generated by $S_r(f)$.
For $x,y\in S(f)$ satisfying $\ind x=\ind y+1$
define am element
$n(x,y)\in\ZZZ((t))$
called
{\it Novikov incidence coefficient}
as follows:
For $k\in\ZZZ$ consider the set
$\G_k(\bar x, \bar y)$
of all $v$-orbits joining $\bar x$ with $\bar yt^k$.
This is a finite set and for every
$\g\in\G_k(\bar x, \bar y)$
the orientation sign $\ve(\g)=\pm 1$
is defined.
(Indeed, choose $\l,\mu\in\RRR$ so that
$\bar x, \bar y t^k\in F^{-1}\([\l,\mu]\)$;
then all the trajectories of
$\G_k(\bar x, \bar y)$
are in the \cob~
$F^{-1}\([\l,\mu]\)$
and we can apply the argument from the page \pageref{finite}).
Set
$$
n_k(\bar x, \bar y)=\sum_{\g\in \G_k(x,y)}\ve(\g)
$$
and
$$
n(x,y)=\sum_k n_k(\bar x, \bar y)\cdot t^k
$$
Now the boundary operator
is defined on the generators of $C_r$ as follows:
$$
\pr_rx=\sum_{y\in S_{r-1}(f)} y\cdot n(x,y)
$$
Since all $n_k(x,y)$ depend
on $v$
we shall sometimes include $v$ in the notation, writing
$n_k(\bar x,\bar y;v), n(\bar x,\bar y;v)$, and $C_*(v)$.
The dependence on the particular choice of liftings of
critical points to $\bar M$ is less important, and we shall
often write for example
$n(x,y)$ instead of $n(\bar x, \bar y)$.
Note also that the resulting complex
does not depend on the particular choice of the function $f$
for which $v$ is gradient.

The next lemma says that $C_*(v)$ is indeed a chain complex.

\bele\label{l:difnov}
For every $r$ we have:
$\pr_r\circ\pr_{r+1}=0$
\enle
\Prf
\wat~for every $x\in S(f)$
the lifting $\bar x$ is in $W_0$.
Then for every $x,y\in S(f)$ the
Novikov incidence coefficient
$n(\bar x, \bar y; v)$
is actually an element of $\ZZZ[[t]]$.
Denote by $\pr_r^{(n)}$
the quotient \ho~
$$
C_r/t^nC_r\to C_{r-1}/t^{n} C_{r-1}
$$
It suffices to prove that
$\pr_r^{(n)}\circ\pr_{r+1}^{(n)}=0$
for every $n$.
But $\pr_r^{(n)}$
is easily identified with the boundary operator in the Morse complex
corresponding to the \cob~
$F^{-1}\([-n,0]\)$. $\qs$

The main aim of the present paper is to investigate the
properties of the incidence coefficients $n(x,y;v)$,
but not the chain \hot~ properties of
$C_*(v)$.
We hope to give the systematic treatment of these
properties in the second part
of this work.
So we shall leave the next proposition without  proof.
(This proposition is proved in \cite{patou}
for the case when the $f$-gradient $v$ is
an $f$-Gradient. The proof is carried over to the general case without
difficulties.)

\bepr\label{p:homnc}
$$
H_*(C_*(v))\approx
H_*(\bar M)\tens{\ZZZ[t, t^{-1}]}\ZZZ((t))
$$\enpr
$\qs$

\section[Incidence coefficients in $\ZZZ((t))$]
{Condition $(\gC\CC)$, and $C^0$-generic
rationality of Novikov incidence coefficients}
\lb{s:cccgr}

We assume the terminology of
\ref{su:mmbt}.

We  say that $v$ satisfies condition
$(\gC\CC)$
($\gC$ for {\it cellular} and $\CC$ for {\it circle})
if the $(F\vert W)$-gradient $v$ satisfies the condition $(\gC\YY)$
from \sub~\ref{su:cyc}.
The set of $f$-gradients $v$ satisfying $(\gC\CC)$
will be denoted by
$\GG\gC\CC(f)$.
First of all we establish that the properties $(\gC\CC)$
and 
$(\gC\CC\TT)$ are $C^0$-generic.

\bepr\lb{p:cgener}
\been\item
The set $\GCC(f)$
is open and dense in $\GG(f)$ \wrt~ $C^0$-topology.
The set
$\GCCT(f)$ 
is open and dense in $\GT(f)$ \wrt~ $C^0$-topology.
\item
Let $v\in \GG(f)$, let $\VV$ be a $C^0$-\nei~  of $v$
in $\GG(f)$.
Then there is $v_0\in \GCC(f)$, \sut
~$v_0\in\VV$
and
$v=v_0$ in some \nei~  of $S(f)$.
\item
Let $v\in \GT(f)$, let $\VV$ be a $C^0$-\nei~  of $v$
in $\GG(f)$.
Ten there is $v_0\in \GCC(f)$, \sut
~$v_0\in\VV$
and
$v=v_0$ in some \nei~  of $S(f)$.
\enen
\enpr

\Prf It suffices to prove 2) and 3).
2) follows immediately from Theorem \ref{t:cyc}. To prove 3) let $v_0$
be an $f$-\gr~\sut~
$v_0\in \GCC(f)\cap\VV$
and $v_0=v$
in a \nei~of $S(f)$
(such $v_0$ exists by 2)).
Now to obtain an $f$-\gr~satisfying the assertion 3),
 apply Lemma \ref{l:kupsm}, and recall that
$\GCC(f)\cap\VV $ is open. $\qs$

Therefore
$\GCCT(f)
=
\GCC(f)\cap\GT(f)$
is $C^0$-open-and-dense in $\GT(f)$.

\beth\lb{t:cccgr}
Let $v\in \GCCT(f)$. Let $x,y\in S(f)$,
$\ind x=\ind y+1$. Then $n(x,y;v)\in \wi L$,
that is
$n(x,y;v)=\frac {P(t)}{t^mQ(t)}$
where $m\in \NNN$ and
$P,Q$ are polynomials with integer coefficients
and $Q(0)=1$.
\enth
\Prf
Let $v\in\GCC(f)$.
The condition $(\gC\YY)$ provides a Morse function
$\phi_1$ on $V_1$
together with its gradient $u_1$,
and  a Morse function
$\phi_0=\phi_1\circ t^{-1}$ on $V_0$
together with its gradient $u_0=(t)_*(u_1)$.
For every $k\in\ZZZ$ we obtain also a Morse function
$\phi_0\circ t^{-k}$
on $V_{-k}$ and a $\phi_0\circ t^{-k}$-gradient
$u_k=(t^k)_*(u_0)$.
We have also an ordered Morse function
$\phi:W\to[0,1]$, adjusted to
$(f,v)$, with an ordered
sequence
$(a_0,..., a_{n+1})$.
Note that $v\mid W_k$ satisfies $(\gC)$
\wrt~ the string
$\SS_k=\(\d, (t^k)_*(u_0),
(t^k)_*(u_1),
\phi_0\circ t^{-k},
 \phi_1\circ t^{-k}  \)
$.

Let $x,y\in S(f), \ind x=s+1, \ind y=s$.
Set $V'=\phi^{-1}(a_{s+1})$.
We can assume that $\bar x\in W, \bar y\in W_1$.
Then $n_k(x,y;v)=0$
for $k\leq 0$.
Consider the oriented submanifold
$X=D(\bar x,v)\cap V_0$
of $V_0$.

Note that $X\sbs V_0^{\{\leq s \}}$ and that
$X\sm
 \Int V_0^{\{\leq s-1 \}}$
is compact. (Indeed, let
$\Sigma=D(\bar x, v)\cap V'$.
Then $\Sigma$ is an embedded $s$-dimensional sphere in $V'$.
Using the terminology of
\sub~ \ref{su:midp}
write $X=\stexp {(-v0)} (\Sigma)$.
Then the compactness of
$X\sm \Int V_0^{\{\leq s-1\}}$
follows from  Lemma \ref{l:compa}).
Denote by $[X]$
the \hog~class of $X$ modulo
$V_0^{\{\leq s-1\}}$.

Applying Theorem \ref{t:consh} 1) several times
we deduce by induction that for every $k$ the submanifold
$X_k=\stind {(-v)}0{-k} (X)$ of $V_{-k}$ is in
$V_{-k}^{\{\leq s\}}$
and
$X_k\sm \Int (V_{-k})^{\{\leq s-1\}}$
is compact.
Set
$V'_{-(k+1)}=V'\cdot t^{k+1}$, and
$X'_k=\stind {(-v)}{-k}{-k-1+a_{s+1}}(X_k)$.

Set $Y_k=D(\bar yt^k, -v)\cap V'_{-k}$;
this is a cooriented embedded $(n-s-1)$-sphere.
Since $v$ satisfies the \TA~
we have:
$X'_k\pitchfork Y_k$.
The proposition \ref{p:couple} implies that
\begin{equation}
n_k(x,y;v)=
X'_k
\krest Y_k=
\langle \bar y t^k, [X_k]\rangle_v\lbl{f:inck}
\end{equation}

Recall from Section \ref{s:hgd} the \ho~
$$
\HH_s(-v):H_s(V_1^{\{\leq s\}}, V_1^{\{\leq s-1\}})
\to
H_s(V_{0}^{\{\leq s\}}, V_{0}^{\{\leq s-1\}})
$$

associated with the \co~ $W$ and
and set $h=\HH_s(-v)\circ (t_*^{-1})$.
Then $h$ is an endomorphism
of $
H_s(V_0^{\{\leq s\}}, V_0^{\{\leq s-1\}})$.

Apply Theorem \ref{t:consh}
to obtain:
$$
n_k(x,y;v)= \langle \bar y , h^k\([X]\)\rangle_v
$$

Now the proof is finished by the application of the next lemma.

\begin{lemm}
\lb{l:cramer}
Let $G$ be a finitely generated abelian group,
 $A$ be an endomorphism of $G$ and
$\lambda :G\to\ZZZ$
 be a homomorphism.
Then for every $p\in G$ the series
 $\sum_{k\geq 0} \lambda(A^k p)t^k\in \ZZZ[[t]]$
is a rational function of $t$ of
 the form $\frac {P(t)}{Q(t)}$,
where $P,Q$ are polynomials
with integral coefficients
and $Q(0)=1$.
\end{lemm}
{\it   \quad   Proof.}
It is sufficient to consider the
case of free finitely generated abelian group $G$.
Consider  a free finitely generated
$\ZZZ [[ t]]$-module
$R=G[[t]]$ and a homomorphism $\phi : R\to R$,
given by
$\phi = 1-At$. Then $\phi$ is invertible, the inverse homomorphism
is given by the formula
$\phi^{-1}=\sum_{k\geq 0} A^k t^k$.
On the other hand the inverse of $\phi$
is given by Cramer
formulas,  which imply
that each matrix entry of the homomorphism, inverse to $\phi$
is the ratio of two integer polynomials. Moreover the denominator
equals to
$\text{det}~(1-At). \qs\qs$

The next result deals with continuity properties of $n(x,y;v)$
\wrt~ small perturbations of $v$.

\beth\label{t:continc}
Let $v\in\GCCT(f)$. Let $U$ be an open \nei~
of $S(f)$.
Then there is $\d>0$ \sut~ for every
$w\in \GCCT(f)$
 with $\Vert w-v\Vert<\d$ and $v\mid U=w\mid U$,
and for every $x,y\in S(f)$ with $\ind x=\ind y+1$ we have:
$n(x,y;v)=n(x,y;w)$
\enth
\Prf
We shall use the terminology of the proof of \ref{t:cccgr},
but we modify it a bit.
 Namely,
we denote $D(\bar x,v)\cap V_0$
by
$X(v)$ in order to stress the dependence on the vector field $v$.
We assume the similar convention for the sets
$X_k, Y_k$,
 the endomorphism $h$, the sphere $\Sigma$ etc.
Thus we have:
$$
n_k(x,y;v)=
\big\langle
\bar y,
\(h(v)\)^k\([X(v)]\)\big\rangle_v\lb{f:inccv}
$$

Recall from Theorem \ref{t:cyc}
that if $w\in \GT(f)$ is sufficiently $C^0$ close to $v$, then
$w$ satisfies $(\gC\CC)$
\wrt~
to the same string
$S=(\d, u_0, u_1, \phi_0, \phi_1)$.
Therefore
the constructions done in the proof of
Theorem \ref{t:cccgr}
apply also to $w$; thus
we obtain the corresponding
objects
$X(w), X_k(w), Y_k(v)$ etc.
Thus the formula (\ref{f:inccv})
holds with $w$ substituted instead of $v$.

The endomorphism
$h(v)=(t^{-1})_*\circ \HH_s(-v)$
is invariant \wrt~ $C^0$-small perturbations of $v$,
see Theorem \ref{t:consh}.
Therefore  we have only to show that
if $v\vert U = w\vert U$ and $\Vert u-w\Vert$
is small enough, we have:
$\langle\cdot, \cdot\rangle_v=
\langle\cdot, \cdot\rangle_w$
and
$[X(v)]=[X(w)]$.
The first property follows from
Proposition \ref{p:couple}.

Proceed to the second one.
We assume the terminology of
\sub~ \ref{su:poh}. In particular, we
 use for the definition of
$\HH, \HH 0$ the triple
$(A,B,C)$ from
(\ref{f:tritri}).
Note that $\Sigma(v)\sbs \Int A$.
Denote by $s(v)$ the fundamental class of $\Sigma(v)$
reduced modulo
$B$; then 
$[X(v)]=\HH 0_s(-v)(s(v))$.
If $w$ is sufficiently $C^0$ close to $v$, then
$\Sigma(w)\sbs\Int A$
and $s(v)=s(w)$ by Proposition \ref{p:discpert} and \ref{p:dischom}.

Lemma \ref{l:h0}
applies to show that $\HH 0_s(-v)=
\HH 0_s(-w)$
and that finishes the proof. $\qs$
\section[Incidence coefficients 
in completions of group rings]
{Incidence coefficients with values in
completions of group rings}
\lb{s:icvc}

\subs{Definitions}
\lb{su:icdef}

In this section we deal with
"equivariant" analogs of Novikov complex.
Let
$f:M\to S^1$
be a Morse map and
$\PP:\wh M\to M$
be a regular covering satisfying $f\circ\PP=0$.
We assume here the terminology of Section \ref{s:bt}.

Let $v\in\GT(f)$.
For every $x\in S(f)$ we choose a lifting
$\wh x$ of $x$ to $\wh M$, and an orientation
of the descending disc $D(x,v)$.
Set $\bar x=p(\wh x)$.
Let $x,y\in S(f), \ind x=\ind y+1$.
The set $\G(x,y)$ of
all $(-v)$-trajectories joining $x$ to $y$
splits  naturally into disjoint parts
defined as follows:
Let $g\in G$. We say that
$\g\in\G_g(\wh x, \wh y)$ if
the lifting of $\g$ to $\wh M$
starting at $\wh x$ finishes at $\wh y\cdot g$.
Then $\{\G_g(\wh x, \wh y)\}_{g\in G}$
is a family of disjoint subsets
of $\G(x,y)$ and
their union equals to $\G(x,y)$.
Note also that the union of
$\G_g(\wh x, \wh y)$ over all $g\in\xi^{-1}(k)$
is in a natural bijective correspondence
with
$\G_k(\bar x, \bar y)$.
In particular, the set
$\cup_{\xi(g)>C}\G_g(\wh x, \wh y)$
is finite for
every
$C\in\RRR$.
Set
\begin{gather}
\wh n_g (\wh x, \wh y)=\sum_{\g\in \G_g(\wh x,\wh y)}\ve(\g)\cdot g\\
\wh n (\wh x, \wh y)=\sum_{g\in G}\wh n_g(\wh x, \wh y)
\end{gather}

Then
$\wh n_g (\wh x, \wh y)$
is an element of $\ZZZ G$
and
$\wh n (\wh x, \wh y)$
is a priori
an element of $\wh{\wh\L}$.
Since for a given $C$ there is only a finite number of
 $(-v)$-trajectories joining $\wh x$ to $\wh y g$
with $\xi(g)>C$, we obtain that
$\wh n (\wh x, \wh y)$
is actually an element of $\Lxi$.
Note that the Novikov incidence coefficient depends on
the particular choice
of the liftings of the critical points to $\wh M$:
it is easy to obtain the following formula:
\begin{equation}
\wh n(\wh xg_1, \wh yg_2)= g_2^{-1} \wh n(\wh x, \wh y) g_1
\lbl{f:chbase}
\end{equation}

\bere
Let $\wh C_r(f)$ be the free right
$\Lxi$-module generated by $S_r(f)$.
Similarly to the \sub~ \ref{su:ncd}
one can prove that
setting for
$x\in S_r(f)$
$$
\wh\pr_r(x)=\sum_{y\in S_{r-1}(f)} x\cdot \wh n(\wh x, \wh y)
$$
we obtain a \ho~ of
$\wh C_r(v)$ to
$\wh C_{r-1}(v)$
satisfying
$\wh\pr_{r+1}
\circ
\wh\pr_r=0$.
The  \hog~  $H_*(\wh C_*(v))$ of the resulting complex is
isomorphic to the \hog~of the completed
chain complex
$C_*^\D(\wh M)\tens{\L}\Lxi$.
We shall not give the proofs  here. See \cite{patou} where these
theorems
are proved in the case when $v$ is an $f$-Gradient.
\enre

\subs{Some non-commutative algebra}
\lb{su:ncalg}
In order to formulate an "equivariant" analog of theorem
\ref{t:cccgr}
we need first of all an analog of the localization ring $\wi L$.
This analog will be introduced via the following non-commutative
 localization procedure.
\pa
{\bf Recollection} (\cite{cohn}, p. 255).
Let $R$ be a (non-commutative) ring.
Assume that for every $n$ a subset $S_n$ of the ring
$\Mat(n\times n, R)$ is given; set $S=\sqcup_n S_n$. For
each matrix $\a\in S_n$
and each pair $1\leq i, j\leq n$ add to the ring $R$
the generators $\bar\a_{i,j}$ subject to the relations
$A\cdot\bar A=\bar A\cdot A=1$, where
$\bar A=(\bar\a_{i,j})$.
We obtain a ring, denoted by $S^{-1} R$, together with a \ho~
$\ell:R\rTo S^{-1} R$.
If $R$ is a subring of a ring $L$, and every $\a\in S$ is
invertible over $L$,
then there is
a commutative diagram:

\begin{diagram}[LaTeXeqno]
R & \rInto       &   L         \\
      &  \rdTo &     \uTo_\ell         \quad\lb{f:incloc} \\
       &                 &   S^{-1} R          \\
\end{diagram}

Apply this procedure to our ring $\L$, setting
$$
S_n=\{1+A\mid A\in\Mat(n\times n, \ZZZ G_{(-1)})\}
$$
We obtain a ring which will be denoted by
$\wi \L_\xi$.
Since every matrix in $S_n$
is invertible in $\Lxi$
(the inverse of $1+A$ being given by the formula
$\sum_{i\geq 0} (-1)^i A^i$)
we have a commutative diagram:

\begin{diagram}[LaTeXeqno]
\L & \rInto       &   \Lxi         \\
      &  \rdTo &     \uTo_\ell         \quad\lb{f:novloc} \\
       &                 &   \wi \L_\xi          \\
\end{diagram}

We shall see that
the \nics~
of $C^0$ generic gradients of Morse maps $M\to S^1$
are  between the
simplest possible
elements of $\LLxi$.
They will be all of the form
$g_1Y(1+A)^{-1}Xg_2$, where
$A\in\Mat(m\times m, \ZZZ G_{(-1)}),\quad X,Y\in
(\ZZZ H)^m,\quad g_1, g_2\in G$.
We shall give a special name for this class of elements of $S^{-1}R$.

\bede\lb{d:typeL}
An element $a\in\Lxi$
{\it is of type $(\LL)$}
if there are
vectors
$X=(X_i)_{1\leq i\leq n},
Y=(Y_i)_{1\leq i\leq n}$ with $X_i, Y_j\in\ZZZ H$,
elements $g_1, g_2\in G$
and an $(m\times m)$-matrix
$A=(a_{ij})$ with entries in $\ZZZ G_{(-1)}$
\sut~
$$
a=g_1\sum_{s\geq 0}\bigg(\sum_{i,j} Y_i a_{i j}^{(s)} X_i\bigg)g_2
$$
where $a^{(s)}_{ij}$
are the matrix entries of $A^s$.
\end{defi}

The elements of type $(\LL)$
are obviously in $\Im\ell$.

In the rest of this subsection we develop some basic
non-commutative linear algebra over the ring $R=\ZZZ H$. This
non-commutative linear-algebraic constructions will be used in
the following.
Let $F$ be a right $R$-module.
Recall that we have chosen an element $\theta\in G_{(-1)}$.
We shall say that an abelian
group \ho~ $\xi:F\to F$
is $\theta$-semilinear, if
$\xi(x\l)=\xi(x)\cdot \theta\l\theta^{-1}$
for every $\l\in R$.
Assume that $F$ is a free right $R$-module with free
base $(e_i)_{1\leq i\leq m}$.
Associate to each $\theta$-semilinear \ho~
$\xi:F\to F$ the matrix
$\wh\xi=(\xi_{ij})_{1\leq i,j\leq m}$ as follows:
$\xi(e_i)=\sum_{j=1}^m e_j\xi_{ji}$.

\bele\label{l:pow}
Let
$x=\sum_{i=1}^m e_ix_i$
be a vector in $F$, and $n\in\NNN$.
Then
\begin{equation}
\xi^n(x)=
\sum_{i,j}e_j\eta_{ji}^{(n)} x_i\theta^{-n}\lbl{f:pow}
\end{equation}
where $\eta_{ji}^{(n)}$
stands for the matrix entries of the matrix $(\wh\xi\theta)^n$.
\enle
\Prf Standard computation by induction in $n$ $\qs$

\begin{coro}\lb{c:summa}
Let $\l:F\to R$
be a \ho~ of right $R$-modules. Then
the element
$$
\sum_{k\geq 0} \l(\xi^k(x))\cdot\theta^k$$
of $\Lxi$ is of type $(\LL)$.
\end{coro}
\Prf
Just apply
$\l$ to
the right hand side of
(\ref{f:pow}). $\qs$

The homomorphisms $F\to R$
to which we shall apply
our Corollary, will be constructed with the help of the
following procedure.
Let $l:F\to\ZZZ$
be a \ho~  of abelian groups \sut~  for every $x\in F$ the set
$\{h\in H\mid l(xh)\not= 0\}$
is finite.
Define a \ho~
$\underline l :F\to R$
as follows:
$$
\underline{ l}(x)=\sum_{h\in H} l(xh^{-1})\cdot h
$$

Then $\underline{ l}$ is obviously a right $\ZZZ H$-module \ho.

\subs{Elements of exponential growth in Novikov rings}
\lb{su:eeg}

Let $x\in\ZZZ((t)), x=\sum_nx_nt^n$. We say,
that $x$ is {\it of exponential growth},
if there are $A,B>0$ \sut~
$\vert x_n\vert\leq Ae^{nB}$ for every $n$.
The elements of exponential growth form obviously 
a subring of $\ZZZ((t))$.
In the present subsection we suggest an analog of 
this notion for elements
of Novikov completion of group rings.

We assume the terminology of
\sub~ \ref{su:nrbt}.
In particular, $\L=\ZZZ G$. For an element
$\l\in\L, \l=\sum_{g\in G} n_g g$
set $\vert\l\vert= \sum\vert n_g\vert$.
Thus $\vert\cdot\vert$ is the restriction to $\ZZZ G$
of the $L^1$-norm on the Banach algebra $L^1(G)$,
see \cite{bousp}, Ch 1, \S 2.2.
In particular, we have: $\vert\l\m\vert\leq \vert\l\vert
\cdot\vert\mu\vert$.
Recall that $\Lxi$ is the Novikov ring associated with a homomorphism
~$\xi:G\to\RRR$.
For $\l\in\Lxi, \l=\sum_{g\in G} n_g g$ and  $c\in\RRR$
set
$\l_{\ug{c}} = \sum_{\xi(g)\geq c} n_g g$.
It follows from the definition of $\Lxi$ that
$\l_{\ug c} \in\ZZZ G$ for every $c$.

\bede\lb{d:eg}
We say that $\l\in\Lxi$ is {\it of exponential growth}
if there are $A,B>0$, \sut~ for every
$c\in\RRR$ we have:
$\vert\l_{\ug c}\vert \leq A e^{-Bc}$.
\footnote{One should probably say
{\it not more than exponential growth}.
We have chosen our present terminology for brevity.}
The subset of all $\l\in\Lxi$ of \eg
~will be denoted by
$\L_\xi^{e.g.}$.
\end{defi}

It is easy to see that
$\L_\xi^{e.g.}$
is a subring of $\Lxi$.

\bele
Elements of type $(\LL)$ are of \eg.
\enle
\Prf
Let $\l\in\Lxi$
be a element of type
$(\LL)$, i.e.
$$\l=g_1\bigg(\sum_{s\geq 0}
X_i a_{ij}^{(s)} Y_j\bigg) g_2,
 \quad A^s=(a^s_{ij} ),
\quad
a_{ij}\in\ZZZ G_{(-1)},
\quad
g_i\in G
$$
We  can assume
$g_1=g_2=1,
X_i=\delta_{ik},
Y_i=\delta_{js}$.
Thus we must show that each  matrix
entry of $\sum_{s\geq 0} A^s$
is in
$\L_\xi^{e.g.}$.
Set
$\Vert A\Vert=\max_{i,j} \vert a_{ij}\vert$.
Then for every
$A,B\in\Mat(m\times m, \ZZZ H)$
we have
$\Vert AB\Vert\leq \Vert A\Vert \Vert B\Vert m$,
and for every $s\geq 0$
we have
$\Vert A^s\Vert \leq
\Vert A\Vert^s m^{s-1}$. This implies easily that if
$N$ is a natural number greater than
$\max(2, \Vert A\Vert m)$, then
$\Vert \sum_{s=0}^k A^s\Vert\leq N e^{k\ln N}$
and our result is proved. $\qs$

\subs{Rationality of \nics}
\lb{su:rnic}

\beth\lb{t:eqcccgr}
Let $v\in\GCCT(f)$. Then for every $x,y\in S(f)$
with $\ind x=\ind y+1$
the Novikov incidence coefficient
$\wh n(\wh x,\wh y;v)$ is of type $(\LL)$.
Therefore
$\wh n(\wh x,\wh y;v)$ is in $\Im\ell$, and is of exponential growth.
\end{theo}
\Prf
The proof follows the lines of the proof of Theorem
\ref{t:cccgr}, and we shall use the terminology from there.
Choose furthermore an element
$\theta\in G_{(-1)}$.
In view of (\ref{f:chbase})
we can choose liftings $\wh x$ of $x$ and $\wh y$ of $y$ to
$\wh M$ as we like.
Choose them so that their projections
$p(\wh x), p(\wh y)$ to $\bar M$ satisfy
$p(\wh x)\sbs W, p(\wh y)\sbs W_1$.
Set $\bar x=p(\wh x), \bar y=p(\wh y)$.

Thus  the descending disc $D(\bar x, v)$ lifts to
$\wh M$ in a natural way, and the manifolds $X_k\sbs V_{-k}$
are also lifted to some manifolds
$\wh X_k\sbs \wh V_{-k}$.
Note that for every $\r\in H$ the point
$\wh y\r\theta^k$ is a lifting of
$\bar yt^k$.
The Novikov incidence coefficient
$\wh n(\wh x, \wh y; v)$
equals to
$\sum_{k\geq 0, \r\in H} \wh n_{\r\theta^k}
(\wh x, \wh y; v)$.

To compute the generic term of this sum we
apply Proposition \ref{p:eqcoup}
to obtain a formula analogous to
(\ref{f:inck}):
\begin{equation}
\wh n_{\r\theta^k}(\wh x, \wh y)=\langle
\wh y\r\theta^k, [\wh X_k]\rangle_v\cdot \rho\theta^k
\lbl{f:eqinck}
\end{equation}
Define a $\theta$-semilinear endomorphism
$\wh h$ of $\wh C_s(u_0)$ setting
$\wh h=\theta^{-1}_*\circ \wh\HH_s(-v\mid W_1)$.
Applying Theorem
\ref{t:eqhgd}
it is not difficult to see that
$[\wh X_k]=\theta^k_*\circ\wh h^k$.
Therefore
$$
\wh n_{\r\theta^k}(\wh x, \wh y)
=
\langle \wh y\r, \wh h^k([\wh X])\rangle=
\langle\wh y, \wh h^k\([\wh X]\)\rho^{-1}\rangle\cdot\rho
$$
Denote by
$\Upsilon$ the \ho~
$\langle\wh y, \cdot\rangle:
\wh C_l(u_0)\to\ZZZ$
and
set $\l=\underline{\Upsilon}$, see the \sub~
\ref{su:ncalg}
for the meaning of $\underline{~}$).
Then
$$
\sum_\rho\wh n_{\r\theta^k}(\wh x, \wh y)
\cdot\r\theta^k
=
      \l \( \wh h^k([\wh X])\)\theta^k
$$

Thus
$$\wh n(\wh x, \wh y; v)=
\sum_{k\geq 0}\l\( \wh h^k([\wh X])\)\theta^k$$
and it suffices to apply Corollary
\ref{c:summa}
to finish the proof.
$\qs$

\newpage
\chapter*{Figures}\lb{ch:fig}

\newpage
\psset{xunit=.7cm,yunit=.7cm}
\begin{pspicture}(1,-3)(12,12)

\rput(0,0){
\psframe[fillstyle=solid, fillcolor=yellow](0,0)(6,6)

\psline{->}(3,-1)(3,7)
\psline{->}(-1,3)(7,3)

\rput(6.5,2.7){$x_+$}
\rput(3.3,7){$x_-$}
\rput(4,5){$Q_r$}
}

\rput(15,0){
\psframe[fillstyle=solid, fillcolor=yellow](0,0)(6,6)
\psbezier[linewidth=0.8mm](3,3)(3,4)(3,5)(4.5,6.5)
\psbezier[linewidth=0.8mm](3,3)(3,2)(3,1)(1.5,-0.5)
\psbezier[linewidth=0.8mm](3,3)(4,3)(5,3)(6.5,4.5)
\psbezier[linewidth=0.8mm](3,3)(2,3)(1,3)(-0.5,3.5)
\psline{->}(3,-1)(3,7)
\psline{->}(-1,3)(7,3)
\rput(6.8,4.6){$W_+$}
\rput(4.8,6.6){$W_-$}
\rput(6.5,2.7){$ x_+$}
\rput(3.3,7){$x_-$}

\psline[linewidth=0.01mm](-1,-1)(7,7)
\psline[linewidth=0.01mm](-1,7)(7,-1)
}

\rput(0,-15){
\psframe(0,0)(6,6)

\psline[fillstyle=solid, fillcolor=yellow](0,0)(3,3)(6,0)(0,0)
\psline[fillstyle=solid, fillcolor=yellow](0,6)(3,3)(6,6)(0,6)
\psline{->}(3,-1)(3,7)
\psline{->}(-1,3)(7,3)
\psline(-1,-1)(7,7)
\psline(-1,7)(7,-1)
\rput(4,5){$Q_r^-$}
\rput(6.5,2.7){$x_+$}
\rput(3.3,7){$x_-$}
}

\rput(15,-15){
\psframe(0,0)(6,6)

\psline[fillstyle=solid, fillcolor=yellow](0,0)(3,3)(0,6)(0,0)
\psline[fillstyle=solid, fillcolor=yellow](6,6)(3,3)(6,0)(6,6)
\psline{->}(3,-1)(3,7)
\psline{->}(-1,3)(7,3)
\psline(-1,-1)(7,7)
\psline(-1,7)(7,-1)
\rput(5,4){$Q_r^+$}
\rput(6.5,2.7){$x_+$}
\rput(3.3,7){$x_-$}
}

\end{pspicture}
\label{fig:morse}

\newpage

\psset{xunit=1cm,yunit=1cm}
\begin{pspicture}(0,0)(16,16)

\psline{->}(-1,0)(16,0)
\psline{->}(0,-1)(0,16)

\psline[linewidth=0.8mm](0,3)(3,3)
\psbezier[linewidth=.8mm](3,3)(3.2,3)(3.4,3.2)(3.4,5)
\psbezier[linewidth=.8mm](3.4,5)(3.4,6.4)(3.4,6.6)(4,6)

\rput(2.7,-.5){$\mu_p$}
\rput(4.3,-.5){$\mu_p+\delta$}
\rput(8.7,-.5){$r_p-\delta$}
\rput(10.3,-.5){$r_p$}
\rput(15,-.5){$t$}

\rput(13,2){$y(t)=\frac C{2t}$}

\rput(-.5,3){$1$}
\rput(-.5,6){$C/B$}
\rput(-.5,15){$\lambda(t)$}

\psline(0,6)(12,6)
\psline(0,3)(12,3)

\psline(3,0)(3,10)
\psline(4,0)(4,10)
\psline(9,0)(9,10)
\psline(10,0)(10,10)

\pscurve(2,12)(2.4,10)(3,8)(4,6)(6,4)(8,3)(10,2.4)(12,2)
\pscurve[linewidth=.8mm](4,6)(4.8,5)(5,4.8)(6,4)(8,3)(9,2.67)

\psbezier[linewidth=.8mm](9.6,4)(9.6,5.8)(9.8,6)(10,6)
\psbezier[linewidth=.8mm](9,2.67)(9.4,2.6)(9.6,2.6)(9.6,4)

\psline[linewidth=0.8mm](10,6)(12,6)

\end{pspicture}
\label{fig:lambda}

\newpage
\psset{xunit=.45cm,yunit=.45cm}
\begin{pspicture}(0,0)(12,20)
\rput(0,0){
\psellipse(5.5,18)(1.5,0.6)
\pscurve[linestyle=dashed](4,2)(4.1,2.2)(4.4,2.4)
(5.5,2.6)(6.6,2.4)(6.9,2.2)(7,2)
\pscurve(4,2)(4.1,1.8)(4.4,1.6)(5.5,1.4)
(6.6,1.6)(6.9,1.8)(7,2)
\psline(4,2)(4,18)
\psline(7,2)(7,18)
\rput(3.2,2){$V_i$}
\rput(3.1,18){$V_{i+1}$}
\rput(8,10){$ W_i$}

\rput(5,-2){ A}

}
\rput(10,-12){\psellipse(5.5,18)(1.5,0.6)
\pscurve[linestyle=dashed](4,2)(4.1,2.2)
(4.4,2.4)(5.5,2.6)
(6.6,2.4)(6.9,2.2)(7,2)
\pscurve(4,2)(4.1,1.8)(4.4,1.6)(5.5,1.4)
(6.6,1.6)(6.9,1.8)(7,2)
\psline(4,2)(4,7)
\psline[linestyle=dashed](4,10.6)(4,11.2)
\psline(4,11.2)(4,18)
\psline(7,2)(7,7)
\psline[linestyle=dashed](7,10.6)(7,11.2)
\psline(7,11.2)(7,18)
\pscurve(4,7)(4,7.4)(3.7,8)(3,9)(2,10.3)(1.2,11)
(1,11.8)(1,12)
\pscurve(7,7)(7,7.4)(7.3,8)(8,9)(9,10.3)(9.8,11)
(10,11.8)(10,12)
\pscurve(2,12)(3,11.4)(4,11.2)(5.5,11)
\pscurve(9,12)(8,11.4)(7,11.2)(5.5,11)
\pscurve(1,12)(1,12.1)(2,12.8)(3,13)(4,13.1)
\pscurve(7,13.1)(8,13)(9,12.8)(10,12.1)(10,12)
\pscurve[linestyle=dashed](4,13.1)(5.5,13.2)(7,13.1)
\pscurve[linestyle=dashed](2,12)(3,11)(3.6,10.5)
(3.8,10.4)(4,10.6)(4,10.8)
\pscurve[linestyle=dashed](9,12)(8,11)(7.4,10.5)
(7.2,10.4)(7,10.6)(7,10.8)

\rput(5,-2){ B}

}

\rput(20,-24){\psellipse(5.5,18)(1.5,0.6)
\pscurve[linestyle=dashed](4,2)(4.1,2.2)(4.4,2.4)
(5.5,2.6)(6.6,2.4)(6.9,2.2)(7,2)
\pscurve(4,2)(4.1,1.8)(4.4,1.6)(5.5,1.4)(6.6,1.6)
(6.9,1.8)(7,2)
\psline(4,2)(4,7)
\psline(4,11.2)(4,18)
\psline(7,2)(7,7)
\psline(7,11.2)(7,18)
\pscurve(4,7)(4,7.4)(3.7,8)(3,9)(2,10.3)(1.2,11)
(1,11.8)(1,12)
\pscurve(7,7)(7,7.4)(7.3,8)(8,9)(9,10.3)(9.8,11)
(10,11.8)(10,12)
\pscurve(2,12)(2.4,11.6)(3,11.6)(3.4,11)(4,11.2)
(4.6,11.4)
(5.5,10.8)(6.4,11.4)(7,11.2)(7.4,11)(8,11.6)
(8.6,11.6)(9,12)
\pscurve(1,12)(1,12.1)(1.6,12.4)(2,12.6)
(2.6,13.1)(3,13)(3.4,12.8)(4,13.2)
\pscurve(7,13.2)(7.6,12.8)(8,13)(8.4,13.1)
(9,12.6)(9.4,12.4)(10,12.1)(10,12)

\rput(5,-2){ C}

}
\end{pspicture}
\label{fig:skl}

\newpage

\psset{xunit=1cm,yunit=1cm}

\begin{pspicture}(0,0)(16,16)

\psline{->}(-1,0)(16,0)
\psline{->}(0,-1)(0,16)

\psline(-.5,-.5)(12.5,12.5)
\psline[linewidth=0.01cm](-.5,12.5)(12.5,-.5)

\psline(0.5,0)(0.5,12.5)
\psline(4,0)(4,12.5)
\psline(5.5,0)(5.5,12.5)
\psline(6.5,0)(6.5,12.5)
\psline(6,0)(6,12.5)
\psline(8,0)(8,12.5)
\psline(11.5,0)(11.5,12.5)
\psline(12,0)(12,12.5)

\psline(0,4)(12.5,4)
\psline(0,8)(12.5,8)
\psline(0,12)(12.5,12)

\psline[linewidth=.8mm](0,0)(0.5,0.5)
\psbezier[linewidth=.8mm](.5,.5)(1.5,1.5)(2,2)(2.5,4)
\psbezier[linewidth=.8mm](2.5,4)(3.5,8)(3.6,8)(4,8)

\psbezier[linewidth=.8mm](4,8)(4.4,8)(4.5,8)(5,7.2)
\psbezier[linewidth=.8mm](5,7.2)(5.3,6.7)(5.4,6.6)(5.5,6.5)

\psline[linewidth=.8mm](5.5,6.5)(6,6)

\rput{-180}(12,12)
{
\psline[linewidth=.8mm](0,0)(0.5,0.5)
\psbezier[linewidth=.8mm](.5,.5)(1.5,1.5)(2,2)(2.5,4)
\psbezier[linewidth=.8mm](2.5,4)(3.5,8)(3.6,8)(4,8)

\psbezier[linewidth=.8mm](4,8)(4.4,8)(4.5,8)(5,7.2)
\psbezier[linewidth=.8mm](5,7.2)(5.3,6.7)(5.4,6.6)(5.5,6.5)
\psline[linewidth=.8mm](5.5,6.5)(6,6)

 }

\rput(0.5,-.5){$\epsilon$}
\rput(4,-.5){$\frac 13$}
\rput(5.2,-.5){$\frac 12 - \! \epsilon$}
\rput(6.8,-.5){$\frac 12 + \! \epsilon$}
\rput(6,-.5){$\frac 12$}
\rput(8,-.5){$\frac 23$}
\rput(11.2,-.5){$1-\!\epsilon$}
\rput(12.2,-.5){$1$}
\rput(15,-.5){$x$}

\rput(-.5,4){$\frac 13$}
\rput(-.5,8){$\frac 23$}
\rput(-.5,12){$1$}
\rput(-.5,15){$\chi(x)$}

\end{pspicture}

\label{fig:chi}

\newpage

\begin{pspicture}(0,0)(12,20)

\rput(8,13){$\Bigg\}$}
\rput(8,15){$\Bigg\}$}

\rput(8.4,13){$W_1$}
\rput(8.4,15){$W$}

\psline{->}(9,17)(11,17)
\rput(10,17.5){$F$}

\psdots[dotscale=1](12,2)(12,14)(12,16)(12,12)
\rput(12.5,2){$-n$}
\rput(12.5,12){$-1$}
\rput(12.5,14){$0$}
\rput(12.5,16){$1$}

\rput(0,10){
\pscurve(4,2)(4.1,1.8)(4.4,1.6)(5.5,1.4)
(6.6,1.6)(6.9,1.8)(7,2)
\rput(7.5,2){$V_{-1}$}
            }

\rput(0,12){
\pscurve(4,2)(4.1,1.8)(4.4,1.6)(5.5,1.4)
(6.6,1.6)(6.9,1.8)(7,2)
\rput(7.5,2){$V_0$}
            }

\rput(0,14){
\pscurve(4,2)(4.1,1.8)(4.4,1.6)(5.5,1.4)
(6.6,1.6)(6.9,1.8)(7,2)
\rput(7.5,2){$V_1$}
            }

\pscurve(4,2)(4.1,1.8)(4.4,1.6)(5.5,1.4)
(6.6,1.6)(6.9,1.8)(7,2)
\psline(4,1)(4,18)
\psline(7,1)(7,18)
\rput(7.5,2){$V_{-n}$}

\psline(12,1)(12,18)
\psline[linewidth=0.01cm](0,5)(0,14)
\psline[linewidth=0.01cm](0,14)(2,14)

\rput(.5,13.5){$V^-$}

\end{pspicture}
\label{fig:cyccov}


\newpage
\begin{thebibliography}{99}\label{refer}
\bibitem{abrob}
\abrob


\bibitem{Arno}
\Arno


\bibitem{Arnoprob}
\Arnoprob

\bibitem{biro}
\biro


\bibitem{bousp}
\bousp

\bibitem{cerf}
\cerf


\bibitem{cohn}
\cohn


\bibitem{dold}

\dold


\bibitem{farrell}
\farrell



\bibitem{farber}
\farber


\bibitem{fried}
\fried


\bibitem{hirsch}
\hirsch


\bibitem{huli}
\huli


\bibitem{kling}
\kling

\bibitem{lam}
\lam


\bibitem{latour}
\latour




\bibitem{laudsiko}
\laudsiko

\bibitem{mengt}
\mengt



\bibitem{milnhcob}
\milnhcob



\bibitem{milnstash}
\milnstash




\bibitem{morse}
\morse


\bibitem{milnWT}
\milnWT


\bibitem{novidok}
\novidok



\bibitem{novshme}
\novshme



\bibitem{novikirh}
\novikirh



\bibitem{noviuspe}
\noviuspe




\bibitem{noviquasi}
\noviquasi






\bibitem{patou}
\patou



\bibitem{pasur}
\pasur




\bibitem{pamrl}
\pamrl



\bibitem{pator}
\pator





\bibitem{paodense}
\paodense






\bibitem{paura}
\paura




\bibitem{paepr}
\paepr


\bibitem{paepri}
\paepri



\bibitem{padok}
\padok



\bibitem{pasbor}
\pasbor



\bibitem{pastpet}
\pastpet



\bibitem{peix}
\peix


\bibitem{pozniak}
\pozniak


\bibitem{ranprepr}
\ranprepr


\bibitem{rani}
\rani

\bibitem{sikothese}
\sikothese



\bibitem{sikoens}
\sikoens



\bibitem{smale}
\smale


\bibitem{smapoi}
\smapoi





\bibitem{thom}
\thom





         \bibitem{tur}
\tur


\bibitem{turtur}
\turtur



\bibitem{witt}
\witt





\end{thebibliography}
\end{document}